\documentclass[a4paper,openany,11pt]{memoir}

\maxtocdepth{subsection}
\settocdepth{subsection}

\usepackage{amsmath, amsthm, amssymb, amsfonts}
\usepackage[labelfont={sf}]{caption}
\usepackage{cite}
\usepackage{colonequals}
\usepackage[autostyle]{csquotes}
\usepackage{enumerate}
\usepackage{enumitem}
\usepackage{euscript}
\usepackage{fancyhdr}
\usepackage{float}
\usepackage[T1]{fontenc}
\usepackage{framed}
\usepackage[utf8]{inputenc}
\usepackage[letterpaper,hmargin=1in,vmargin={0.75in,0.2in}]{geometry}
\usepackage{lmodern}
\usepackage{hyperref}
\usepackage{mathrsfs}
\usepackage{mathtools}
\usepackage{multicol}
\usepackage{parskip}
\usepackage{relsize}
\usepackage{stmaryrd}
\usepackage[svgnames]{xcolor}
\usepackage{thmtools}
\usepackage{tikz}
\usepackage{titlesec}
\usepackage{upgreek}
\usepackage{xfrac}
\usepackage{yfonts}

\newcommand{\arr}{\mathbb{R}}
\newcommand{\cee}{\mathbb{C}}

\newcommand{\dee}{\mathrm{d}}
\newcommand{\deee}{\,\mathrm{d}}

\newcommand{\enn}{\mathbb{N}}
\newcommand{\ennplus}{\enn_+}
\newcommand{\rzee}{{\mathbb{R}/\mathbb{Z}}}

\newcommand{\zee}{\mathbb{Z}}

\DeclareMathOperator{\ach}{\mathbf{H}}
\DeclareMathOperator{\achn}{\mathbf{H}_\ast}

\DeclareMathOperator{\einf}{\mathbf{EInf}}

\DeclareMathOperator{\bbr}{\mathrm{BBor}}

\DeclareMathOperator{\fou}{\mathrm{Fou}}

\DeclareMathOperator{\prob}{\mathrm{Prob}}

\renewcommand{\dotfill}{%
  \leavevmode
  \cleaders\hbox to 0.65em{\hss\textbf{.}\hss}\hfill\kern0pt
}
\renewcommand{\epsilon}{\varepsilon}
\renewcommand{\hat}{\widehat}
\renewcommand{\phi}{\varphi}

\newcommand{\ddotfill}{%
  \leavevmode
  \cleaders\hbox to 0.65em{\hss.\hss}\hfill\kern0pt
}

\newcommand{\lvertt}{\left \vert}
\newcommand{\rvertt}{\right \vert}

\newcommand{\dlvertt}{\lvertt \lvertt}
\newcommand{\drvertt}{\rvertt \rvertt}

\newcommand{\eusc}[1]{\EuScript{#1}}

\newcommand{\ov}[1]{\overline{#1}}

\newsavebox{\depbox}
\newlength{\deplabelwidth}
\newlength{\depboxwidth}
\newcommand{\depsep}{\,;\,}
\newcommand{\dep}[2]{\par\noindent\hangindent=\deplabelwidth\hangafter=1%
  \makebox[\deplabelwidth][l]{\underline{#1:}}#2}
\newenvironment{deps}
  {\par\addvspace{6pt}\centering\setlength{\fboxsep}{8pt}%
   \begin{lrbox}{\depbox}\begin{minipage}{\depboxwidth}%
   \setlength{\parskip}{0pt}\setlength{\parindent}{0pt}}
  {\end{minipage}\end{lrbox}\fbox{\usebox{\depbox}}\par\addvspace{6pt}}
\newcommand{\depitem}[2]{[#1, \ref{#2}]}
\newcommand{\depitemc}[3]{[#1, \ref{#2}\ref{#3}]}
\newcommand{\depcont}[1]{\par\noindent\hangindent=\deplabelwidth\hangafter=1\makebox[\deplabelwidth][l]{}#1}

\makeatletter \def\vhrulefill#1{\leavevmode\leaders\hrule\@height#1\hfill \kern\z@}
\makeatother

\makeatletter
\newcommand{\hooktwoheadright}[2][]{%
  \lhook\joinrel
  \ext@arrow 0359\rightarrowfill@ {#1}{#2}%
  \mathrel{\mspace{-15mu}}\rightarrow
}
\makeatother

\titleformat{\part}[display]
  {\normalfont\sffamily\HUGE\bfseries}
  {Part \thepart}{20pt}{\HUGE}

\titleformat{\chapter}[display]
  {\normalfont\sffamily\LARGE\bfseries}
  {\chaptertitlename\ \thechapter}{20pt}{\Huge}
  
\titleformat{\section}
  {\normalfont\sffamily\LARGE\bfseries}
  {\thesection}{1em}{}
  
\titleformat{\subsection}
  {\normalfont\sffamily\Large\bfseries}
  {\thesubsection}{1em}{}
  
\titleformat{\subsubsection}
  {\normalfont\sffamily\Large\bfseries}
  {\thesubsubsection}{1em}{}

\definecolor{gloryred}{HTML}{BB133E}
\definecolor{gloryblue}{HTML}{2C3ABC}
\definecolor{glorygreen}{HTML}{24802A}
\definecolor{gloryorange}{HTML}{FF6C0C}
\definecolor{topgrey}{HTML}{CFCFCF}
\definecolor{midgrey}{HTML}{AAAAAA}
\definecolor{downgrey}{HTML}{828282}

\colorlet{headcolor1}{Black}
\colorlet{rulecolor1.0}{gloryred}
\colorlet{rulecolor1.1}{gloryblue}
\colorlet{rulecolor1.2}{glorygreen}

\newenvironment{leftbar1.0}{
  
  \MakeFramed {\advance\hsize-\width \FrameRestore}}
 {\endMakeFramed}
 
 \newenvironment{leftbar1.1}{
  
  \MakeFramed {\advance\hsize-\width \FrameRestore}}
 {\endMakeFramed}
 
 \newenvironment{leftbar1.2}{
  
  \MakeFramed {\advance\hsize-\width \FrameRestore}}
 {\endMakeFramed}

 \declaretheoremstyle[headfont=\sffamily\bfseries,
 notefont=\sffamily\bfseries,
 notebraces={}{},
 headpunct=,
 bodyfont=\sffamily,
 headformat=\color{headcolor1}\NAME~\NUMBER:~\!\!\NOTE\hfill\smallskip\linebreak,
 preheadhook=\begin{leftbar1.0},
 postfoothook=\end{leftbar1.0},
 ]{customTheorem}
 \declaretheorem[style=customTheorem, numberwithin=chapter]{Theorem}

 \declaretheoremstyle[headfont=\sffamily\bfseries,
 notefont=\sffamily\bfseries,
 notebraces={}{},
 headpunct=,
 bodyfont=\sffamily,
 headformat=\color{headcolor1}\NAME~\NUMBER:~\!\!\NOTE\hfill\smallskip\linebreak,
 preheadhook=\begin{leftbar1.1},
 postfoothook=\end{leftbar1.1},
 ]{customLemma}

 \declaretheoremstyle[headfont=\sffamily\bfseries,
 notefont=\sffamily\bfseries,
 notebraces={}{},
 headpunct=,
 bodyfont=\sffamily,
 headformat=\color{headcolor1}\NAME~\NUMBER\hfill\smallskip\linebreak,
 preheadhook=\begin{leftbar1.2},
 postfoothook=\end{leftbar1.2},
 ]{customProp}

 \declaretheoremstyle[headfont=\sffamily\bfseries,
 notefont=\sffamily\bfseries,
 notebraces={}{},
 headpunct=,
 bodyfont=\sffamily,
 headformat=\color{headcolor1}Proposition~\NUMBER:~\!\!\NOTE\hfill\smallskip\linebreak,
 preheadhook=\begin{leftbar1.2},
 postfoothook=\end{leftbar1.2},
 ]{customPropn}
 \declaretheorem[style=customPropn, numberwithin=section]{Propositionn}

\declaretheorem[style=customLemma, sibling=Propositionn]{Lemma}

\colorlet{headcolor2}{Black}
\colorlet{rulecolor2.0}{topgrey}
\colorlet{rulecolor2.1}{midgrey}
\colorlet{rulecolor2.2}{downgrey}

\newenvironment{leftbar2.0}{%
  \MakeFramed {\advance\hsize-\width \FrameRestore}}%
 {\endMakeFramed}
 
 \newenvironment{leftbar2.1}{%
  \MakeFramed {\advance\hsize-\width \FrameRestore}}%
 {\endMakeFramed}
 
 \newenvironment{leftbar2.2}{%
  \MakeFramed {\advance\hsize-\width \FrameRestore}}%
 {\endMakeFramed}

  \declaretheoremstyle[headfont=\sffamily\bfseries,%
 notefont=\sffamily\bfseries,%
 notebraces={}{},%
 headpunct=,%
 bodyfont=\sffamily,%
 headformat=\color{headcolor2}\NAME~\NUMBER:~\!\!\NOTE{}\hfill\smallskip\linebreak,%
 preheadhook=\begin{leftbar2.0},%
 postfoothook=\end{leftbar2.0},%
 ]{customDefinition}
 
\declaretheorem[style=customDefinition, sibling=Propositionn]{Definition}

 \declaretheoremstyle[headfont=\sffamily\bfseries,%
 notefont=\sffamily\bfseries,%
 notebraces={}{},%
 headpunct=,%
 bodyfont=\sffamily,%
 headformat=\color{headcolor2}\NAME~\NUMBER:~\!\!\NOTE{}\hfill\smallskip\linebreak,%
 preheadhook=\begin{leftbar2.1},%
 postfoothook=\end{leftbar2.1},%
 ]{customCorollary}

 \declaretheorem[style=customCorollary, sibling=Theorem]{Corollary}

 \declaretheoremstyle[headfont=\sffamily\bfseries,%
 notefont=\sffamily\bfseries,%
 notebraces={}{},%
 headpunct=,%
 bodyfont=\sffamily,%
 headformat=\color{headcolor2}\NAME~\NUMBER:~\!\!\NOTE{}\hfill\smallskip\linebreak,%
 preheadhook=\begin{leftbar2.2},%
 postfoothook=\end{leftbar2.2},%
 ]{customRemark}

\declaretheorem[style=customRemark, sibling=Propositionn]{Remark}

\colorlet{headcolor3}{Black}
\colorlet{rulecolor3}{Black}

\newenvironment{leftbar3}{%
  \MakeFramed {\advance\hsize-\width \FrameRestore}}%
 {\endMakeFramed}

 \declaretheoremstyle[headfont=\sffamily\bfseries,%
 notefont=\sffamily\bfseries,%
 notebraces={}{},%
 headpunct=,%
 bodyfont=\sffamily,%
 headformat=\color{headcolor3}\NAME~\!\!\NOTE\hfill\smallskip\linebreak,%
 preheadhook=\begin{leftbar3},%
 postfoothook=\end{leftbar3},%
 qed=$\blacksquare$,
 ]{prf1}

\theoremstyle{Proof}
 
\declaretheorem[style=prf1, numberwithin=chapter]{Proof}

\colorlet{headcolor4}{Black}
\colorlet{rulecolor4}{gloryorange}

\newenvironment{leftbar4}{
  
  \MakeFramed {\advance\hsize-\width \FrameRestore}}
 {\endMakeFramed}

 \declaretheoremstyle[headfont=\sffamily\bfseries,
 notefont=\sffamily\bfseries,
 notebraces={}{},
 headpunct=,
 bodyfont=\sffamily,
 headformat=\color{headcolor4}Theorem~\NUMBER:~\!\!\NOTE\hfill\smallskip\linebreak,
 preheadhook=\begin{leftbar4},
 postfoothook=\end{leftbar4},
 ]{customoTheorem}

 \declaretheoremstyle[headfont=\sffamily\bfseries,
 notefont=\sffamily\bfseries,
 notebraces={}{},
 headpunct=,
 bodyfont=\sffamily,
 headformat=\color{headcolor4}Lemma~\NUMBER:~\!\!\NOTE\hfill\smallskip\linebreak,
 preheadhook=\begin{leftbar4},
 postfoothook=\end{leftbar4},
 ]{customoLemma}
 \declaretheorem[style=customoLemma, sibling=Propositionn]{oLemma}

 \declaretheoremstyle[headfont=\sffamily\bfseries,
 notefont=\sffamily\bfseries,
 notebraces={}{},
 headpunct=,
 bodyfont=\sffamily,
 headformat=\color{headcolor4}Proposition~\NUMBER:~\!\!\NOTE\hfill\smallskip\linebreak,
 preheadhook=\begin{leftbar4},
 postfoothook=\end{leftbar4},
 ]{customoProp}

\numberwithin{equation}{section}

\fancypagestyle{plain}{%
\fancyhf{} % clear all header and footer fields
\fancyfoot[C]{\sffamily \bfseries \fontsize{11pt}{11pt}\selectfont\thepage} % except the center

}
\DeclareMathAlphabet{\mathpzc}{OT1}{pzc}{m}{it}

\DeclareQuoteStyle{single}{\lq}[]{\rq}{}[]{}
\DeclareQuoteStyle{double}{\lq\lq}[]{\rq\rq}{}[]{}

\thanksmarkseries{arabic}

\begin{document}

\sffamily

\setquotestyle{single}

\title{\Huge \bfseries Rudolph-Johnson without dynamics: endomorphic rigidity on the circle}
\author{\LARGE Peter Burton\thanks{\url{peterburton1728@gmail.com}} and Kate Juschenko\thanks{\url{k.juschenko@austin.utexas.edu}}}
\date{\LARGE\today}

\maketitle

\renewcommand{\abstractnamefont}{\Large \sffamily\bfseries}
\renewcommand{\abstracttextfont}{\sffamily}

\begin{abstract} We prove that entropy alone, with no invariance hypothesis, forces a Borel probability measure on the circle toward Lebesgue measure under every sublacunary set of multiplicative endomorphisms $x \mapsto mx$. Here, we refer to an infinite set of multipliers $m$ as sublacunary if its consecutive ratios tend to $1$. Examples of such sets include the primes following perfect cubes, the range of the partition function, the integers $\lfloor n^{\log n} \rfloor$ and the semigroup generated by $2$ and $3$ which is the classical case of Furstenberg's conjecture.\\

Specifically, we show that every measure $\mu$ has a weak-star limit of endomorphs which dominates Lebesgue measure scaled by the upper entropy dimension of $\mu$. These theorems are effective at finite entropy resolutions, and we provide the associated bounds for the primes and the perfect powers. The semigroup case may be interpreted as an invariance-free form of the Rudolph-Johnson theorem. We also show that an analog of Lyons' conjecture at full entropy dimension for every sublacunary multiplier set while failing by an arbitrarily large factor for the lacunary semigroup of powers of $3$. The proofs combine Fourier analysis with the algebraic structure of the multiplicative endomorphisms. \end{abstract}

\renewcommand\printtoctitle[1]{\LARGE\bfseries #1}

\renewcommand*{\cftsectionfont}{\sffamily}
\renewcommand*{\cftsectionpagefont}{\sffamily}

\renewcommand*{\cftsubsectionfont}{\sffamily}
\renewcommand*{\cftsubsectionpagefont}{\sffamily}

\tableofcontents*

\newpage

\section{Invitation to article} \label{sec.ad}

The strain of mathematics in this article originally derives from Furstenberg's $\times 2 \times 3$ conjecture, which we take the time to state although it will play no further role in the article. Reifying the unit circle as $\rzee$, for $n \in \ennplus$ we consider the endomorphism $x \mapsto nx$ of $\rzee$. We also define a Borel probability measure $\mu$ on $\rzee$ to be $\times n$\emph{-invariant} if the following holds for all continuous functions $\psi:\rzee \to \cee$.
\begin{equation*} \int_{\rzee} \psi(nx) \deee \mu(x) = \int_{\rzee}  \psi(x) \deee \mu(x)  \end{equation*}
Furstenberg's $\times 2 \times 3$ conjecture asserts that Lebesgue measure is the unique atomless probability measure on $\rzee$ which is both $\times 2$ and $\times 3$ invariant. The conjecture developed from the material in \cite{MR0213508}, and has remained central in ergodic theory since then. It also appears as Conjecture 3 of \cite{MR2827904}, where Lindenstrauss provides extensive context situating it in homogeneous dynamics and number theory.\\

Perhaps the most significant progress to date on Furstenberg's conjecture itself came with the work of D.~Rudolph in \cite{MR1062766} and A.~Johnson in \cite{MR1194793}, which established the conjecture under the assumption that $\mu$ has positive Kolmogorov-Sinai entropy in the sense of \cite{MR0648108}. Continuing from there, J.~Bourgain, E.~Lindenstrauss, P.~Michel and A.~Venkatesh in their article \cite{MR2563089} obtained an effective version of Rudolph-Johnson's results and moreover situated them in a context which is robust beyond the asymptotic regime of ergodic theory. We take the liberty of abbreviating the four authors of \cite{MR2563089} to BLMV in the sequel. \\

The purpose of the present article is to prove theorems which make a case that the underlying phenomenon in Rudolph-Johnson arises from the interaction of Fourier analysis and the rigid algebra of circle endomorphisms, and ultimately the invariance hypotheses of ergodic theory should be thought of as inessential. What the phenomenon requires is entropy in the measure, and nothing arithmetic of the multipliers beyond sublacunarity, the condition that consecutive multipliers have ratios tending to $1$. We term this phenomenon `endomorphic rigidity', to be interpreted as follows:

\begin{center} \emph{Endomorphic rigidity occurs when a set of endomorphisms of a compact group in some sense forces an arbitrary probability measure toward Haar measure.} \end{center}

This is emphatically in distinction to the phenomenon of \emph{equidistribution} captured by ergodic theorems, whereby a single measure-preserving system forces time-averages of test functions to the space-mean.\\

We now briefly introduce some notation crucial for our context, which will be used to formulate two initial theorems in support of our case.

\begin{itemize}
\item We write $\lambda$ for Lebesgue (Haar) probability measure on $\rzee$, and $e(t) \coloneqq e^{2\pi i t}$ for $t \in \arr$.
\item We identify $\rzee$ with its coset transversal $[0,1)$, and for $n \in \ennplus$ we write $\ach(\mu,n)$ for the Shannon entropy of the distribution induced by $\mu$ on the partition of $\rzee$ into $n$ arcs of equal length starting from $0$. We define the \emph{upper entropy dimension} of $\mu$ as follows; it always lies in $[0,1]$ and is a well-known quantity in information theory. (See, for example, \cite{MR0107575}.)
\begin{equation*} \ach_\ast^\uparrow(\mu) \coloneqq \limsup_{n \to \infty} \frac{\ach(\mu,n)}{\log(n)} \end{equation*}
We refer the reader to Definition \ref{def.entropy} below for entropy-related formalisms.
\item Finally, for $c \in [0,1]_\arr$ we write $\nu \succeq c\lambda$ to mean that the following holds for every nonnegative continuous function $\psi:\rzee \to [0,\infty)$.
\begin{align*} \int \psi \deee \nu \geq c \int \psi \deee \lambda \end{align*}
\end{itemize}

\subsection{Qualitative endomorphic rigidity} \label{sec.equisparse}

Techniques in the vein of Weyl's method \cite{MR1511862} for the squares and Vinogradov's bilinear method \cite{MR62183} for the primes show that many sets $\mathscr{M}$ of natural numbers admit an increasing sequence $(\mathscr{F}_n)_{n \in \enn}$ of nonempty finite subsets with union $\mathscr{M}$ such that the following individual Lebesgue-equidistribution statement holds for every irrational $\theta \in \arr$.
\begin{align} \label{sec.equisparse.1} \lim_{n \to \infty} \frac{1}{\lvert \mathscr{F}_n \rvert} \sum_{\ell \in \mathscr{F}_n} e(\ell \theta) =0 \end{align}

We refer the reader to Chapter 1 of \cite{MR419394} for the classical theory and to \cite{MR4720212} for a recent treatment covering the primes, almost primes, squarefree numbers and practical numbers at once. A Fubini argument over $\mu\times\mu$ together with dominated convergence shows that if (\ref{sec.equisparse.1}) holds then for every atomless Borel probability measure $\mu$ on $\rzee$ there exists a sequence $(\ell_k)_{k \in \enn}$ of elements of $\mathscr{M}$ such that we have the following for all continuous $\psi:\rzee \to \cee$.
\begin{align} \label{sec.equisparse.2}  \lim_{k \to \infty} \int_{\rzee} \psi(\ell_kx) \deee \mu(x) = \int_{\rzee} \psi(x) \deee \lambda(x) \end{align}

Where (\ref{sec.equisparse.1}) is available, the qualitative result obtained from the classical route is stronger than anything the present article gives for the same set of multipliers. This can be obtained by applying the Fubini argument to the atomless part of an arbitrary $\mu$, thus producing a limit measure which dominates $c\lambda$ where $c$ is the atomless mass of $\mu$. Since atomic measures have upper entropy dimension zero we find $c \geq \ach_\ast^\uparrow(\mu)$.\\

The value of an entropy-based route to qualitative endomorphic rigidity lies with sets for which (\ref{sec.equisparse.1}) is unknown or false, and the four examples below are chosen on that principle. We also emphasize that the effective endomorphic rigidity results illustrated in Section \ref{sec.vino} below appear to be new even for sets where (\ref{sec.equisparse.1}) is known, such as the primes.\\

Our first illustrative theorem offers an entropy-based route to the type of convergence in (\ref{sec.equisparse.2}) which does not rely on individual equidistribution results of the form in (\ref{sec.equisparse.1}). We will give four examples of such sets.

\begin{itemize} \item The set consisting of the first prime after each perfect cube. Here (\ref{sec.equisparse.1}) is not known, since the sequence selects a single prime from each window $(n^3,(n+1)^3)$. In contrast, the available results on $p\theta$ for primes in short intervals, beginning with \cite{MR858400}, control averages over all primes in a window as opposed to individual selections.
\item The range of the partition function, defined as usual to be the map $\mathtt{p}:\enn \to \enn$ given as follows.
\begin{align} \mathtt{p}(n) \coloneqq \Bigl\lvert \Bigl\{ (a_1,\ldots,a_n) \in \enn^n \, : \, a_1 \geq \cdots \geq a_n \mbox{ and } a_1 + \cdots + a_n = n \Bigr\} \Bigr\rvert \label{def.canonsets.c4.1} \end{align}
Here (\ref{sec.equisparse.1}) is again unknown. Hardy and Ramanujan \cite{MR2280879} give $\log \mathtt{p}(n) \sim \pi\sqrt{2n/3}$, so $\mathtt{p}(\enn)$ has the growth of the sequences $\lfloor e^{n^b} \rfloor$ with $b \in (0,1)$, for which Problem 26 of \cite{MR3613710} records that no exponential sum estimates are available.
\item The set of natural numbers of the form $\lfloor n^{\log(n)} \rfloor$. Here (\ref{sec.equisparse.1}) is again unknown. Since $n^{\log(n)} = \exp((\log n)^2)$, this is the sequence which \cite{MR2585398} describes as ``out of reach at the moment''. The estimates of Karatsuba \cite{MR330068} handle $\lfloor n^{(\log n)^a} \rfloor$ only for $a < 1/2$, as recorded in Problem 26 of \cite{MR3613710}, and the recent \cite{FMK25} treats exactly that range.
\item The semigroup $\mathbb{S}_{2,3} \coloneqq \{2^c3^d:c,d \in \enn\}$ at the core of Furstenberg's conjecture. Here (\ref{sec.equisparse.1}) fails along the natural exhaustions by size and by boxes: Furstenberg pointed out at the end of \cite{MR0213508} that it fails for suitable Liouville $\theta$, and \cite{MR4746037} constructs such $\theta$ explicitly and demonstrates that the set of exceptional $\theta$ has positive Hausdorff dimension. See also \cite{MR4741329}, which shows there is a set of $\theta$ with Hausdorff dimension one whose $\times 2 \times 3$ empirical measures, averaged over the boxes $\{2^c3^d : c,d \leq N\}$ rather than by size, fail to converge. \end{itemize}

The four sets above share a qualitative property: the ratio of consecutive elements tends to $1$. Following Boshernitzan \cite{Bosh94} we call an infinite set of natural numbers with this property \emph{sublacunary}; the formal definition is Definition \ref{def.sublac} below. Our first theorem recovers the full convergence (\ref{sec.equisparse.2}) along every sublacunary set for measures of upper entropy dimension $1$, and otherwise provides a lower bound by $\lambda$ in terms of the upper entropy dimension.

\begin{Theorem}[Endomorphic rigidity along sublacunary sets] \label{thm.ad.sublac} Let $\mathscr{M}$ be an infinite set of positive integers whose increasing enumeration $(a_j)_{j \in \enn}$ satisfies $\lim_{j \to \infty} a_{j+1}/a_j = 1$, and let $\mu$ be a Borel probability measure on $\rzee$. Then there exist a sequence $(m_k)_{k \in \enn}$ of elements of $\mathscr{M}$ and a Borel probability measure $\nu$ on $\rzee$ satisfying $\nu \succeq \ach_\ast^\uparrow(\mu) \lambda$ such that the following holds for every continuous $\psi: \rzee \to \cee$.
\begin{equation*} \lim_{k \to \infty} \int_{\rzee} \psi(m_kx) \deee \mu(x) = \int_{\rzee} \psi(x) \deee \nu(x) \end{equation*} \end{Theorem}

 Taking into account metrizability of the weak-star topology (see Theorem 3.16 of \cite{MR1157815}), Theorem \ref{thm.ad.sublac} is the sequential form of Theorem \ref{thm.entdim} in Section \ref{sec.intro.app.special.state} below. Theorem \ref{thm.equisparse} now records the instances of Theorem \ref{thm.ad.sublac} for the four sets above.

\begin{Theorem}[Endomorphic rigidity along sparse sets] \label{thm.equisparse} Let $\mu$ be a Borel probability measure on $\rzee$. \begin{description}
\item[\textsf{(Endomorphic rigidity along cube-succeeding primes)}] \label{thm.equisparse.c1} There exist a sequence $(j_k)_{k \in \enn}$ of natural numbers and a Borel probability measure $\nu$ on $\rzee$ satisfying $\nu \succeq \ach_\ast^\uparrow(\mu) \lambda$ such that if we write $p_k$ for the first prime after $j_k^3$ then the following holds for every continuous $\psi: \rzee \to \cee$.
\begin{equation*} \lim_{k \to \infty} \int_{\rzee} \psi(p_kx) \deee \mu(x) = \int_{\rzee} \psi(x) \deee \nu(x) \end{equation*}
\item[\textsf{(Endomorphic rigidity along partition range)}] \label{thm.equisparse.c2} There exist a sequence $(m_k)_{k \in \enn}$ of natural numbers and a Borel probability measure $\eta$ on $\rzee$ satisfying $\eta \succeq \ach_\ast^\uparrow(\mu) \lambda$ such that the following holds for every continuous $\psi: \rzee \to \cee$.
\begin{equation*} \lim_{k \to \infty} \int_{\rzee} \psi(\mathtt{p}(m_k)x) \deee \mu(x) = \int_{\rzee} \psi(x) \deee \eta(x) \end{equation*}
\item[\textsf{(Endomorphic rigidity along logarithmic powers)}] \label{thm.equisparse.c3} There exist a sequence $(n_k)_{k \in \enn}$ of positive natural numbers and a Borel probability measure $\varpi$ on $\rzee$ satisfying $\varpi \succeq \ach_\ast^\uparrow(\mu) \lambda$ such that the following holds for every continuous $\psi: \rzee \to \cee$.
\begin{equation*} \lim_{k \to \infty} \int_{\rzee} \psi\left(\left\lfloor n_k^{\log(n_k)} \right \rfloor x\right) \deee \mu(x) = \int_{\rzee} \psi(x) \deee \varpi(x) \end{equation*}
\item[\textsf{(Endomorphic rigidity along multiplicative semigroup)}] \label{thm.equisparse.c4} There exist sequences $(c_k)_{k \in \enn}$ and $(d_k)_{k \in \enn}$ of natural numbers and a Borel probability measure $\varsigma$ on $\rzee$ satisfying $\varsigma \succeq \ach_\ast^\uparrow(\mu) \lambda$ such that the following holds for every continuous $\psi: \rzee \to \cee$.
\begin{equation*} \lim_{k \to \infty} \int_{\rzee} \psi(2^{c_k}3^{d_k}x) \deee \mu(x) = \int_{\rzee} \psi(x) \deee \varsigma(x) \end{equation*}
 \end{description} \end{Theorem}
 
The first clause in Theorem \ref{thm.equisparse} follows from Corollary \ref{cor.rjprime}, while the remaining three clauses follow from Clauses \ref{thm.rjcanon.c2}, \ref{thm.rjcanon.c3} and \ref{thm.rjcanon.c4} of Theorem \ref{thm.rjcanon} below.\\
 
 The semigroup statement in the fourth clause in Theorem \ref{thm.equisparse} may be thought of as an `invariance-unnecessary Rudolph-Johnson' theorem. We discuss this perspective further after Corollary \ref{cor.rjclass}, again in Section \ref{sec.intro.app.special.state} below. However this fourth clause is less original to the current article, as it can be obtained with controlled modifications to the proof of Theorem 1.4 in the BLMV article \cite{MR2563089}.\\
 
\subsection{Entropy-effective endomorphic rigidity} \label{sec.vino}

We now present an effective theorem, in two clauses, for sets of great classical interest. Here the individual equidistribution in (\ref{sec.equisparse.1}) is known, but the dominated convergence method provides no rate. Again, we obtain an entropy minorization in general and full endomorphic rigidity to $\lambda$ for measures of upper entropy dimension $1$.

\begin{Theorem}[Entropy-effective endomorphic rigidity]\label{thm.ad.primeff}
Let $\epsilon \in (0,1]$, let $r \in \enn$ and let $\mu$ be a Borel probability measure on $\rzee$ and adopt the notation below.
\begin{equation*} \ell \coloneqq \left\lceil \frac{256}{\epsilon}\log_2\left(\frac{8(r+2)}{\epsilon}\right) \right\rceil \end{equation*}
Then we have the following.
\begin{description}\item[\textsf{(Entropy-effective endomorphic rigidity along primes)}] For every integer $N$ satisfying
\begin{align*} \log_2(N) \geq \frac{9\max(2^{51} + \ell, 6\ell + 15)}{\epsilon} \end{align*}
there exists a prime $p \leq 2N$ such that for all $c_0, \ldots, c_r \in \cee$ we have the following.
\begin{equation*} \int_{\rzee} \lvertt \sum_{j=0}^{r} c_j\, e(pjx) \rvertt^2 \! \dee \mu(x) \;\geq\; \left(\frac{\ach(\mu,N)}{\log(N)} - \epsilon\right)\sum_{j=0}^{r} \lvert c_j \rvert^2 \end{equation*}
\item[\textsf{(Entropy-effective endomorphic rigidity along powers)}] For every $v \in \ennplus$ and every integer $N$ satisfying
\begin{align*} \log_2(N) \geq \frac{9((2v-1)(\ell+v+1) + \ell)}{\epsilon}\end{align*}
there exists $w \in \ennplus$ with $w^v \leq 2N$ such that for all $c_0, \ldots, c_r \in \cee$ we have the following.
\begin{equation*} \int_{\rzee} \lvertt \sum_{j=0}^{r} c_j\, e(w^vjx) \rvertt^2 \! \dee \mu(x) \;\geq\; \left(\frac{\ach(\mu,N)}{\log(N)} - \epsilon\right)\sum_{j=0}^{r} \lvert c_j \rvert^2 \end{equation*} \end{description}
\end{Theorem}

The two clauses of Theorem \ref{thm.ad.primeff} are Clauses \ref{thm.blprime.c1} and \ref{thm.blprime.c2} of Theorem \ref{thm.blprime} below, read with $\sigma \coloneqq \ach(\mu,N)/\log(N)$. In that context the parameter $\ell$ of Theorem \ref{thm.blprime} carries a factor $1 - \sigma + \epsilon \leq 2$, so it is bounded by the $\ell$ displayed above. Moreover, the thresholds appearing in Theorem \ref{thm.ad.primeff} imply those in Theorem \ref{thm.blprime} since the latter increase with $\ell$.

\subsection{Related to Lyons' conjecture}

In \cite{MR941238}, Lyons conjectured the conclusion of Theorem \ref{thm.lyons} below with no entropy hypothesis for every continuous $\mu$ and every nonlacunary semigroup $\mathscr{M}$, which would imply Furstenberg's conjecture. Badea and Grivaux disproved Lyons' conjecture in \cite{MR4082895}, constructing atomless $\mu$ with $\inf_{s \in \mathbb{S}_{2,3}} \lvert \hat{\mu}(s) \rvert > 0$, where as usual for $s \in \zee$ we introduce the frequency-$s$ Fourier coefficient of $\mu$ as follows.
\begin{align} \label{eq.intro.x1} \hat{\mu}(s) \coloneqq \int_{\rzee} e(-sx) \deee \mu(x) \end{align}
It was further shown in \cite{MR4822284} that generic $\times p$-invariant continuous measures carry large Fourier coefficients along $(q^n)_{n \in \enn}$, in contrast to the theorem of Johnson and Rudolph \cite{MR1351328} that for ergodic $\times p$-invariant measures with positive entropy we have convergence to $\lambda$ along a full-density set of multipliers of the form $q^n$. \\

Positive entropy is thus already known to be the dividing line in the invariant setting, and Theorem \ref{thm.lyons} below exhibits entropy dimension as the dividing line with no invariance at all. The obstruction is that every such measure has an entropy-dimension deficit at least as large as its Fourier floor. Thus Lyons' conclusion survives for measures of full entropy dimension, uniformly over all sublacunary multiplier sets, with no multiplicative structure required.

\begin{Theorem}[Lyons' conjecture at full entropy dimension for sublacunary multipliers] \label{thm.lyons} Let $\mathscr{M}$ be an infinite set of positive integers whose increasing enumeration $(a_j)_{j \in \enn}$ satisfies $\lim_{j \to \infty} a_{j+1}/a_j = 1$, and let $\mu$ be a Borel probability measure on $\rzee$. Then we have:
\begin{equation*} \inf_{m \in \mathscr{M}} \lvert \hat{\mu}(m) \rvert \leq 1 - \ach_\ast^\uparrow(\mu) \end{equation*}
In particular, if $\ach_\ast^\uparrow(\mu) = 1$ then $\hat{\mu}(m_k) \to 0$ along some sequence $(m_k)_{k \in \enn}$ in $\mathscr{M}$. \end{Theorem}

We now describe how Theorem \ref{thm.lyons} may be obtained from Theorem \ref{thm.entdim} below. Since $\mathscr{M}$ is sublacunary, that theorem returns $\nu$ in the weak-star trajectory closure of $\mu$ along $\mathscr{M}$ with $\nu \succeq \ach_\ast^\uparrow(\mu)\lambda$. Since $\hat{\nu}(1)$ is a limit of the numbers $\hat{\mu}(m)$ with $m \in \mathscr{M}$ by the definition of the weak-star topology, the infimum in the theorem is at most $\lvert \hat{\nu}(1) \rvert$.\\

\begin{Remark}[Rigidity sequences] \label{rem.rigid} Bergelson, del Junco, Lema\'nczyk and Rosenblatt \cite{BdJLR14} and Eisner and Grivaux \cite{EG11} define a strictly increasing sequence $(n_k)_{k \in \enn}$ of positive integers to be a \emph{rigidity sequence} if some weakly mixing measure-preserving system is rigid along it. (See Chapter 9 of \cite{MR3410920}, in particular Section 9.2 and the concluding remarks on p.~175, for these ergodic-theoretic notions.) These authors prove that this holds if and only if some atomless Borel probability measure $\mu$ on $\rzee$ satisfies $\hat{\mu}(n_k) \to 1$. Eisner and Grivaux construct in \cite{EG11} a rigidity sequence with $n_{k+1}/n_k \to 1$, that is, a sublacunary rigidity sequence.\\

Since the tails of a sublacunary sequence are sublacunary, by applying Theorem \ref{thm.lyons} to them we find that every Borel probability measure $\mu$ on $\rzee$ with $\hat{\mu}(n_k) \to 1$ along a sublacunary sequence satisfies $\ach_\ast^\uparrow(\mu) = 0$. Thus the measures witnessing rigidity along sublacunary sequences have entropy dimension zero. In \cite{BdJLR14} it is asked whether the increasing enumeration of $\mathbb{S}_{2,3}$ is a rigidity sequence; Theorem \ref{thm.lyons} leaves this open, but forces any witnessing measure to have entropy dimension zero. \end{Remark}

The sublacunarity hypothesis in Theorem \ref{thm.lyons} cannot be dropped. For $p \in (0,1/2]$ we write $\nu_p$ for the bias-$p$ ternary Cantor measure on $\rzee$, in the sense of \cite{MR0625600}, and we note that $\nu_p$ is $\times 3$-invariant.

\begin{Theorem}[Rigidity failure for a lacunary semigroup] \label{thm.ad.cantor} For all $p \in (0,1/2]$ we have the inequality below.
\begin{equation*} \inf_{j \in \enn} \lvert \hat{\nu}_p(3^j) \rvert \geq 1 - \frac{8}{\log(1/p)} \ach_\ast^\uparrow(\nu_p) \end{equation*}\end{Theorem}

Since $8/\log(1/p) \to 0$ as $p \to 0$, we find that the conclusion of Theorem \ref{thm.lyons} fails for $\mathscr{M} \coloneqq \{3^j : j \in \enn\}$ by an arbitrarily large factor. Theorem \ref{thm.ad.cantor} is presented for contrast and it does not involve the core method of the article. Its proof is a short computation given in Section \ref{sec.pfapp.cantor} below.\\

In the other direction, Gao and Wen show in \cite{MR4176691} that the Fourier coefficients of the analogous $\times 2$-invariant coin-tossing measures vanish along $(b^n)_{n \in \enn}$ for every $b$ multiplicatively independent of $2$, and along the semigroup generated by $3$ and $5$, using lower bounds for linear forms in logarithms. The fourth clause of Theorem \ref{thm.equisparse} recovers  a weaker statement of the same kind from entropy alone: since $\hat{\nu}_p(3^c h) = \hat{\nu}_p(h)$ for all $c \in \enn$ and $h \in \zee$, that clause returns a sequence of powers of $2$ along which the Fourier coefficients of $\nu_p$ converge in modulus to a limit at most $1 - \ach_\ast^\uparrow(\nu_p)$.

\subsection{Further directions}

The authors are optimistic that theory developed here may be carried across the bridges established in \cite{MR2358379, MR2648695} to obtain results in the vein of \cite{MR2247967}.\\

We also remark for contextual purposes that this article is an outgrowth of longstanding efforts on the part of the authors to resolve new cases of Furstenberg's $\times 2 \times 3$ conjecture. The companion article \cite{EntRates} is in that sense a direct successor to the present one: it takes Clause \ref{thm.blprime.c3} of Theorem \ref{thm.blprime} as its only input and derives rates of decay of resolution entropy for $\times 2 \times 3$-invariant measures whose coordinate function adheres to eigenfunctions.

\section{Discussion of AI use}

The artificial intelligence models Claude Opus 4.8, Opus 5, Fable 5 and Fable 5.1 from Anthropic PBC were extensively used in the composition of this article, as follows.
\begin{description}
	\item[\textsf{Mathematical falsification}] Most importantly, to accelerate exclusion of incorrect ideas on the part of the authors.
	\item[\textsf{Mathematical bookkeeping}] Secondarily, to accelerate the accurate management of constants and parameters in proof design.
	\item[\textsf{Latex typography}] Lastly, to accelerate the identification of typographical errors in the Latex manuscript. \end{description}

\hrulefill
\vspace{12 pt}

The specific formal style in this manuscript is not a marker of AI-produced text. Indeed, quite the opposite is true: the Claude models used cannot reliably reproduce this style. The authors wrote the manuscript line-by-line in order to verify it themselves. Claude was called for additional verification and error-cleaning after every line.\\

The authors are optimistic that this style is a functioning compromise between human readability and AI verifiability. It should be reasonably possible for a graduate student with a course in measure theory to formally check each line to their own satisfaction. It should also be possible for a frontier AI system to verify the manuscript without the need to impute meaning or prove facts independently.

\chapter{Introduction} \label{cha.intro}

\section{Formulation of main theorem} \label{sec.intro.main}

In light of Section \ref{sec.ad} we next choose to formulate our main theorem in the strongest possible way upfront, as Theorem \ref{thm.main} in Section \ref{sec.intro.main.state} below. This will first require us to go through some preliminaries in Section \ref{sec.intro.main.defs}.

\subsection{Preliminaries to statement of main theorem} \label{sec.intro.main.defs}

These preliminaries are categorized into five headings:

\begin{itemize}
	\item Elementary analysis
	\item Measure endomorphs
	\item Resolution entropy theory
	\item The exponential infimum
	\item Progression-shadowing
\end{itemize}

Of these, only the last is original to this article as opposed to a declaration of standard material and the reader may prefer initially to skip much of the first four headings.

\subsubsection{Elementary analysis}

\label{sec.intro.main.defs.elemen}

We adopt the following notational conventions.

\begin{enumerate}[label=\upshape(\thesection:\alph*)]
\item \label{conv.exp} For $x \in \arr$ we adopt the notation $e(x) \coloneqq e^{2 \pi i x}$. For $n \in \zee$, it will be convenient to have a symbol to denote the function $x \mapsto e(nx)$ itself. To this end we define a function $\mathfrak{e}^n: \arr \to \cee$ by $\mathfrak{e}^n(x) \coloneqq e(nx)$. 
\item \label{conv.period} We will identify $\zee$-periodic functions on $\arr$ (such as $\mathfrak{e}^n$) with the corresponding functions on $\rzee$. We will sometimes identify the group $\rzee$ with its coset transversal $[0,1)_{\arr}$.
\item \label{conv.setdil} If $\mathscr{E}$ and $\mathscr{F}$ are subsets of $\zee$ we write $\mathscr{E}\mathscr{F} \coloneqq \{jk : j \in \mathscr{E} \mbox{ and } k \in \mathscr{F}\}$. In particular, for $m \in \zee$ and $\mathscr{F}$ a subset of $\zee$ we have $\{m\}\mathscr{F} = \{mj : j \in \mathscr{F}\}$.
\item \label{conv.zeeinterval} For $x,y \in \arr_{\pm \infty}$ (where $\arr_{\pm \infty} \coloneqq \arr \cup \{-\infty,+\infty\}$) it will be convenient to adopt the notation $[x,y]_{\zee} \coloneqq \{n \in \zee: x \leq n \leq y \}$. For distinctness we will decorate the standard notations denoting real intervals with a subscript $\mathbb{R}$. Thus, for example, the closed unit interval in the real line will be denoted $[0,1]_{\arr}$ and the strictly positive numbers will be denoted $(0,\infty)_\arr$.
\item \label{conv.prob} We adopt the notation $\operatorname{FMea}(\rzee)$ to denote the set of Borel measures $\mu$ on $\rzee$ such that $\mu(\rzee) < \infty$, and we adopt the notation $\prob(\rzee)$ to denote the set of all Borel probability measures on $\rzee$. Thus for $\mu \in \operatorname{FMea}(\rzee)$ we have $\mu \in \prob(\rzee)$ if and only if $\mu(\rzee) = 1$. Also, for a nonempty finite set $S$ we write $\prob(S)$ for the set of functions from $S$ to $[0,1]_\arr$ summing to $1$.
\item \label{conv.integral} If $\mu \in \operatorname{FMea}(\rzee)$ and $\psi:\rzee \to \cee$ is a $\mu$-integrable function it will sometimes be convenient to adopt the notation $\mu(\psi)$ to denote the complex number corresponding to the integral of $\psi$ over $\mu$.
\item \label{conv.bbor} We adopt the notation $\bbr(\rzee)$ for the space of bounded Borel measurable functions from $\rzee$ to $\cee$, equipped with the supremum norm. Every $\mu \in \operatorname{FMea}(\rzee)$ evaluates on $\bbr(\rzee)$ by way of \ref{conv.integral}.
\item \label{conv.indicator} For a subset $E$ of $\rzee$ we write $\mathbf{1}_E: \rzee \to \{0,1\}$ for its indicator function. Assuming $E$ is Borel, we have that $\mathbf{1}_E \in \bbr(\rzee)$.
\item \label{conv.continuous} We adopt the notation $C(\rzee)$ for the Banach space of continuous functions from $\rzee$ to $\cee$, equipped with the supremum norm.
\item \label{conv.dual} We write $C(\rzee)^\dagger$ to denote the space of continuous linear functionals on $C(\rzee)$, and refer to the elements of $C(\rzee)^\dagger$ as \emph{finite complex measures}. It follows from Theorem 7.17 in \cite{MR1681462} that $\operatorname{FMea}(\rzee)$ may be identified with the set of elements $\zeta$ of $C(\rzee)^\dagger$ satisfying $\zeta(\psi) \geq 0$ for all nonnegative $\psi \in C(\rzee)$, and we will use this identification without further comment.
\item \label{conv.submeasure} If $\mu,\nu \in \operatorname{FMea}(\rzee)$, we define \emph{$\nu$ to be a submeasure of $\mu$} if $\nu(E) \leq \mu(E)$ for all Borel subsets $E$ of $\rzee$. We denote this in symbols by $\nu \preceq \mu$ or $\mu \succeq \nu$.
\item \label{conv.Tmap} For $n \in \ennplus$ we write $\mathtt{T}_n$ to denote the map from $\rzee$ to $\rzee$ given by $\mathtt{T}_n(x) \coloneqq nx$.
\item \label{conv.leb} We write $\lambda$ to denote the Haar probability measure on $\rzee$, which may be identified with the elementary Lebesgue measure on the coset transversal $[0,1)_\arr$.
\end{enumerate}

\subsubsection{Measure endomorphs}

The following definition is the core connection between the algebraic endomorphy of the circle and harmonic analysis.

\begin{Definition}[Endomorphs of complex measures on $\rzee$] \label{def.endomorph} For $n \in \ennplus$ and $\zeta \in C(\rzee)^\dagger$, we define the \emph{$\times n$-endomorph of $\zeta$} to be the element $n_\circledcirc \zeta$ of $C(\rzee)^\dagger$ given by stipulating that $[n_\circledcirc\zeta](\psi) \coloneqq \zeta(\psi \circ \mathtt{T}_n)$ for all $\psi \in C(\rzee)$.\\

We define $\zeta$ to be \emph{$\times n$-invariant} if $n_\circledcirc \zeta = \zeta$ and for $a,b \in \ennplus$ we define $\zeta$ to be \emph{$\times a \times b$-invariant} if it is both $\times a$-invariant and $\times b$-invariant. \end{Definition}

We immediately find that this endomorphy is in fact an action of $\ennplus$, as the following remark records.

\begin{Remark}[Composition of endomorphs] \label{rem.endocomp}
Let $m,n \in \ennplus$ and let $\zeta \in C(\rzee)^\dagger$. Directly from \ref{conv.Tmap} we have $\mathtt{T}_m \circ \mathtt{T}_n = \mathtt{T}_{mn}$, and so unwrapping Definition \ref{def.endomorph} twice provides the following for all $\psi \in C(\rzee)$.
\begin{align*} \bigl[ m_{\circledcirc} ( n_{\circledcirc} \zeta ) \bigr](\psi) = [n_{\circledcirc}\zeta](\psi \circ \mathtt{T}_m) = \zeta(\psi \circ \mathtt{T}_m \circ \mathtt{T}_n) = \zeta(\psi \circ \mathtt{T}_{mn}) = [(mn)_{\circledcirc}\zeta](\psi) \end{align*}
Thus $m_{\circledcirc}(n_{\circledcirc}\zeta) = (mn)_{\circledcirc}\zeta$.\\

In particular, if $\zeta$ is $\times a \times b$-invariant for some $a,b \in \ennplus$, then induction provides $(a^j b^k)_{\circledcirc}\, \zeta = \zeta$ for all $j,k \in \enn$. \end{Remark}

\subsubsection{Resolution entropy theory} \label{sec.intro.main.defs.ent}

The following definition is the heart of effectivizing Rudolph-Johnson type results, transferring the asymptotic regime of Kolmogorov-Sinai entropy to a quantitative, finitary regime.

\begin{Definition}[Resolution entropy] \label{def.entropy}
Let $n \in \ennplus$.
\begin{itemize} \item For $\ell \in [0,n-1]_\zee$, we define the \emph{index-$\ell$ resolution-$n$ cell in $\rzee$} to be the subset $\mathscr{P}_n(\ell)$ of $\rzee$ corresponding to the subset $[\ell/n,(\ell+1)/n)_{\arr}$ of the coset transversal $[0,1)_\arr$. \end{itemize} 
Now, also let $\mu \in \prob(\rzee)$.
\begin{itemize} \item We define the \emph{resolution-$n$ entropy of $\mu$} as follows.
\begin{align}\label{def.entropy.0}\ach(\mu , n) \coloneqq - \sum_{\ell=0}^{n-1} \mu\bigl(\mathscr{P}_n(\ell)\bigr) \log\Bigl(\mu\bigl(\mathscr{P}_n(\ell)\bigr)\Bigr)\end{align} 
Here, we understand the conventions that $\log$ denotes the natural logarithm and that $0 \log 0 \coloneqq 0$. 
\item We also define the \emph{normalized resolution-$n$ entropy of $\mu$} as follows.
\begin{align} \achn(\mu , n) \coloneqq \begin{dcases}
\hspace{0.5 cm} 1 & \mbox{ if }n=1
\\ \frac{\ach(\mu , n)}{\log(n)} & \mbox{ else }\end{dcases} \label{def.entropy-1} \end{align}
\item Finally, we define the \emph{upper entropy dimension of $\mu$} by $\ach_\ast^{\uparrow}(\mu) \coloneqq \limsup_{n \to \infty} \achn(\mu,n)$. \end{itemize} \end{Definition}

We record the following standard result which follows, for example, from Lemma 2.1.1 and Theorem 2.6.4 in \cite{MR2239987}.

\begin{oLemma}[Range of resolution entropy] \label{olem.entrange} For all $n \in \ennplus$ and all $\mu \in \prob(\rzee)$ we have $\ach(\mu,n) \in [0,\log(n)]_{\arr}$. \end{oLemma}

The remark below emphasizes the status of resolution entropy as an intensive quantity.

\begin{Remark}[Normalization of resolution entropy] \label{rem.achnunit} Let $\mu \in \prob(\rzee)$. Combining Lemma \ref{olem.entrange} with the definition of the normalized resolution-$n$ entropy in Definition \ref{def.entropy} we find that $\achn(\mu , n)$ is an element of the unit interval $[0,1]_{\arr}$ for all $n \in \ennplus$. Thus it follows that $\ach_\ast^{\uparrow}(\mu)$ is also an element of $[0,1]_{\arr}$.\end{Remark}

\subsubsection{The exponential infimum} \label{sec.intro.main.defs.einf}

The definition below is typically the output of our machinery: it provides a measurement of measure-uniformity that is far stronger than the resolution entropy as in Definition \ref{def.entropy}.

\begin{Definition}[Exponential infimum] \label{def.foustat} Let $\mu \in \prob(\rzee)$ and let $\mathscr{F}$ be a nonempty finite subset of $\zee$. We define the \emph{bandprofile-$\mathscr{F}$ exponential infimum of $\mu$} as follows.
\begin{align*} \einf(\mu , \mathscr{F}) \coloneqq \inf \left \{ \mu \left ( \left \vert \sum_{j \in \mathscr{F}} \alpha(j) \,\mathfrak{e}^j \right \vert^2  \right): \alpha \in \cee^\mathscr{F} \mbox{ and } \sum_{j \in \mathscr{F}} |\alpha(j)|^2 = 1  \right \} \end{align*} \end{Definition}

We collect certain straightforward facts about Definition \ref{def.foustat}.

\begin{Remark}[Normalization of exponential infimum] \label{rem.unitint} It follows directly from Definition \ref{def.foustat} that for all $\mu \in \prob(\rzee)$ and all nonempty finite subsets $\mathscr{F}$ of $\zee$ the quantity $ \einf(\mu , \mathscr{F})$ is an element of the unit interval $[0,1]_{\arr}$. \end{Remark}

\begin{Remark}[Exponential infimum and singleton sets] \label{rem.einfsing}  It follows directly from the definition of the exponential infimum in Definition \ref{def.foustat} that we have $\einf(\mu,\{n\})=1$ for all $\mu \in \prob(\rzee)$ and all $n \in \zee$.\end{Remark}

\begin{Remark}[Monotonicity of exponential infimum] \label{rem.statmono}
It follows directly from the definition of the exponential infimum in Definition \ref{def.foustat} that for all $\mu \in \prob(\rzee)$ and all nonempty finite subsets $\mathscr{E},\mathscr{F}$ of $\zee$ with $\mathscr{E} \subseteq \mathscr{F}$ we have $\einf(\mu,\mathscr{F}) \leq \einf(\mu,\mathscr{E})$. \end{Remark}

\subsubsection{Progression-shadowing}

The quantitative definition below, together with the depth function of Definition \ref{def.shadep} built on it, is the article's vehicle for carrying Rudolph-Johnson type theory beyond semigroups, and the authors believe it to be new. Its qualitative content is classical: Lemma \ref{lem.ratio} below shows that a set has divergent depth exactly when it is sublacunary.

\begin{Definition}[Progression-shadowing sets: quantitative] \label{def.progshad} 
Let $s,\ell \in \enn$ and let $\mathscr{M}$ be a subset of $\enn$. We say $\mathscr{M}$ is \emph{$(s,\ell)$ progression-shadowing} if for every $m \in [1,2^\ell]_{\zee}$ there exists $q \in \mathscr{M}$ such that:
\begin{align} \vert q - m 2^s \vert \leq 2^{s-\ell} \label{def.progshad.1} \end{align} \end{Definition}

\begin{Remark}[Elementary properties of progression-shadowing] \label{rem.progshad}
Let $s,\ell \in \enn$. The following assertions follow directly from Definition \ref{def.progshad}.
\begin{enumerate}[label=\upshape(\arabic*)]
\item \label{rem.progshad.c1} If $\mathscr{M}$ and $\mathscr{N}$ are subsets of $\enn$ with $\mathscr{M} \subseteq \mathscr{N}$, and if $\mathscr{M}$ is $(s,\ell)$ progression-shadowing, then $\mathscr{N}$ is $(s,\ell)$ progression-shadowing.
\item \label{rem.progshad.c2} If $\mathscr{M}$ is a subset of $\enn$ which is $(s,\ell)$ progression-shadowing, then $\mathscr{M}$ is $(s,\ell_0)$ progression-shadowing for all $\ell_0 \in [0,\ell]_{\zee}$.
\item \label{rem.progshad.c3} Taking $q \coloneqq m2^s$ for each $m \in [1,2^\ell]_{\zee}$ shows that $\enn$ is $(s,\ell)$ progression-shadowing.
\end{enumerate} \end{Remark}

It will be essential for our effective results to have a function quantifying the progression-shadowing properties of a given set, which we introduce below.

\begin{Definition}[Progression-shadowing depth] \label{def.shadep}
Let $\mathscr{M}$ be a subset of $\enn$ and let $s \in \enn$. We define the \emph{scale-$s$ progression-shadowing depth of $\mathscr{M}$} to be the element $\mathtt{Dep}_{\mathscr{M}}(s)$ of $\enn \cup \{\infty\}$ given as follows.
\begin{align} \mathtt{Dep}_{\mathscr{M}}(s) \coloneqq \sup \bigl\{ \ell \in \enn : \mathscr{M} \mbox{ is } (s,\ell) \mbox{ progression-shadowing} \bigr\} \label{def.shadep.1} \end{align}
Here we understand a convention that the supremum is taken in $\enn \cup \{\infty\}$, so that in particular it returns $0$ when no $\ell$ satisfies the stated condition, and $\infty$ when every $\ell$ does.\end{Definition}

Naturally Definitions \ref{def.progshad} and \ref{def.shadep} are synergetic, as the following remark clarifies.

\begin{Remark}[Shadowing depth and progression-shadowing] \label{rem.shadep}
Let $\mathscr{M}$ be a subset of $\enn$ and let $s,\ell \in \enn$.
\begin{enumerate}[label=\upshape(\arabic*)]
\item \label{rem.shadep.c1} Clause \ref{rem.progshad.c2} of Remark \ref{rem.progshad} asserts that the supremum in (\ref{def.shadep.1}) is taken over an initial segment of $\enn$, and consequently for all $\ell \in \ennplus$ we have $\mathtt{Dep}_{\mathscr{M}}(s) \geq \ell$ if and only if $\mathscr{M}$ is $(s,\ell)$ progression-shadowing.
\item \label{rem.shadep.c2} Combining Clause \ref{rem.shadep.c1} of the current remark with Clause \ref{rem.progshad.c1} of Remark \ref{rem.progshad} we find that $\mathtt{Dep}_{\mathscr{M}}(s) \leq \mathtt{Dep}_{\mathscr{N}}(s)$ for all subsets $\mathscr{N}$ of $\enn$ with $\mathscr{M} \subseteq \mathscr{N}$.
\item \label{rem.shadep.c3} Combining Clause \ref{rem.shadep.c1} of the current remark with Clause \ref{rem.progshad.c3} of Remark \ref{rem.progshad} we find that $\mathtt{Dep}_{\enn}(s) = \infty$.
\end{enumerate} \end{Remark}

We have chosen simplicity over generality by using the fixed base $2$ in Definition \ref{def.progshad}. The remark below discharges this.

\begin{Remark}[Independence of the base] \label{rem.progshadbase}
The arguments in the present article use no property of the base $2$ appearing in Definition \ref{def.progshad} beyond the fact that it is an integer greater than $1$. Thus for any $b \in [2,\infty)_{\zee}$ one may define $\mathscr{M}$ to be \emph{$(b,s,\ell)$ progression-shadowing} if for every $m \in [1,b^\ell]_{\zee}$ there exists $q \in \mathscr{M}$ with $\lvert q - mb^s \rvert \leq b^{s-\ell}$, define the associated depth function $\mathtt{Dep}^{(b)}_{\mathscr{M}}$ exactly as in Definition \ref{def.shadep}, and replace the base $2$ by $b$ throughout. All results below remain valid, at the cost of allowing the constants involved to depend on $b$.\\

We fix the base $2$ purely for notational simplicity in an article which is already heavy with notation, and we emphasize that this base is unrelated to the multiplicatively independent pair of integers appearing in the classical Rudolph-Johnson theorem.\\

We also note that the depth functions $\mathtt{Dep}^{(b)}_{\mathscr{M}}$ arising from distinct bases need not be comparable, since Definition \ref{def.progshad} constrains $\mathscr{M}$ near specified points rather than throughout an interval, and so the $2 \mapsto b$ generalization does contain new cases. In every application appearing in the present article, however, the hypothesis is verified by way of an upper bound for the gaps of $\mathscr{M}$ as in Lemma \ref{lem.gapshad}, a criterion which delivers every base at once.\\

It may even be possible to extend the article's methods to replace the sequence $(b^s)_{s \in \enn}$ with a sequence of integers $(b_s)_{s \in \enn}$ such that $b_s$ divides $b_{s+1}$.\end{Remark}

\subsection{Statement of main theorem} \label{sec.intro.main.state}

\subsubsection{Statement of main theorem: quantitative}

We now present our main theorem, which will be proved in Chapter \ref{cha.pfmain} below.

\begin{Theorem}[Quantitative endomorphic rigidity via entropy] \label{thm.main} \quad
\vspace{-1 \baselineskip}
\begin{enumerate}[label=\upshape(\arabic*)] \item \label{thm.main.c1} For all:
\begin{enumerate}[label=\upshape(\arabic*)] \item \label{thm.main.c1.1} $N \in (2^{100},\infty)_{\zee}$, for which we adopt the notation $n \coloneqq \lfloor \log_2(N) \rfloor$,
\item \label{thm.main.c1.2} and all $\sigma \in (0,1]_{\arr}$ and all $\beta \in (0,1]_{\arr}$ and all $\gamma \in (0,1/4]_{\arr}$ such that $\sigma > \beta + 2\gamma$,
\item \label{thm.main.c1.3} and all $\ell \in [2/\gamma, \beta n]_{\zee}$,
\item \label{thm.main.c1.4} and all $\mu \in \prob(\rzee)$ such that $\achn(\mu,N) \geq \sigma$,
\item \label{thm.main.c1.5} and all subsets $\mathscr{M}$ of $\enn$ such that $\mathtt{Dep}_{\mathscr{M}}(s) \geq \ell$ for all $s \in [\beta n - \ell, n - \ell]_{\zee}$,
\end{enumerate}
\item \label{thm.main.c2} there exists a subset $\mathscr{G}$ of $\mathscr{M} \cap [1,2^{n+1}]_{\zee}$ satisfying $\lvert \mathscr{G} \rvert \geq 2^{\ell}/5$,
\end{enumerate}
such that for all $m \in \mathscr{G}$ and all nonempty finite subsets $\mathscr{F}$ of $\enn$ we have the inequality below.
\begin{align} \label{thm.main.conc1}  \einf\bigl(\mu,\{m\}\mathscr{F} \bigr) \geq \frac{1}{1-2\gamma}\left(\sigma - \beta - 2\gamma - \frac{2}{n}\right) - 160\left( \frac{\max(\mathscr{F})^{7/2}}{2^{\gamma \ell/2}} + \frac{\max(\mathscr{F})^{4}}{2^{\ell}} \right) \end{align} \end{Theorem}

We note that the hypothesis on $\mathtt{Dep}_{\mathscr{M}}$ is not exotic: Lemma \ref{lem.depbound} below records lower bounds on the progression-shadowing depth of several interesting sets, and Section \ref{sec.intro.app.ramp} converts each of these bounds into a specialization of the present theorem.

\subsubsection{Statement of main theorem: existential}

\label{sec.intro.app.special.exist}

The next Theorem \ref{thm.exist} presents the appropriate regime to dissolve the subtracted factor at the right of (\ref{thm.main.conc1}) into an $\epsilon$. \\

 Theorem \ref{thm.exist} will be derived from Theorem \ref{thm.main} in Section \ref{sec.pfapp.special.exist} below. 

\begin{Theorem}[Existential endomorphic rigidity via entropy] \label{thm.exist} \quad
\vspace{-1 \baselineskip}
\begin{enumerate}[label=\upshape(\arabic*)] \item \label{thm.exist.c1} For all:
\begin{enumerate}[label=\upshape(\arabic*)] \item \label{thm.exist.c1.1} $\epsilon \in (0,1]_{\arr}$ and all $r \in \enn$ and all $\sigma \in (0,1]_{\arr}$,
\item \label{thm.exist.c1.2} and all $\ell \in \enn$ such that $\ell \geq (128(1-\sigma+\epsilon)/\epsilon)\log_2\bigl(8(r+2)/\epsilon\bigr)$,
\item \label{thm.exist.c1.3} and all subsets $\mathscr{M}$ of $\enn$,
\item \label{thm.exist.c1.4} and all $t \in \enn$ such that $\mathtt{Dep}_{\mathscr{M}}(s) \geq \ell$ for all $s \in [t,\infty)_{\zee}$,
\item \label{thm.exist.c1.5} and all $N \in \enn$ such that $N \geq 2^{9(t+\ell)/\epsilon}$,
\item \label{thm.exist.c1.6} and all $\mu \in \prob(\rzee)$ such that $\achn(\mu,N) \geq \sigma$,
\end{enumerate}
\item \label{thm.exist.c2} there exists $m \in \mathscr{M} \cap [1,2N]_{\zee}$, \end{enumerate}
such that the inequality below holds.
\begin{align} \einf\bigl(\mu,\{m\}[0,r]_{\zee}\bigr) \geq \sigma - \epsilon \label{thm.exist.conc1} \end{align} \end{Theorem}

We pause to record the simplest instance of Theorem \ref{thm.exist}. Clause \ref{rem.shadep.c3} of Remark \ref{rem.shadep} provides $\mathtt{Dep}_{\enn}(s) = \infty$ for every $s \in \enn$, and so Theorem \ref{thm.exist} may be instantiated with $\mathscr{M} \coloneqq \enn$, with $t \coloneqq 0$, and with any resolution $N$ at least $2^{9\ell/\epsilon}$.\\

The threshold on $\ell$ in Clause \ref{thm.exist.c1.2} is largest when $\sigma$ is smallest, and since $\epsilon \leq 1$ it never exceeds $(256/\epsilon)\log_2\bigl(8(r+2)/\epsilon\bigr)$. Taking $\ell$ to be the least integer above that quantity and $N$ to be the least admissible resolution, the resulting statement carries no arithmetic hypothesis on the multiplier whatsoever. Thus both $\ell$ and $N$ depend only on $\epsilon$ and $r$; in particular they are independent of $\mu$.\\

Instantiating with $\sigma \coloneqq \achn(\mu,N)$ when that quantity is positive, and noting that the conclusion is immediate from Definition \ref{def.foustat} when it vanishes, we find that every element $\mu$ of $\prob(\rzee)$ admits a positive integer $m$ bounded in terms of $\epsilon$ and $r$ alone, such that $\einf(\mu,\{m\}[0,r]_{\zee})$ underestimates the entropy of $\mu$ at a resolution fixed in advance by at most $\epsilon$.\\

Beyond the free-multiplier instance just described, the hypothesis of Theorem \ref{thm.exist} is supplied with explicit thresholds by the canonical multiplier sets of Lemma \ref{lem.depbound} below, and the resulting specializations are recorded in Section \ref{sec.intro.app.ramp}.

\section{Formulation of application theorems} \label{sec.intro.app}

\subsection{Further preliminary material on progression-shadowing}

\subsubsection{Additional definitions} \label{sec.intro.app.special.defs}

The limiting theorems of Section \ref{sec.intro.app.special.state} below are stated for the sets described by the following classical condition. The terminology follows Boshernitzan \cite{Bosh94}, and is chosen by contrast with the lacunary sets, for which the ratios in question stay bounded away from $1$.

\begin{Definition}[Sublacunary sets] \label{def.sublac}
Let $\mathscr{M}$ be an infinite subset of $\enn$, and write $(m_j)_{j \in \enn}$ for the increasing enumeration of $\mathscr{M} \setminus \{0\}$. We define $\mathscr{M}$ to be \emph{sublacunary} if the following holds.
\begin{align} \lim_{j \to \infty} \frac{m_{j+1}}{m_j} = 1 \label{def.sublac.1} \end{align} \end{Definition}

Lemma \ref{lem.ratio} below identifies the sublacunary sets with the sets whose progression-shadowing depth diverges. Thus the hypothesis of the limiting theorems is the condition (\ref{def.sublac.1}), while Definition \ref{def.progshad} and the depth function of Definition \ref{def.shadep} remain the bookkeeping which makes the effective theorems possible. Lemma \ref{lem.ratio} will be proved in Section \ref{sec.pfapp.ratio} below.

\begin{Lemma}[Sublacunary sets and progression-shadowing depth] \label{lem.ratio}
For all infinite subsets $\mathscr{M}$ of $\enn$, the following assertions are equivalent.
\begin{enumerate}[label=\upshape(\arabic*)]
\item \label{lem.ratio.c1} $\mathscr{M}$ is sublacunary.
\item \label{lem.ratio.c2} $\lim_{s \to \infty} \mathtt{Dep}_{\mathscr{M}}(s) = \infty$.
\item \label{lem.ratio.c3} For every $\ell \in \ennplus$ there exists $s_0 \in \enn$ such that $\mathscr{M}$ is $(s,\ell)$ progression-shadowing for all $s \in [s_0,\infty)_\zee$.
\end{enumerate} \end{Lemma}

\subsubsection{Bounds on progression-shadowing depth for canonical sets} \label{sec.intro.app.exhibit}

We now recall the basic criterion on semigroup generators from Johnson's work in \cite{MR1194793}.

\begin{Definition}[Multiplicative independence] \label{def.multindep}
We define $a,b \in [2,\infty)_\zee$ to be \emph{multiplicatively independent} if the only $j,k \in \zee$ satisfying $a^j b^k = 1$ are $j = k = 0$. \end{Definition}

The sets defined below are those for which we record bounds on the progression-shadowing depth of Definition \ref{def.shadep}.

\begin{Definition}[Canonical multiplier sets] \label{def.canonsets} \quad \vspace{-1 \baselineskip}
\begin{enumerate}[label=\upshape(\arabic*)]
\item \label{def.canonsets.c1} We write $\mathbb{P}$ to denote the set of prime numbers.
\item \label{def.canonsets.c2} Let $v \in \ennplus$. We define the \emph{set of $v$-th powers} to be the subset $\mathbb{W}_v$ of $\enn$ given as follows.
\begin{align} \mathbb{W}_v \coloneqq \bigl\{ n^v : n \in \enn \bigr\} \label{def.canonsets.c2.1} \end{align}
\item \label{def.canonsets.c3} We define the \emph{set of logarithmic powers} to be the subset $\mathbb{L}$ of $\enn$ given as follows.
\begin{align} \mathbb{L} \coloneqq \bigl\{ \lfloor n^{\log(n)} \rfloor : n \in \ennplus \bigr\} \label{def.canonsets.c3.1} \end{align}
\item \label{def.canonsets.c4} Recall the partition function $\mathtt{p}:\enn \to \enn$ from (\ref{def.canonsets.c4.1}). Write $\mathtt{p}(\enn)$ to denote the \emph{set of partition numbers}, that is the subset $\bigl\{ \mathtt{p}(n) : n \in \enn \bigr\}$ of $\ennplus$.
\item \label{def.canonsets.c5} Let $a,b \in [2,\infty)_\zee$. We define the \emph{multiplicative semigroup generated by $a$ and $b$} to be the subset $\mathbb{S}_{a,b}$ of $\ennplus$ given as follows.
\begin{align} \mathbb{S}_{a,b} \coloneqq \bigl\{ a^j b^k : j,k \in \enn \bigr\} \label{def.canonsets.c5.1} \end{align}
\end{enumerate} \end{Definition}

We emphasize that Definition \ref{def.shadep} defines $\mathtt{Dep}_{\mathscr{M}}(s)$ as a supremum, and that every statement we make about a canonical set is a lower bound on this quantity rather than an evaluation of it. For the primes, an exact evaluation would require a level of control on prime gaps in dyadic windows which seems far beyond the current state of the art.\\

We now provide depth bounds for three of the canonical sets of Definition \ref{def.canonsets}. We choose to involve the logarithmic powers and the partition numbers only in the qualitative Theorem \ref{thm.rjcanon} below, where their sublacunarity is verified directly, and hence no depth bound is recorded for them.

\begin{Lemma}[Bounds on progression-shadowing depth for canonical sets] \label{lem.depbound} \quad \vspace{-1 \baselineskip}
\begin{enumerate}[label=\upshape(\arabic*)]
\item \label{lem.depbound.c1} Let $n_0 \in \ennplus$ and let $\mathscr{M}$ be a subset of $\enn$ such that $\mathscr{M} \cap [n^3+1,(n+1)^3-1]_{\zee} \neq \emptyset$ for all $n \in [n_0,\infty)_{\zee}$. Then for all $s \in [3\log_2(n_0) \vee 12,\infty)_\zee$ we have:
\begin{align} \mathtt{Dep}_{\mathscr{M}}(s) \geq \frac{s}{5} - 3 \label{lem.depbound.c1.1} \end{align}
\item \label{lem.depbound.c2} Let $v \in \ennplus$. For all $s \in \enn$ we have:
\begin{align} \mathtt{Dep}_{\mathbb{W}_v}(s) \geq \frac{s}{2v-1} - v - 1 \label{lem.depbound.c2.1} \end{align}
\item \label{lem.depbound.c5} Let $a,b \in [2,\infty)_\zee$ be multiplicatively independent. Then there exist $\mathtt{c}_{a,b} \in (0,\infty)_\arr$ and $\mathtt{s}_{a,b} \in \enn$, both depending only on $a$ and $b$, such that the following holds for all $s \in [\mathtt{s}_{a,b},\infty)_\zee$.
\begin{align} \mathtt{Dep}_{\mathbb{S}_{a,b}}(s) \geq \mathtt{c}_{a,b} \log(s) \label{lem.depbound.c5.1} \end{align}
\end{enumerate} \end{Lemma}

\begin{Lemma}[Progression-shadowing depth of the primes] \label{lem.primedep}
For all $s \in [2^{51},\infty)_\zee$ we have:
\begin{align} \mathtt{Dep}_{\mathbb{P}}(s) \geq \frac{s}{5} - 3 \label{lem.primedep.1} \end{align} \end{Lemma}

The constants $\mathtt{c}_{a,b}$ and $\mathtt{s}_{a,b}$ of Clause \ref{lem.depbound.c5} are produced explicitly in the proof from the effective irrationality exponent $\kappa_{a,b}$, which we obtain from the deep work of A.~Baker and G.~W\"ustholz in \cite{MR1234835} cited as Lemma \ref{olem.bw} below. Both are therefore effectively computable from that reference, and they are the only constants appearing in the present article which we do not evaluate numerically.  These constants enter only the effective semigroup clause of Theorem \ref{thm.blprime}; the qualitative semigroup clause of Theorem \ref{thm.rjcanon} uses no transcendence theory, resting instead on Furstenberg's lemma, cited as Lemma \ref{olem.furst} below.\\

The bounds we obtain for $\mathbb{W}_v$ and for $\mathbb{P}$ are worth comparing. The bound for $\mathbb{W}_v$ is elementary, resting on nothing beyond the estimate $(n+1)^v - n^v \leq v(n+1)^{v-1}$, whereas the bound $s/5 - 3$ for $\mathbb{P}$ rests on our citations to work of A.W. Dudek in \cite{MR3584567} on the existence of primes in short intervals.\\

Clause \ref{lem.depbound.c2} of Lemma \ref{lem.depbound} provides that $\mathtt{Dep}_{\mathbb{W}_v}(s)$ is at least $s/(2v-1) - v - 1$, while Lemma \ref{lem.primedep} provides that $\mathtt{Dep}_{\mathbb{P}}(s)$ is at least $s/5 - 3$ once $s$ is at least $2^{51}$. Thus the obvious expression $s/3 - 3$ for the squares is superior to the bound we obtain for primes, which is matched to leading order by the bound $s/5 - 4$ which Clause \ref{lem.depbound.c2} returns for the cubes.\\

We also note that our method can treat the range of any integer polynomial with positive leading coefficient as we treat $\mathbb{W}_v$, where $v$ is taken to be the degree.

\subsection{Fixed-resolution endomorphic rigidity for canonical multiplier sets} \label{sec.intro.app.ramp}

We now specialize Theorem \ref{thm.exist} to three of our canonical sets. Theorem \ref{thm.blprime} will be derived from Theorem \ref{thm.exist} in Section \ref{sec.pfapp.special.ramp} below.

\begin{Theorem}[Rigidity after a canonical endomorph] \label{thm.blprime}
Let $\epsilon \in (0,1]_\arr$, let $r \in \enn$, let $\sigma \in (0,1]_\arr$ and let $\mu \in \prob(\rzee)$, and adopt the notation below.
\begin{align} \ell \coloneqq \left\lceil \frac{128(1-\sigma+\epsilon)}{\epsilon}\log_2\left(\frac{8(r+2)}{\epsilon}\right) \right\rceil \label{thm.blprime.h1} \end{align}
\begin{enumerate}[label=\upshape(\arabic*)]
\item \label{thm.blprime.c1} Define $t \coloneqq 2^{51} \vee 5(\ell+3)$. Then for all $N \in [2^{9(t+\ell)/\epsilon},\infty)_\zee$ such that $\achn(\mu,N) \geq \sigma$ there exists $m \in \mathbb{P} \cap [1,2N]_\zee$ such that $\einf\bigl(m_\circledcirc\mu,[0,r]_\zee\bigr) \geq \sigma - \epsilon$.
\item \label{thm.blprime.c2} Let $v \in \ennplus$ and define $t \coloneqq (2v-1)(\ell+v+1)$. Then for all $N \in [2^{9(t+\ell)/\epsilon},\infty)_\zee$ such that $\achn(\mu,N) \geq \sigma$ there exists $m \in \mathbb{W}_v \cap [1,2N]_\zee$ such that $\einf\bigl(m_\circledcirc\mu,[0,r]_\zee\bigr) \geq \sigma - \epsilon$.
\item \label{thm.blprime.c3} Let $a,b \in [2,\infty)_\zee$ be multiplicatively independent, let $\mathtt{c}_{a,b}$ and $\mathtt{s}_{a,b}$ be as provided by Clause \ref{lem.depbound.c5} of Lemma \ref{lem.depbound}, and define $t \coloneqq \mathtt{s}_{a,b} \vee \bigl\lceil e^{\ell/\mathtt{c}_{a,b}} \bigr\rceil$. Then for all $N \in [2^{9(t+\ell)/\epsilon},\infty)_\zee$ such that $\achn(\mu,N) \geq \sigma$ there exists $m \in \mathbb{S}_{a,b} \cap [1,2N]_\zee$ such that $\einf\bigl(m_\circledcirc\mu,[0,r]_\zee\bigr) \geq \sigma - \epsilon$.
\end{enumerate} \end{Theorem}

We also note that if $\mu$ is $\times a \times b$-invariant, then iterating the invariance gives $m_\circledcirc \mu = \mu$ for every $m \in \mathbb{S}_{a,b}$, so that the conclusion of Clause \ref{thm.blprime.c3} becomes $\einf(\mu,[0,r]_\zee) \geq \achn(\mu,N) - \epsilon$. This is an effective Rudolph-Johnson statement at a fixed resolution. The subsequent Section \ref{sec.intro.app.special} which follows dissolves the effectivity in this statement to arrive at Rudolph-Johnson statements of classical type.

\subsection{Rudolph-Johnson type theorems} \label{sec.intro.app.special}

\subsubsection{Rudolph-Johnson type theorems: preliminary definitions} \label{sec.intro.app.special.prelim}

\begin{enumerate}[label=\upshape(\thesection:\alph*)] \item \label{conv.weakstar} We recall that the \emph{weak-star topology} on $C(\rzee)^\dagger$ is defined to be the coarsest topology for which the map $\zeta \mapsto \zeta(\psi)$ is continuous for every $\psi \in C(\rzee)$. Standard theory entails that a sequence $(\zeta_n)_{n \in \enn}$ of elements of $C(\rzee)^\dagger$ converges to $\zeta \in C(\rzee)^\dagger$ if and only if the sequence of complex numbers $(\zeta_n(\psi))_{n \in \enn}$ converges to $\zeta(\psi)$ for every $\psi \in C(\rzee)$.\end{enumerate}

We remark that it is useful for us to introduce the weak-star abstraction, since we define by duality; the representation of measures through duality goes back to F.~Riesz \cite{JFM40.0388.03}. This topology is often called \emph{vague} in mathematics --- the name is due to Bourbaki (Ch.~III \S 1 of \cite{MR2018901}); see also \S 4.1 of \cite{MR3642325} and \S 1.10 of \cite{MR2760403}. In applied contexts sometimes this topology is called simply \emph{weak}, although the latter terminology becomes formally erroneous if extended to noncompact underlying spaces.

\begin{Definition}[Endomorphic trajectory] \label{def.endotraj} Let $\mathscr{M}$ be a subset of $\enn$ and let $\mu \in \prob(\rzee)$.
\begin{itemize}
	\item We define the \emph{endomorphic $\mathscr{M}$-trajectory of $\mu$} to be the subset $\Sigma(\mu,\mathscr{M})$ of $\prob(\rzee)$ given as follows. 
\begin{align*} \Sigma(\mu,\mathscr{M}) \coloneqq \bigl \{ m_\circledcirc \mu : m \in \mathscr{M} \cap \ennplus \bigr \} \end{align*}
\item We define the \emph{endomorphic $\mathscr{M}$-trajectory closure of $\mu$} to be the subset $\overline{\Sigma(\mu,\mathscr{M})}$ of $\prob(\rzee)$ given by the closure of $\Sigma(\mu,\mathscr{M})$ in the weak-star topology.
\end{itemize}
\end{Definition}

\subsubsection{Rudolph-Johnson type theorems: statements} \label{sec.intro.app.special.state}

We now soften Theorem \ref{thm.exist} into a statement about weak-star trajectory closures.

\begin{Theorem}[Limiting endomorphic rigidity via upper entropy dimension] \label{thm.entdim}
For all sublacunary subsets $\mathscr{M}$ of $\enn$ and all $\mu \in \prob(\rzee)$ there exists $\nu \in \overline{\Sigma(\mu,\mathscr{M})}$ such that $\nu \succeq \ach_\ast^\uparrow(\mu)\lambda$.\end{Theorem}

\begin{Remark}[Sharpness of the coefficient in Theorem \ref{thm.entdim}] \label{rem.sharp}
The coefficient $\ach_\ast^\uparrow(\mu)$ in Theorem \ref{thm.entdim} cannot be increased while retaining the theorem's scope. To see this, let $\delta_0$ denote the point mass at $0$ and let $c \in [0,1]_\arr$ be arbitrary. Then the measure $c\lambda + (1-c)\delta_0$ is $\times n$-invariant for all $n \in \ennplus$, has upper entropy dimension $c$, and dominates no multiple $d\lambda$ with $d > c$. The same family of measures is noted after Theorem 1.3 of \cite{MR2563089} as showing the sharpness of the invariant statement there. \end{Remark}

We specialize to our canonical sets.

\begin{Theorem}[Rudolph-Johnson for canonical multiplier sets] \label{thm.rjcanon}
Let $\mu \in \prob(\rzee)$.
\begin{enumerate}[label=\upshape(\arabic*)]
\item \label{thm.rjcanon.c1} Let $n_0 \in \ennplus$ and let $\mathscr{M}$ be a subset of $\enn$ such that $\mathscr{M} \cap [n^3+1,(n+1)^3-1]_{\zee} \neq \emptyset$ for all $n \in [n_0,\infty)_{\zee}$. Then there exists $\nu \in \overline{\Sigma(\mu,\mathscr{M})}$ such that $\nu \succeq \ach_\ast^\uparrow(\mu)\lambda$.
\item \label{thm.rjcanon.c2} There exists $\nu \in \overline{\Sigma(\mu,\mathbb{L})}$ such that $\nu \succeq \ach_\ast^\uparrow(\mu)\lambda$.
\item \label{thm.rjcanon.c3} There exists $\nu \in \overline{\Sigma(\mu,\mathtt{p}(\enn))}$ such that $\nu \succeq \ach_\ast^\uparrow(\mu)\lambda$.
\item \label{thm.rjcanon.c4} Let $a,b \in [2,\infty)_\zee$ be multiplicatively independent. Then there exists $\nu \in \overline{\Sigma(\mu,\mathbb{S}_{a,b})}$ such that $\nu \succeq \ach_\ast^\uparrow(\mu)\lambda$.
\end{enumerate} \end{Theorem}

Towards providing an interesting example for Clause \ref{thm.rjcanon.c1} above, we record the following heavy result due to A.W. Dudek.

\begin{oLemma}[Primes between consecutive cubes, Theorem 1.1 in \cite{MR3584567}] \label{olem.primegap}
For all $n \in \enn$ with $n \geq \exp(\exp(33.3))$ there exists $p \in \mathbb{P}$ such that $n^3 < p < (n+1)^3$. \end{oLemma}

Combining Clause \ref{thm.rjcanon.c1} of Theorem \ref{thm.rjcanon} with the above lemma we obtain the corollary below.

\begin{Corollary}[Rudolph-Johnson for cube-succeeding primes] \label{cor.rjprime} Let $\mathscr{M}$ denote the set consisting of the first prime after each perfect cube. Then for all $\mu \in \prob(\rzee)$ there exists $\nu \in \overline{\Sigma(\mu,\mathscr{M})}$ such that $\nu \succeq \ach_\ast^\uparrow(\mu)\lambda$. \end{Corollary}

Clause \ref{thm.rjcanon.c4} of Theorem \ref{thm.rjcanon} may be thought of as an `invariance-unnecessary' Rudolph-Johnson theorem. If we further assume that $\mu$ is $\times a \times b$-invariant then by iterating the invariance assumption we find $\Sigma(\mu,\mathbb{S}_{a,b})=\{\mu\}$. Since the weak-star topology is Hausdorff we obtain $\ov{\Sigma(\mu,\mathbb{S}_{a,b})}=\{\mu\}$ and so the following corollary.

\begin{Corollary}[Classical Rudolph-Johnson] \label{cor.rjclass} For all multiplicatively independent $a,b \in [2,\infty)_{\zee}$ and all $\times a \times b$-invariant $\mu \in \prob(\rzee)$ we have that $\mu \succeq \ach_\ast^\uparrow(\mu)\lambda$. \end{Corollary}

We pause to record how Corollary \ref{cor.rjclass} contains the original conclusion from the work of Rudolph in \cite{MR1062766} and Johnson in \cite{MR1194793} that $\lambda$ is the unique ergodic $\times a \times b$-invariant element of $\prob(\rzee)$ with positive entropy. The present article makes no other use of asymptotic ergodic theory, and these details are well-known in that field.  We also note that Corollary \ref{cor.rjclass} uses no transcendence theory: the sublacunarity of $\mathbb{S}_{a,b}$ which Theorem \ref{thm.entdim} consumes is Furstenberg's lemma, cited as Lemma \ref{olem.furst} below, while the Baker--W\"ustholz input of Lemma \ref{olem.bw} is spent only in Clause \ref{thm.blprime.c3} of Theorem \ref{thm.blprime}.\\

Write $\mathscr{P}_a$ for the partition of $\rzee$ into the resolution-$a$ cells $\mathscr{P}_a(0),\ldots,\mathscr{P}_a(a-1)$ of Definition \ref{def.entropy}, and note that the resolution-$a^j$ cells are the elements of the $j$-fold refinement of $\mathscr{P}_a$ under $\mathtt{T}_a$. Since $\mathscr{P}_a$ is generating for $\mathtt{T}_a$ in the sense of Theorem 4.18 from \cite{MR0648108}, if $\mu$ is $\times a$-invariant and we write $h_\mu(\mathtt{T}_a)$ for the Kolmogorov-Sinai dynamical entropy of $\mu$ with respect to $\mathtt{T}_a$ then that theorem implies:
\begin{align*} \lim_{j \to \infty} \achn(\mu,a^j) =  \frac{h_\mu(\mathtt{T}_a)}{\log(a)}  \end{align*}
Since $\ach_\ast^\uparrow(\mu)$ is a limit superior over all resolutions we find in particular that:
\begin{align*} \ach_\ast^\uparrow(\mu) \geq \frac{h_\mu(\mathtt{T}_a)}{\log(a)} \end{align*}
 Granting this, suppose that $\mu$ is ergodic with respect to the action generated by $\mathtt{T}_a$ and $\mathtt{T}_b$ and satisfies $h_\mu(\mathtt{T}_a)>0$. Corollary \ref{cor.rjclass} then supplies a strictly positive $c$ with $\mu \succeq c\lambda$, so that $\mu$ and $\lambda$ are not mutually singular. We also note that $\lambda$ is ergodic for $\mathtt{T}_a$ alone and hence for the joint action, while Lemma 4.6 in \cite{MR2723325}, applied at the level of the joint action, provides that distinct ergodic invariant measures are mutually singular. Thus we conclude that $\mu = \lambda$.
 
\chapter{Proof of Theorem \ref{thm.main}} \label{cha.pfmain}

\section{Proof of Theorem \ref{thm.main}: declarations} \label{cha.pfmain.defs}

The proof of Theorem \ref{thm.main} occupies the present chapter, and we organize it into four parts whose names we borrow from the structure of a computer program. 
\begin{itemize} \item The \emph{declarations} fix notation and elementary conventions on which the argument depends. 
\item The \emph{imports} are those results which we quote and use without proof. 
\item The \emph{subroutines} form a chain of statements, each of which is instantiated with explicit data, returns explicit data, and is called in turn by its successor. 
\item The \emph{main} argument instantiates the last subroutine and assembles the conclusion of Theorem \ref{thm.main}. \end{itemize}
We feel this terminology is germane as the proofs below proceed formally by instantiation and return.

\subsection{Further elementary conventions} \label{cha.pfmain.defs.conv}

We fix the following further conventions for the remainder of the article.
\begin{enumerate}[label=\upshape(\thesection:\alph*)]
\item \label{conv.einfint} We will simplify notation in the second coordinate of the exponential infimum from Definition \ref{def.foustat} when the bandprofile is an initial interval: for $n \in \enn$ we will write $\einf(\mu,n)$ instead of $\einf(\mu,[0,n]_{\zee})$.
\item \label{conv.min} For $x,y \in \arr$ we adopt infix notations for maximum and minimum: 
\begin{align*} x \wedge y &\coloneqq \min(x,y) & x \vee y &\coloneqq \max(x,y) \end{align*}
It will sometimes be convenient to use large prefix versions of these operators. Thus, for example, if $m,n \in \zee$ satisfy $m \leq n$ and $x_m,\ldots,x_n$ is a list of real numbers we define:
\begin{align*} \bigwedge_{j=m}^n x_j & \coloneqq \min  \bigl \{x_j: j \in [m,n]_{\zee}  \bigr \} &  \bigvee_{j=m}^n x_j & \coloneqq \max \bigl \{x_j: j \in [m,n]_{\zee} \bigr \}  \end{align*}
\end{enumerate}

\begin{Remark}[Lipschitz estimate for exponentials] \label{rem.expdiff} For all $x,y \in \arr$ we have:
\begin{align} \label{rem.expdiff.1} |e(x) - e(y)| = \bigl \vert e(x-y) - 1 \bigr \vert \leq 2 \pi |x-y| \end{align}
Here, the equality follows from the unimodularity of $e(y)$, while the inequality follows from the elementary estimate $|e(t)-1| = 2\lvert\sin(\pi t)\rvert \leq 2\pi |t|$ for $t \in \arr$. \end{Remark}

\subsection{Harmonic analysis on the circle} \label{cha.pfmain.defs.fou}

We now declare and remark upon some tools from basic analysis on the circle.

\subsubsection{Fourier quantities}

The following convention generalizes that in (\ref{eq.intro.x1}).

\begin{Definition}[Fourier coefficients] \label{def.foucoef} For $n \in \zee$ and $\zeta \in C(\rzee)^\dagger$ we define the \emph{frequency-$n$ Fourier coefficient of $\zeta$} to be the complex number given by $\hat{\zeta}(n) \coloneqq \zeta(\mathfrak{e}^{-n})$.\\

Sometimes the variables $\zeta$ and $n$ as above may need to be replaced with more complicated symbolic expressions having the same logical type, and in these cases we may adopt the more cumbersome notation $\fou(\zeta,n)$ to denote $\hat{\zeta}(n)$. \end{Definition}

The following weak gauge of measure proximity will be useful.

\begin{Definition}[Finite-bandprofile Fourier max-seminorm] \label{def.founorm} Let $\mathscr{F}$ be a finite subset of $\zee$ and let $\zeta \in C(\rzee)^\dagger$. We define the \emph{bandprofile-$\mathscr{F}$ max-seminorm of $\zeta$} as follows.
\begin{align*} ||\zeta||_{\fou(\mathscr{F})} \coloneqq \bigvee_{j \in \mathscr{F}} |\hat{\zeta}(j)|  \end{align*}
Here we understand the convention that the maximum is taken in $[0,\infty)_{\arr}$, so that in particular $||\zeta||_{\fou(\emptyset)} = 0$. For $n \in \enn$ we simplify notation by writing $||\cdot||_{\fou(n)}$ instead of $||\cdot||_{\fou([0,n]_{\zee})}$. \end{Definition}

We also note the following symmetry.

\begin{Remark}[Conjugate symmetry of Fourier coefficients] \label{rem.fournompos} It follows directly from Definition \ref{def.foucoef} that if $\mu,\nu \in \operatorname{FMea}(\rzee)$ then for all $n \in \zee$ we have $\fou(\mu - \nu, -n) = \ov{\fou(\mu-\nu,n)}$. \end{Remark}

\subsubsection{Connections to endomorphy}

We now collect certain further remarks about how harmonic analytic quantities behave under circle endomorphy.

\begin{Remark}[Submeasures and endomorphs] \label{rem.pushsub} It follows directly from Definition \ref{def.endomorph} that if $\mu,\nu \in \operatorname{FMea}(\rzee)$ satisfy $\nu \preceq \mu$ then for all $n \in \ennplus$ we have $n_\circledcirc \nu \preceq n_{\circledcirc} \mu$.\\

We note that if $\zeta \in \prob(\rzee)$ then also $n_\circledcirc \zeta \in \prob(\rzee)$, directly from \ref{conv.dual}.\end{Remark}

\begin{Remark}[Endomorphs and change-of-variable formulas] \label{rem.endochangeofvar} Let $n \in \ennplus$ and let $\mu \in \operatorname{FMea}(\rzee)$. 
\begin{enumerate}[label=\upshape(\arabic*)] \item \label{rem.endochangeofvar.c1} It follows directly from Definition \ref{def.endomorph} that $n_\circledcirc \mu$ satisfies $\widehat{\,n_\circledcirc \mu\,}(\ell) = \hat{\mu}(n\ell)$ for all $\ell \in \zee$.
\item \label{rem.endochangeofvar.c2} It is similarly direct that $n_\circledcirc \mu$ is the unique element of $\operatorname{FMea}(\rzee)$ which satisfies the following for all $\psi \in C(\rzee)$.
\begin{align*} \int_{\rzee} \psi(x) \deee [n_\circledcirc\mu](x) = \int_{\rzee} \psi(nx) \deee \mu(x) \end{align*}
\item \label{rem.endochangeofvar.c3} Combining Definition \ref{def.endomorph} with \ref{conv.dual} and Proposition 10.1 in \cite{MR1681462} we find that for all Borel subsets $E$ of $\rzee$ we have $[n_\circledcirc \mu](E) = \mu\bigl(\mathtt{T}_n^{-1}(E)\bigr)$.\end{enumerate}\end{Remark}

\begin{Remark}[Endomorphy of exponential infimum] \label{rem.einfendo} It follows directly from the definition of the exponential infimum in Definition \ref{def.foustat} that for all $\mu \in \prob(\rzee)$, all nonempty finite subsets $\mathscr{F}$ of $\zee$ and all $p \in \ennplus$ we have $\einf(\mu,\{p\}\mathscr{F}) = \einf(p_\circledcirc \mu,\mathscr{F})$. In particular, if $\mu$ is $\times p$-invariant then we have $\einf(\mu,\{p\}\mathscr{F}) = \einf(\mu,\mathscr{F})$.  \end{Remark}

\begin{Remark}[Endomorphs and the Fourier max-seminorm] \label{rem.endofou} Let $\mu,\nu \in \prob(\rzee)$, let $n \in \ennplus$ and let $\mathscr{F}$ be a nonempty finite subset of $\zee$. Combining the definition of the Fourier max-seminorm in Definition \ref{def.founorm} with Clause \ref{rem.endochangeofvar.c1} in Remark \ref{rem.endochangeofvar} we find $||n_\circledcirc\mu-n_\circledcirc\nu||_{\fou(\mathscr{F})} = ||\mu-\nu||_{\fou(\{n\}\mathscr{F})}$.\end{Remark}

\subsubsection{Resolution norms}

We now introduce certain norms which measure uniformity of a measure in a stronger way than the resolution entropy from Definition \ref{def.entropy}.

\begin{Definition}[Resolution norms of measures] \label{def.mupnorm} Let $n \in \ennplus$, let $\mu \in \prob(\rzee)$ and let $p\in [1,\infty)_{\arr}$. Recalling the resolution-$n$ cells $\mathscr{P}_n(0),\ldots,\mathscr{P}_n(n-1)$ from Definition \ref{def.entropy}, we define the \emph{resolution-$n$ $\ell^p$-norm of $\mu$} to be the element $||\mu(\mathscr{P}_n)||_p$ of $[0,\infty)_{\arr}$ given as follows.
\begin{align} ||\mu(\mathscr{P}_n)||_p \coloneqq \left(\sum_{j=0}^{n-1}\mu(\mathscr{P}_n(j))^p\right)^{1/p} \label{eq.mupnorm} \end{align} \end{Definition}

\begin{Remark}[Resolution entropy and $\ell^p$-norms] \label{rem.mupnorm} Let $n \in \ennplus$, let $\mu \in \prob(\rzee)$ and observe that the formula in (\ref{eq.mupnorm}) returns the somewhat degenerate result $||\mu(\mathscr{P}_n)||_1 = 1$ for the resolution-$n$ $\ell^1$-norm of $\mu$. However, the resolution-$n$ entropy $\ach(\mu,n)$ originally introduced in Definition \ref{def.entropy} can be expressed as a certain renormalized limit of $||\mu(\mathscr{P}_n)||_p$ as $p \downarrow 1$. More explicitly, for $p \in (1,\infty)_{\arr}$ we define the \emph{resolution}-$n$ \emph{exponent-}$p$ \emph{R\'{e}nyi entropy} of $\mu$ as follows. 
\begin{align*} \ach^p(\mu,n) \coloneqq -\frac{p}{p-1}\log\bigl(||\mu(\mathscr{P}_n)||_p\bigr)\end{align*} 
Then according to \cite{MR0132570}, we have that $\ach^p(\mu,n)$ is nonincreasing with $p$ and we recover the entropy $\ach(\mu,n)$ as in (\ref{def.entropy.0}) by the limit below.
\begin{align*} \lim_{p \downarrow 1} \ach^p(\mu,n) = \ach(\mu,n) \end{align*} \end{Remark}

\begin{Remark}[Refinement monotonicity of resolution norms] \label{rem.mupnormref} Let $\mu \in \prob(\rzee)$ and let $a,b \in \ennplus$ be such that $a$ divides $b$. It follows directly from Definition \ref{def.entropy} that each cell of $\mathscr{P}_{a}$ is the disjoint union of $b/a$ cells of $\mathscr{P}_{b}$. Combining this observation with the superadditivity of the map $t \mapsto t^2$ on $[0,\infty)_{\arr}$ and Definition \ref{def.mupnorm} we find $\bigl \vert \bigl \vert \mu(\mathscr{P}_{b}) \bigr \vert \bigr \vert_2 \leq \bigl \vert \bigl \vert \mu(\mathscr{P}_{a}) \bigr \vert \bigr \vert_2$. \end{Remark}

\subsection{Gram matrices} \label{cha.pfmain.defs.gram}

We fix the following conventions regarding linear algebra.

\begin{enumerate}[resume, label=\upshape(\thesection:\alph*)]
\item \label{conv.matrix} Given a nonempty finite set $F$, we adopt the convention that a \emph{matrix indexed by $F$} refers to a function $\eusc{A}: F \times F \to \cee$. We define such a matrix $\eusc{A}$ to be \emph{self-adjoint} if $\eusc{A}(\omega,\tau) = \ov{\eusc{A}(\tau,\omega)}$ for all $\tau,\omega \in F$.
\item \label{conv.leasteigen} Let $F$ be a nonempty finite set and let $\eusc{A}$ be a self-adjoint matrix indexed by $F$. The spectral theorem (see, for example, \S 4.1 in \cite{MR2978290}) provides that $\cee^F$ admits an orthonormal basis of eigenvectors of $\eusc{A}$ with real eigenvalues, and we adopt the notation $\rho_\downarrow(\eusc{A})$ to denote the least eigenvalue of $\eusc{A}$.
\end{enumerate}

\begin{Definition}[Gram matrices] \label{def.gramabstract} Let $\mathcal{H}$ be a Hilbert space, let $F$ be a nonempty finite set and $\mathcal{X}:F \to \mathcal{H}$ be a function. 
\begin{itemize} \item We define the \emph{Gram matrix of $\mathcal{X}$} to be the $|F| \times |F|$ matrix $\eusc{G}_{\mathcal{X}}: F \times F \to \cee$ given by stipulating that $\eusc{G}_{\mathcal{X}}(\tau,\omega) \coloneqq \langle \mathcal{X}(\tau) , \mathcal{X}(\omega) \rangle$ for all $\tau,\omega \in F$. \end{itemize}
Assume now that there exists $n \in \enn$ such that $F = [0,n]_\zee$.  
\begin{itemize}  \item Adopting the notation $\xi_j \coloneqq \mathcal{X}(j)$ for $j \in [0,n]_\zee$, we may alternatively denote the Gram matrix $\eusc{G}_\mathcal{X}$ as $\eusc{G}_{\xi_0,\ldots,\xi_n}$. \end{itemize} \end{Definition}

Keeping the notation of Definition \ref{def.gramabstract}, we offer the following illustration.
\begin{align*} \eusc{G}_{\xi_0,\ldots,\xi_{n}} &= \begin{pmatrix} \langle \xi_0, \xi_0 \rangle & \cdots & \langle \xi_0, \xi_{n} \rangle  
\\ \vdots & \ddots & \vdots 
\\  \langle \xi_{n}, \xi_0 \rangle & \cdots & \langle \xi_{n}, \xi_{n} \rangle  
 \end{pmatrix} \end{align*}

 \begin{Remark}[Gram matrices and norm] \label{rem.gramnorm} It follows directly from the definition of the Gram matrix in Definition \ref{def.gramabstract} that for all Hilbert spaces $\mathcal{H}$, all $n \in \enn$, all lists $\xi_0,\ldots,\xi_{n}$ of vectors in $\mathcal{H}$ and all $\alpha \in \cee^{[0,n]_{\zee}}$ we have:
 \begin{align*} \dlvertt \sum_{j=0}^n \alpha(j) \xi_j \drvertt^2 = \sum_{j,k=0}^n \alpha(j) \ov{\alpha(k)} \, \eusc{G}_{\xi_0,\ldots,\xi_n}(j,k) \end{align*}
 \end{Remark}
 
 \begin{Remark}[Self-adjointness of Gram matrices] \label{rem.gramsa} Directly from the conjugate symmetry of the inner product we have $\eusc{G}_{\mathcal{X}}(\tau,\omega) = \ov{\eusc{G}_{\mathcal{X}}(\omega,\tau)}$, so that every Gram matrix is self-adjoint. In particular $\rho_\downarrow(\eusc{G}_{\mathcal{X}})$ is well-defined by way of \ref{conv.leasteigen}, and every quadratic form $\sum_{j,k} \alpha(j)\ov{\alpha(k)}\, \eusc{G}_{\mathcal{X}}(j,k)$ is real-valued. \end{Remark}

The following standard result appears, for example, as Theorem 4.2.6 in \cite{MR2978290}.

\begin{oLemma}[Courant-Fischer minimax principle at the bottom] \label{olem.courfish} Let $F$ be a nonempty finite set and let $\eusc{A}:F \times F \to \cee$ be self-adjoint. Then we have:
\begin{align*} \rho_\downarrow(\eusc{A}) = \inf \left\{ \sum_{\tau,\omega \in F} \alpha(\tau) \ov{\alpha(\omega)} \,\eusc{A}(\tau,\omega) : \alpha \in \cee^{F} \mbox{ and } \sum_{\tau \in F}  |\alpha(\tau)|^2 =1 \right\} \end{align*} \end{oLemma}

The remark below provides certain useful reinterpretations of $\rho_\downarrow$.

\begin{Remark}[Variational formula for least eigenvalue of Gram matrix] \label{rem.gramvar} Let $\mathcal{H}$ be a Hilbert space, let $n \in \enn$ and let $\xi_0,\ldots,\xi_n$ be a list of vectors in $\mathcal{H}$.
\begin{enumerate}[label=\upshape(\arabic*)] \item \label{rem.gramvar.c0} Since Remark \ref{rem.gramsa} provides self-adjointness, it follows from Lemma \ref{olem.courfish} that we have:
\begin{align*} \rho_\downarrow(\eusc{G}_{\xi_0,\ldots,\xi_n}) = \inf \left\{ \sum_{j,k=0}^n \alpha(j) \ov{\alpha(k)} \,\eusc{G}_{\xi_0,\ldots,\xi_n}(j,k) : \alpha \in \cee^{[0,n]_{\zee}} \mbox{ and } \sum_{j=0}^n  |\alpha(j)|^2 =1 \right\} \end{align*}
\item \label{rem.gramvar.c1} Combining Clause \ref{rem.gramvar.c0} with Remark \ref{rem.gramnorm} we find that the following holds for all $\alpha \in \cee^{[0,n]_\zee}$.
\begin{align*}  \dlvertt \sum_{j=0}^n \alpha(j) \xi_j \drvertt^2 \geq \rho_\downarrow(\eusc{G}_{\xi_0,\ldots,\xi_n}) \sum_{j=0}^n  |\alpha(j)|^2  \end{align*} \end{enumerate} \end{Remark}

\subsection{Bandprofiled Toeplitz matrices} \label{cha.pfmain.defs.toe}

\subsubsection{Bandprofiled Toeplitz matrices: definitions}

We continue the running list of conventions.

\begin{enumerate}[resume, label=\upshape(\thesection:\alph*)]
\item \label{conv.setarith} If $\mathscr{E}$ and $\mathscr{F}$ are subsets of $\zee$ we write $\mathscr{E} + \mathscr{F} \coloneqq \{j + k : j \in \mathscr{E} \mbox{ and } k \in \mathscr{F}\}$ and $-\mathscr{F} \coloneqq \{-j : j \in \mathscr{F}\}$, and we abbreviate $\mathscr{E} - \mathscr{F} \coloneqq \mathscr{E} + (-\mathscr{F})$.
\item \label{conv.finfun} We will sometimes identify a finite subset $\mathscr{F}$ of $\zee$ with the unique increasing function from $[0,|\mathscr{F}|-1]_{\zee}$ to $\mathscr{F}$. Thus, for example, $\mathscr{F}(0)$ should be interpreted as the least element of $\mathscr{F}$ and $\mathscr{F}(|\mathscr{F}|-1)$ should be interpreted as the greatest element of $\mathscr{F}$.
\item \label{conv.ltwo} If $\mu \in \prob(\rzee)$ we equip the space $L^2(\rzee,\mu)$ with the inner product given by $\langle \psi , \phi \rangle \coloneqq \mu(\psi \ov{\phi})$, so that in particular complex inner products are linear in the first coordinate and conjugate-linear in the second coordinate.
\end{enumerate}

The following definition provides a local bridge from Fourier analysis to linear algebra.

\begin{Definition}[Bandprofiled Toeplitz matrices of measures] \label{def.toeplitz} Let $\mu \in \prob(\rzee)$, let $\mathscr{F}$ be a finite subset of $\zee$ and let $n \in [0,|\mathscr{F}|-1]_{\zee}$. We define the \emph{bandprofile-$\mathscr{F}$ index-$n$ Toeplitz matrix of $\mu$} to be the $(n+1) \times (n+1)$ matrix $\eusc{T}_{\mathscr{F},n}^\mu: [0,n]_{\zee} \times [0,n]_{\zee} \to \cee$ given by the stipulation that $\eusc{T}_{\mathscr{F},n}^\mu(j,k) \coloneqq \hat{\mu}(\mathscr{F}(k)-\mathscr{F}(j))$ for all $j,k \in [0,n]_{\zee}$. \end{Definition}

Noting that $\hat{\mu}(-\ell) = \ov{\hat{\mu}(\ell)}$ for all $\mu \in \prob(\rzee)$ and all $\ell \in \zee$ (as in Remark \ref{rem.fournompos}), we offer the following illustration.
\begin{align*} \eusc{T}_{[0,3]_{\zee},3}^{\mu} = \begin{pmatrix} 1 & \hat{\mu}(1) & \hat{\mu}(2) & \hat{\mu}(3) \vspace{3 pt}
\\ \ov{\hat{\mu}(1)} & 1 & \hat{\mu}(1) &  \hat{\mu}(2) \vspace{3 pt}
\\ \ov{\hat{\mu}(2)} & \ov{\hat{\mu}(1)} &  1 &  \hat{\mu}(1) \vspace{3 pt}
\\ \ov{\hat{\mu}(3)} & \ov{\hat{\mu}(2)} & \ov{\hat{\mu}(1)} &  1  \end{pmatrix} \end{align*}

\subsubsection{Bandprofiled Toeplitz matrices: remarks}

We now collect certain remarks related to Definition \ref{def.toeplitz}.

\begin{Remark}[Gram and Toeplitz] \label{rem.gramtoe} Let $\mathscr{F}$ be a finite subset of $\zee$, let $\mu \in \prob(\rzee)$ and for $n \in [0,|\mathscr{F}|-1]_{\zee}$ define $\xi_n \in L^2(\rzee,\mu)$ to be the $\mu$-equivalence class of the function $\mathfrak{e}^{\mathscr{F}(n)}$. By combining the definition of the Gram matrix in Definition \ref{def.gramabstract} with Definition \ref{def.toeplitz} we find that $\eusc{T}^\mu_{\mathscr{F},n} = \eusc{G}_{\xi_0,\ldots,\xi_{n}}$. \end{Remark}

\begin{Remark}[Alternative expressions for the exponential infimum] \label{rem.einf} Let $\mu \in \prob(\rzee)$ and let $\mathscr{F}$ be a nonempty finite subset of $\zee$. 
\begin{enumerate}[label=\upshape(\arabic*)]
	\item \label{rem.einf.c1} The exponential infimum may be expressed more directly in terms of Fourier coefficients as follows.
\begin{align*} \einf(\mu , \mathscr{F}) =\inf \left \{ \sum_{j,k \in \mathscr{F}} \alpha(j) \ov{\alpha(k)}\,\hat{\mu}(k-j) :\alpha \in \cee^\mathscr{F} \mbox{ and } \sum_{j \in \mathscr{F}} |\alpha(j)|^2 = 1 \right \} \end{align*} \end{enumerate}\end{Remark}

\begin{Remark}[Toeplitz matrices and bilinear forms]  \label{rem.toebilin} Let $\mu \in \prob(\rzee)$, let $\mathscr{F}$ be a nonempty finite subset of $\zee$ and let $n \in [0,|\mathscr{F}|-1]_{\zee}$. 
\begin{enumerate}[label=\upshape(\arabic*)] \item Combining the identification of Gram and Toeplitz matrices described in Remark \ref{rem.gramtoe} with the variational formula in Remark \ref{rem.gramvar} we find that the following holds.
\begin{align*} \rho_\downarrow(\eusc{T}^\mu_{\mathscr{F},n}) = \inf \left \{ \sum_{j,k=0}^{n} \alpha(j) \ov{\alpha(k)} \,\hat{\mu}(\mathscr{F}(k)-\mathscr{F}(j)) :\alpha \in \cee^{[0,n]_{\zee}} \mbox{ and } \sum_{j=0}^{n} |\alpha(j)|^2 =1 \right \}\end{align*}
\item \label{rem.toebilin.c2} Since $\{\mathscr{F}(0),\ldots,\mathscr{F}(n)\} \subseteq \mathscr{F}$, we can combine the above display with the expression in Clause \ref{rem.einf.c1} of Remark \ref{rem.einf} to find that $\rho_\downarrow(\eusc{T}^\mu_{\mathscr{F},n}) \geq  \einf(\mu , \mathscr{F})$.
\item \label{rem.toebilin.c3}  Since $\{\mathscr{F}(0),\ldots,\mathscr{F}(|\mathscr{F}|-1)\} = \mathscr{F}$, we also find that $\rho_\downarrow(\eusc{T}^\mu_{\mathscr{F},|\mathscr{F}|-1}) =  \einf(\mu , \mathscr{F})$. \end{enumerate} \end{Remark}
 
\begin{Remark}[Exponential infimum and submeasures] \label{rem.einfmono}
Let $\mathscr{F}$ be a nonempty finite subset of $\zee$, let $t \in [0,1]_{\arr}$ and let $\nu,\mu \in \prob(\rzee)$ be such that $t \nu \preceq \mu$. It follows directly from the definition of the exponential infimum in Definition \ref{def.foustat} that we have $t \einf(\nu,\mathscr{F}) \leq \einf(\mu,\mathscr{F})$. \end{Remark}

\section{Proof of Theorem \ref{thm.main}: imports} \label{cha.pfmain.imports}

The following standard facts appear, for example, respectively in \S 2.3 and as Lemma 2.9 of \cite{MR2563089}.

\begin{oLemma}[Basic estimates for resolution entropy] \label{olem.entropy} Let $n \in \ennplus$.
\begin{enumerate}[label=\upshape(\arabic*)] \item \label{olem.entropy.c2} Let $\mu,\nu,\varsigma \in \prob(\rzee)$ and assume that $\mu = s \nu + (1-s) \varsigma$ for some $s \in [0,1]_{\arr}$. Then we have:
\begin{align}\label{olem.entropy.c2.1} 0 \leq \ach(\mu,n) -s \ach(\nu,n) - (1-s) \ach(\varsigma,n) \leq -s\log(s)-(1-s)\log(1-s)\end{align}
\item \label{olem.entropy.c5} For all $\mu \in \prob(\rzee)$ and all $r \in [n,\infty)_\zee$ we have $\ach(\mu,r) - \log(\lceil r/n \rceil +1) \leq \ach(\mu,n)$.  \end{enumerate}	\end{oLemma} 

The following standard fact appears, for example, in Ch.~I \S 2.5 of \cite{MR2039503}.

\begin{oLemma}[Fej\'er kernel] \label{olem.fejer}
For $K \in \ennplus$ define the \emph{Fej\'er kernel} $\mathtt{F}_K: \arr \to \arr$ as follows.
\begin{align} \mathtt{F}_K(u) \coloneqq \sum_{j=-K+1}^{K-1}\left(1 - \frac{\lvert j \rvert}{K}\right) e(ju) \label{olem.fejer.h1} \end{align}
Also define $\mathfrak{d}: \arr \to [0,1/2]_{\arr}$ by $\mathfrak{d}(u) \coloneqq \inf_{p \in \zee} \lvert u - p \rvert$. Then for all $u \in \arr$ we have $\mathtt{F}_K(u) \in [0,K]_{\arr}$, while for all noninteger $u \in \arr$ we have:
\begin{align} \mathtt{F}_K(u) = \frac{\sin^2(\pi K u)}{K\sin^2(\pi u)} \leq \frac{1}{4K\,\mathfrak{d}(u)^2} \label{olem.fejer.conc1} \end{align} \end{oLemma}

\section{Proof of Theorem \ref{thm.main}: subroutines} \label{cha.pfmain.infra}

\subsection{Least eigenvalue and off-diagonal inequality} \label{cha.pfmain.infra.lsv}

The following Lemma \ref{lem.lsv} is a routine linear-algebraic estimate, converting small off-diagonal entries into a lower bound for the minimal eigenvalue. We will interpret it in the Fourier-analytic context of Clause \ref{rem.toebilin.c3} from Remark \ref{rem.toebilin}, for ultimate consumption by Lemma \ref{lem.unifeinf}.\\

The dependencies in the statement of Lemma \ref{lem.lsv} are as follows.
\begin{deps}
\dep{Conventions}{\depitem{$\rho_\downarrow$}{conv.leasteigen}}
\dep{Definitions}{\depitem{Gram matrices}{def.gramabstract}}
\end{deps}

We now move on to the formal statement.

\begin{Lemma}[Least eigenvalue and off-diagonal maximum inequality] \label{lem.lsv} For all Hilbert spaces $\mathcal{H}$, all $n \in \enn$ and all lists $\xi_0,\ldots,\xi_{n}$ of unit vectors in $\mathcal{H}$ we have:
\begin{align} \rho_\downarrow(\eusc{G}_{\xi_0,\ldots,\xi_n})  \geq 1- n(n+1)\bigvee_{j,k \in [0,n]_{\zee}\,: j \neq k} \bigl \vert \eusc{G}_{\xi_0,\ldots,\xi_n}(j,k) \bigr \vert \label{lem.lsv.conc} \end{align}
  \end{Lemma}

We adopt the convention that when $n = 0$ the maximum appearing in (\ref{lem.lsv.conc}) is empty, and that its product with $n(n+1) = 0$ is to be interpreted as $0$.\\

The dependencies in the proof of Lemma \ref{lem.lsv} beyond those in its statement are as follows.
\begin{deps}
\dep{Remarks}{\depitem{Self-adjointness of Gram matrices}{rem.gramsa}\depsep\depitem{Variational formula for $\rho_\downarrow$}{rem.gramvar}}
\end{deps}

We now move on to the formal proof.

\begin{Proof}[of Lemma \ref{lem.lsv}] Let $\mathcal{H},n$ and $\xi_0,\ldots,\xi_n$ be as in the statement of Lemma \ref{lem.lsv}. It will be convenient to simplify notation in the current proof by writing $\eusc{G}$ instead of $\eusc{G}_{\xi_0,\ldots,\xi_n}$. Since $\xi_0$ was assumed to be a unit vector, if $n =0$ then $\eusc{G} = (1)$ and so the desired inequality (\ref{lem.lsv.conc}) is verified as $1 \geq 1$. Thus we may assume without loss of generality that $n \geq 1$.\\

Let $\alpha: [0,n]_{\zee} \to \cee$ satisfy $\sum_{j=0}^n |\alpha(j)|^2 =1$. We compute:
\begin{align}  \lvertt \sum_{j,k \in [0,n]_{\zee}\,: \,j \neq k}  \alpha(j) \ov{\alpha(k)}\, \eusc{G}(j,k) \rvertt & \leq \sum_{j,k \in [0,n]_{\zee} \,: \,j \neq k} \lvert \eusc{G}(j,k) \rvert \label{lem.lsv.1}
 \\ & \leq n(n+1)\bigvee_{j,k \in [0,n]_{\zee} \,: \,j \neq k} \lvert \eusc{G}(j,k) \rvert  \label{lem.lsv.2} \end{align}
 This normalization of $\alpha$ implies that $|\alpha(j)| \leq 1$ for all $j \in [0,n]_{\zee}$, and thereby justifies the inequality in (\ref{lem.lsv.1}). Again using our normalization of $\alpha$ we find:
\begin{align*}  \sum_{j,k =0}^n  \alpha(j) \ov{\alpha(k)}\, \bigl  \langle \xi_j, \xi_k\bigr \rangle  =  1+\sum_{j,k \in [0,n]_{\zee}\,: \,j \neq k}  \alpha(j) \ov{\alpha(k)}\, \bigl  \langle \xi_j, \xi_k\bigr \rangle \end{align*}
By construction we have $\eusc{G}(j,k) = \langle \xi_j,\xi_k \rangle$ and so both sides of the previous display are real-valued by Remark \ref{rem.gramsa}, and so we may combine it with the computation ending in (\ref{lem.lsv.2}) to find:
\begin{align*}  \sum_{j,k =0}^n  \alpha(j) \ov{\alpha(k)}\, \eusc{G}(j,k) \geq  1- n(n+1)\bigvee_{j,k \in [0,n]_{\zee} \,: \,j \neq k} \lvert \eusc{G}(j,k) \rvert \end{align*}
Combining the previous display with the expression for $\rho_\downarrow(\eusc{G}_{\xi_0,\ldots,\xi_n})$ appearing in Clause \ref{rem.gramvar.c0} in Remark \ref{rem.gramvar} completes the proof of Lemma \ref{lem.lsv}. \end{Proof}
 
\subsection{Fourier decay and exponential infimum} \label{cha.pfmain.infra.unifeinf}

The following Lemma \ref{lem.unifeinf} is the Fourier-analytic translation of Lemma \ref{lem.lsv}. This will be consumed later in Lemma \ref{lem.equildy}, where it converts a Fourier max-seminorm estimate on a submeasure $\nu$ into an exponential infimum bound for $\nu$ and thus obtains leverage from monotonicity of the exponential infimum under the submeasure relation.\\

The dependencies in the statement of Lemma \ref{lem.unifeinf} are as follows.
\begin{deps}
\dep{Definitions}{\depitem{Exponential infimum}{def.foustat}\depsep\depitem{Fourier max-seminorm}{def.founorm}}
\end{deps}

We now move on to the formal statement.

\begin{Lemma}[Exponential infimum forced by Fourier decay] \label{lem.unifeinf} For all $\mu \in \prob(\rzee)$ and all nonempty finite subsets $\mathscr{F}$ of $\zee$ we have:
\begin{align} \einf(\mu,\mathscr{F}) \geq 1- (|\mathscr{F}|^2-|\mathscr{F}|)\, ||\mu||_{\fou((\mathscr{F} -\mathscr{F})\setminus\{0\})}   \label{lem.unifeinf.conc} \end{align} \end{Lemma}

The dependencies in the proof of Lemma \ref{lem.unifeinf} beyond those in the statement are as follows.
\begin{deps}
\dep{Definitions}{\depitem{Bandprofiled Toeplitz matrices}{def.toeplitz}}
\dep{Remarks}{\depitem{Exponential infimum of singletons}{rem.einfsing}\depsep\depitem{Gram and Toeplitz}{rem.gramtoe}\depsep\depitemc{Toeplitz matrices and bilinear forms}{rem.toebilin}{rem.toebilin.c3}}
\dep{Lemmas}{\depitem{Least eigenvalue and off-diagonal maximum}{lem.lsv}}
\end{deps}
We now move on to the formal proof.
\begin{Proof}[of Lemma \ref{lem.unifeinf}] Let $\mu$ and $\mathscr{F}$ be as in the statement of Lemma \ref{lem.unifeinf}, and define $n \in \enn$ by $n \coloneqq |\mathscr{F}|-1$. Observe that if $n =0$ then Remark \ref{rem.einfsing} implies that the left side of the inequality in (\ref{lem.unifeinf.conc}) is equal to $1$, and hence the inequality is satisfied. Thus we may assume without loss of generality that $n>0$. We have:
\begin{align} \bigvee_{j,k \in [0,n]_{\zee}: j \neq k } \lvert \eusc{T}^\mu_{\mathscr{F},n}(j,k) \rvert = \bigvee_{j,k \in \mathscr{F}\,:\, j \neq k} \lvert \hat{\mu}(k-j) \rvert \leq  ||\mu||_{\fou((\mathscr{F} -\mathscr{F})\setminus \{0\})} \label{lem.unifeinf.1} \end{align}
Here, the equality on the left follows from unwrapping the definition of the Toeplitz matrix in Definition \ref{def.toeplitz}, while the inequality on the right follows from unwrapping the definition of the Fourier max-seminorm in Definition \ref{def.founorm}.\\

Now, define a list of vectors $\xi_0,\ldots,\xi_{n} \in L^2(\rzee,\mu)$ by stipulating that $\xi_j \coloneqq \mathfrak{e}^{\mathscr{F}(j)}$ for all $j \in [0,n]_{\zee}$. In the remainder of the current proof, it will be convenient to adopt the notation of the inner product $\langle \cdot, \cdot \rangle$ on $L^2(\rzee,\mu)$ from \ref{conv.ltwo} along with the associated norm $||\cdot||$.\\

Since $\lvert \mathfrak{e}^{\mathscr{F}(j)} \rvert$ is identically equal to $1$ on $\rzee$, we have $||\xi_j|| = 1$ for all $j \in [0,n]_{\zee}$, and so $\xi_0,\ldots,\xi_n$ is a list of unit vectors. Remark \ref{rem.gramtoe} implies that $\eusc{G}_{\xi_0,\ldots,\xi_n} = \eusc{T}^\mu_{\mathscr{F},n}$ and so we may rewrite the expression in (\ref{lem.unifeinf.1}) as follows.
\begin{align*} \bigvee_{j,k \in [0,n]_{\zee}: j \neq k } \lvert \eusc{G}_{\xi_0,\ldots,\xi_n}(j,k) \rvert \leq  ||\mu||_{\fou((\mathscr{F} -\mathscr{F})\setminus\{0\})}  \end{align*}
Combining the previous display with Lemma \ref{lem.lsv} we find:
\begin{align*} \rho_\downarrow(\eusc{G}_{\xi_0,\ldots,\xi_n}) \geq 1 - n(n+1)\, ||\mu||_{\fou((\mathscr{F} -\mathscr{F})\setminus\{0\})}   \end{align*}
Combining the above display with Clause \ref{rem.toebilin.c3} in Remark \ref{rem.toebilin} and again using our prior observation that $\eusc{G}_{\xi_0,\ldots,\xi_n} = \eusc{T}^\mu_{\mathscr{F},n}$ we find:
\begin{align*} \einf(\mu,\mathscr{F}) \geq 1 - n(n+1)\, ||\mu||_{\fou((\mathscr{F} -\mathscr{F})\setminus\{0\})}   \end{align*}
Unwrapping the definition of $n$ in the above display completes the proof of Lemma \ref{lem.unifeinf}. \end{Proof}

\subsection{Averaged Fourier decay at a single resolution} \label{cha.pfmain.infra.fouravg}

The purpose of Lemma \ref{lem.fouravg} below is to extract an effective constant from the proof of Lemma 4.5 in \cite{MR2563089}. We find that the $\alpha_1$ appearing there may be taken as $32$. Lemma \ref{lem.fouravg} will be consumed by Lemma \ref{lem.shadave} below in an averaged form. \\

The dependencies in the statement of Lemma \ref{lem.fouravg} are as follows.
\begin{deps}
\dep{Definitions}{\depitem{Resolution entropy}{def.entropy}\depsep\depitem{Fourier coefficients}{def.foucoef}\depsep\depitem{Resolution norms of measures}{def.mupnorm}}
\end{deps}

We now move on to the formal statement.

\begin{Lemma}[Averaged Fourier decay at a single resolution] \label{lem.fouravg}
For all $\nu \in \prob(\rzee)$, all $M \in \ennplus$ and all $g \in \zee \setminus \{0\}$ we have:
\begin{align} \frac{1}{M} \sum_{m=0}^{M-1} \bigl \vert \hat{\nu}(mg) \bigr \vert^2 \leq 32 \, \lvert g \rvert \, \bigl \vert \bigl \vert \nu(\mathscr{P}_M) \bigr \vert \bigr \vert_2^2 \label{lem.fouravg.conc1} \end{align} \end{Lemma}

An outline of the proof of Lemma \ref{lem.fouravg} is as follows.
\begin{itemize}
\item After accounting for Fourier symmetries, the left side of (\ref{lem.fouravg.conc1}) is the flat average of $\lvert \hat{\varsigma}(m) \rvert^2$ over $m \in [0,M-1]_{\zee}$. This is recorded at (\ref{lem.fouravg.p13.1}) and then we establish the estimate for $\varsigma$ at frequency one. The passage back to $\nu$ costs exactly the factor $g$, by combining the decomposition in (\ref{lem.fouravg.p11.2}) with refinement monotonicity.
\item We then transfer the argument to the two-torus. The flat average is dominated at (\ref{lem.fouravg.p1}) by a Fej\'er-weighted average. On the other hand, the weighted average is the integral of the kernel $G(y-x)$ against $\varsigma$ in each variable separately, so we can bound by the mass that the pair $(x,y)$ places near the diagonal of $(\rzee)^2$ when both coordinates are sampled from $\varsigma$. 
\item By (\ref{lem.fouravg.p2}) and (\ref{lem.fouravg.p3}) the kernel is a ridge along the diagonal, of height $4$ and width comparable to $1/M$, decaying like the inverse square of the distance from the diagonal. Thus mean-square Fourier decay along $[0,M-1]_{\zee}$ is controlled by near-diagonal correlation.
\item The cell partition discretizes the ridge, using the matrix $B$ of (\ref{lem.fouravg.p6}). We obtain a bound on its row sums at (\ref{lem.fouravg.p9.3}) using the Basel evaluation. The mass of the diagonal band itself is $\sum_k \varsigma_k^2$, which is the resolution norm carried by the conclusion, at (\ref{lem.fouravg.p10.3}). This can be interpreted as the probability that two independent samples of $\varsigma$ land in the same cell.
\end{itemize}

The dependencies in the proof of Lemma \ref{lem.fouravg} beyond those in the statement are as follows.
\begin{deps}
\dep{Conventions}{\depitem{Indicator functions}{conv.indicator}}
\dep{Definitions}{\depitem{Endomorphs of complex measures}{def.endomorph}}
\dep{Remarks}{\depitem{Conjugate symmetry of Fourier coefficients}{rem.fournompos}}
\depcont{\depitemc{Fourier coefficients of endomorphs}{rem.endochangeofvar}{rem.endochangeofvar.c1}}
\depcont{\depitemc{Change of variable for endomorphs}{rem.endochangeofvar}{rem.endochangeofvar.c3}}
\depcont{\depitem{Refinement monotonicity of resolution norms}{rem.mupnormref}}
\dep{Imports}{\depitem{Fej\'er kernel}{olem.fejer}}
\end{deps}

We now move on to the formal proof.

\begin{Proof}[of Lemma \ref{lem.fouravg}] Let $\nu, M$ and $g$ be as in the statement of Lemma \ref{lem.fouravg}. Instantiating Remark \ref{rem.fournompos} with the pair $(\nu,0)$ provides $\hat{\nu}(-n) = \ov{\hat{\nu}(n)}$ for all $n \in \zee$, so that both sides of (\ref{lem.fouravg.conc1}) take the same value at $g$ and at $-g$. Thus we may assume without loss of generality that $g \in \ennplus$. We define $\varsigma \coloneqq g_{\,\circledcirc}\, \nu$ according to Definition \ref{def.endomorph}; Clause \ref{rem.endochangeofvar.c3} of Remark \ref{rem.endochangeofvar} expresses $\varsigma$ as the pushforward of $\nu$ along $\mathtt{T}_g$, so that in particular $\varsigma \in \prob(\rzee)$, while Clause \ref{rem.endochangeofvar.c1} of the same remark provides $\hat{\varsigma}(j) = \hat{\nu}(jg)$ for all $j \in \zee$.\\

We define $c: \zee \to [0,\infty)_{\arr}$ as follows.
\begin{align*} c(j) \coloneqq \begin{dcases} \frac{2}{M}\left( 1 - \frac{\lvert j \rvert}{2M} \right) & \mbox{ if } j \in [-2M+1,\, 2M-1]_{\zee} \\ \hspace{0.5 cm} 0 & \mbox{ else} \end{dcases} \end{align*}
For $m \in [0,M-1]_{\zee}$ we have $c(m) = (2M-m)/M^2 \geq (M+1)/M^2 \geq 1/M$, and since every summand on the right below is nonnegative we obtain:
\begin{align} \frac{1}{M} \sum_{m=0}^{M-1} \lvert \hat{\varsigma}(m) \rvert^2 \leq \sum_{j=-2M+1}^{2M-1} c(j)\, \lvert \hat{\varsigma}(j) \rvert^2 \label{lem.fouravg.p1} \end{align}
We define $G: \arr \to \arr$ by the equality on the left below, while the equality on the right follows by unwrapping the assignment of $c$ against the definition of the Fej\'er kernel in (\ref{olem.fejer.h1}).
\begin{align} G(u) \coloneqq \sum_{j=-2M+1}^{2M-1} c(j)\, e(ju) = \frac{2}{M}\, \mathtt{F}_{2M}(u) \label{lem.fouravg.p2} \end{align}
The function $G$ is $\zee$-periodic, and we regard it as a function on $\rzee$ by way of \ref{conv.period}. Combining (\ref{lem.fouravg.p2}) with Lemma \ref{olem.fejer} provides $G(u) \in [0,4]_{\arr}$ for all $u \in \arr$, together with the following bound for all $u \in \arr \setminus \zee$, in which we adopt the notation $\mathfrak{d}$ from Lemma \ref{olem.fejer}.
\begin{align} G(u) \leq \frac{2}{M} \cdot \frac{1}{8M\, \mathfrak{d}(u)^2} = \frac{1}{4M^2\, \mathfrak{d}(u)^2} \label{lem.fouravg.p3} \end{align}
We define $\Psi: \rzee \to \cee$ as follows.
\begin{align*} \Psi(y) \coloneqq \sum_{j=-2M+1}^{2M-1} c(j)\, \hat{\varsigma}(j)\, e(jy) \end{align*}
For $y \in \rzee$ we compute:
\begin{align} \Psi(y) & = \sum_{j=-2M+1}^{2M-1} c(j)\, e(jy)\, \varsigma(\mathfrak{e}^{-j}) \label{lem.fouravg.p4.1}
\\ & = \varsigma\bigl(x \mapsto G(y-x)\bigr) \label{lem.fouravg.p4.2} \end{align}
Here, the equality in (\ref{lem.fouravg.p4.1}) follows by unwrapping the assignment of $\Psi$, while (\ref{lem.fouravg.p4.2}) follows from (\ref{lem.fouravg.p4.1}) since unwrapping (\ref{lem.fouravg.p2}) provides:
\begin{align*} G(y-x) = \sum_{j=-2M+1}^{2M-1} c(j)\, e(jy)\, e(-jx) \end{align*}
Similarly, we compute:
\begin{align} \varsigma(\Psi) & = \sum_{j=-2M+1}^{2M-1} c(j)\, \hat{\varsigma}(j)\, \varsigma(\mathfrak{e}^{j}) \label{lem.fouravg.p5.1}
\\ & = \sum_{j=-2M+1}^{2M-1} c(j)\, \lvert \hat{\varsigma}(j) \rvert^2 \label{lem.fouravg.p5.2} \end{align}
Here, the equality in (\ref{lem.fouravg.p5.1}) follows by unwrapping the assignment of $\Psi$, while (\ref{lem.fouravg.p5.2}) follows from (\ref{lem.fouravg.p5.1}) since $\varsigma(\mathfrak{e}^{j}) = \ov{\hat{\varsigma}(j)}$.\\

Recalling from Definition \ref{def.entropy} and \ref{conv.period} that the cell $\mathscr{P}_M(k)$ corresponds to the interval $[k/M,(k+1)/M)_{\arr}$ of the coset transversal, for $k,\ell \in [0,M-1]_{\zee}$ we define $d_{k,\ell} \in \enn$ as follows.
\begin{align*} d_{k,\ell} \coloneqq \inf_{p \in \zee} \lvert \ell - k - pM \rvert \end{align*}
We also define $B: [0,M-1]_{\zee}^2 \to (0,\infty)_{\arr}$ as follows.
\begin{align*} B(k,\ell) \coloneqq \begin{cases} \hspace{0.35 cm} 4 & \mbox{ if } d_{k,\ell} \leq 1 \\ \hspace{0.2 cm} 1/d_{k,\ell}^{\,2} & \mbox{ if } d_{k,\ell} \geq 2 \end{cases} \end{align*}
Now let $k,\ell \in [0,M-1]_{\zee}$, let $x \in \mathscr{P}_M(k)$ and $y \in \mathscr{P}_M(\ell)$, and let $x_0$ and $y_0$ be the corresponding transversal representatives, so that $y_0 - x_0 = (\ell-k)/M + \theta$ for some $\theta \in (-1/M,1/M)_{\arr}$. We consider two cases according to the value of $d_{k,\ell}$.
\begin{itemize}
\item Suppose $d_{k,\ell} \leq 1$. Since $G(u) \in [0,4]_{\arr}$ for all $u \in \arr$, we find $G(y-x) \leq 4 = B(k,\ell)$.
\item Suppose $d_{k,\ell} \geq 2$. We compute:
\begin{align} \mathfrak{d}(y_0 - x_0) & \geq \frac{d_{k,\ell}}{M} - \lvert \theta \rvert \label{lem.fouravg.p6.1}
\\ & \geq \frac{d_{k,\ell}-1}{M} \label{lem.fouravg.p6.2}
\\ & \geq \frac{d_{k,\ell}}{2M} \label{lem.fouravg.p6.3} \end{align}
This computation may be justified as follows.
\begin{itemize}
\item[$\circ$] The inequality in (\ref{lem.fouravg.p6.1}) follows from the triangle inequality, since unwrapping the assignments of $\mathfrak{d}$ and $d_{k,\ell}$ provides $\mathfrak{d}\bigl((\ell-k)/M\bigr) = d_{k,\ell}/M$.
\item[$\circ$] (\ref{lem.fouravg.p6.2}) follows from (\ref{lem.fouravg.p6.1}) since $\lvert \theta \rvert < 1/M$.
\item[$\circ$] (\ref{lem.fouravg.p6.3}) follows from (\ref{lem.fouravg.p6.2}) since $d_{k,\ell} \geq 2$.
\end{itemize}
In particular $y_0 - x_0 \notin \zee$, and so combining (\ref{lem.fouravg.p3}) with (\ref{lem.fouravg.p6.3}) provides $G(y-x) \leq 1/d_{k,\ell}^{\,2} = B(k,\ell)$.\end{itemize}
In both cases we have obtained the following bound, in which $k$, $\ell$, $x$ and $y$ were arbitrary.
\begin{align} G(y-x) \leq B(k,\ell) \label{lem.fouravg.p6} \end{align}
We abbreviate $\varsigma_k \coloneqq \varsigma(\mathscr{P}_M(k))$ for $k \in [0,M-1]_{\zee}$, and we note that Definition \ref{def.entropy} provides that the cells $\mathscr{P}_M(0),\ldots,\mathscr{P}_M(M-1)$ partition $\rzee$. Temporarily fix $\ell \in [0,M-1]_{\zee}$ and $y \in \mathscr{P}_M(\ell)$. In the indicator notation of \ref{conv.indicator}, the following holds for all $x \in \rzee$.
\begin{align} G(y-x) \leq \sum_{k=0}^{M-1} B(k,\ell)\, \mathbf{1}_{\mathscr{P}_M(k)}(x) \label{lem.fouravg.p7.1} \end{align}
Here, each $x \in \rzee$ lies in exactly one cell, and for $x \in \mathscr{P}_M(k)$ the right side of (\ref{lem.fouravg.p7.1}) equals $B(k,\ell)$, so that (\ref{lem.fouravg.p7.1}) is an instance of (\ref{lem.fouravg.p6}).\\

Evaluating both sides of (\ref{lem.fouravg.p7.1}) against $\varsigma$ we obtain:
\begin{align} \Psi(y) \leq \sum_{k=0}^{M-1} B(k,\ell)\, \varsigma_k \label{lem.fouravg.p7.2} \end{align}
Here, 
\begin{itemize} \item the evaluation preserves the inequality in (\ref{lem.fouravg.p7.1}) since both sides are real-valued elements of $\bbr(\rzee)$; 
\item the left side of (\ref{lem.fouravg.p7.1}) evaluates to $\Psi(y)$ by the computation ending in (\ref{lem.fouravg.p4.2}); 
\item and the right side of (\ref{lem.fouravg.p7.1}) evaluates to the right side of (\ref{lem.fouravg.p7.2}) since $\varsigma(\mathbf{1}_{\mathscr{P}_M(k)}) = \varsigma_k$ for all $k \in [0,M-1]_{\zee}$. \end{itemize}
Since $\ell$ and $y$ were arbitrary, and again since the cells partition $\rzee$, the following holds for all $y \in \rzee$.
\begin{align} \Psi(y) \leq \sum_{\ell=0}^{M-1} \left( \sum_{k=0}^{M-1} B(k,\ell)\, \varsigma_k \right) \mathbf{1}_{\mathscr{P}_M(\ell)}(y) \label{lem.fouravg.p8.1} \end{align}
Here, each $y \in \rzee$ lies in exactly one cell, and for $y \in \mathscr{P}_M(\ell)$ the right side of (\ref{lem.fouravg.p8.1}) equals the right side of (\ref{lem.fouravg.p7.2}), so that (\ref{lem.fouravg.p8.1}) is an instance of (\ref{lem.fouravg.p7.2}).\\

Evaluating both sides of (\ref{lem.fouravg.p8.1}) against $\varsigma$ we obtain:
\begin{align} \sum_{j=-2M+1}^{2M-1} c(j)\, \lvert \hat{\varsigma}(j) \rvert^2 \leq \sum_{\ell=0}^{M-1} \sum_{k=0}^{M-1} B(k,\ell)\, \varsigma_k\, \varsigma_\ell \label{lem.fouravg.p8.2} \end{align}
Here, 
\begin{itemize} \item the left side of (\ref{lem.fouravg.p8.1}) is real-valued since the computation ending in (\ref{lem.fouravg.p4.2}) expresses each value $\Psi(y)$ as the evaluation of $\varsigma$ against a real-valued function, and so the evaluation preserves the inequality;
\item the left side of (\ref{lem.fouravg.p8.1}) evaluates to $\varsigma(\Psi)$, which equals the left side of (\ref{lem.fouravg.p8.2}) by the computation ending in (\ref{lem.fouravg.p5.2}); 
\item and the right side of (\ref{lem.fouravg.p8.1}) evaluates to the right side of (\ref{lem.fouravg.p8.2}) since $\varsigma(\mathbf{1}_{\mathscr{P}_M(\ell)}) = \varsigma_\ell$ for all $\ell \in [0,M-1]_{\zee}$.\end{itemize}
Now fix $k \in [0,M-1]_{\zee}$. For each $j \in \ennplus$, an element $\ell$ of $[0,M-1]_{\zee}$ satisfying $d_{k,\ell} = j$ is congruent to $k+j$ or to $k-j$ modulo $M$, and each congruence class modulo $M$ contains exactly one element of $[0,M-1]_{\zee}$; thus at most two elements $\ell$ of $[0,M-1]_{\zee}$ satisfy $d_{k,\ell} = j$. Similarly, $d_{k,\ell} = 0$ holds only for $\ell = k$. We compute:
\begin{align} \sum_{\ell=0}^{M-1} B(k,\ell) & \leq 12 + 2\sum_{j=2}^{\infty} \frac{1}{j^2} \label{lem.fouravg.p9.1}
\\ & \leq 10 + \frac{\pi^2}{3} \label{lem.fouravg.p9.2}
\\ & \leq 14 \label{lem.fouravg.p9.3} \end{align}
This computation may be justified as follows.
\begin{itemize}
\item[$\circ$] The inequality in (\ref{lem.fouravg.p9.1}) follows from the assignment of $B$, since at most three elements $\ell$ of $[0,M-1]_{\zee}$ satisfy $d_{k,\ell} \leq 1$, each contributing a summand of at most $4$, while for each $j \in [2,\infty)_{\zee}$ at most two elements $\ell$ satisfy $d_{k,\ell} = j$, each contributing $1/j^2$.
\item[$\circ$] (\ref{lem.fouravg.p9.2}) follows from (\ref{lem.fouravg.p9.1}) since the solution to the Basel problem provides:
\begin{align*} \sum_{j=1}^{\infty} \frac{1}{j^2} = \frac{\pi^2}{6} \end{align*}
\item[$\circ$] (\ref{lem.fouravg.p9.3}) follows from (\ref{lem.fouravg.p9.2}) since $\pi^2 \leq 12$.
\end{itemize}
We compute:
\begin{align} \sum_{\ell=0}^{M-1} \sum_{k=0}^{M-1} B(k,\ell)\, \varsigma_k\, \varsigma_\ell & \leq \sum_{k=0}^{M-1} \varsigma_k^2 \sum_{\ell=0}^{M-1} B(k,\ell) \label{lem.fouravg.p10.1}
\\ & \leq 14 \sum_{k=0}^{M-1} \varsigma_k^2 \label{lem.fouravg.p10.2}
\\ & = 14\, \bigl \vert \bigl \vert \varsigma(\mathscr{P}_M) \bigr \vert \bigr \vert_2^2 \label{lem.fouravg.p10.3} \end{align}
This computation may be justified as follows.
\begin{itemize}
\item[$\circ$] The inequality in (\ref{lem.fouravg.p10.1}) follows by applying the elementary AGM inequality $2\varsigma_k\varsigma_\ell \leq \varsigma_k^2 + \varsigma_\ell^2$ to each summand and identifying the two resulting double sums by way of the symmetry $B(k,\ell) = B(\ell,k)$.
\item[$\circ$] (\ref{lem.fouravg.p10.2}) follows from (\ref{lem.fouravg.p10.1}) by the computation ending in (\ref{lem.fouravg.p9.3}), in which $k$ was arbitrary.
\item[$\circ$] (\ref{lem.fouravg.p10.3}) follows from (\ref{lem.fouravg.p10.2}) directly from Definition \ref{def.mupnorm}.
\end{itemize}
It remains to return from $\varsigma$ to $\nu$. For all $\ell \in [0,M-1]_{\zee}$ we compute:
\begin{align} \varsigma_\ell & = \nu\bigl(\mathtt{T}_g^{-1}(\mathscr{P}_M(\ell))\bigr) \label{lem.fouravg.p11.1}
\\ & = \sum_{j=0}^{g-1} \nu\bigl(\mathscr{P}_{gM}(\ell + jM)\bigr) \label{lem.fouravg.p11.2} \end{align}
This computation may be justified as follows.
\begin{itemize}
\item[$\circ$] The equality in (\ref{lem.fouravg.p11.1}) follows from Clause \ref{rem.endochangeofvar.c3} of Remark \ref{rem.endochangeofvar}, applied to the assignment $\varsigma = g_{\,\circledcirc}\, \nu$.
\item[$\circ$] (\ref{lem.fouravg.p11.2}) follows from (\ref{lem.fouravg.p11.1}) since directly from \ref{conv.Tmap} and Definition \ref{def.entropy} we have the disjoint union:
\begin{align*} \mathtt{T}_g^{-1}\bigl(\mathscr{P}_M(\ell)\bigr) = \bigcup_{j=0}^{g-1} \mathscr{P}_{gM}(\ell + jM) \end{align*}
\end{itemize}
Applying the Cauchy-Schwartz inequality to the lists $(1)_{j=0}^{g-1}$ and $\bigl(\nu(\mathscr{P}_{gM}(\ell+jM))\bigr)_{j=0}^{g-1}$, the computation ending in (\ref{lem.fouravg.p11.2}) provides the following for all $\ell \in [0,M-1]_{\zee}$.
\begin{align} \varsigma_\ell^{\,2} \leq g \sum_{j=0}^{g-1} \nu\bigl(\mathscr{P}_{gM}(\ell + jM)\bigr)^2 \label{lem.fouravg.p12.1} \end{align}
Summing (\ref{lem.fouravg.p12.1}) over $\ell \in [0,M-1]_{\zee}$ we obtain:
\begin{align} \bigl \vert \bigl \vert \varsigma(\mathscr{P}_M) \bigr \vert \bigr \vert_2^2 \leq g\, \bigl \vert \bigl \vert \nu(\mathscr{P}_{gM}) \bigr \vert \bigr \vert_2^2 \label{lem.fouravg.p12.2} \end{align}
Here, the left side is $\sum_{\ell=0}^{M-1} \varsigma_\ell^{\,2}$ directly from Definition \ref{def.mupnorm}, while the right side collects the summands by way of the bijection $(j,\ell) \mapsto \ell + jM$ from $[0,g-1]_{\zee} \times [0,M-1]_{\zee}$ to $[0,gM-1]_{\zee}$ together with Definition \ref{def.mupnorm}. \\

We now compute:
\begin{align} \frac{1}{M} \sum_{m=0}^{M-1} \lvert \hat{\nu}(mg) \rvert^2 & = \frac{1}{M} \sum_{m=0}^{M-1} \lvert \hat{\varsigma}(m) \rvert^2 \label{lem.fouravg.p13.1}
\\ & \leq 14\, \bigl \vert \bigl \vert \varsigma(\mathscr{P}_M) \bigr \vert \bigr \vert_2^2 \label{lem.fouravg.p13.2}
\\ & \leq 14\, g\, \bigl \vert \bigl \vert \nu(\mathscr{P}_{gM}) \bigr \vert \bigr \vert_2^2 \label{lem.fouravg.p13.3}
\\ & \leq 14\, g\, \bigl \vert \bigl \vert \nu(\mathscr{P}_{M}) \bigr \vert \bigr \vert_2^2 \label{lem.fouravg.p13.4}
\\ & \leq 32\, \lvert g \rvert\, \bigl \vert \bigl \vert \nu(\mathscr{P}_{M}) \bigr \vert \bigr \vert_2^2 \label{lem.fouravg.p13.5} \end{align}
This computation may be justified as follows.
\begin{itemize}
\item[$\circ$] The equality in (\ref{lem.fouravg.p13.1}) holds since $\hat{\varsigma}(m) = \hat{\nu}(mg)$ for all $m \in \zee$, as established at the opening of the current proof.
\item[$\circ$] (\ref{lem.fouravg.p13.2}) follows from (\ref{lem.fouravg.p13.1}) by combining (\ref{lem.fouravg.p1}) with (\ref{lem.fouravg.p8.2}) and with the computation ending in (\ref{lem.fouravg.p10.3}).
\item[$\circ$] (\ref{lem.fouravg.p13.3}) follows from (\ref{lem.fouravg.p13.2}) by (\ref{lem.fouravg.p12.2}).
\item[$\circ$] (\ref{lem.fouravg.p13.4}) follows from (\ref{lem.fouravg.p13.3}) by Remark \ref{rem.mupnormref}, applied with the resolutions $M$ and $gM$, where we note that $M$ divides $gM$.
\item[$\circ$] (\ref{lem.fouravg.p13.5}) follows from (\ref{lem.fouravg.p13.4}) since we have assumed $g \in \ennplus$, and where we choose to replace the constant $14$ by the larger constant $32$ for notational convenience in the sequel.
\end{itemize}
The computation ending in (\ref{lem.fouravg.p13.5}) verifies the desired inequality (\ref{lem.fouravg.conc1}) and thereby completes the proof of Lemma \ref{lem.fouravg}.
\end{Proof}

\subsection{Flat components at a good scale}

Lemma \ref{lem.flatscale} is a version of Lemma 4.7 in \cite{MR2563089}, which those authors deduce from their Lemma 3.6. Their statement carries a single parameter, and that parameter does three jobs at once: it is the fraction of the available digits which the argument declines to look at, it is the mass surrendered in passing from the measure to the flat component extracted from it, and it is the exponent governing how flat that component is.\\

For the purposes in BLMV the three parameters may be identified without loss, but doing so would cause substantial deterioration in our quantitative estimates. We therefore carry the block fraction and the surrendered mass as $\beta$ and the flatness exponent as $\gamma$, obtaining independent control of the two, at the cost of one further hypothesis relating them to the entropy level.\\

Since we need the lemma in a form their statement does not provide, we give the proof in full rather than importing it.\\

The dependencies in the statement of Lemma \ref{lem.flatscale} are as follows.
\begin{deps}
\dep{Conventions}{\depitem{Submeasures}{conv.submeasure}}
\dep{Definitions}{\depitem{Resolution entropy}{def.entropy}}
\depcont{\depitem{Resolution norms of measures}{def.mupnorm}}\end{deps}

We now move on to the formal statement.

\begin{Lemma}[Flat components at a good scale] \label{lem.flatscale} \quad
\vspace{-1 \baselineskip}
\begin{enumerate}[label=\upshape(\arabic*)] \item \label{lem.flatscale.c1} For all:
\begin{enumerate}[label=\upshape(\arabic*)] \item \label{lem.flatscale.c1.1} $n \in \enn$ and all $\sigma \in (0,1]_{\arr}$ and all $\beta \in (0,1]_{\arr}$ and all $\gamma \in (0,1/4]_{\arr}$ such that $\sigma > \beta + 2\gamma$,
\item \label{lem.flatscale.c1.2} and all $\ell \in [2/\gamma, \beta n]_{\zee}$,
\item \label{lem.flatscale.c1.3} and all $\mu \in \prob(\rzee)$ such that $\achn(\mu,2^n) \geq \sigma$,
\end{enumerate}
\item \label{lem.flatscale.c2} there exist:
\begin{enumerate}[label=\upshape(\arabic*)] \item \label{lem.flatscale.c2.1} $s \in [\beta n - \ell , n - \ell ]_{\zee}$,
\item \label{lem.flatscale.c2.2} and elements $w_0,\ldots,w_{2^s-1}$ of $[0,\infty)_{\arr}$ and elements $\nu_0,\ldots,\nu_{2^s-1}$ of $\prob(\rzee)$, \end{enumerate}
\end{enumerate}
such that:
\begin{enumerate}[resume, label=\upshape(\arabic*)] \item \label{lem.flatscale.c3} we have $\sum_{k=0}^{2^s-1}w_k \geq \frac{\sigma - \beta - 2\gamma}{1-2\gamma}$,
 \item \label{lem.flatscale.c4} and we have $\sum_{k=0}^{2^s-1}w_k \nu_k \preceq \mu$,
\item \label{lem.flatscale.c5} and for all $k \in [0,2^s-1]_{\zee}$ we have $\nu_k(\mathscr{P}_{2^s}(k)) = 1$,
\item \label{lem.flatscale.c6} and we have:
\begin{align*}  \sum_{k=0}^{2^s-1} w_k \, \bigl \vert \bigl \vert \nu_k(\mathscr{P}_{2^{s+\ell}}) \bigr \vert \bigr \vert _2^2 \leq \frac{1}{2^{\gamma \ell}}\end{align*}
  \end{enumerate}\end{Lemma}
  
An outline of the proof of Lemma \ref{lem.flatscale} is as follows.
\begin{itemize}
\item Splitting $\mu$ over the cells of one resolution exhibits its entropy a block of $\ell$ digits lower as the entropy at the coarser resolution plus the cell masses weighting the fibre entropies. This is the identity ending at (\ref{lem.flatscale.3.5}), and it plays the role of the conditional-entropy expansion opening the proof of Lemma 3.6 in \cite{MR2563089}.
\item Telescoping downward from $n$ in blocks of length $\ell$, stopping above the bottom $\beta n$ digits, leaves entropy $(\sigma - \beta)n\log(2)$ across the $M$ blocks. This is the choice of scale in the first half of Lemma 3.6 of \cite{MR2563089}; telescoping downward keeps the block boundaries integral.
\item Each fibre at that scale is cut at $2^{-\gamma\ell}$. Below the cut, comparison against the uniform distribution of the same mass bounds the contribution by $y_k \ell \log(2) + 1 - y_k$; above it the cut itself gives $\gamma\ell\log(2)(1-y_k)$. Averaging against the cell masses converts the entropy bound into one on the mass above the cut. The cut is that of Lemma 2.7 of \cite{MR2563089}, which their Lemma 3.6 invokes, and $\mathscr{Y}_k$ and $\mathscr{Z}_k$ are the classes it separates.
\item Rearranging spends $\ell \geq 2/\gamma$ at (\ref{lem.flatscale.19.0}), converting the additive $1$ into a second $\gamma$ so that the loss is proportional to the entropy deficit. This is why $\gamma$ is carried apart from $\beta$, and has no counterpart in \cite{MR2563089}, where the two coincide and the retained mass is bounded by a difference rather than a ratio.
\item The measures demanded are the cells below the cut, renormalized. Their normalizing factors cancel against the weights at (\ref{lem.flatscale.22.1}) and (\ref{lem.flatscale.23.1}), leaving sums of cell masses equal to $1$.\end{itemize}
The dependencies in the proof of Lemma \ref{lem.flatscale} beyond those in the statement are as follows.
\begin{deps}
\dep{Imports}{\depitem{Range of resolution entropy}{olem.entrange}}\end{deps}

We now move on to the formal proof.
  
\begin{Proof}[of Lemma \ref{lem.flatscale}] Let $n,\sigma,\beta,\gamma,\ell$ and $\mu$ be as in Clause \ref{lem.flatscale.c1} of Lemma \ref{lem.flatscale}. Clause \ref{lem.flatscale.c1}\ref{lem.flatscale.c1.2} provides $2/\gamma \leq \ell \leq \beta n$, so that $n \geq 2/(\beta\gamma)$ and in particular $n \in \ennplus$. \\

For $s \in [0,n]_{\zee}$ and $k \in [0,2^s-1]_{\zee}$ we define $w^s_k \in [0,1]_{\arr}$ by $w^s_k \coloneqq \mu(\mathscr{P}_{2^s}(k))$, and we define $\mathscr{X}_s$ to be the subset of $[0,2^s-1]_{\zee}$ given as follows.
\begin{align*} \mathscr{X}_s \coloneqq \{k \in [0,2^s-1]_\zee : w_k^s >0 \} \end{align*}
We note that the following is immediate from the definition of resolution entropy in Definition \ref{def.entropy}.
\begin{align} \label{lem.flatscale.0} \ach\bigl(\mu,2^{s}\bigr) = - \sum_{k \in \mathscr{X}_s} w^s_k \log\bigl(w^s_k\bigr)  \end{align}
For $k \in \mathscr{X}_s$ we define $\nu^s_k \in \prob(\rzee)$ by stipulating that $\nu^s_k(E) \coloneqq \mu\bigl(E \cap \mathscr{P}_{2^s}(k)\bigr)/w^s_k$ for every Borel subset $E$ of $\rzee$.\\

Let $s \in [0,n-\ell]_{\zee}$. Definition \ref{def.entropy} identifies $\mathscr{P}_{2^{s+\ell}}(j)$ with $[j2^{-s-\ell},(j+1)2^{-s-\ell})_{\arr}$ and identifies $\mathscr{P}_{2^{s}}(k)$ with $[k2^{-s},(k+1)2^{-s})_{\arr}$, so that for all $j \in [0,2^{s+\ell}-1]_{\zee}$ and all $k \in [0,2^s-1]_{\zee}$ we have the following.
\begin{align} \label{lem.flatscale.1} \mathscr{P}_{2^{s+\ell}}(j) \cap \mathscr{P}_{2^{s}}(k) = \begin{dcases} \mathscr{P}_{2^{s+\ell}}(j) & \mbox{ if } j \in [k2^\ell,(k+1)2^\ell-1]_{\zee} 
\\ \emptyset & \mbox{ if } j \in [0,2^{s+\ell}-1]_{\zee} \setminus  [k2^\ell,(k+1)2^\ell-1]_{\zee}\end{dcases} \end{align}
Combining the above dichotomy with the assignment of $\nu^s_k$ we obtain the next display for all $k \in \mathscr{X}_s$.
\begin{align} \label{lem.flatscale.1.1} \nu^s_k\bigl(\mathscr{P}_{2^{s+\ell}}(j)\bigr)  = \begin{dcases} \mu\bigl(\mathscr{P}_{2^{s+\ell}}(j)\bigr)/ w^s_k & \mbox{ if } j \in [k2^\ell,(k+1)2^\ell-1]_{\zee} 
\\ 0 & \mbox{ if } j \in [0,2^{s+\ell}-1]_{\zee} \setminus  [k2^\ell,(k+1)2^\ell-1]_{\zee}\end{dcases} \end{align}
For each fixed $k \in \mathscr{X}_s$ the cells $\mathscr{P}_{2^{s+\ell}}(j)$ for $j \in [k2^\ell,(k+1)2^\ell-1]_{\zee}$ are pairwise disjoint with union $\mathscr{P}_{2^{s}}(k)$. Combining this observation with the assignment of $\nu^s_k$ we find:
\begin{align} \sum_{j=k2^\ell}^{(k+1)2^\ell-1} \nu^s_k\bigl(\mathscr{P}_{2^{s+\ell}}(j)\bigr) = 1 \label{lem.flatscale.2} \end{align}
Combining the lower case of (\ref{lem.flatscale.1.1}) with the definition of resolution entropy in Definition \ref{def.entropy} we obtain the following for all $k \in \mathscr{X}_s$.
\begin{align} \label{lem.flatscale.2.1} \ach\bigl(\nu^s_k , 2^{s+\ell}\bigr)=- \sum_{j=k2^\ell}^{(k+1)2^\ell-1} \nu^s_k\bigl(\mathscr{P}_{2^{s+\ell}}(j)\bigr) \log\Bigl(\nu^s_k\bigl(\mathscr{P}_{2^{s+\ell}}(j)\bigr)\Bigr) \end{align}
We compute:
\begin{align} \ach\bigl(\mu,2^{s+\ell}\bigr) & = - \sum_{j=0}^{2^{s+\ell}-1} \mu\bigl(\mathscr{P}_{2^{s+\ell}}(j)\bigr) \log\Bigl(\mu\bigl(\mathscr{P}_{2^{s+\ell}}(j)\bigr)\Bigr) \label{lem.flatscale.3.1}
\\ & = - \sum_{k =0}^{2^s-1} \sum_{j=k2^\ell}^{(k+1)2^\ell-1} \mu\bigl(\mathscr{P}_{2^{s+\ell}}(j)\bigr) \log\Bigl(\mu\bigl(\mathscr{P}_{2^{s+\ell}}(j)\bigr)\Bigr) \label{lem.flatscale.3.1.1}
\\ & = - \sum_{k \in \mathscr{X}_s} \sum_{j=k2^\ell}^{(k+1)2^\ell-1} \mu\bigl(\mathscr{P}_{2^{s+\ell}}(j)\bigr) \log\Bigl(\mu\bigl(\mathscr{P}_{2^{s+\ell}}(j)\bigr)\Bigr) \label{lem.flatscale.3.2}
\\ & = - \sum_{k \in \mathscr{X}_s} \sum_{j=k2^\ell}^{(k+1)2^\ell-1} w^s_k \, \nu^s_k\bigl(\mathscr{P}_{2^{s+\ell}}(j)\bigr) \Bigl[ \log\bigl(w^s_k\bigr) + \log\Bigl(\nu^s_k\bigl(\mathscr{P}_{2^{s+\ell}}(j)\bigr)\Bigr) \Bigr] \label{lem.flatscale.3.3}
\\ & = - \sum_{k \in \mathscr{X}_s} w^s_k \log\bigl(w^s_k\bigr) - \sum_{k \in \mathscr{X}_s} w^s_k \sum_{j=k2^\ell}^{(k+1)2^\ell-1} \nu^s_k\bigl(\mathscr{P}_{2^{s+\ell}}(j)\bigr) \log\Bigl(\nu^s_k\bigl(\mathscr{P}_{2^{s+\ell}}(j)\bigr)\Bigr) \label{lem.flatscale.3.4}
\\ & = \ach\bigl(\mu,2^{s}\bigr) + \sum_{k \in \mathscr{X}_s} w^s_k \, \ach\bigl(\nu^s_k , 2^{s+\ell}\bigr) \label{lem.flatscale.3.5} \end{align}
This computation may be justified as follows.
\begin{itemize}
	\item The equality in (\ref{lem.flatscale.3.1}) follows from the definition of the resolution entropy in (\ref{def.entropy.0}), applied at resolution $2^{s+\ell}$.
	\item (\ref{lem.flatscale.3.1.1}) follows from (\ref{lem.flatscale.3.1}) since the sets $[k2^\ell,(k+1)2^\ell-1]_{\zee}$ for $k \in [0,2^s-1]_{\zee}$ are pairwise disjoint with union $[0,2^{s+\ell}-1]_{\zee}$.
		\item (\ref{lem.flatscale.3.2}) follows from (\ref{lem.flatscale.3.1.1}) since for $k \in [0,2^s-1]_{\zee} \setminus \mathscr{X}_s$ and $j \in [k2^\ell,(k+1)2^\ell-1]_{\zee}$ the upper case of (\ref{lem.flatscale.1}) exhibits $\mathscr{P}_{2^{s+\ell}}(j)$ as a subset of the $\mu$-null set $\mathscr{P}_{2^s}(k)$.
	\item (\ref{lem.flatscale.3.3}) follows from (\ref{lem.flatscale.3.2}) by (\ref{lem.flatscale.1.1}).	
	\item (\ref{lem.flatscale.3.4}) follows from (\ref{lem.flatscale.3.3}) by (\ref{lem.flatscale.2}).
	\item (\ref{lem.flatscale.3.5}) follows from (\ref{lem.flatscale.3.4}) by combining (\ref{lem.flatscale.0}) and (\ref{lem.flatscale.2.1}).
\end{itemize}
Since $s \in [0,n-\ell]_{\zee}$ was arbitrary, the computation ending in (\ref{lem.flatscale.3.5}) holds for every $s \in [0,n-\ell]_{\zee}$. \\

We define $M \in \ennplus$ and, for $u \in [1,M]_{\zee}$, an element $s_u$ of $\zee$ as follows.
\begin{align} M & \coloneqq \left\lfloor \frac{(1-\beta)n}{\ell} \right\rfloor + 1 & s_u & \coloneqq n - u\ell \label{lem.flatscale.4} \end{align}
We compute:
\begin{align} n - M\ell & = n - \ell \left\lfloor \frac{(1-\beta)n}{\ell} \right\rfloor - \ell \label{lem.flatscale.5.1}
\\ & \geq n - (1-\beta)n - \ell \label{lem.flatscale.5.2}
\\ & = \beta n - \ell \label{lem.flatscale.5.3}
\\ & \geq 0 \label{lem.flatscale.5.4} \end{align}
This computation may be justified as follows.
\begin{itemize}
	\item The equality in (\ref{lem.flatscale.5.1}) follows from the assignment of $M$ at the left of (\ref{lem.flatscale.4}).
	\item (\ref{lem.flatscale.5.2}) follows from (\ref{lem.flatscale.5.1}) since $\lfloor x \rfloor \leq x$ for all $x \in \arr$ and since $\ell \in \ennplus$.
	\item (\ref{lem.flatscale.5.4}) follows from (\ref{lem.flatscale.5.3}) since Clause \ref{lem.flatscale.c1}\ref{lem.flatscale.c1.2} of our current hypotheses ensures $\ell \leq \beta n$.
\end{itemize}
Since the assignment of $s_u$ at the right of (\ref{lem.flatscale.4}) is decreasing in $u$ with $s_1 = n - \ell$ and $s_M = n - M\ell$, the computation ending in (\ref{lem.flatscale.5.4}) provides the following for all $u \in [1,M]_\zee$.
\begin{align}  \label{lem.flatscale.6.0} s_u \in [\beta n - \ell, n - \ell]_{\zee} \end{align}
We compute again:
\begin{align} M & = \left\lfloor \frac{(1-\beta)n}{\ell} \right\rfloor + 1 \label{lem.flatscale.6.1}
\\ & \leq \frac{(1-\beta)n + \ell}{\ell} \label{lem.flatscale.6.2}
\\ & \leq \frac{n}{\ell} \label{lem.flatscale.6.3} \end{align}
This computation may be justified as follows.
\begin{itemize}
	\item The equality in (\ref{lem.flatscale.6.1}) follows from the assignment of $M$ at the left of (\ref{lem.flatscale.4}).
	\item (\ref{lem.flatscale.6.2}) follows from (\ref{lem.flatscale.6.1}) since $\lfloor x \rfloor \leq x$ for all $x \in \arr$.
	\item (\ref{lem.flatscale.6.3}) follows from (\ref{lem.flatscale.6.2}) since Clause \ref{lem.flatscale.c1}\ref{lem.flatscale.c1.2} provides $\ell \leq \beta n$.
\end{itemize}
We compute again:
\begin{align} \sum_{u=1}^{M} \sum_{k \in \mathscr{X}_{s_u}} w^{s_u}_k \, \ach\bigl(\nu^{s_u}_k , 2^{s_u+\ell}\bigr) & = \sum_{u=1}^{M} \Bigl( \ach\bigl(\mu,2^{s_u+\ell}\bigr) - \ach\bigl(\mu,2^{s_u}\bigr) \Bigr) \label{lem.flatscale.7.1}
\\ & = \ach\bigl(\mu,2^{n}\bigr) - \ach\bigl(\mu,2^{n-M\ell}\bigr) \label{lem.flatscale.7.2}
\\ & \geq \sigma n \log(2) - \ach\bigl(\mu,2^{n-M\ell}\bigr) \label{lem.flatscale.7.3}
\\ & \geq \sigma n \log(2) - (n - M\ell)\log(2) \label{lem.flatscale.7.4}
\\ & \geq \sigma n \log(2) - \beta n \log(2) \label{lem.flatscale.7.5}
\\ & = (\sigma - \beta) \, n \log(2) \label{lem.flatscale.7.6} \end{align}
This computation may be justified as follows.
\begin{itemize}
	\item The equality in (\ref{lem.flatscale.7.1}) follows by instantiating the computation ending in (\ref{lem.flatscale.3.5}) with the scale $s_u$ for each $u \in [1,M]_{\zee}$, which is legitimate since (\ref{lem.flatscale.6.0}) implies $s_u \in [0,n-\ell]_{\zee}$.
	\item (\ref{lem.flatscale.7.2}) follows from (\ref{lem.flatscale.7.1}) since the assignment of $s_u$ at the right of (\ref{lem.flatscale.4}) provides $s_u + \ell = n - (u-1)\ell$, so that the sum telescopes between the resolutions $2^n$ and $2^{n-M\ell}$.
	\item (\ref{lem.flatscale.7.3}) follows from (\ref{lem.flatscale.7.2}) by combining Clause \ref{lem.flatscale.c1}\ref{lem.flatscale.c1.3} of our current Lemma \ref{lem.flatscale} with the definition of the normalized resolution entropy in (\ref{def.entropy-1}), which provides $\ach(\mu,2^n) = \achn(\mu,2^n)\, n \log(2)$.
	\item (\ref{lem.flatscale.7.4}) follows from (\ref{lem.flatscale.7.3}) by instantiating Lemma \ref{olem.entrange} with the resolution $2^{n-M\ell}$, which is an element of $\ennplus$ by the computation ending in (\ref{lem.flatscale.5.4}).
	\item (\ref{lem.flatscale.7.5}) follows from (\ref{lem.flatscale.7.4}) since the assignment of $M$ at the left of (\ref{lem.flatscale.4}) provides $M\ell > (1-\beta)n$ and hence $n - M\ell \leq \beta n$.
\end{itemize}
Since the outer sum at the left of (\ref{lem.flatscale.7.1}) has $M$ terms, the computation ending in (\ref{lem.flatscale.7.6}) provides $u \in [1,M]_{\zee}$ such that the following holds.
\begin{align} \sum_{k \in \mathscr{X}_{s_u}} w^{s_u}_k \, \ach\bigl(\nu^{s_u}_k , 2^{s_u+\ell}\bigr) \geq \frac{(\sigma - \beta)\, n \log(2)}{M} \label{lem.flatscale.8} \end{align}
Fixing such a $u$ and noting the bounds in (\ref{lem.flatscale.6.0}), we assign the parameter $s$ demanded by Clause \ref{lem.flatscale.c2}\ref{lem.flatscale.c2.1} of our current Lemma \ref{lem.flatscale} as $s \coloneqq s_u$. We compute:
\begin{align} \sum_{k \in \mathscr{X}_{s}} w^{s}_k \, \ach\bigl(\nu^{s}_k , 2^{s+\ell}\bigr) & \geq \frac{(\sigma - \beta)\, n \log(2)}{M} \label{lem.flatscale.9.1}
\\ & \geq (\sigma - \beta)\, \ell \log(2) \label{lem.flatscale.9.3} \end{align}
This computation may be justified as follows.
\begin{itemize}
	\item The inequality in (\ref{lem.flatscale.9.1}) is (\ref{lem.flatscale.8}) under the assignment $s \coloneqq s_u$.
	\item (\ref{lem.flatscale.9.3}) follows from (\ref{lem.flatscale.9.1}) by the computation ending in (\ref{lem.flatscale.6.3}), since Clause \ref{lem.flatscale.c1}\ref{lem.flatscale.c1.1} provides $\sigma - \beta > 0$.\end{itemize}
For $k \in \mathscr{X}_s$ we define $\mathscr{Y}_k$ to be the subset of $[k2^\ell,(k+1)2^\ell-1]_{\zee}$ given as follows.
\begin{align} \mathscr{Y}_k \coloneqq \Bigl\{ j \in [k2^\ell,(k+1)2^\ell-1]_{\zee} \; : \; 0 < \nu^s_k\bigl(\mathscr{P}_{2^{s+\ell}}(j)\bigr) < 2^{-\gamma \ell} \Bigr\} \label{lem.flatscale.10} \end{align}
We also define $y_k \in [0,1]_{\arr}$ as follows, the upper bound holding by (\ref{lem.flatscale.2}).
\begin{align} y_k \coloneqq \sum_{j \in \mathscr{Y}_k} \nu^s_k\bigl(\mathscr{P}_{2^{s+\ell}}(j)\bigr) \label{lem.flatscale.11} \end{align}
Let $k \in \mathscr{X}_s$ and assume for the moment that $\mathscr{Y}_k$ is nonempty, so that the assignment of $y_k$ in (\ref{lem.flatscale.11}) provides $y_k > 0$. For $j \in \mathscr{Y}_k$, we can instantiate the elementary inequality $\log(x) \leq x-1$ with the assignment $x \coloneqq y_k 2^{-\ell}/\nu^s_k\bigl(\mathscr{P}_{2^{s+\ell}}(j)\bigr)$ to obtain the following.
\begin{align} - \nu^s_k\bigl(\mathscr{P}_{2^{s+\ell}}(j)\bigr) \log\left( \frac{2^{\ell}\, \nu^s_k\bigl(\mathscr{P}_{2^{s+\ell}}(j)\bigr)}{y_k} \right) \leq \frac{y_k}{2^{\ell}} - \nu^s_k\bigl(\mathscr{P}_{2^{s+\ell}}(j)\bigr) \label{lem.flatscale.12.0} \end{align}
We compute:
\begin{align} \begin{split} - \sum_{j \in \mathscr{Y}_k} \nu^s_k\bigl(\mathscr{P}_{2^{s+\ell}}(j)\bigr)\log\Bigl(\nu^s_k\bigl(\mathscr{P}_{2^{s+\ell}}(j)\bigr)\Bigr) & = y_k \ell \log(2) + y_k \log\left(\frac{1}{y_k}\right)
\\ & \hspace{0.4 in} - \sum_{j \in \mathscr{Y}_k} \nu^s_k\bigl(\mathscr{P}_{2^{s+\ell}}(j)\bigr) \log\left( \frac{2^{\ell}\, \nu^s_k\bigl(\mathscr{P}_{2^{s+\ell}}(j)\bigr)}{y_k} \right) \end{split} \label{lem.flatscale.12.1}
\\ & \leq y_k \ell \log(2) + y_k \log\left(\frac{1}{y_k}\right) + \sum_{j \in \mathscr{Y}_k} \left( \frac{y_k}{2^{\ell}} - \nu^s_k\bigl(\mathscr{P}_{2^{s+\ell}}(j)\bigr) \right) \label{lem.flatscale.12.2}
\\ & = y_k \ell \log(2) + y_k \log\left(\frac{1}{y_k}\right) + \frac{\lvert \mathscr{Y}_k \rvert \, y_k}{2^{\ell}} - y_k \label{lem.flatscale.12.3}
\\ & \leq y_k \ell \log(2) + y_k \log\left(\frac{1}{y_k}\right) \label{lem.flatscale.12.4}
\\ & \leq y_k \ell \log(2) + 1 - y_k \label{lem.flatscale.12.5} \end{align}
This computation may be justified as follows.
\begin{itemize}
	\item The equality in (\ref{lem.flatscale.12.1}) follows by combining the assignment of $y_k$ in (\ref{lem.flatscale.11}) with the observation that $\log\bigl(2^\ell x / y_k \bigr) = \log(x) + \ell\log(2) - \log(y_k)$ for all $x \in (0,\infty)_{\arr}$.
	\item (\ref{lem.flatscale.12.2}) follows from (\ref{lem.flatscale.12.1}) by summing (\ref{lem.flatscale.12.0}) over $j \in \mathscr{Y}_k$.
	\item (\ref{lem.flatscale.12.3}) follows from (\ref{lem.flatscale.12.2}) by the assignment of $y_k$ in (\ref{lem.flatscale.11}).
	\item (\ref{lem.flatscale.12.4}) follows from (\ref{lem.flatscale.12.3}) since (\ref{lem.flatscale.10}) exhibits $\mathscr{Y}_k$ as a subset of $[k2^\ell,(k+1)2^\ell-1]_{\zee}$ and hence provides $\lvert \mathscr{Y}_k \rvert \leq 2^\ell$.
	\item (\ref{lem.flatscale.12.5}) follows from (\ref{lem.flatscale.12.4}) by instantiating the inequality $\log(x) \leq x-1$ with $x \coloneqq 1/y_k$.\end{itemize}
If instead $\mathscr{Y}_k$ is empty then $y_k = 0$, so that the left side of (\ref{lem.flatscale.12.1}) is an empty sum while the right side of (\ref{lem.flatscale.12.5}) is $1$. In either case we therefore have the following for all $k \in \mathscr{X}_s$.
\begin{align} - \sum_{j \in \mathscr{Y}_k} \nu^s_k\bigl(\mathscr{P}_{2^{s+\ell}}(j)\bigr)\log\Bigl(\nu^s_k\bigl(\mathscr{P}_{2^{s+\ell}}(j)\bigr)\Bigr) \leq y_k \ell \log(2) + 1 - y_k \label{lem.flatscale.12.6} \end{align}
We define $\mathscr{Z}_k$ to be the subset of $[k2^\ell,(k+1)2^\ell-1]_{\zee}$ given as follows.
\begin{align} \mathscr{Z}_k \coloneqq [k2^\ell,(k+1)2^\ell-1]_{\zee} \setminus \mathscr{Y}_k \label{lem.flatscale.13} \end{align}
By the assignment of $\mathscr{Y}_k$ in (\ref{lem.flatscale.10}), every $j \in \mathscr{Z}_k$ with $\nu^s_k\bigl(\mathscr{P}_{2^{s+\ell}}(j)\bigr) > 0$ satisfies $\nu^s_k\bigl(\mathscr{P}_{2^{s+\ell}}(j)\bigr) \geq 2^{-\gamma\ell}$ and hence $-\log\Bigl(\nu^s_k\bigl(\mathscr{P}_{2^{s+\ell}}(j)\bigr)\Bigr) \leq \gamma \ell \log(2)$. Multiplying by $\nu^s_k\bigl(\mathscr{P}_{2^{s+\ell}}(j)\bigr)$, we obtain the following for all $j \in \mathscr{Z}_k$.
\begin{align} - \nu^s_k\bigl(\mathscr{P}_{2^{s+\ell}}(j)\bigr) \log\Bigl( \nu^s_k\bigl(\mathscr{P}_{2^{s+\ell}}(j)\bigr) \Bigr) \leq \gamma \ell \log(2) \, \nu^s_k\bigl(\mathscr{P}_{2^{s+\ell}}(j)\bigr) \label{lem.flatscale.14.0} \end{align}
We compute:
\begin{align} - \sum_{j \in \mathscr{Z}_k} \nu^s_k\bigl(\mathscr{P}_{2^{s+\ell}}(j)\bigr) \log\Bigl( \nu^s_k\bigl(\mathscr{P}_{2^{s+\ell}}(j)\bigr) \Bigr) & \leq \gamma \ell \log(2) \sum_{j \in \mathscr{Z}_k} \nu^s_k\bigl(\mathscr{P}_{2^{s+\ell}}(j)\bigr) \label{lem.flatscale.14.1}
\\ & = \gamma \ell \log(2) \left( \sum_{j=k2^\ell}^{(k+1)2^\ell-1} \nu^s_k\bigl(\mathscr{P}_{2^{s+\ell}}(j)\bigr) - \sum_{j \in \mathscr{Y}_k} \nu^s_k\bigl(\mathscr{P}_{2^{s+\ell}}(j)\bigr) \right) \label{lem.flatscale.14.2}
\\ & = \gamma \ell \log(2) \, (1 - y_k) \label{lem.flatscale.14.3} \end{align}
This computation may be justified as follows.
\begin{itemize}
	\item The inequality in (\ref{lem.flatscale.14.1}) follows by summing (\ref{lem.flatscale.14.0}) over $j \in \mathscr{Z}_k$.
		\item (\ref{lem.flatscale.14.2}) follows from (\ref{lem.flatscale.14.1}) by the assignment of $\mathscr{Z}_k$ in (\ref{lem.flatscale.13}).
	\item (\ref{lem.flatscale.14.3}) follows from (\ref{lem.flatscale.14.2}) by combining (\ref{lem.flatscale.2}) with the assignment of $y_k$ in (\ref{lem.flatscale.11}).
\end{itemize}
For $k \in \mathscr{X}_s$ we compute:
\begin{align} \ach\bigl(\nu^s_k , 2^{s+\ell}\bigr) & = - \sum_{j=k2^\ell}^{(k+1)2^\ell-1} \nu^s_k\bigl(\mathscr{P}_{2^{s+\ell}}(j)\bigr) \log\Bigl( \nu^s_k\bigl(\mathscr{P}_{2^{s+\ell}}(j)\bigr) \Bigr) \label{lem.flatscale.15.1}
\\ & = - \sum_{j \in \mathscr{Y}_k} \nu^s_k\bigl(\mathscr{P}_{2^{s+\ell}}(j)\bigr) \log\Bigl( \nu^s_k\bigl(\mathscr{P}_{2^{s+\ell}}(j)\bigr) \Bigr) - \sum_{j \in \mathscr{Z}_k} \nu^s_k\bigl(\mathscr{P}_{2^{s+\ell}}(j)\bigr) \log\Bigl( \nu^s_k\bigl(\mathscr{P}_{2^{s+\ell}}(j)\bigr) \Bigr) \label{lem.flatscale.15.2}
\\ & \leq y_k \ell \log(2) + 1 - y_k - \sum_{j \in \mathscr{Z}_k} \nu^s_k\bigl(\mathscr{P}_{2^{s+\ell}}(j)\bigr) \log\Bigl( \nu^s_k\bigl(\mathscr{P}_{2^{s+\ell}}(j)\bigr) \Bigr) \label{lem.flatscale.15.3}
\\ & \leq y_k \ell \log(2) + 1 - y_k + \gamma \ell \log(2) \, (1 - y_k) \label{lem.flatscale.15.4}
\\ & = y_k \ell \log(2) + (1 - y_k)\bigl( 1 + \gamma \ell \log(2) \bigr) \label{lem.flatscale.15.5} \end{align}
This computation may be justified as follows.
\begin{itemize}
	\item The equality in (\ref{lem.flatscale.15.1}) is (\ref{lem.flatscale.2.1}).
	\item (\ref{lem.flatscale.15.2}) follows from (\ref{lem.flatscale.15.1}) since the assignment of $\mathscr{Z}_k$ in (\ref{lem.flatscale.13}) exhibits $[k2^\ell,(k+1)2^\ell-1]_{\zee}$ as the disjoint union of $\mathscr{Y}_k$ and $\mathscr{Z}_k$.
		\item (\ref{lem.flatscale.15.3}) follows from (\ref{lem.flatscale.15.2}) by (\ref{lem.flatscale.12.6}).
	\item (\ref{lem.flatscale.15.4}) follows from (\ref{lem.flatscale.15.3}) by the computation ending in (\ref{lem.flatscale.14.3}).\end{itemize}
Since the cells $\mathscr{P}_{2^s}(k)$ for $k \in [0,2^s-1]_{\zee}$ are pairwise disjoint with union $\rzee$, and since the assignment of $\mathscr{X}_s$ annihilates the remaining terms, we have the following.
\begin{align} \sum_{k \in \mathscr{X}_s} w^s_k = 1 \label{lem.flatscale.16} \end{align}
We define $Y \in [0,1]_{\arr}$ as follows, the bounds holding by (\ref{lem.flatscale.16}) together with $y_k \in [0,1]_{\arr}$.
\begin{align} Y \coloneqq \sum_{k \in \mathscr{X}_s} w^s_k \, y_k \label{lem.flatscale.17} \end{align}
We compute:
\begin{align} (\sigma - \beta)\, \ell \log(2) & \leq \sum_{k \in \mathscr{X}_s} w^s_k \, \ach\bigl(\nu^s_k , 2^{s+\ell}\bigr) \label{lem.flatscale.18.1}
\\ & \leq \sum_{k \in \mathscr{X}_s} w^s_k \Bigl( y_k \ell \log(2) + (1-y_k)\bigl(1 + \gamma \ell \log(2)\bigr) \Bigr) \label{lem.flatscale.18.2}
\\ & = Y \ell \log(2) + (1 - Y)\bigl( 1 + \gamma \ell \log(2) \bigr) \label{lem.flatscale.18.3} \end{align}
This computation may be justified as follows.
\begin{itemize}
	\item The inequality in (\ref{lem.flatscale.18.1}) is (\ref{lem.flatscale.9.3}).
	\item (\ref{lem.flatscale.18.2}) follows from (\ref{lem.flatscale.18.1}) by multiplying the computation ending in (\ref{lem.flatscale.15.5}) by $w^s_k$, which is nonnegative, and summing over $k \in \mathscr{X}_s$.
	\item (\ref{lem.flatscale.18.3}) follows from (\ref{lem.flatscale.18.2}) by combining (\ref{lem.flatscale.16}) with the assignment of $Y$ in (\ref{lem.flatscale.17}).
\end{itemize}
Since Clause \ref{lem.flatscale.c1}\ref{lem.flatscale.c1.2} of our current Lemma \ref{lem.flatscale} provides $\ell \geq 2/\gamma$, so that $\gamma \ell \geq 2$, and since $2\log(2) \geq 1$ we have $\gamma \ell \log(2) \geq 1$ or equivalently:
\begin{align} \label{lem.flatscale.19.0} (1-2\gamma)\ell\log(2) \leq \ell\log(2) - 1 - \gamma\ell\log(2) \end{align}
We compute:
\begin{align} (1 - Y)(1 - 2\gamma)\, \ell \log(2) & \leq  (1 - Y)\, \ell \log(2) - (1-Y)\bigl( 1 + \gamma \ell \log(2) \bigr) \label{lem.flatscale.19.2}
\\ & \leq (1 - Y)\, \ell \log(2) + Y \ell \log(2) - (\sigma - \beta)\, \ell \log(2) \label{lem.flatscale.19.3}
\\ & = (1 - \sigma + \beta)\, \ell \log(2) \label{lem.flatscale.19.4} \end{align}
Here, the inequality in (\ref{lem.flatscale.19.2}) follows from multiplying (\ref{lem.flatscale.19.0}) by the nonnegative quantity $1 - Y$, while (\ref{lem.flatscale.19.3}) follows from (\ref{lem.flatscale.19.2}) by the computation ending in (\ref{lem.flatscale.18.3}).\\

Since Clause \ref{lem.flatscale.c1}\ref{lem.flatscale.c1.1} provides $\gamma \leq 1/4$ we have $1 - 2\gamma > 0$, so that dividing the computation ending in (\ref{lem.flatscale.19.4}) by $(1-2\gamma)\ell\log(2)$ and rearranging provides the following.
\begin{align} Y \geq \frac{\sigma - \beta - 2\gamma}{1 - 2\gamma} \label{lem.flatscale.20} \end{align}
For $k \in \mathscr{X}_s$ with $y_k > 0$ we assign the parameters demanded by Clause \ref{lem.flatscale.c2}\ref{lem.flatscale.c2.2} of our current Lemma \ref{lem.flatscale} as $w_k \coloneqq w^s_k \, y_k$ and by stipulating that the following holds for every Borel subset $E$ of $\rzee$.
\begin{align} \nu_k(E) \coloneqq \frac{1}{y_k} \sum_{j \in \mathscr{Y}_k} \nu^s_k\bigl( E \cap \mathscr{P}_{2^{s+\ell}}(j) \bigr) \label{lem.flatscale.21} \end{align}
For every other $k \in [0,2^s-1]_{\zee}$ we use the dummy assignments $w_k \coloneqq 0$ and $\nu_k \coloneqq 2^s \lambda(\,\cdot \cap \mathscr{P}_{2^s}(k))$. \\

We first verify Clause \ref{lem.flatscale.c5}. For $k \in \mathscr{X}_s$ with $y_k > 0$, the upper case of (\ref{lem.flatscale.1}) provides $\mathscr{P}_{2^{s+\ell}}(j) \cap \mathscr{P}_{2^s}(k) = \mathscr{P}_{2^{s+\ell}}(j)$ for every $j \in \mathscr{Y}_k$, so that the assignment of $\nu_k$ in (\ref{lem.flatscale.21}) together with the assignment of $y_k$ in (\ref{lem.flatscale.11}) provides $\nu_k(\mathscr{P}_{2^s}(k)) = 1$; the same assignment with $E \coloneqq \rzee$ shows that $\nu_k$ is an element of $\prob(\rzee)$. For every other $k$ the assignment of $\nu_k$ provides $\nu_k(\mathscr{P}_{2^s}(k)) = 1$ directly. \\

We next verify Clause \ref{lem.flatscale.c3}. Since $w_k = 0$ for $k \in [0,2^s-1]_{\zee} \setminus \mathscr{X}_s$, and since $w^s_k y_k = 0$ for those $k \in \mathscr{X}_s$ with $y_k = 0$, the assignment of $w_k$ together with the assignment of $Y$ in (\ref{lem.flatscale.17}) provides $\sum_{k=0}^{2^s-1} w_k = Y$, and (\ref{lem.flatscale.20}) is the desired inequality. \\

It remains to verify Clauses \ref{lem.flatscale.c4} and \ref{lem.flatscale.c6}. To that end we fix a Borel subset $E$ of $\rzee$ and compute:
\begin{align} \sum_{k=0}^{2^s-1} w_k \, \nu_k(E) & = \sum_{k \in \mathscr{X}_s} \sum_{j \in \mathscr{Y}_k} w^s_k \, \nu^s_k\bigl( E \cap \mathscr{P}_{2^{s+\ell}}(j) \bigr) \label{lem.flatscale.22.1}
\\ & = \sum_{k \in \mathscr{X}_s} \sum_{j \in \mathscr{Y}_k} \mu\bigl( E \cap \mathscr{P}_{2^{s+\ell}}(j) \bigr) \label{lem.flatscale.22.2}
\\ & \leq \mu(E) \label{lem.flatscale.22.3} \end{align}
This computation may be justified as follows.
\begin{itemize}
	\item The equality in (\ref{lem.flatscale.22.1}) follows by combining the assignment of $w_k$ with the assignment of $\nu_k$ in (\ref{lem.flatscale.21}), the terms indexed by the remaining $k$ vanishing since $w_k = 0$ for those $k$.
	\item (\ref{lem.flatscale.22.2}) follows from (\ref{lem.flatscale.22.1}) by combining the assignment of $\nu^s_k$ with the upper case of (\ref{lem.flatscale.1}), which provides $\mathscr{P}_{2^{s+\ell}}(j) \cap \mathscr{P}_{2^s}(k) = \mathscr{P}_{2^{s+\ell}}(j)$ for every $j \in \mathscr{Y}_k$.
	\item (\ref{lem.flatscale.22.3}) follows from (\ref{lem.flatscale.22.2}) since the cells $\mathscr{P}_{2^{s+\ell}}(j)$ indexed by the pairs $(k,j)$ with $k \in \mathscr{X}_s$ and $j \in \mathscr{Y}_k$ are pairwise disjoint subsets of $\rzee$.
\end{itemize}
Since $E$ was an arbitrary Borel subset of $\rzee$, the computation ending in (\ref{lem.flatscale.22.3}) verifies Clause \ref{lem.flatscale.c4} by way of \ref{conv.submeasure}.
We now turn to Clause \ref{lem.flatscale.c6}, and to that end we fix $k \in \mathscr{X}_s$ with $y_k > 0$ and compute:
\begin{align} w_k \, \bigl\lvert\bigl\lvert \nu_k(\mathscr{P}_{2^{s+\ell}}) \bigr\rvert\bigr\rvert_2^2 & = \frac{w^s_k}{y_k} \sum_{j \in \mathscr{Y}_k} \nu^s_k\bigl(\mathscr{P}_{2^{s+\ell}}(j)\bigr)^2 \label{lem.flatscale.23.1}
\\ & \leq \frac{w^s_k}{y_k \, 2^{\gamma \ell}} \sum_{j \in \mathscr{Y}_k} \nu^s_k\bigl(\mathscr{P}_{2^{s+\ell}}(j)\bigr) \label{lem.flatscale.23.2}
\\ & = \frac{w^s_k}{2^{\gamma \ell}} \label{lem.flatscale.23.3} \end{align}
This computation may be justified as follows.
\begin{itemize}
	\item The equality in (\ref{lem.flatscale.23.1}) follows by combining the assignment of $w_k$ with the definition of the resolution norm in (\ref{eq.mupnorm}), where the assignment of $\nu_k$ in (\ref{lem.flatscale.21}) provides $\nu_k(\mathscr{P}_{2^{s+\ell}}(j)) = \nu^s_k(\mathscr{P}_{2^{s+\ell}}(j))/y_k$ for $j \in \mathscr{Y}_k$ and $\nu_k(\mathscr{P}_{2^{s+\ell}}(j)) = 0$ for every other $j \in [0,2^{s+\ell}-1]_{\zee}$.
	\item (\ref{lem.flatscale.23.2}) follows from (\ref{lem.flatscale.23.1}) since the assignment of $\mathscr{Y}_k$ in (\ref{lem.flatscale.10}) provides $\nu^s_k(\mathscr{P}_{2^{s+\ell}}(j)) < 2^{-\gamma\ell}$ for every $j \in \mathscr{Y}_k$.
	\item (\ref{lem.flatscale.23.3}) follows from (\ref{lem.flatscale.23.2}) by the assignment of $y_k$ in (\ref{lem.flatscale.11}).
\end{itemize}
Since $k$ was an arbitrary element of $\mathscr{X}_s$ with $y_k > 0$, we compute again:
\begin{align} \sum_{k=0}^{2^s-1} w_k \, \bigl\lvert\bigl\lvert \nu_k(\mathscr{P}_{2^{s+\ell}}) \bigr\rvert\bigr\rvert_2^2 & \leq \frac{1}{2^{\gamma \ell}} \sum_{k \in \mathscr{X}_s} w^s_k \label{lem.flatscale.24.1}
\\ & = \frac{1}{2^{\gamma \ell}} \label{lem.flatscale.24.2} \end{align}
This computation may be justified as follows.
\begin{itemize}
	\item The inequality in (\ref{lem.flatscale.24.1}) follows by summing the computation ending in (\ref{lem.flatscale.23.3}) over those $k \in \mathscr{X}_s$ with $y_k > 0$, the terms indexed by every other $k \in [0,2^s-1]_{\zee}$ vanishing since $w_k = 0$ for those $k$ and the terms so introduced on the right being nonnegative.
	\item (\ref{lem.flatscale.24.2}) follows from (\ref{lem.flatscale.24.1}) by (\ref{lem.flatscale.16}).
\end{itemize}
The computation ending in (\ref{lem.flatscale.24.2}) verifies Clause \ref{lem.flatscale.c6} and thereby completes the proof of Lemma \ref{lem.flatscale}.
\end{Proof}

\subsection{Averaging along shadow sets} \label{cha.pfmain.infra.shadave}

The following lemma is the root of the current article's technology, transferring $\ell^2$-resolution norms to an averaged Fourier estimate. The analog in BLMV is Lemma 4.6 from \cite{MR2563089}, but here we take the averaging set $\mathscr{Q}$ from the progression-shadowing hypothesis of Definition \ref{def.progshad} instead of a gap estimate for $\mathbb{S}_{a,b}$. The constraint on $\ell$ thereby leaves the lemma altogether, becoming a hypothesis supplied from outside, and what appears in its place is the explicit shadowing error appearing at the right of (\ref{lem.shadave.conc1}). The only analytic input is Lemma \ref{lem.fouravg}, as the only analytic input to BLMV's Lemma 4.6 is their Lemma 4.5.\\

That a gap estimate for a multiplicative semigroup can be exchanged for a hypothesis on an arbitrary set of multipliers is a core thesis of the present article, and the exchange is made here.\\

The dependencies in the statement of Lemma \ref{lem.shadave} are as follows.
\begin{deps}
\dep{Definitions}{\depitem{Progression-shadowing sets}{def.progshad}}
\depcont{\depitem{Resolution entropy}{def.entropy}}
\depcont{\depitem{Fourier coefficients}{def.foucoef}}
\depcont{\depitem{Resolution norms of measures}{def.mupnorm}}\end{deps}

We now move on to the formal statement.

\begin{Lemma}[Shadow-averaging lemma] \label{lem.shadave} \quad
\vspace{-1 \baselineskip}
\begin{enumerate}[label=\upshape(\arabic*)] \item \label{lem.shadave.c1} For all:
\begin{enumerate}[label=\upshape(\arabic*)] \item \label{lem.shadave.c1.1} $s \in \enn$ and all $\ell \in [2,\infty)_{\zee}$,
\item \label{lem.shadave.c1.2} and all subsets $\mathscr{M}$ of $\enn$ which are $(s,\ell)$ progression-shadowing,
\end{enumerate}
\item \label{lem.shadave.c2} there exists a subset $\mathscr{Q}$ of $\mathscr{M} \cap [1,2^{s+\ell+1}]_{\zee}$ satisfying $\lvert \mathscr{Q} \rvert = 2^\ell$,
\item \label{lem.shadave.c3} such that for all:
\begin{enumerate}[label=\upshape(\arabic*)] \item \label{lem.shadave.c3.1} $k \in [0,2^s-1]_{\zee}$ and all $\nu \in \prob(\rzee)$ satisfying $\nu(\mathscr{P}_{2^s}(k)) = 1$,
\item \label{lem.shadave.c3.2} and all $g \in \zee \setminus \{0\}$,
\end{enumerate}
\end{enumerate}
we have:
\begin{align} \frac{1}{|\mathscr{Q}|} \sum_{q \in \mathscr{Q}} \lvert \hat{\nu}(qg) \rvert^2 \leq 8 \left(16 \lvert g \rvert \, \bigl \vert \bigl \vert \nu(\mathscr{P}_{2^{s+\ell}}) \bigr \vert \bigr \vert_2^2 + \frac{\pi^2 g^2}{2^{2\ell}} \right) \label{lem.shadave.conc1} \end{align} \end{Lemma}

An outline of the proof of Lemma \ref{lem.shadave} is as follows.

\begin{itemize}
\item The shadowing hypothesis is harvested once and for all, around (\ref{lem.shadave.p1}): for each $m \in [1,2^\ell]_{\zee}$ it supplies an element $q_m$ of $\mathscr{M}$ lying within $2^{s-\ell}$ of $m2^s$, and $\mathscr{Q}$ is the set of these. Consecutive multiples $m2^s$ are separated by $2^s$ while the perturbations have size at most $2^{s-\ell}$, so the $q_m$ stay distinct and $\mathscr{Q}$ has the cardinality and range demanded by Clause \ref{lem.shadave.c2}.
\item Because $\nu$ is carried by a single cell of resolution $2^s$, the characters $\mathfrak{e}^{-q_mg}$ and $\mathfrak{e}^{-m2^sg}$ agree across that cell up to a unimodular constant and an error of $2\pi \lvert g \rvert 2^{-\ell}$. This is the computation around (\ref{lem.shadave.p4.1}), and it is where the shadowing error is spent. What it buys, recorded in (\ref{lem.shadave.p5}), is the replacement of the irregular frequencies $q_mg$ by the arithmetic progression $m2^sg$.
\item Passing to the endomorph $\varsigma \coloneqq 2^s_{\,\circledcirc}\nu$ converts $\hat{\nu}(m2^sg)$ into $\hat{\varsigma}(mg)$, so that after enlarging the range of summation to $[0,2^{\ell+1}-1]_{\zee}$ the average is of exactly the shape bounded by Lemma \ref{lem.fouravg}, which is applied at (\ref{lem.shadave.p9.5}). 
\item Two identifications of scale return the estimate to $\nu$: the support hypothesis matches the resolution norm of $\varsigma$ at scale $2^{\ell+1}$ with that of $\nu$ at scale $2^{s+\ell+1}$. This is established at (\ref{lem.shadave.p8}), and refinement monotonicity relaxes the latter to scale $2^{s+\ell}$ at (\ref{lem.shadave.p9.7}).
\end{itemize}

The dependencies in the proof of Lemma \ref{lem.shadave} beyond those in the statement are as follows.
\begin{deps}
\dep{Conventions}{\depitem{Periodic functions on $\rzee$}{conv.period}\depsep\depitem{$\mathtt{T}_n$}{conv.Tmap}}
\dep{Definitions}{\depitem{Endomorphs of complex measures}{def.endomorph}}
\dep{Remarks}{\depitem{Lipschitz estimate for exponentials}{rem.expdiff}}
\depcont{\depitemc{Fourier coefficients of endomorphs}{rem.endochangeofvar}{rem.endochangeofvar.c1}}
\depcont{\depitemc{Change of variable for endomorphs}{rem.endochangeofvar}{rem.endochangeofvar.c3}}
\depcont{\depitem{Refinement monotonicity of resolution norms}{rem.mupnormref}}
\dep{Lemmas}{\depitem{Averaged Fourier decay at a single resolution}{lem.fouravg}}
\end{deps}

We now move on to the formal proof.

\begin{Proof}[of Lemma \ref{lem.shadave}] Let $s, \ell$ and $\mathscr{M}$ be as in Clause \ref{lem.shadave.c1} of Lemma \ref{lem.shadave}. For each $m \in [1,2^\ell]_{\zee}$, instantiating Definition \ref{def.progshad} with our currently-local $\mathscr{M}$ returns an element $q_m$ of $\mathscr{M}$ such that:
\begin{align} \lvert q_m - m 2^s \rvert \leq 2^{s-\ell} \label{lem.shadave.p1} \end{align}
We assign the set demanded by Clause \ref{lem.shadave.c2} as $\mathscr{Q} \coloneqq \{ q_m : m \in [1,2^\ell]_{\zee} \}$, so that $\mathscr{Q} \subseteq \mathscr{M}$ by construction.\\

We now verify the remaining assertions of Clause \ref{lem.shadave.c2}. For all $m \in [1,2^\ell]_{\zee}$ we have:
\begin{align*} q_m \leq m2^s + 2^{s-\ell} \leq 2^{s+\ell} + 2^{s-\ell} \leq 2^{s+\ell+1} \end{align*}
Here, the first inequality follows from (\ref{lem.shadave.p1}) while the second follows since $m \leq 2^\ell$. Also combining (\ref{lem.shadave.p1}) with the bound $m \geq 1$ we find $q_m \geq 2^s(1 - 2^{-\ell}) > 0$, and since $q_m \in \enn$ this yields $q_m \geq 1$. Thus $\mathscr{Q} \subseteq \mathscr{M} \cap [1,2^{s+\ell+1}]_{\zee}$. Moreover, for all distinct $m,m' \in [1,2^\ell]_{\zee}$ we have:
\begin{align*} \lvert q_m - q_{m'} \rvert \geq 2^s - 2^{s-\ell+1} \geq 2^{s-1} \end{align*}
Here, the first inequality follows by combining the triangle inequality with (\ref{lem.shadave.p1}) and the bound $\lvert m - m' \rvert \geq 1$, while the second follows since $\ell \geq 2$. In particular the assignment $m \mapsto q_m$ is injective, and so $\lvert \mathscr{Q} \rvert = 2^\ell$. This completes the verification of Clause \ref{lem.shadave.c2}, and so we are challenged with data $(k,\nu,g)$ as in Clause \ref{lem.shadave.c3}.\\

Toward verifying the conclusion (\ref{lem.shadave.conc1}), we temporarily fix $m \in [1,2^\ell]_{\zee}$ and define $\theta_m \in \cee$ as follows.
\begin{align*} \theta_m \coloneqq e\bigl((m2^s - q_m)\,g\,k2^{-s}\bigr) \end{align*}
Recalling from Definition \ref{def.entropy} and \ref{conv.period} that the cell $\mathscr{P}_{2^s}(k)$ corresponds to the interval $[k2^{-s}, (k+1)2^{-s})_{\arr}$ of the coset transversal, we compute:
\begin{align} \bigl \lvert \hat{\nu}(q_m g) - \theta_m \hat{\nu}(m2^sg) \bigr \rvert & = \lvertt \nu\bigl(\mathfrak{e}^{-q_mg} - \theta_m \, \mathfrak{e}^{-m2^sg}\bigr) \rvertt \label{lem.shadave.p4.1}
\\ & \leq \sup_{t \in [k2^{-s},\,(k+1)2^{-s})_{\arr}} \bigl \lvert e(-q_mgt) - \theta_m \, e(-m2^sgt) \bigr \rvert \label{lem.shadave.p4.2}
\\ & = \sup_{t \in [k2^{-s},\,(k+1)2^{-s})_{\arr}} \bigl \lvert e\bigl((m2^s - q_m)gt\bigr) - e\bigl((m2^s-q_m)\,g\,k2^{-s}\bigr) \bigr \rvert \label{lem.shadave.p4.3}
\\ & \leq 2\pi \, \lvert m2^s - q_m \rvert \, \lvert g \rvert \, 2^{-s} \label{lem.shadave.p4.4}
\\ & \leq 2\pi \, \lvert g \rvert \, 2^{-\ell} \label{lem.shadave.p4.5} \end{align}
This computation may be justified as follows.
\begin{itemize}
\item The equality in (\ref{lem.shadave.p4.1}) follows by combining Definition \ref{def.foucoef} with the linearity of the functional $\nu$.
\item (\ref{lem.shadave.p4.2}) follows from (\ref{lem.shadave.p4.1}) since $\nu \in \prob(\rzee)$ satisfies $\nu(\mathscr{P}_{2^s}(k)) = 1$, so that the following holds for all $\psi \in \bbr(\rzee)$.
\begin{align*} \lvert \nu(\psi) \rvert \leq \sup_{x \in \mathscr{P}_{2^s}(k)} |\psi(x)| \end{align*}
\item (\ref{lem.shadave.p4.3}) follows from (\ref{lem.shadave.p4.2}) by factoring out the unimodular quantity $e(-m2^sgt)$ and unwrapping the assignment of $\theta_m$.
\item (\ref{lem.shadave.p4.4}) follows from (\ref{lem.shadave.p4.3}) by Remark \ref{rem.expdiff}, where we note that $\lvert t - k2^{-s} \rvert \leq 2^{-s}$ for all $t$ in the relevant interval.
\item (\ref{lem.shadave.p4.5}) follows from (\ref{lem.shadave.p4.4}) by (\ref{lem.shadave.p1}).
\end{itemize}
Since $\lvert \theta_m \rvert = 1$, combining the computation ending in (\ref{lem.shadave.p4.5}) with the triangle inequality gives:
\begin{align*} \lvert \hat{\nu}(q_mg) \rvert \leq \lvert \hat{\nu}(m2^sg) \rvert + 2\pi \lvert g \rvert 2^{-\ell} \end{align*} 
The elementary bound $(u+v)^2 \leq 2(u^2 + v^2)$ holds for all $u,v \in [0,\infty)_{\arr}$, and combining with the previous display we obtain:
\begin{align} \lvert \hat{\nu}(q_mg) \rvert^2 \leq 2\, \lvert \hat{\nu}(m2^sg) \rvert^2 + \frac{8\pi^2 g^2}{2^{2\ell}} \label{lem.shadave.p5} \end{align}
Since $m \in [1,2^\ell]_{\zee}$ was arbitrary, the inequality in (\ref{lem.shadave.p5}) holds for all such $m$. Continuing, we define $\varsigma \coloneqq 2^s_{\,\circledcirc}\, \nu$ according to Definition \ref{def.endomorph}. Clause \ref{rem.endochangeofvar.c3} of Remark \ref{rem.endochangeofvar} expresses $\varsigma$ as the pushforward of $\nu$ along $\mathtt{T}_{2^s}$, and in particular provides $\varsigma \in \prob(\rzee)$.\\

Directly from \ref{conv.Tmap} and Definition \ref{def.entropy}, we obtain the following for all $j \in [0,2^{\ell+1}-1]_{\zee}$.
\begin{align} \mathtt{T}_{2^s}^{-1}\bigl(\mathscr{P}_{2^{\ell+1}}(j)\bigr) = \bigcup_{r=0}^{2^s-1} \mathscr{P}_{2^{s+\ell+1}}\bigl(r 2^{\ell+1} + j\bigr) \label{lem.shadave.p6} \end{align}
Also directly from Definition \ref{def.entropy}, we have the following for all $u \in [0,2^{s+\ell+1}-1]_{\zee}$.
\begin{align*} \mathscr{P}_{2^{s+\ell+1}}(u) \subseteq \mathscr{P}_{2^s}\bigl(\bigl\lfloor u \, 2^{-\ell-1} \bigr\rfloor\bigr)  \end{align*}
Since $\nu \in \prob(\rzee)$ satisfies $\nu(\mathscr{P}_{2^s}(k)) = 1$, we have $\nu(\mathscr{P}_{2^s}(r)) = 0$ for all $r \in [0,2^s-1]_{\zee} \setminus \{k\}$, and so the previous display implies that for all $u \in [0,2^{s+\ell+1}-1]_{\zee} \setminus [k2^{\ell+1},(k+1)2^{\ell+1}-1]_{\zee}$ we have:
\begin{align} \nu\bigl(\mathscr{P}_{2^{s+\ell+1}}(u)\bigr) = 0 \label{lem.shadave.p6.2} \end{align}

Now, temporarily fix $j \in [0,2^{\ell+1}-1]_{\zee}$ and $r \in [0,2^s-1]_{\zee} \setminus \{k\}$. Since $r \leq 2^s-1$ and $j \leq 2^{\ell+1}-1$ we have $r2^{\ell+1}+j \in [0,2^{s+\ell+1}-1]_{\zee}$. We now consider the two cases $r<k$ and $r>k$ separately.
\begin{itemize}
\item If $r<k$ then $r+1 \leq k$ and so we have:
\begin{align*} r2^{\ell+1}+j \leq (r+1)2^{\ell+1}-1 \leq k2^{\ell+1}-1 \end{align*}
\item If $r>k$ then $r \geq k+1$ and so we have:
\begin{align*} r2^{\ell+1}+j \geq r2^{\ell+1} \geq (k+1)2^{\ell+1} \end{align*}
\end{itemize}
In either case we find $r2^{\ell+1}+j \notin [k2^{\ell+1},(k+1)2^{\ell+1}-1]_{\zee}$. Thus (\ref{lem.shadave.p6.2}) implies that the following holds for all $j \in [0,2^{\ell+1}-1]_{\zee}$ and all $r \in [0,2^s-1]_{\zee} \setminus \{k\}$.
\begin{align} \nu\Bigl(\mathscr{P}_{2^{s+\ell+1}}\bigl(r2^{\ell+1}+j\bigr)\Bigr) = 0 \label{lem.shadave.p6.3} \end{align}
Again let $j \in [0,2^{\ell+1}-1]_{\zee}$ and calculate:
\begin{align} \varsigma\bigl(\mathscr{P}_{2^{\ell+1}}(j)\bigr) & = \nu\left(\mathtt{T}_{2^s}^{-1}\bigl(\mathscr{P}_{2^{\ell+1}}(j)\bigr)\right) \label{lem.shadave.p7.1}
\\ & = \sum_{r=0}^{2^s-1} \nu\Bigl(\mathscr{P}_{2^{s+\ell+1}}\bigl(r2^{\ell+1}+j\bigr)\Bigr) \label{lem.shadave.p7.2}
\\ & = \nu\Bigl(\mathscr{P}_{2^{s+\ell+1}}\bigl(k2^{\ell+1}+j\bigr)\Bigr) \label{lem.shadave.p7} \end{align}
This calculation may be justified as follows.
\begin{itemize}
\item The equality in (\ref{lem.shadave.p7.1}) follows from Clause \ref{rem.endochangeofvar.c3} in Remark \ref{rem.endochangeofvar}, applied to the assignment $\varsigma = 2^s_{\,\circledcirc} \nu$.
\item (\ref{lem.shadave.p7.2}) follows from (\ref{lem.shadave.p7.1}) by (\ref{lem.shadave.p6}), where the additivity of $\nu$ is justified by the disjointness of the union appearing there.
\item (\ref{lem.shadave.p7}) follows from (\ref{lem.shadave.p7.2}) by (\ref{lem.shadave.p6.3}).
\end{itemize}
We may combine (\ref{lem.shadave.p6.2}) with (\ref{lem.shadave.p7}) to obtain the following for all $u \in [0,2^{s+\ell+1}-1]_{\zee}$.
\begin{align*} \nu\bigl(\mathscr{P}_{2^{s+\ell+1}}(u)\bigr) = \begin{dcases} \varsigma\Bigl(\mathscr{P}_{2^{\ell+1}}\bigl(u - k2^{\ell+1}\bigr)\Bigr) & \mbox{ if } u \in [k2^{\ell+1}, (k+1)2^{\ell+1}-1]_{\zee} \\ \hspace{0.5 cm} 0 & \mbox{ else} \end{dcases} \end{align*}
Since the assignment $u \mapsto u - k2^{\ell+1}$ is a bijection from $[k2^{\ell+1},(k+1)2^{\ell+1}-1]_{\zee}$ to $[0,2^{\ell+1}-1]_{\zee}$, combining the previous display with Definition \ref{def.mupnorm} we obtain:
\begin{align} \bigl \vert \bigl \vert \varsigma(\mathscr{P}_{2^{\ell+1}}) \bigr \vert \bigr \vert_2^2 = \bigl \vert \bigl \vert \nu(\mathscr{P}_{2^{s+\ell+1}}) \bigr \vert \bigr \vert_2^2 \label{lem.shadave.p8} \end{align}
We are now ready to assemble the desired inequality. We compute:
\begin{align} \frac{1}{\lvert \mathscr{Q} \rvert} \sum_{q \in \mathscr{Q}} \lvert \hat{\nu}(qg) \rvert^2 & = \frac{1}{2^\ell}\sum_{m = 1}^{2^\ell} \lvert \hat{\nu}(q_m g) \rvert^2 \label{lem.shadave.p9.1}
\\ & \leq \frac{2}{2^\ell}\sum_{m = 1}^{2^\ell} \lvert \hat{\nu}(m2^s g) \rvert^2 + \frac{8\pi^2g^2}{2^{2\ell}} \label{lem.shadave.p9.2}
\\ & = \frac{2}{2^\ell}\sum_{m = 1}^{2^\ell} \lvert \hat{\varsigma}(m g) \rvert^2 + \frac{8\pi^2g^2}{2^{2\ell}} \label{lem.shadave.p9.3}
\\ & \leq \frac{2}{2^\ell}\sum_{m = 0}^{2^{\ell+1}-1} \lvert \hat{\varsigma}(m g) \rvert^2 + \frac{8\pi^2g^2}{2^{2\ell}} \label{lem.shadave.p9.4}
\\ & \leq 128 \lvert g \rvert \, \bigl \vert \bigl \vert \varsigma(\mathscr{P}_{2^{\ell+1}}) \bigr \vert \bigr \vert_2^2 + \frac{8\pi^2g^2}{2^{2\ell}} \label{lem.shadave.p9.5}
\\ & = 128 \lvert g \rvert \, \bigl \vert \bigl \vert \nu(\mathscr{P}_{2^{s+\ell+1}}) \bigr \vert \bigr \vert_2^2 + \frac{8\pi^2g^2}{2^{2\ell}} \label{lem.shadave.p9.6}
\\ & \leq 128 \lvert g \rvert \, \bigl \vert \bigl \vert \nu(\mathscr{P}_{2^{s+\ell}}) \bigr \vert \bigr \vert_2^2 + \frac{8\pi^2g^2}{2^{2\ell}} \label{lem.shadave.p9.7} \end{align}
This computation may be justified as follows.
\begin{itemize}
\item The equality in (\ref{lem.shadave.p9.1}) follows since by construction we have $\mathscr{Q} = \{q_m: m \in [1,2^\ell]_\zee \}$.
\item (\ref{lem.shadave.p9.2}) follows from (\ref{lem.shadave.p9.1}) by averaging (\ref{lem.shadave.p5}) over $m \in [1,2^\ell]_{\zee}$.
\item (\ref{lem.shadave.p9.3}) follows from (\ref{lem.shadave.p9.2}) by Clause \ref{rem.endochangeofvar.c1} in Remark \ref{rem.endochangeofvar}.
\item (\ref{lem.shadave.p9.4}) follows from (\ref{lem.shadave.p9.3}) since the additional summands are nonnegative.
\item (\ref{lem.shadave.p9.5}) follows from (\ref{lem.shadave.p9.4}) by instantiating Lemma \ref{lem.fouravg} with the assignment $(\nu,M,g) \coloneqq (\varsigma,2^{\ell+1},g)$, where the frequency $g$ is nonzero by Clause \ref{lem.shadave.c3}\ref{lem.shadave.c3.2} in the current lemma, and where we note that $2 \cdot 2^{-\ell} = 4 \cdot 2^{-\ell-1}$.
\item (\ref{lem.shadave.p9.6}) follows from (\ref{lem.shadave.p9.5}) by (\ref{lem.shadave.p8}).
\item (\ref{lem.shadave.p9.7}) follows from (\ref{lem.shadave.p9.6}) by Remark \ref{rem.mupnormref}.
\end{itemize}
The computation ending in (\ref{lem.shadave.p9.7}) verifies the desired inequality (\ref{lem.shadave.conc1}) and thereby completes the proof of Lemma \ref{lem.shadave}.
\end{Proof}

\subsection{Averaged Fourier decay for entropic measures} \label{cha.pfmain.infra.decompave}

The purpose of Lemma \ref{lem.decompave} below is to trade the resolution-norm factor appearing in the conclusion of Lemma \ref{lem.shadave} for an explicit rate of Fourier decay, at the price of a hypothesis on the resolution entropy of $\mu$ and of passing to a submeasure of $\mu$ produced along the way. The mechanism effecting this trade is the device built for exactly this purpose by BLMV, which appears above as Lemma \ref{lem.flatscale}. The submeasure and the multiplier set returned here are consumed below in Lemma \ref{lem.extract}.\\

The dependencies in the statement of Lemma \ref{lem.decompave} are as follows.
\begin{deps}
\dep{Conventions}{\depitem{Submeasures}{conv.submeasure}}
\dep{Definitions}{\depitem{Progression-shadowing depth}{def.shadep}\depsep\depitem{Resolution entropy}{def.entropy}\depsep\depitem{Fourier coefficients}{def.foucoef}}
\end{deps}

We now move on to the formal statement.

\begin{Lemma}[Averaged Fourier decay for entropic measures] \label{lem.decompave} \quad
\vspace{-1 \baselineskip}
\begin{enumerate}[label=\upshape(\arabic*)] \item \label{lem.decompave.c1} For all:
\begin{enumerate}[label=\upshape(\arabic*)] \item \label{lem.decompave.c1.1} $n \in \enn$ and all $\sigma \in (0,1]_{\arr}$ and all $\beta \in (0,1]_{\arr}$ and all $\gamma \in (0,1/4]_{\arr}$ such that $\sigma > \beta + 2\gamma$,
\item \label{lem.decompave.c1.2} and all $\ell \in [2/\gamma, \beta n]_{\zee}$,
\item \label{lem.decompave.c1.3} and all $\mu \in \prob(\rzee)$ such that $\achn(\mu,2^n) \geq \sigma$,
\item \label{lem.decompave.c1.4} and all subsets $\mathscr{M}$ of $\enn$ such that $\mathtt{Dep}_{\mathscr{M}}(s_0) \geq \ell$ for all $s_0 \in [\beta n - \ell, n - \ell]_{\zee}$,
\end{enumerate}
\item \label{lem.decompave.c2} there exist:
\begin{enumerate}[label=\upshape(\arabic*)] \item \label{lem.decompave.c2.1} $s \in [\beta n - \ell, n - \ell]_{\zee}$,
\item \label{lem.decompave.c2.2} and $w \in [(\sigma - \beta - 2\gamma)/(1-2\gamma), 1]_{\arr}$,
\item \label{lem.decompave.c2.3} and $\nu \in \prob(\rzee)$ satisfying $w \nu \preceq \mu$,
\item \label{lem.decompave.c2.4} and a subset $\mathscr{Q}$ of $\mathscr{M} \cap [1,2^{s+\ell+1}]_{\zee}$ satisfying $\lvert \mathscr{Q} \rvert = 2^\ell$, \end{enumerate}
\end{enumerate}
such that for all $g \in \zee \setminus \{0\}$ we have the inequality below.
\begin{align} \frac{1}{\lvert \mathscr{Q} \rvert} \sum_{q \in \mathscr{Q}} \lvert \hat{\nu}(qg) \rvert^2 \leq \frac{8}{w^2}\left( \frac{32 \lvert g \rvert}{2^{\gamma \ell}} + \frac{\pi^2 g^2}{2^{2\ell}} \right) \label{lem.decompave.conc1} \end{align} \end{Lemma}

An outline of the proof of Lemma \ref{lem.decompave} is as follows.
\begin{itemize}
\item The decomposition from Lemma \ref{lem.flatscale} is imported wholesale, so we are returned a scale $s$ together with a weighted family of measures carried by the individual cells of resolution $2^s$. At (\ref{lem.decompave.1}) and (\ref{lem.decompave.1.1}) we assemble these into a single total mass $w$ and a single measure $\nu$. That $w$ is at least $(\sigma - \beta - 2\gamma)/(1-2\gamma)$ and that $w\nu \preceq \mu$ are read off the same import, discharging Clauses \ref{lem.decompave.c2.1} through \ref{lem.decompave.c2.3}.
\item The returned scale is then checked against the hypotheses of Lemma \ref{lem.shadave} and the depth hypothesis converts to progression-shadowing at the scale $s$ by Clause \ref{rem.shadep.c1} of Remark \ref{rem.shadep}. Instantiating that lemma returns the multiplier set $\mathscr{Q}$ demanded by Clause \ref{lem.decompave.c2.4}.
\item Because each measure in the family is carried by a single cell of resolution $2^s$, the shadow-averaging estimate applies to each of them separately, which is the content of (\ref{lem.decompave.2}).
\item What remains is to average those estimates against the weights. The Cauchy-Schwartz inequality at (\ref{lem.decompave.3.2}) trades the square of a weighted sum for the weighted sum of squares at the cost of one factor of $w$, the per-cell bounds are inserted, and the entropy hypothesis enters only at (\ref{lem.decompave.3.7}), where Lemma \ref{lem.flatscale} converts the weighted sum of resolution norms into the decay $2^{-\gamma \ell}$.
\end{itemize}

The dependencies in the proof of Lemma \ref{lem.decompave} beyond those in the statement are as follows.
\begin{deps}
\dep{Definitions}{\depitem{Resolution norms of measures}{def.mupnorm}}
\dep{Remarks}{\depitemc{Shadowing depth and progression-shadowing}{rem.shadep}{rem.shadep.c1}}
\dep{Lemmas}{\depitem{Flat components at a good scale}{lem.flatscale}\depsep\depitem{Shadow-averaging lemma}{lem.shadave}}
\end{deps}

We now move on to the formal proof.

\begin{Proof}[of Lemma \ref{lem.decompave}] Let $n,\sigma,\beta,\gamma,\ell,\mu$ and $\mathscr{M}$ be as in Clause \ref{lem.decompave.c1} of Lemma \ref{lem.decompave}. Instantiating Clause \ref{lem.flatscale.c1} of Lemma \ref{lem.flatscale} with the data $(n,\sigma,\beta,\gamma,\ell,\mu)$, we are returned by Clause \ref{lem.flatscale.c2} of that lemma an integer $s \in [\beta n - \ell, n - \ell]_{\zee}$ together with elements $w_0,\ldots,w_{2^s-1}$ of $[0,\infty)_{\arr}$ and elements $\nu_0,\ldots,\nu_{2^s-1}$ of $\prob(\rzee)$. We assign the parameter $s$ demanded by Clause \ref{lem.decompave.c2}\ref{lem.decompave.c2.1} of Lemma \ref{lem.decompave} as currently-local, and we define $w \in [0,\infty)_{\arr}$ as follows.
\begin{align} \label{lem.decompave.1} w \coloneqq \sum_{k=0}^{2^s-1} w_k \end{align}
Clause \ref{lem.flatscale.c3} of Lemma \ref{lem.flatscale} asserts that $w \geq (\sigma - \beta - 2\gamma)/(1-2\gamma)$. Clause \ref{lem.decompave.c1}\ref{lem.decompave.c1.1} of Lemma \ref{lem.decompave} provides $\sigma > \beta + 2\gamma$, so that the numerator of this ratio is positive, and provides $\gamma \leq 1/4$, so that its denominator is at least $1/2$. Hence $w > 0$, and in particular division by $w$ is legal in the sequel.\\

Clause \ref{lem.flatscale.c4} of Lemma \ref{lem.flatscale} asserts that $\sum_{k=0}^{2^s-1}w_k\nu_k \preceq \mu$, and so instantiating \ref{conv.submeasure} with the Borel set $\rzee$ we obtain $w \leq 1$. Thus $w$ lies in the interval demanded by Clause \ref{lem.decompave.c2}\ref{lem.decompave.c2.2} of Lemma \ref{lem.decompave}, and so we may assign this $w$ as the parameter demanded there. We now define $\nu$ as follows.
\begin{align} \label{lem.decompave.1.1} \nu \coloneqq \frac{1}{w}\sum_{k=0}^{2^s-1} w_k \nu_k \end{align}
Since each $\nu_k \in \prob(\rzee)$, it follows immediately from (\ref{lem.decompave.1}) that $\nu \in \prob(\rzee)$. Finally, unwrapping the assignment of $\nu$ and again citing Clause \ref{lem.flatscale.c4} of Lemma \ref{lem.flatscale} we find $w\nu \preceq \mu$. Thus we may assign this $\nu$ as the parameter demanded by Clause \ref{lem.decompave.c2}\ref{lem.decompave.c2.3}.\\

In order to instantiate Lemma \ref{lem.shadave} we verify that the data $(s,\ell,\mathscr{M})$ satisfies the hypotheses of Clause \ref{lem.shadave.c1} of that lemma.
\begin{itemize}
\item Clause \ref{lem.decompave.c1}\ref{lem.decompave.c1.2} of the current Lemma \ref{lem.decompave} provides $\ell \leq \beta n$, whence $\beta n - \ell \geq 0$ and so $s \in \enn$.
\item The same clause provides $\ell \geq 2/\gamma$, while Clause \ref{lem.decompave.c1}\ref{lem.decompave.c1.1} of Lemma \ref{lem.decompave} provides $\gamma \leq 1/4$ and hence $1/\gamma \geq 4$. Thus $\ell \geq 8$, so that in particular we have the required inequality $\ell \geq 2$.
\item Since $s \in [\beta n - \ell, n - \ell]_{\zee}$, Clause \ref{lem.decompave.c1}\ref{lem.decompave.c1.4} of Lemma \ref{lem.decompave} provides $\mathtt{Dep}_{\mathscr{M}}(s) \geq \ell$, and so Clause \ref{rem.shadep.c1} of Remark \ref{rem.shadep} asserts that $\mathscr{M}$ is $(s,\ell)$ progression-shadowing.
\end{itemize}
Instantiating Clause \ref{lem.shadave.c1} of Lemma \ref{lem.shadave} with the data $(s,\ell,\mathscr{M})$, we are returned by Clause \ref{lem.shadave.c2} of that lemma a subset $\mathscr{Q}$ of $\mathscr{M} \cap [1,2^{s+\ell+1}]_{\zee}$ satisfying $\lvert \mathscr{Q} \rvert = 2^\ell$. We assign this set as the parameter demanded by Clause \ref{lem.decompave.c2}\ref{lem.decompave.c2.4} of Lemma \ref{lem.decompave}. All parameters demanded by Clause \ref{lem.decompave.c2} of Lemma \ref{lem.decompave} have been assigned, and so in order to complete the current proof it remains to verify the desired inequality (\ref{lem.decompave.conc1}). We also note that Lemma \ref{lem.shadave} remains pending at the entry to Clause \ref{lem.shadave.c3}.\\

Fix $g \in \zee \setminus \{0\}$ throughout the remainder of the current proof, and temporarily fix $k \in [0,2^s-1]_{\zee}$. Clause \ref{lem.flatscale.c5} of Lemma \ref{lem.flatscale} asserts that $\nu_k(\mathscr{P}_{2^s}(k)) = 1$, and so we may instantiate Clause \ref{lem.shadave.c3} of Lemma \ref{lem.shadave} with the data $(k,\nu_k,g)$. Since $k$ was arbitrary, the conclusion (\ref{lem.shadave.conc1}) of that lemma provides the following for all $k \in [0,2^s-1]_{\zee}$.
\begin{align} \frac{1}{\lvert \mathscr{Q} \rvert} \sum_{q \in \mathscr{Q}} \lvert \hat{\nu}_k(qg) \rvert^2 \leq 8 \left( 16 \lvert g \rvert \, \bigl \vert \bigl \vert \nu_k(\mathscr{P}_{2^{s+\ell}}) \bigr \vert \bigr \vert_2^2 + \frac{\pi^2 g^2}{2^{2\ell}} \right) \label{lem.decompave.2} \end{align}
We are now ready to assemble the desired inequality. We compute:
\begin{align} \frac{1}{\lvert \mathscr{Q} \rvert} \sum_{q \in \mathscr{Q}} \lvert \hat{\nu}(qg) \rvert^2 & = \frac{1}{w^2} \cdot \frac{1}{\lvert \mathscr{Q} \rvert} \sum_{q \in \mathscr{Q}} \lvertt \sum_{k=0}^{2^s-1} w_k \hat{\nu}_k(qg) \rvertt^2 \label{lem.decompave.3.1}
\\ & \leq \frac{1}{w^2} \cdot \frac{1}{\lvert \mathscr{Q} \rvert} \sum_{q \in \mathscr{Q}} \left( \sum_{k=0}^{2^s-1} w_k \right)\left( \sum_{k=0}^{2^s-1} w_k \lvert \hat{\nu}_k(qg) \rvert^2 \right) \label{lem.decompave.3.2}
\\ & = \frac{1}{w} \cdot  \sum_{k=0}^{2^s-1} \frac{w_k}{\lvert \mathscr{Q} \rvert} \sum_{q \in \mathscr{Q}}   \lvert \hat{\nu}_k(qg) \rvert^2 \label{lem.decompave.3.3}
\\ & \leq \frac{8}{w} \sum_{k=0}^{2^s-1} w_k \left( 16 \lvert g \rvert \, \bigl \vert \bigl \vert \nu_k(\mathscr{P}_{2^{s+\ell}}) \bigr \vert \bigr \vert_2^2 + \frac{\pi^2 g^2}{2^{2\ell}} \right) \label{lem.decompave.3.5}
\\ & = \frac{128 \lvert g \rvert}{w}\left( \sum_{k=0}^{2^s-1} w_k \, \bigl \vert \bigl \vert \nu_k(\mathscr{P}_{2^{s+\ell}}) \bigr \vert \bigr \vert_2^2\right) +\frac{8\pi^2 g^2}{2^{2\ell}}\label{lem.decompave.3.6}
\\ & \leq \frac{256 \lvert g \rvert}{w \, 2^{\gamma \ell}} + \frac{8 \pi^2 g^2}{2^{2\ell}} \label{lem.decompave.3.7}
\\ & \leq \frac{8}{w^2}\left( \frac{32 \lvert g \rvert}{2^{\gamma \ell}} + \frac{\pi^2 g^2}{2^{2\ell}} \right) \label{lem.decompave.3.8} \end{align}
This computation may be justified as follows.
\begin{itemize}
\item The equality in (\ref{lem.decompave.3.1}) follows directly from the assignment of $\nu$ in (\ref{lem.decompave.1.1}).
\item (\ref{lem.decompave.3.2}) follows from (\ref{lem.decompave.3.1}) since for each $q \in \mathscr{Q}$ we can apply the Cauchy-Schwartz inequality to the lists $\bigl(\sqrt{w_k}\bigr)_{k=0}^{2^s-1}$ and $\bigl(\sqrt{w_k} \, \hat{\nu}_k(qg)\bigr)_{k=0}^{2^s-1}$.
\item (\ref{lem.decompave.3.3}) follows from (\ref{lem.decompave.3.2}) by using (\ref{lem.decompave.1}) and exchanging the sums.
\item (\ref{lem.decompave.3.5}) follows from (\ref{lem.decompave.3.3}) by instantiating (\ref{lem.decompave.2}) with each $k \in [0,2^s-1]_{\zee}$.
\item (\ref{lem.decompave.3.6}) follows from (\ref{lem.decompave.3.5}) by using (\ref{lem.decompave.1}) and factoring.
\item (\ref{lem.decompave.3.7}) follows from (\ref{lem.decompave.3.6}) by Clause \ref{lem.flatscale.c6} of Lemma \ref{lem.flatscale}.
\item (\ref{lem.decompave.3.8}) follows from (\ref{lem.decompave.3.7}) since $w \leq 1$. \end{itemize}
The computation ending in (\ref{lem.decompave.3.8}) verifies the desired inequality (\ref{lem.decompave.conc1}) and thereby completes the proof of Lemma \ref{lem.decompave}.\end{Proof}

\subsection{Extraction of spreading multipliers} \label{cha.pfmain.infra.extract}

Lemma \ref{lem.extract} below converts the averaged-over-$\mathscr{Q}$ estimate from Lemma \ref{lem.decompave} into an estimate holding pointwise on a large subset, using probabilistic concentration over the finite set $\mathscr{Q}$. The pointwise form is what the sequel demands: in Lemma \ref{lem.equildy} the estimate is applied simultaneously at every frequency of the band $(\{m\}\mathscr{F} - \{m\}\mathscr{F}) \setminus \{0\}$, where it is the maximum over that band and not any average across it. The resulting bound on the Fourier max-seminorm of $\nu$ is what Lemma \ref{lem.unifeinf} converts into an exponential infimum.\\

The dependencies in the statement of Lemma \ref{lem.extract} are as follows.
\begin{deps}
\dep{Conventions}{\depitem{Submeasures}{conv.submeasure}}
\dep{Definitions}{\depitem{Progression-shadowing depth}{def.shadep}}
\depcont{\depitem{Resolution entropy}{def.entropy}}
\depcont{\depitem{Fourier coefficients}{def.foucoef}}
\end{deps}

We now move on to the formal statement.

\begin{Lemma}[Extraction of spreading multipliers] \label{lem.extract} \quad
\vspace{-1 \baselineskip}
\begin{enumerate}[label=\upshape(\arabic*)] \item \label{lem.extract.c1} For all:
\begin{enumerate}[label=\upshape(\arabic*)] \item \label{lem.extract.c1.1} $n \in \enn$ and all $\sigma \in (0,1]_{\arr}$ and all $\beta \in (0,1]_{\arr}$ and all $\gamma \in (0,1/4]_{\arr}$ such that $\sigma > \beta + 2\gamma$,
\item \label{lem.extract.c1.2} and all $\ell \in [2/\gamma, \beta n]_{\zee}$,
\item \label{lem.extract.c1.3} and all $\mu \in \prob(\rzee)$ such that $\achn(\mu,2^n) \geq \sigma$,
\item \label{lem.extract.c1.4} and all subsets $\mathscr{M}$ of $\enn$ such that $\mathtt{Dep}_{\mathscr{M}}(s_0) \geq \ell$ for all $s_0 \in [\beta n - \ell, n - \ell]_{\zee}$,
\end{enumerate}
\item \label{lem.extract.c2} there exist:
\begin{enumerate}[label=\upshape(\arabic*)] \item \label{lem.extract.c2.1} $s \in [\beta n - \ell, n - \ell]_{\zee}$,
\item \label{lem.extract.c2.2} and $w \in [(\sigma - \beta - 2\gamma)/(1-2\gamma), 1]_{\arr}$,
\item \label{lem.extract.c2.3} and $\nu \in \prob(\rzee)$ satisfying $w \nu \preceq \mu$,
\item \label{lem.extract.c2.4} and a subset $\mathscr{G}$ of $\mathscr{M} \cap [1,2^{s+\ell+1}]_{\zee}$ satisfying $\lvert \mathscr{G} \rvert \geq 2^\ell/5$, \end{enumerate}
\end{enumerate}
such that for all $m \in \mathscr{G}$ and all $g \in \zee \setminus \{0\}$ we have the inequality below.
\begin{align} \lvert \hat{\nu}(mg) \rvert \leq \frac{20}{w}\left( \frac{4 \lvert g \rvert^{3/2}}{2^{\gamma \ell/2}} + \frac{g^2}{2^{\ell}} \right) \label{lem.extract.conc1} \end{align} \end{Lemma}

An outline of the proof of Lemma \ref{lem.extract} is as follows.
\begin{itemize}
\item The hypotheses in Clause \ref{lem.extract.c1} of the current lemma are verbatim those of Lemma \ref{lem.decompave}, so that lemma is instantiated at once and the parameters $s$, $w$ and $\nu$ pass through unchanged. Everything that follows is a refinement of the multiplier set returned by Lemma \ref{lem.decompave}, and the conclusion of the current lemma differs in trading an average over that set for a pointwise bound on a large subset of it.
\item A threshold function $\Xi$ is calibrated at (\ref{lem.extract.2}) so that the averaged bound supplied by Lemma \ref{lem.decompave} reads as a Chebyshev inequality. We find that for each frequency $g$, the multipliers at which $\lvert \hat{\nu}(qg) \rvert$ exceeds $\Xi(g)$ form a set of size at most $\lvert \mathscr{Q} \rvert / 5g^2$, which is (\ref{lem.extract.5}). The factor $5$ and the decay in $g$ are chosen in advance so that the exceptional sets will later sum to less than the whole.
\item High frequencies may be dismissed since once $\lvert g \rvert \geq 2^{\ell}$ the threshold exceeds $1$ by (\ref{lem.extract.6}), while a Fourier coefficient of a probability measure never does. Thus the exceptional set is empty there and only finitely many frequencies need to be discarded.
\item Discarding these is a union bound, and the sum $\sum_{g =1}^\infty 1/g^2$ keeps the total below four fifths of $\mathscr{Q}$ at (\ref{lem.extract.8.5}). We obtain the set $\mathscr{G}$ demanded by Clause \ref{lem.extract.c2.4}, and on it the threshold bound is put into the stated shape by the computation ending in (\ref{lem.extract.9.3}).
\end{itemize}

The dependencies in the proof of Lemma \ref{lem.extract} beyond those in the statement are as follows.
\begin{deps}
\dep{Lemmas}{\depitem{Averaged Fourier decay for entropic measures}{lem.decompave}}
\end{deps}

We now move on to the formal proof.

\begin{Proof}[of Lemma \ref{lem.extract}] Let $n,\sigma,\beta,\gamma,\ell,\mu$ and $\mathscr{M}$ be as in Clause \ref{lem.extract.c1} of Lemma \ref{lem.extract}. The hypotheses gathered in Clause \ref{lem.decompave.c1} of Lemma \ref{lem.decompave} coincide with those gathered in Clause \ref{lem.extract.c1} of the current lemma, and so we may instantiate the former with the data $(n,\sigma,\beta,\gamma,\ell,\mu,\mathscr{M})$.\\

From Clause \ref{lem.decompave.c2} of Lemma \ref{lem.decompave} we are returned an integer $s \in [\beta n - \ell, n - \ell]_{\zee}$, an element $w$ of the interval named in Clause \ref{lem.decompave.c2}\ref{lem.decompave.c2.2} of Lemma \ref{lem.decompave}, an element $\nu$ of $\prob(\rzee)$ satisfying $w\nu \preceq \mu$, and a subset $\mathscr{Q}$ of $\mathscr{M} \cap [1,2^{s+\ell+1}]_{\zee}$ satisfying $\lvert \mathscr{Q} \rvert = 2^\ell$.\\

We assign $s$, $w$ and $\nu$ as the parameters demanded by Clauses \ref{lem.extract.c2}\ref{lem.extract.c2.1}, \ref{lem.extract.c2}\ref{lem.extract.c2.2} and \ref{lem.extract.c2}\ref{lem.extract.c2.3} of the current lemma respectively, deferring the assignment of $\mathscr{G}$. We note that the bounds in Clause \ref{lem.extract.c1}\ref{lem.extract.c1.1} of Lemma \ref{lem.extract} ensure that $(\sigma - \beta - 2\gamma)/(1-2\gamma) > 0$, and the bound assumed in Clause \ref{lem.decompave.c2}\ref{lem.decompave.c2.2} implies $w > 0$.\\

Now, we define a function $\Xi$ from $\zee \setminus \{0\}$ to $(0,\infty)_{\arr}$ as follows.
\begin{align} \Xi(g) \coloneqq \frac{\lvert g \rvert}{w} \left( \frac{1280 \lvert g \rvert}{2^{\gamma \ell}} + \frac{40 \pi^2 g^2}{2^{2\ell}} \right)^{1/2} \label{lem.extract.1} \end{align}
The assignment of $\Xi$ is arranged so that for all $g \in \zee \setminus \{0\}$ we have the following:
\begin{align} \frac{\Xi(g)^2}{5g^2} = \frac{8}{w^2}\left( \frac{32 \lvert g \rvert}{2^{\gamma \ell}} + \frac{\pi^2 g^2}{2^{2\ell}} \right) \label{lem.extract.2} \end{align}
Temporarily fix $g \in \zee \setminus \{0\}$ and define a subset $\mathscr{A}_g$ of $\mathscr{Q}$ as follows.
\begin{align} \mathscr{A}_g \coloneqq \bigl\{ q \in \mathscr{Q} : \lvert \hat{\nu}(qg) \rvert > \Xi(g) \bigr\} \label{lem.extract.3} \end{align}
We compute:
\begin{align} \lvert \mathscr{A}_g \rvert \, \Xi(g)^2 & \leq \sum_{q \in \mathscr{Q}} \lvert \hat{\nu}(qg) \rvert^2 \label{lem.extract.4.1}
\\ & \leq \lvert \mathscr{Q} \rvert \cdot \frac{8}{w^2}\left( \frac{32 \lvert g \rvert}{2^{\gamma \ell}} + \frac{\pi^2 g^2}{2^{2\ell}} \right) \label{lem.extract.4.2}
\\ & = \frac{\lvert \mathscr{Q} \rvert \, \Xi(g)^2}{5g^2} \label{lem.extract.4.3} \end{align}
This computation may be justified as follows.
\begin{itemize}
\item The inequality in (\ref{lem.extract.4.1}) follows from the definition of $\mathscr{A}_g$ in (\ref{lem.extract.3}), since the summands are nonnegative and since each $q \in \mathscr{A}_g$ contributes a summand exceeding $\Xi(g)^2$.
\item (\ref{lem.extract.4.2}) follows from (\ref{lem.extract.4.1}) by the conclusion (\ref{lem.decompave.conc1}) of Lemma \ref{lem.decompave}.
\item (\ref{lem.extract.4.3}) follows from (\ref{lem.extract.4.2}) by (\ref{lem.extract.2}).
\end{itemize}
Since $\Xi(g) > 0$ we may divide the computation ending in (\ref{lem.extract.4.3}) by $\Xi(g)^2$, and since $g$ was arbitrary we obtain the following for all $g \in \zee \setminus \{0\}$.
\begin{align} \lvert \mathscr{A}_g \rvert \leq \frac{\lvert \mathscr{Q} \rvert}{5g^2} \label{lem.extract.5} \end{align}
We next observe that $\mathscr{A}_g$ is empty whenever $\lvert g \rvert \geq 2^\ell$. Indeed, discarding the first summand in the assignment of $\Xi$ from (\ref{lem.extract.1}) and using $w \leq 1$ we find that the following holds for all $g \in \zee \setminus \{0\}$.
\begin{align} \Xi(g) \geq \frac{\lvert g \rvert}{w} \left( \frac{40 \pi^2 g^2}{2^{2\ell}} \right)^{1/2} = \frac{\sqrt{40}\, \pi \, g^2}{w \, 2^{\ell}} \geq \frac{19 g^2}{2^{\ell}} \label{lem.extract.6} \end{align}
Thus if $\lvert g \rvert \geq 2^\ell$ then $\Xi(g) \geq 19 \cdot 2^{\ell} > 1$, while $\lvert \hat{\nu}(qg) \rvert \leq 1$ for all $q \in \mathscr{Q}$ since $\nu \in \prob(\rzee)$; the claim then follows from the definition of $\mathscr{A}_g$ in (\ref{lem.extract.3}).\\

We now assign the subset $\mathscr{G}$ demanded by Clause \ref{lem.extract.c2}\ref{lem.extract.c2.4} of Lemma \ref{lem.extract} as follows.
\begin{align} \mathscr{G} \coloneqq \mathscr{Q} \setminus \bigcup_{g \in [-2^\ell,2^\ell]_{\zee} \setminus \{0\}} \mathscr{A}_g \label{lem.extract.7} \end{align}
Since $\mathscr{G}$ is a subset of $\mathscr{Q}$, it is a subset of $\mathscr{M} \cap [1,2^{s+\ell+1}]_{\zee}$ as required. In order to bound the cardinality of $\mathscr{G}$ from below, we compute:
\begin{align} \lvertt \bigcup_{g \in [-2^\ell,2^\ell]_{\zee} \setminus \{0\}} \mathscr{A}_g \rvertt & \leq \sum_{g \in [-2^\ell,2^\ell]_{\zee} \setminus \{0\}} \lvert \mathscr{A}_g \rvert \label{lem.extract.8.1}
\\ & \leq \frac{\lvert \mathscr{Q} \rvert}{5} \sum_{g \in [-2^\ell,2^\ell]_{\zee} \setminus \{0\}} \frac{1}{g^2} \label{lem.extract.8.2}
\\ & = \frac{2\lvert \mathscr{Q} \rvert}{5} \sum_{g=1}^{2^\ell} \frac{1}{g^2} \label{lem.extract.8.3}
\\ & \leq \frac{\pi^2 \lvert \mathscr{Q} \rvert}{15} \label{lem.extract.8.4}
\\ & \leq \frac{4 \lvert \mathscr{Q} \rvert}{5} \label{lem.extract.8.5} \end{align}
This computation may be justified as follows.
\begin{itemize}
\item (\ref{lem.extract.8.2}) follows from (\ref{lem.extract.8.1}) by (\ref{lem.extract.5}).
\item (\ref{lem.extract.8.3}) follows from (\ref{lem.extract.8.2}) since $g \mapsto g^{-2}$ is an even function on $\zee$.
\item (\ref{lem.extract.8.4}) follows from (\ref{lem.extract.8.3}) since the solution to the Basel problem provides:
\begin{align*} \sum_{g=1}^{\infty} \frac{1}{g^2} = \frac{\pi^2}{6} \end{align*}
\item (\ref{lem.extract.8.5}) follows from (\ref{lem.extract.8.4}) since $\pi^2 \leq 12$.
\end{itemize}
Combining the computation ending in (\ref{lem.extract.8.5}) with the assignment of $\mathscr{G}$ in (\ref{lem.extract.7}) and with the equality $\lvert \mathscr{Q} \rvert = 2^\ell$ we find $\lvert \mathscr{G} \rvert \geq \lvert \mathscr{Q} \rvert/5 = 2^{\ell}/5$, and so the assignment demanded by Clause \ref{lem.extract.c2}\ref{lem.extract.c2.4} of Lemma \ref{lem.extract} is complete.\\

It remains to verify the conclusion (\ref{lem.extract.conc1}) of Lemma \ref{lem.extract}, and to that end we fix $m \in \mathscr{G}$ and $g \in \zee \setminus \{0\}$. We compute:
\begin{align} \Xi(g) & \leq \frac{\lvert g \rvert}{w} \left( \left( \frac{1280 \lvert g \rvert}{2^{\gamma \ell}} \right)^{1/2} + \left( \frac{40 \pi^2 g^2}{2^{2\ell}} \right)^{1/2} \right) \label{lem.extract.9.1}
\\ & = \frac{1}{w} \left( \frac{16\sqrt{5}\, \lvert g \rvert^{3/2}}{2^{\gamma \ell/2}} + \frac{\sqrt{40}\, \pi \, g^2}{2^{\ell}} \right) \label{lem.extract.9.2}
\\ & \leq \frac{20}{w} \left( \frac{4 \lvert g \rvert^{3/2}}{2^{\gamma \ell/2}} + \frac{g^2}{2^{\ell}} \right) \label{lem.extract.9.3} \end{align}
This computation may be justified as follows.
\begin{itemize}
\item The inequality in (\ref{lem.extract.9.1}) follows by combining the assignment of $\Xi$ in (\ref{lem.extract.1}) with the elementary bound $(a+b)^{1/2} \leq a^{1/2}+b^{1/2}$, valid for all $a,b \in [0,\infty)_{\arr}$.
\item (\ref{lem.extract.9.2}) follows from (\ref{lem.extract.9.1}) by distributing the factor $\lvert g \rvert / w$ across the bracket and evaluating the two square roots.
\item (\ref{lem.extract.9.3}) follows from (\ref{lem.extract.9.2}) since $16\sqrt{5} \leq 80$ and $\sqrt{40}\,\pi \leq 20$.
\end{itemize}
 We claim that $\lvert \hat{\nu}(mg) \rvert \leq \Xi(g)$. Indeed, if $\lvert g \rvert \leq 2^\ell$ then the assignment of $\mathscr{G}$ in (\ref{lem.extract.7}) provides $m \notin \mathscr{A}_g$, so that the claim follows from the definition of $\mathscr{A}_g$ in (\ref{lem.extract.3}); while if $\lvert g \rvert > 2^\ell$ then (\ref{lem.extract.6}) provides $\Xi(g) \geq 19 \cdot 2^{\ell} \geq 1$, so that the claim follows since $\lvert \hat{\nu}(mg) \rvert \leq 1$. Combining the computation ending in (\ref{lem.extract.9.3}) with the last-verified claim we obtain the desired inequality (\ref{lem.extract.conc1}), and thereby complete the proof of Lemma \ref{lem.extract}. \end{Proof}

\subsection{Endomorphic rigidity at dyadic resolution} \label{cha.pfmain.infra.equildy}

Lemma \ref{lem.equildy} below carries the content of Theorem \ref{thm.main} at dyadic scales. What remains in the proof of that theorem is the passage to an arbitrary resolution $N$, effected by Clause \ref{olem.entropy.c5} of Lemma \ref{olem.entropy}. This costs an additive $2/n$ in the entropy level, which is where it appears in the conclusion of that theorem.

The dependencies in the statement of Lemma \ref{lem.equildy} are as follows.
\begin{deps}
\dep{Conventions}{\depitem{Dilations of sets of integers}{conv.setdil}}
\dep{Definitions}{\depitem{Progression-shadowing depth}{def.shadep}}
\depcont{\depitem{Resolution entropy}{def.entropy}}
\depcont{\depitem{Exponential infimum}{def.foustat}}
\end{deps}

We now move on to the formal statement.

\begin{Lemma}[Endomorphic rigidity at dyadic resolution] \label{lem.equildy} \quad
\vspace{-1 \baselineskip}
\begin{enumerate}[label=\upshape(\arabic*)] \item \label{lem.equildy.c1} For all:
\begin{enumerate}[label=\upshape(\arabic*)] \item \label{lem.equildy.c1.1} $n \in \enn$ and all $\sigma \in (0,1]_{\arr}$ and all $\beta \in (0,1]_{\arr}$ and all $\gamma \in (0,1/4]_{\arr}$ such that $\sigma > \beta + 2\gamma$,
\item \label{lem.equildy.c1.2} and all $\ell \in [2/\gamma, \beta n]_{\zee}$,
\item \label{lem.equildy.c1.3} and all $\mu \in \prob(\rzee)$ such that $\achn(\mu,2^n) \geq \sigma$,
\item \label{lem.equildy.c1.4} and all subsets $\mathscr{M}$ of $\enn$ such that $\mathtt{Dep}_{\mathscr{M}}(s_0) \geq \ell$ for all $s_0 \in [\beta n - \ell, n - \ell]_{\zee}$,
\end{enumerate}
\item \label{lem.equildy.c2} there exists a subset $\mathscr{G}$ of $\mathscr{M} \cap [1,2^{n+1}]_{\zee}$ satisfying $\lvert \mathscr{G} \rvert \geq 2^\ell/5$,
\end{enumerate}
such that for all $m \in \mathscr{G}$ and all nonempty finite subsets $\mathscr{F}$ of $\enn$ we have the inequality below.
\begin{align} \einf\bigl(\mu,\{m\}\mathscr{F}\bigr) \geq \frac{\sigma - \beta - 2\gamma}{1-2\gamma} - 160\left( \frac{\max(\mathscr{F})^{7/2}}{2^{\gamma \ell/2}} + \frac{\max(\mathscr{F})^{4}}{2^{\ell}} \right) \label{lem.equildy.conc1} \end{align} \end{Lemma}

An outline of the proof of Lemma \ref{lem.equildy} is as follows.
\begin{itemize}
\item As with Lemma \ref{lem.extract}, the hypotheses in Clause \ref{lem.equildy.c1} of the current lemma are verbatim those of the lemma preceding it, so Lemma \ref{lem.extract} is instantiated at once. The multiplier set it returns already has the required cardinality, and since $s \leq n - \ell$ its ambient interval $[1,2^{s+\ell+1}]_{\zee}$ is contained in $[1,2^{n+1}]_{\zee}$, which discharges Clause \ref{lem.equildy.c2}.
\item The degenerate bandprofile $\mathscr{F} = \{0\}$ is disposed of first, since Remark \ref{rem.einfsing} makes the left side of the conclusion equal to $1$ there.
\item Every frequency of the band $(\{m\}\mathscr{F} - \{m\}\mathscr{F}) \setminus \{0\}$ is of the form $mg$ with $0 < \lvert g \rvert \leq \max(\mathscr{F})$, and the bound supplied by Lemma \ref{lem.extract} is nondecreasing in $\lvert g \rvert$. A single evaluation at $\max(\mathscr{F})$ therefore dominates the whole band, which is the Fourier max-seminorm estimate (\ref{lem.equildy.1}).
\item The assembly returns from $\nu$ to $\mu$ across the submeasure relation, converts by Lemma \ref{lem.unifeinf}, and pays a crude count of $\{m\}\mathscr{F}$ for the powers of $\max(\mathscr{F})$. At (\ref{lem.equildy.2.4}) the factor $w$ gained on the way back cancels the factor $1/w$ carried by Lemma \ref{lem.extract}.
\end{itemize}

The dependencies in the proof of Lemma \ref{lem.equildy} beyond those in the statement are as follows.
\begin{deps}
\dep{Conventions}{\depitem{Submeasures}{conv.submeasure}\depsep\depitem{Sums and differences of sets of integers}{conv.setarith}}
\dep{Definitions}{\depitem{Finite-bandprofile Fourier max-seminorm}{def.founorm}}
\dep{Remarks}{\depitem{Exponential infimum and singleton sets}{rem.einfsing}\depsep\depitem{Exponential infimum and submeasures}{rem.einfmono}}
\dep{Lemmas}{\depitem{Exponential infimum forced by Fourier decay}{lem.unifeinf}\depsep\depitem{Extraction of spreading multipliers}{lem.extract}}
\end{deps}

We now move on to the formal proof.

\begin{Proof}[of Lemma \ref{lem.equildy}] Let $n,\sigma,\beta,\gamma,\ell,\mu$ and $\mathscr{M}$ be as in Clause \ref{lem.equildy.c1} of Lemma \ref{lem.equildy}. The hypotheses gathered in Clause \ref{lem.extract.c1} of Lemma \ref{lem.extract} coincide with those gathered in Clause \ref{lem.equildy.c1} of the current lemma, and so we may instantiate the former with the data $(n,\sigma,\beta,\gamma,\ell,\mu,\mathscr{M})$.\\

From Clause \ref{lem.extract.c2} of Lemma \ref{lem.extract} we are returned an integer $s \in [\beta n - \ell, n - \ell]_{\zee}$, an element $w$ of the interval named in Clause \ref{lem.extract.c2}\ref{lem.extract.c2.2} of Lemma \ref{lem.extract}, an element $\nu$ of $\prob(\rzee)$ satisfying $w \nu \preceq \mu$, and a subset $\mathscr{G}$ of $\mathscr{M} \cap [1,2^{s+\ell+1}]_{\zee}$ satisfying $\lvert \mathscr{G} \rvert \geq 2^\ell/5$. Since $s \leq n - \ell$ we have $2^{s+\ell+1} \leq 2^{n+1}$, so that $\mathscr{G}$ is a subset of $\mathscr{M} \cap [1,2^{n+1}]_{\zee}$, and we assign $\mathscr{G}$ as the set demanded by Clause \ref{lem.equildy.c2} of the current lemma.\\

It remains to verify the conclusion (\ref{lem.equildy.conc1}) of Lemma \ref{lem.equildy}, and to that end we fix $m \in \mathscr{G}$ together with a nonempty finite subset $\mathscr{F}$ of $\enn$, and we adopt the notation $P \coloneqq \max(\mathscr{F})$. We now dispose of the degenerate case when $P = 0$. In this case $\mathscr{F} = \{0\}$ and hence $\{m\}\mathscr{F} = \{0\}$, so that Remark \ref{rem.einfsing} provides $\einf(\mu,\{m\}\mathscr{F}) = 1$. Since the right side of (\ref{lem.equildy.conc1}) is then equal to $(\sigma - \beta - 2\gamma)/(1-2\gamma)$, which is at most $1$ because $\sigma \leq 1$ and $\beta > 0$, the desired inequality holds. We may therefore assume that $P \geq 1$.\\

Since $m \geq 1$, the set $(\{m\}\mathscr{F} - \{m\}\mathscr{F}) \setminus \{0\}$ consists of the integers $mg$ with $g \in (\mathscr{F}-\mathscr{F}) \setminus \{0\}$, and each such $g$ satisfies $0 < \lvert g \rvert \leq P$ because $\mathscr{F}$ is a subset of $[0,P]_{\zee}$. Since the right side of the conclusion (\ref{lem.extract.conc1}) of Lemma \ref{lem.extract} is nondecreasing in $\lvert g \rvert$, combining that conclusion with the definition of the Fourier max-seminorm in Definition \ref{def.founorm} we obtain:
\begin{align} \bigl \vert \bigl \vert \nu \bigr \vert \bigr \vert_{\fou\left(\left(\{m\}\mathscr{F} - \{m\}\mathscr{F}\right) \setminus \{0\}\right)} \leq \frac{20}{w}\left( \frac{4 P^{3/2}}{2^{\gamma \ell/2}} + \frac{P^{2}}{2^{\ell}} \right) \label{lem.equildy.1} \end{align}
We now compute:
\begin{align} \einf\bigl(\mu,\{m\}\mathscr{F}\bigr) & \geq w \,\einf\bigl(\nu,\{m\}\mathscr{F}\bigr) \label{lem.equildy.2.1}
\\ & \geq w \left( 1 - \left( \lvert \{m\}\mathscr{F} \rvert^2 - \lvert \{m\}\mathscr{F} \rvert \right) \bigl \vert \bigl \vert \nu \bigr \vert \bigr \vert_{\fou\left(\left(\{m\}\mathscr{F} - \{m\}\mathscr{F}\right) \setminus \{0\}\right)} \right) \label{lem.equildy.2.2}
\\ & \geq w \left( 1 - \frac{40P^2}{w}\left( \frac{4 P^{3/2}}{2^{\gamma \ell/2}} + \frac{P^{2}}{2^{\ell}} \right) \right) \label{lem.equildy.2.3}
\\ & = w - 40\left( \frac{4 P^{7/2}}{2^{\gamma \ell/2}} + \frac{P^{4}}{2^{\ell}} \right) \label{lem.equildy.2.4}
\\ & \geq \frac{\sigma - \beta - 2\gamma}{1-2\gamma} - 160 \left( \frac{P^{7/2}}{2^{\gamma \ell/2}} + \frac{P^{4}}{2^{\ell}} \right) \label{lem.equildy.2.5} \end{align}
This computation may be justified as follows.
\begin{itemize}
\item The inequality in (\ref{lem.equildy.2.1}) follows from Remark \ref{rem.einfmono}, applied with the submeasure relation $w\nu \preceq \mu$ returned above.
\item (\ref{lem.equildy.2.2}) follows from (\ref{lem.equildy.2.1}) by the conclusion (\ref{lem.unifeinf.conc}) of Lemma \ref{lem.unifeinf}, applied to $\nu$ and $\{m\}\mathscr{F}$, where multiplication by $w$ preserves the inequality since $w \geq 0$.
\item (\ref{lem.equildy.2.3}) follows from (\ref{lem.equildy.2.2}) by combining (\ref{lem.equildy.1}) with the observation that $\lvert \{m\}\mathscr{F} \rvert \leq \lvert \{m\}[0,P]_{\zee} \rvert = P+1$ and so:
\begin{align*} \lvert \{m\}\mathscr{F} \rvert^2 - \lvert \{m\}\mathscr{F} \rvert \leq (P+1)P \leq 2P^2 \end{align*}
\item (\ref{lem.equildy.2.5}) follows from (\ref{lem.equildy.2.4}) since Clause \ref{lem.extract.c2}\ref{lem.extract.c2.2} of Lemma \ref{lem.extract} provides $w \geq (\sigma - \beta - 2\gamma)/(1-2\gamma)$, and since $40 \cdot 4 \leq 160$ and $40 \leq 160$.
\end{itemize}
Unwrapping the definition of $P$ from the computation ending in (\ref{lem.equildy.2.5}) verifies the desired inequality (\ref{lem.equildy.conc1}) and thereby completes the proof of Lemma \ref{lem.equildy}.\end{Proof}

\section{Proof of Theorem \ref{thm.main}: main} \label{cha.pfmain.main}

An outline of the proof of Theorem \ref{thm.main} is as follows.
\begin{itemize}
\item The theorem is stated at an arbitrary resolution $N$ whereas Lemma \ref{lem.equildy} works at dyadic resolutions. Thus we descend to $2^n$ for $n = \lfloor \log_2(N) \rfloor$ at a cost bounded by Clause \ref{olem.entropy.c5} of Lemma \ref{olem.entropy}. The computation ending in (\ref{thm.main.p2.5}) records the outcome as the reduced entropy level $\sigma_0$ assigned in (\ref{thm.main.p1}).
\item  To produce the multiplier set we use Lemma \ref{lem.shadave} directly with the depth hypothesis at the single scale $n-\ell$, which is the upper endpoint of the window over which that hypothesis is assumed. The remainder of the proof runs under (\ref{thm.main.p3}).
\item Under (\ref{thm.main.p3}) the reduced data $(n,\sigma_0,\beta,\gamma,\ell,\mu,\mathscr{M})$ satisfies Clause \ref{lem.equildy.c1} of Lemma \ref{lem.equildy}, and the multiplier set returned by that lemma is the one demanded here. Because the loss $2/n$ was charged to the entropy level rather than to the conclusion, unwrapping $\sigma_0$ in that lemma's conclusion reproduces the desired inequality (\ref{thm.main.conc1}).\end{itemize}

The dependencies in the proof of Theorem \ref{thm.main} beyond those in the statement are as follows.
\begin{deps}
\dep{Remarks}{\depitemc{Shadowing depth and progression-shadowing}{rem.shadep}{rem.shadep.c1}}
\dep{Lemmas}{\depitem{Shadow-averaging lemma}{lem.shadave}\depsep\depitem{Endomorphic rigidity at dyadic resolution}{lem.equildy}}
\dep{Imports}{\depitemc{Basic estimates for resolution entropy}{olem.entropy}{olem.entropy.c5}}
\end{deps}

We now move on to the formal proof.

\begin{Proof}[of Theorem \ref{thm.main}] Let $N,\sigma,\beta,\gamma,\ell,\mu$ and $\mathscr{M}$ be as in Clause \ref{thm.main.c1} of Theorem \ref{thm.main}, and recall from that clause the notation $n = \lfloor \log_2(N) \rfloor$. Since $N > 2^{100}$ we have $n \geq 100$, and by the definition of $n$ we have $2^n \leq N < 2^{n+1}$. We define $\sigma_0 \in \arr$ as follows.
\begin{align} \sigma_0 \coloneqq \sigma - \frac{2}{n} \label{thm.main.p1} \end{align}
We compute:
\begin{align} \achn(\mu,2^n) & = \frac{\ach(\mu,2^n)}{n \log(2)} \label{thm.main.p2.1}
\\ & \geq \frac{\ach(\mu,N) - \log(3)}{n \log(2)} \label{thm.main.p2.2}
\\ & \geq \frac{\sigma \log(N) - \log(3)}{n \log(2)} \label{thm.main.p2.3}
\\ & \geq \sigma - \frac{\log(3)}{n \log(2)} \label{thm.main.p2.4}
\\ & \geq \sigma_0 \label{thm.main.p2.5} \end{align}
This computation may be justified as follows.
\begin{itemize}
\item The equality in (\ref{thm.main.p2.1}) follows from the definition of the normalized resolution entropy in Definition \ref{def.entropy}, the case $n \geq 100$ ensuring that $2^n \neq 1$.
\item (\ref{thm.main.p2.2}) follows from (\ref{thm.main.p2.1}) by instantiating Clause \ref{olem.entropy.c5} of Lemma \ref{olem.entropy} with the resolutions $2^n$ and $N$, where we note that $N < 2^{n+1}$ provides $\lceil N 2^{-n} \rceil \leq 2$.
\item (\ref{thm.main.p2.3}) follows from (\ref{thm.main.p2.2}) by combining the hypothesis in Clause \ref{thm.main.c1}\ref{thm.main.c1.4} with the definition of the normalized resolution entropy in Definition \ref{def.entropy}.
\item (\ref{thm.main.p2.4}) follows from (\ref{thm.main.p2.3}) since $\log(N) \geq n \log(2)$.
\item (\ref{thm.main.p2.5}) follows from (\ref{thm.main.p2.4}) by combining the assignment of $\sigma_0$ in (\ref{thm.main.p1}) with the inequality $\log(3) \leq 2\log(2)$.
\end{itemize}
We now dispose of a degenerate case when $\sigma_0 \leq \beta + 2\gamma$. Toward producing the set demanded by Clause \ref{thm.main.c2} of Theorem \ref{thm.main} in this case, we verify that the data $(n-\ell,\ell,\mathscr{M})$ satisfies the hypotheses gathered in Clause \ref{lem.shadave.c1} of Lemma \ref{lem.shadave}.
\begin{itemize}
\item Clause \ref{thm.main.c1}\ref{thm.main.c1.3} of Theorem \ref{thm.main} provides $\ell \leq \beta n$, while Clause \ref{thm.main.c1}\ref{thm.main.c1.2} provides $\beta \leq 1$. Thus $\ell \leq n$, whence $n - \ell \in \enn$.
\item Clause \ref{thm.main.c1}\ref{thm.main.c1.3} also provides $\ell \geq 2/\gamma$, while Clause \ref{thm.main.c1}\ref{thm.main.c1.2} provides $\gamma \leq 1/4$ and hence $1/\gamma \geq 4$. Thus $\ell \geq 8$, so that in particular we have the required inequality $\ell \geq 2$.
\item Since $n - \ell$ lies in $[\beta n - \ell, n - \ell]_{\zee}$, Clause \ref{thm.main.c1}\ref{thm.main.c1.5} of Theorem \ref{thm.main} provides $\mathtt{Dep}_{\mathscr{M}}(n-\ell) \geq \ell$, and so Clause \ref{rem.shadep.c1} of Remark \ref{rem.shadep} asserts that $\mathscr{M}$ is $(n-\ell,\ell)$ progression-shadowing.
\end{itemize}
Instantiating Clause \ref{lem.shadave.c1} of Lemma \ref{lem.shadave} with the data $(n-\ell,\ell,\mathscr{M})$, we are returned by Clause \ref{lem.shadave.c2} of that lemma a subset $\mathscr{Q}$ of $\mathscr{M} \cap [1,2^{n+1}]_{\zee}$ satisfying $\lvert \mathscr{Q} \rvert = 2^\ell$. Assigning $\mathscr{G} \coloneqq \mathscr{Q}$ therefore discharges Clause \ref{thm.main.c2} of Theorem \ref{thm.main}. Moreover Clause \ref{thm.main.c1}\ref{thm.main.c1.2} of Theorem \ref{thm.main} provides $\gamma \leq 1/4$, so that $1-2\gamma \geq 1/2 > 0$, and unwrapping the assignment of $\sigma_0$ in (\ref{thm.main.p1}) we obtain:
\begin{align*} \frac{1}{1-2\gamma}\left(\sigma - \beta - 2\gamma - \frac{2}{n}\right) = \frac{\sigma_0 - \beta - 2\gamma}{1-2\gamma} \leq 0 \end{align*}
Here the inequality follows from the hypothesis $\sigma_0 \leq \beta + 2\gamma$ of the current case. Since the term subtracted in (\ref{thm.main.conc1}) is nonnegative, the previous display bounds the right side of (\ref{thm.main.conc1}) above by $0$, while it is immediate from Definition \ref{def.foustat} that $\einf\bigl(\mu,\{m\}\mathscr{F}\bigr) \geq 0$ for every $m \in \mathscr{G}$ and every nonempty finite subset $\mathscr{F}$ of $\enn$. Therefore in this case the desired inequality (\ref{thm.main.conc1}) is verified, and we may assume for the remainder of the current proof that the following holds.
\begin{align} \sigma_0 > \beta + 2\gamma \label{thm.main.p3} \end{align}
In order to instantiate Lemma \ref{lem.equildy} we now verify that the data $(n,\sigma_0,\beta,\gamma,\ell,\mu,\mathscr{M})$ satisfies the hypotheses gathered in Clause \ref{lem.equildy.c1} of that lemma.
\begin{itemize}
\item Since $n \geq 100$ we have $n \in \enn$. Clause \ref{thm.main.c1}\ref{thm.main.c1.2} of Theorem \ref{thm.main} provides $\beta > 0$ and $\gamma > 0$, so combining these with (\ref{thm.main.p3}) we find $\sigma_0 > 0$, while combining the assignment of $\sigma_0$ in (\ref{thm.main.p1}) with the bound $\sigma \leq 1$ from the same clause we find $\sigma_0 \leq 1$. Thus $\sigma_0 \in (0,1]_{\arr}$.
\item Clause \ref{thm.main.c1}\ref{thm.main.c1.2} of Theorem \ref{thm.main} further provides $\beta \in (0,1]_{\arr}$ and $\gamma \in (0,1/4]_{\arr}$, while (\ref{thm.main.p3}) provides the inequality $\sigma_0 > \beta + 2\gamma$, all as demanded by Clause \ref{lem.equildy.c1}\ref{lem.equildy.c1.1} of Lemma \ref{lem.equildy}.
\item Theorem \ref{thm.main} hypothesizes $\ell \in [2/\gamma, \beta n]_{\zee}$, as demanded by Clause \ref{lem.equildy.c1}\ref{lem.equildy.c1.2}.
\item The computation ending in (\ref{thm.main.p2.5}) provides $\achn(\mu,2^n) \geq \sigma_0$, as demanded by Clause \ref{lem.equildy.c1}\ref{lem.equildy.c1.3}.
\item Theorem \ref{thm.main} hypothesizes that $\mathtt{Dep}_{\mathscr{M}}(s_0) \geq \ell$ for all $s_0 \in [\beta n - \ell, n - \ell]_{\zee}$, as demanded by Clause \ref{lem.equildy.c1}\ref{lem.equildy.c1.4} of Lemma \ref{lem.equildy}.
\end{itemize}
Instantiating Clause \ref{lem.equildy.c1} of Lemma \ref{lem.equildy} with the data $(n,\sigma_0,\beta,\gamma,\ell,\mu,\mathscr{M})$, we are returned by Clause \ref{lem.equildy.c2} of that lemma a subset $\mathscr{G}$ of $\mathscr{M} \cap [1,2^{n+1}]_{\zee}$ satisfying $\lvert \mathscr{G} \rvert \geq 2^\ell/5$, and we assign $\mathscr{G}$ as the set demanded by Clause \ref{thm.main.c2} of Theorem \ref{thm.main}. Fixing $m \in \mathscr{G}$ together with a nonempty finite subset $\mathscr{F}$ of $\enn$, the conclusion (\ref{lem.equildy.conc1}) of Lemma \ref{lem.equildy} provides:
\begin{align*} \einf\bigl(\mu,\{m\}\mathscr{F}\bigr) \geq \frac{\sigma_0 - \beta - 2\gamma}{1-2\gamma} - 160\left( \frac{\max(\mathscr{F})^{7/2}}{2^{\gamma \ell/2}} + \frac{\max(\mathscr{F})^{4}}{2^{\ell}} \right) \end{align*}
Unwrapping the assignment of $\sigma_0$ in (\ref{thm.main.p1}) from the previous display verifies the desired inequality (\ref{thm.main.conc1}) and thereby completes the proof of Theorem \ref{thm.main}.\end{Proof}

\chapter{Proofs of subsequent theorems} \label{cha.pfapp}

\section{Proof of Theorem \ref{thm.exist}} \label{sec.pfapp.special.exist}

An outline of the proof of Theorem \ref{thm.exist} is as follows.
\begin{itemize}
\item We initially dispose of the degenerate case $\sigma \leq \epsilon$, where the conclusion demands nothing, and the proof thereafter runs under $\sigma > \epsilon$.
\item Otherwise the theorem is a change of parameters into Theorem \ref{thm.main}, with $\beta \coloneqq \epsilon/8$ and $\gamma \coloneqq \epsilon/\bigl(16(1-\sigma+\epsilon)\bigr)$. The thresholds in Clauses \ref{thm.exist.c1.2} and \ref{thm.exist.c1.5} are engineered backwards from this instantiation. Thus the threshold on $N$ forces $N > 2^{100}$ together with the window bound (\ref{thm.exist.p3}), while the threshold on $\ell$ forces the floor $\ell \geq 2/\gamma$ at (\ref{thm.exist.p1}) and the exponential clearance (\ref{thm.exist.p4}). Each hypothesis of Theorem \ref{thm.main} is then read off from exactly one of these facts. 
\item Any member of the returned multiplier set serves, and what remains is error accounting at (\ref{thm.exist.p5}). The rearrangement at (\ref{thm.exist.p6.1}) separates the flat cost $\beta$ from the cost $2\gamma(1-\sigma)$ charged against the entropy deficit; those two together with $2/n$ are each held to $\epsilon/8$ before the factor $8/7$ is paid, while the two terms carrying $160$ are held to $2\epsilon/7$ jointly at (\ref{thm.exist.p7.3}). This is the step the numerals $128$, $8$ and $9$ in Clauses \ref{thm.exist.c1.2} and \ref{thm.exist.c1.5} exist to certify.
\end{itemize}

The dependencies in the proof of Theorem \ref{thm.exist} beyond those in the statement are as follows.
\begin{deps}
\dep{Definitions}{\depitem{Progression-shadowing sets}{def.progshad}}
\dep{Remarks}{\depitemc{Shadowing depth and progression-shadowing}{rem.shadep}{rem.shadep.c1}}
\dep{Theorems}{\depitem{Quantitative endomorphic rigidity via entropy}{thm.main}}
\end{deps}

We now move on to the formal proof.

\begin{Proof}[of Theorem \ref{thm.exist}] Let $\epsilon,r,\sigma,\ell,\mathscr{M},t,N$ and $\mu$ be as in Clause \ref{thm.exist.c1} of Theorem \ref{thm.exist}. Clause \ref{thm.exist.c1}\ref{thm.exist.c1.1} provides $r \in \enn$ and $\epsilon \leq 1$, so that $8(r+2)/\epsilon \geq 16$ and hence $\log_2\bigl(8(r+2)/\epsilon\bigr) \geq 4$. Combining this bound with the bound $\sigma \leq 1$ from the same clause, the hypothesis of Clause \ref{thm.exist.c1}\ref{thm.exist.c1.2} gives:
\begin{align} \ell \geq \frac{512(1-\sigma+\epsilon)}{\epsilon} \geq 512 \label{thm.exist.p1} \end{align}
In particular we have $\ell \geq 1$, and so Clause \ref{thm.exist.c1}\ref{thm.exist.c1.4} of Theorem \ref{thm.exist} at the scale $t$ provides $\mathtt{Dep}_{\mathscr{M}}(t) \geq \ell$. Thus Clause \ref{rem.shadep.c1} of Remark \ref{rem.shadep} asserts that $\mathscr{M}$ is $(t,\ell)$ progression-shadowing, and instantiating Definition \ref{def.progshad} with $m \coloneqq 1$ returns $q \in \mathscr{M}$ satisfying $\lvert q - 2^t \rvert \leq 2^{t-\ell}$. Consequently $0 < 2^t - 2^{t-\ell} \leq q \leq 2^t + 2^{t-\ell} \leq 2^{t+1}$. Clause \ref{thm.exist.c1}\ref{thm.exist.c1.5} of Theorem \ref{thm.exist} provides $N \geq 2^{9(t+\ell)/\epsilon}$, and since $\epsilon \leq 1$ the exponent is at least $t$, so that $N \geq 2^t$. We conclude that $q \in \mathscr{M} \cap [1,2N]_{\zee}$.\\

We now dispose of a degenerate case when $\sigma \leq \epsilon$. Assigning $m \coloneqq q$, it is immediate from Definition \ref{def.foustat} that $\einf\bigl(\mu,\{m\}[0,r]_{\zee}\bigr) \geq 0$. Therefore in this case the desired inequality (\ref{thm.exist.conc1}) is verified as $\sigma - \epsilon \leq 0$, and we may assume for the remainder of the current proof that $\sigma > \epsilon$.\\

We now assign the data with which Theorem \ref{thm.main} will be instantiated as follows.
\begin{align*} \beta &\coloneqq \frac{\epsilon}{8} & \gamma &\coloneqq \frac{\epsilon}{16(1-\sigma+\epsilon)} & \mathscr{F} &\coloneqq [0,r]_{\zee} \end{align*}
Since Clause \ref{thm.exist.c1}\ref{thm.exist.c1.1} of Theorem \ref{thm.exist} provides $\sigma \leq 1$ we have $1-\sigma+\epsilon \geq \epsilon > 0$, so that $\gamma$ is well-defined and satisfies $\gamma \leq 1/16$.\\

We now verify that $\beta + 2\gamma < \sigma/4$. Combining the bound $\sigma \leq 1$ with the assumption $\sigma > \epsilon$ we find $(1-\sigma)(\sigma-\epsilon) \geq 0$, so that $\epsilon \leq \sigma(1-\sigma+\epsilon)$ and therefore $\epsilon/(1-\sigma+\epsilon) \leq \sigma$. Unwrapping the assignment of $\gamma$ we conclude $2\gamma \leq \sigma/8$, while unwrapping the assignment of $\beta$ and again citing $\sigma > \epsilon$ we conclude $\beta < \sigma/8$ and so $\beta + 2\gamma < \sigma/4$ as required.\\

Adopting the notation $n \coloneqq \lfloor \log_2(N) \rfloor$ from Clause \ref{thm.main.c1}\ref{thm.main.c1.1} of Theorem \ref{thm.main}, Clause \ref{thm.exist.c1}\ref{thm.exist.c1.5} of Theorem \ref{thm.exist} provides $N \geq 2^{9(t+\ell)/\epsilon}$ and hence $n \geq 9(t+\ell)/\epsilon - 1$. Since (\ref{thm.exist.p1}) together with $\epsilon \leq 1$ provides $(t+\ell)/\epsilon \geq \ell \geq 512$, we obtain:
\begin{align} n \geq \frac{8(t+\ell)}{\epsilon} \label{thm.exist.p3} \end{align}
In particular $n \geq 8\ell/\epsilon \geq 4096/\epsilon$. Moreover, unwrapping the assignment of $\gamma$ and citing Clause \ref{thm.exist.c1}\ref{thm.exist.c1.2} of Theorem \ref{thm.exist} we have $\gamma \ell/2 \geq 4\log_2\bigl(8(r+2)/\epsilon\bigr)$, which gives $\gamma \ell/2 \geq \log_2\bigl(4096(r+2)^4/\epsilon^4\bigr)$ and hence:
\begin{align} 2^{\gamma \ell/2} \geq \frac{4096(r+2)^4}{\epsilon^4} \label{thm.exist.p4} \end{align}
We now verify that the data assigned above satisfies the hypotheses of Theorem \ref{thm.main}.
\begin{itemize} \item The bound $n \geq 4096/\epsilon$ recorded above together with $\epsilon \leq 1$ gives $n \geq 4096$, and since $n = \lfloor \log_2(N) \rfloor$ we have $N \geq 2^n$. Thus $N \in (2^{100},\infty)_{\zee}$, as demanded by Clause \ref{thm.main.c1}\ref{thm.main.c1.1} of Theorem \ref{thm.main}.
\item Clause \ref{thm.exist.c1}\ref{thm.exist.c1.1} of Theorem \ref{thm.exist} provides $\sigma \in (0,1]_{\arr}$ and $\epsilon \in (0,1]_{\arr}$, so that the assignment of $\beta$ gives $\beta \in (0,1]_{\arr}$, while the bound $\gamma \leq 1/16$ established above gives $\gamma \in (0,1/4]_{\arr}$. Together with the inequality $\beta + 2\gamma < \sigma/4$ established above, these are the requirements of Clause \ref{thm.main.c1}\ref{thm.main.c1.2}.
\item Clause \ref{thm.exist.c1}\ref{thm.exist.c1.2} of Theorem \ref{thm.exist} provides $\ell \in \enn$, while unwrapping the assignment of $\gamma$ gives $2/\gamma = 32(1-\sigma+\epsilon)/\epsilon$, which (\ref{thm.exist.p1}) bounds above by $\ell$. Moreover (\ref{thm.exist.p3}) together with the assignment of $\beta$ gives $\beta n = \epsilon n/8 \geq t+\ell \geq \ell$. Thus $\ell \in [2/\gamma, \beta n]_{\zee}$, as demanded by Clause \ref{thm.main.c1}\ref{thm.main.c1.3}.
\item Clause \ref{thm.exist.c1}\ref{thm.exist.c1.6} of Theorem \ref{thm.exist} provides $\achn(\mu,N) \geq \sigma$, as demanded by Clause \ref{thm.main.c1}\ref{thm.main.c1.4}.
\item Combining (\ref{thm.exist.p3}) with the assignment of $\beta$ provides $\beta n \geq t+\ell$, so that every $s \in [\beta n - \ell, n - \ell]_{\zee}$ satisfies $s \geq t$. Clause \ref{thm.exist.c1}\ref{thm.exist.c1.4} of Theorem \ref{thm.exist} therefore provides $\mathtt{Dep}_{\mathscr{M}}(s) \geq \ell$ for all such $s$, as demanded by Clause \ref{thm.main.c1}\ref{thm.main.c1.5}.\end{itemize}
Instantiating Clause \ref{thm.main.c1} of Theorem \ref{thm.main} with the data $(N,\sigma,\beta,\gamma,\ell,\mu,\mathscr{M})$, we are returned by Clause \ref{thm.main.c2} of that theorem a subset $\mathscr{G}$ of $\mathscr{M} \cap [1,2^{n+1}]_{\zee}$ satisfying $\lvert \mathscr{G} \rvert \geq 2^\ell/5$. Since $2^\ell/5 > 0$ we have $\mathscr{G} \neq \emptyset$, and so we may fix $m \in \mathscr{G}$. Since $2^n \leq N$ we have $2^{n+1} \leq 2N$, so that $m \in \mathscr{M} \cap [1,2N]_{\zee}$.\\

We assign the currently-local $m$ as the multiplier demanded by Clause \ref{thm.exist.c2} of Theorem \ref{thm.exist}. Moreover $\mathscr{F} = [0,r]_{\zee}$ is a nonempty finite subset of $\enn$ satisfying $\max(\mathscr{F}) = r$, and so the conclusion (\ref{thm.main.conc1}) of Theorem \ref{thm.main} provides:
\begin{align} \einf\bigl(\mu,\{m\}[0,r]_{\zee}\bigr) \geq \frac{1}{1-2\gamma}\left(\sigma - \beta - 2\gamma - \frac{2}{n}\right) - 160\left( \frac{r^{7/2}}{2^{\gamma \ell/2}} + \frac{r^{4}}{2^{\ell}} \right) \label{thm.exist.p5} \end{align}
We bound the two parts of the right side of (\ref{thm.exist.p5}) separately, beginning with the entropy part. We compute:
\begin{align} \frac{1}{1-2\gamma}\left(\sigma - \beta - 2\gamma - \frac{2}{n}\right) & = \sigma - \frac{1}{1-2\gamma}\left(\beta + 2\gamma(1-\sigma) + \frac{2}{n}\right) \label{thm.exist.p6.1}
\\ & \geq \sigma - \frac{8}{7}\left(\beta + 2\gamma(1-\sigma) + \frac{2}{n}\right) \label{thm.exist.p6.2}
\\ & \geq \sigma - \frac{8}{7}\left(\frac{\epsilon}{8} + \frac{\epsilon}{8} + \frac{2}{n}\right) \label{thm.exist.p6.3}
\\ & \geq \sigma - \frac{8}{7} \cdot \frac{3\epsilon}{8} \label{thm.exist.p6.4}
\\ & = \sigma - \frac{3\epsilon}{7} \label{thm.exist.p6.5} \end{align}
This computation may be justified as follows.
\begin{itemize}
\item The equality in (\ref{thm.exist.p6.1}) follows since $\sigma(1-2\gamma) - \beta - 2\gamma(1-\sigma) - 2/n = \sigma - \beta - 2\gamma - 2/n$.
\item (\ref{thm.exist.p6.2}) follows from (\ref{thm.exist.p6.1}) since the bound $\gamma \leq 1/16$ established above provides $1-2\gamma \geq 7/8$, and since the bracketed quantity is nonnegative.
\item (\ref{thm.exist.p6.3}) follows from (\ref{thm.exist.p6.2}) by unwrapping the assignments of $\beta$ and $\gamma$, the latter providing $2\gamma(1-\sigma) = \epsilon(1-\sigma)/\bigl(8(1-\sigma+\epsilon)\bigr) \leq \epsilon/8$.
\item (\ref{thm.exist.p6.4}) follows from (\ref{thm.exist.p6.3}) since the bound $n \geq 4096/\epsilon$ recorded above provides $2/n \leq \epsilon/8$.
\end{itemize}
We now bound the Fourier part of the right side of (\ref{thm.exist.p5}). We compute:
\begin{align} 160\left( \frac{r^{7/2}}{2^{\gamma \ell/2}} + \frac{r^{4}}{2^{\ell}} \right) & \leq \frac{320\,(r+2)^4}{2^{\gamma \ell /2}} \label{thm.exist.p7.1}
\\ & \leq \frac{320\,\epsilon^4}{4096} \label{thm.exist.p7.2}
\\ & \leq \frac{2\epsilon}{7} \label{thm.exist.p7.3} \end{align}
This computation may be justified as follows.
\begin{itemize}
\item (\ref{thm.exist.p7.1}) follows since $r^{7/2} \leq (r+2)^4$ and $r^4 \leq (r+2)^4$, and since the bound $\gamma \leq 1/16$ established above provides $\gamma \ell /2 \leq \ell$ and hence $2^{\ell} \geq 2^{\gamma \ell/2}$.
\item (\ref{thm.exist.p7.2}) follows from (\ref{thm.exist.p7.1}) by (\ref{thm.exist.p4}).
\item (\ref{thm.exist.p7.3}) follows from (\ref{thm.exist.p7.2}) since $\epsilon \leq 1$ provides $\epsilon^4 \leq \epsilon$, and since $320/4096 \leq 2/7$.
\end{itemize}
Combining (\ref{thm.exist.p5}) with the computations ending in (\ref{thm.exist.p6.5}) and (\ref{thm.exist.p7.3}) we obtain $\einf\bigl(\mu,\{m\}[0,r]_{\zee}\bigr) \geq \sigma - 3\epsilon/7 - 2\epsilon/7$, and since $5/7 \leq 1$ this verifies the desired inequality (\ref{thm.exist.conc1}) and thereby completes the proof of Theorem \ref{thm.exist}.
\end{Proof}

\section{Proof of Lemma \ref{lem.depbound}} \label{sec.pfapp.exhibit}

\subsection{Proof of Lemma \ref{lem.depbound}: preliminaries}

\subsubsection{Proof of Lemma \ref{lem.depbound}: gap criterion}

The following small lemma provides a sufficient condition for progression-shadowing in the form our applications to canonical sets can verify.

\begin{Lemma}[Gap criterion for progression-shadowing] \label{lem.gapshad}
Let $s,\ell \in \enn$ and let $\mathscr{M}$ be a subset of $\enn$. Suppose that every closed interval of length $2^{s-\ell}$ whose left endpoint lies in $[2^s, 2^{s+\ell}]_\arr$ contains an element of $\mathscr{M}$. Then $\mathscr{M}$ is $(s,\ell)$ progression-shadowing, and consequently $\mathtt{Dep}_{\mathscr{M}}(s) \geq \ell$. \end{Lemma} 

\begin{Proof}[of Lemma \ref{lem.gapshad}] Let $s,\ell$ and $\mathscr{M}$ be as in the statement of Lemma \ref{lem.gapshad} and let $m \in [1,2^\ell]_{\zee}$. Then $m 2^s$ lies in $[2^s,2^{s+\ell}]_\arr$, and so the hypothesis applied to the interval $[m2^s, m2^s + 2^{s-\ell}]_\arr$ provides an element $q$ of $\mathscr{M}$ lying in that interval. In particular $\lvert q - m2^s \rvert \leq 2^{s-\ell}$, which is the condition (\ref{def.progshad.1}) demanded by Definition \ref{def.progshad}.\\

Since $m$ was arbitrary, $\mathscr{M}$ is $(s,\ell)$ progression-shadowing. If $\ell \in \ennplus$ then Clause \ref{rem.shadep.c1} of Remark \ref{rem.shadep} provides $\mathtt{Dep}_{\mathscr{M}}(s) \geq \ell$, while if $\ell = 0$ then the same conclusion holds directly, since Definition \ref{def.shadep} defines $\mathtt{Dep}_{\mathscr{M}}(s)$ as an element of $\enn \cup \{\infty\}$. \end{Proof}

\subsubsection{Proof of Lemma \ref{lem.depbound}: number theory imports}

We now record certain external results from number theory which will be needed in the proof of Lemma \ref{lem.depbound}.\\

The following classical result follows, for example, from Theorem 36 of \cite{MR2445243}.

\begin{oLemma}[Dirichlet approximation] \label{olem.dirichlet}
For all $\theta \in \arr$ and all $N \in \ennplus$ there exist $q \in [1,N]_\zee$ and $p \in \zee$ such that $\lvert q\theta - p \rvert \leq 1/N$. \end{oLemma}

The following result is a standard specialization of the celebrated theorem of A.~Baker and G.~W\"ustholz on linear forms in logarithms, which appears as the Theorem on p.~20 of \cite{MR1234835}. The two-logarithm form from which ours follows is recorded as Theorem 4.3 of \cite{MR2563089}.

\begin{oLemma}[Effective irrationality measure for ratios of logarithms] \label{olem.bw}
Let $a,b \in [2,\infty)_\zee$ be multiplicatively independent. Then there exists an effectively computable $\kappa_{a,b} \in [2,\infty)_\arr$, depending only on $a$ and $b$, such that the following holds for all $p \in \zee$ and all $q \in [2,\infty)_\zee$.
\begin{align} \left\lvert \frac{\log(a)}{\log(b)} - \frac{p}{q} \right\rvert \geq q^{-\kappa_{a,b}} \label{olem.bw.conc1} \end{align} \end{oLemma}
 
\subsubsection{Proof of Lemma \ref{lem.depbound}: orbit filling}

The proof of Lemma \ref{lem.rotfill} is a writing-out, in the quantitative case we need, of the two most fundamental ideas in Diophantine approximation. These are that the pigeonhole principle produces a good rational approximation, and that repeated small rotations fill the circle.
\begin{itemize}
\item Dirichlet approximation supplies an increment, which is the $q \in [1,N]_{\zee}$ with $\eta = \lvert q\theta - p_0 \rvert \leq 1/N$, assigned at (\ref{lem.rotfill.1}).
\item The irrationality-measure hypothesis bounds the increment from below, $\eta \geq N^{1-\kappa}$, so it is small but not degenerate. 
\item Iterating the increment $u \leq \lceil 1/\eta \rceil \leq N^{\kappa-1}+1$ times, in the direction dictated by the sign of $\sigma$, lands within $\eta \leq 1/N$ of the target $y$. The multiplier $m = uq$ then costs at most $ 2N^{\kappa}$.
\end{itemize}

\begin{Lemma}[Orbit filling from an irrationality measure] \label{lem.rotfill}
Let $\theta \in \arr$ and $\kappa \in [2,\infty)_\arr$ be such that $\lvert \theta - p/q \rvert \geq q^{-\kappa}$ for all $p \in \zee$ and all $q \in [2,\infty)_\zee$. Then for all $N \in [4,\infty)_\zee$ and all $y \in \arr$ there exist $m \in [0,2N^{\kappa}]_\zee$ and $p \in \zee$ such that $\lvert m\theta - p - y \rvert \leq 1/N$. \end{Lemma}

\begin{Proof}[of Lemma \ref{lem.rotfill}] Let $\theta, \kappa, N$ and $y$ be as in the statement of Lemma \ref{lem.rotfill}. Since a pair $(m,p)$ certifying the conclusion for $y - \lfloor y \rfloor$ yields the pair $(m, p - \lfloor y \rfloor)$ for $y$, we may assume without loss of generality that $y \in [0,1)_\arr$. \\

We now instantiate Lemma \ref{olem.dirichlet} with $\theta$ and $N$ as currently-local, so we are returned $q \in [1,N]_\zee$ and $p_0 \in \zee$ satisfying $\lvert q\theta - p_0 \rvert \leq 1/N$. We define $\sigma \in [-1/N,1/N]_\arr$ and $\eta \in [0,1/N]_\arr$ as follows.
\begin{align} \label{lem.rotfill.1} \sigma &\coloneqq q\theta - p_0 & \eta & \coloneqq \lvert \sigma \rvert \end{align}
We now verify the lower bound $\eta \geq N^{1-\kappa}$ by considering two cases.
\begin{itemize}
\item \textbf{(Main case)} Suppose $q \in [2,N]_\zee$. We obtain:  
\begin{align*} \eta = q\left\lvert \theta - \frac{p_0}{q} \right\rvert \geq q^{1-\kappa} \geq N^{1-\kappa} \end{align*}
Here, the left equality follows directly from the definition of $\eta$ in (\ref{lem.rotfill.1}), the middle inequality follows from the hypothesis of the current Lemma \ref{lem.rotfill}, and the right inequality holds since we have assumed that $q \leq N$ and $1-\kappa \leq 0$.
 
\item \textbf{(Edge case)} Suppose $q = 1$. We obtain:
\begin{align*} \eta = \left\lvert \theta - \frac{2p_0}{2} \right\rvert \geq 2^{-\kappa} \geq N^{1-\kappa} \end{align*}
Here, the left equality again follows directly from the definition of $\eta$ in (\ref{lem.rotfill.1}), the middle inequality follows by applying the hypothesis of the current Lemma \ref{lem.rotfill} to the pair $(2p_0,2)$, and the right inequality holds since we have assumed that $N \geq 4$ and $\kappa \geq 2$.\end{itemize}

Since Lemma \ref{olem.dirichlet} returns $q \in [1,N]_\zee$ these cases are exhaustive, and so we indeed find $\eta \geq N^{1-\kappa}$ and in particular $\sigma \neq 0$. We now divide cases again by the sign of $\sigma$.
\begin{itemize}
\item \textbf{(Positive case)} Suppose $\sigma > 0$. We define $u \coloneqq \lceil y/\eta \rceil$, and assign the parameters demanded by the current conclusion as $m \coloneqq uq$ and $p \coloneqq up_0$. Then $u\sigma = u\eta \in [y,y+\eta]_\arr$, while $m\theta - p = u(q\theta - p_0) = u\sigma$, and so $\lvert m\theta - p - y \rvert \leq \eta$. Since $y \leq 1$ we also have $u \leq \lceil 1/\eta \rceil$.
\item \textbf{(Negative case)} Suppose $\sigma < 0$. We define $u \coloneqq \lceil (1-y)/\eta \rceil$, and we assign the parameters demanded by the current conclusion as $m \coloneqq uq$ and $p \coloneqq up_0 - 1$. Then $u\sigma = -u\eta \in [-(1-y)-\eta, -(1-y)]_\arr$, so that $u\sigma + 1 \in [y-\eta,y]_\arr$, while $m\theta - p = u(q\theta - p_0) + 1 = u\sigma + 1$, and so $\lvert m\theta - p - y \rvert \leq \eta$. Since $1-y \leq 1$ we also have $u \leq \lceil 1/\eta \rceil$.
\end{itemize}
In both cases $\lvert m\theta - p - y \rvert \leq \eta \leq 1/N$, and we compute:
\begin{align*} m = uq \leq \left \lceil \frac{1}{\eta} \right \rceil N \leq \left( N^{\kappa-1}+1 \right)N \leq 2N^{\kappa} \end{align*}
Here, the left inequality holds since $q \leq N$, the middle inequality holds since $\eta \geq N^{1-\kappa}$ and the right inequality holds since $N \leq N^{\kappa}$. As $m$ is moreover a nonnegative integer, this completes the proof of Lemma \ref{lem.rotfill}. \end{Proof}

\subsection{Proof of Lemma \ref{lem.depbound}: main}

Everything in the proof below is a direct and elementary computation. Each clause is verified through Lemma \ref{lem.gapshad}; the first two consume nothing further, and the third consumes Lemma \ref{olem.bw} together with Lemma \ref{lem.rotfill}.

\begin{Proof}[of Lemma \ref{lem.depbound}] Each clause is verified by exhibiting a value of $\ell$ for which the hypothesis of Lemma \ref{lem.gapshad} holds for the multiplier set named in that clause.\\

\underline{Verification of Clause \ref{lem.depbound.c1}}:\\

Let $n_0$ and $\mathscr{M}$ be as in Clause \ref{lem.depbound.c1} of Lemma \ref{lem.depbound}, let $s \in [3\log_2(n_0) \vee 12,\infty)_\zee$ and define $\ell \in \enn$ as follows.
\begin{align*} \ell \coloneqq \left\lfloor \frac{s-9}{5} \right\rfloor \end{align*}
Now let $x \in [2^s, 2^{s+\ell}]_\arr$ and define $n \coloneqq \lceil x^{1/3} \rceil$. We have $n \geq 2^{s/3} \geq n_0$, where the right inequality holds since $s \geq 3\log_2(n_0)$. Thus the hypothesis of Clause \ref{lem.depbound.c1} provides $m \in \mathscr{M}$ with $n^3 < m < (n+1)^3$.\\

We compute:
\begin{align} m & \leq (n+1)^3 \label{lem.depbound.p1.1}
\\ & \leq x + 6x^{2/3} + 12x^{1/3} + 8 \label{lem.depbound.p1.2}
\\ & \leq x + 7x^{2/3} \label{lem.depbound.p1.3}
\\ & \leq x + 7 \cdot 2^{2(s+\ell)/3} \label{lem.depbound.p1.4}
\\ & \leq x + 2^{s-\ell} \label{lem.depbound.p1.5} \end{align}
This computation may be justified as follows.
\begin{itemize}
\item The inequality in (\ref{lem.depbound.p1.1}) holds by the assignment of $m$.
\item (\ref{lem.depbound.p1.2}) follows from (\ref{lem.depbound.p1.1}) by expanding $(n+1)^3 \leq (x^{1/3}+2)^3$, where we note that $n \leq x^{1/3}+1$.
\item (\ref{lem.depbound.p1.3}) follows from (\ref{lem.depbound.p1.2}) since $12x^{1/3} + 8 \leq x^{2/3}$, which holds because $s \geq 12$ provides $x^{1/3} \geq 2^{s/3} \geq 16$.
\item (\ref{lem.depbound.p1.4}) follows from (\ref{lem.depbound.p1.3}) since $x \leq 2^{s+\ell}$.
\item (\ref{lem.depbound.p1.5}) follows from (\ref{lem.depbound.p1.4}) since $5\ell \leq s - 9$ and $3\log_2(7) \leq 9$, so that $\log_2(7) + 2(s+\ell)/3 \leq s-\ell$.
\end{itemize}
Since also $m > n^3 \geq x$, the closed interval of length $2^{s-\ell}$ with left endpoint $x$ contains $m$. Since $x$ was arbitrary in $[2^s, 2^{s+\ell}]_\arr$, Lemma \ref{lem.gapshad} provides the following.
\begin{align*} \mathtt{Dep}_{\mathscr{M}}(s) \geq \ell \geq \frac{s-9}{5} - 1 \geq \frac{s}{5} - 3 \end{align*}
Here, the middle inequality holds since the floor function decreases its argument by less than $1$, while the right inequality holds since $9/5 + 1 \leq 3$. This completes the verification of Clause \ref{lem.depbound.c1}.\\

\underline{Verification of Clause \ref{lem.depbound.c2}}:\\

Let $v \in \ennplus$ and let $s \in \enn$. Since $\mathtt{Dep}_{\mathbb{W}_v}(s)$ is nonnegative by Definition \ref{def.shadep}, we may assume that $s \geq v^2$, as otherwise the right side of (\ref{lem.depbound.c2.1}) is negative and the conclusion is immediate. Define $\ell \in \enn$ as follows.
\begin{align*} \ell \coloneqq \left\lfloor \frac{s - v^2}{2v-1} \right\rfloor \end{align*}
Now let $x \in [2^s, 2^{s+\ell}]_\arr$ and define $n \coloneqq \lceil x^{1/v} \rceil$, so that $n^v$ is an element of $\mathbb{W}_v \cap [x,\infty)_\arr$.\\

We compute:
\begin{align} n^v & \leq \bigl(x^{1/v} + 1\bigr)^v \label{lem.depbound.p2.1}
\\ & \leq x + 2^v x^{(v-1)/v} \label{lem.depbound.p2.2}
\\ & \leq x + 2^{v + (s+\ell)(v-1)/v} \label{lem.depbound.p2.3}
\\ & \leq x + 2^{s-\ell} \label{lem.depbound.p2.4} \end{align}
This computation may be justified as follows.
\begin{itemize}
\item The inequality in (\ref{lem.depbound.p2.1}) follows since $n \leq x^{1/v} + 1$.
\item (\ref{lem.depbound.p2.2}) follows from (\ref{lem.depbound.p2.1}) by expanding the binomial and bounding each of the $2^v - 1$ terms other than $x$ by $x^{(v-1)/v}$, where we note that $x \geq 1$.
\item (\ref{lem.depbound.p2.3}) follows from (\ref{lem.depbound.p2.2}) since $x \leq 2^{s+\ell}$.
\item (\ref{lem.depbound.p2.4}) follows from (\ref{lem.depbound.p2.3}) since the assignment of $\ell$ provides $\ell(2v-1) \leq s - v^2$, which rearranges to $v + (s+\ell)(v-1)/v \leq s - \ell$.
\end{itemize}
Since also $n^v \geq x$, the closed interval of length $2^{s-\ell}$ with left endpoint $x$ contains $n^v$.  Since $x$ was arbitrary in $[2^s, 2^{s+\ell}]_\arr$, Lemma \ref{lem.gapshad} provides the following.
\begin{align*} \mathtt{Dep}_{\mathbb{W}_v}(s) \geq \ell \geq \frac{s-v^2}{2v-1} - 1 \geq \frac{s}{2v-1} - v - 1 \end{align*}
Here, the middle inequality holds since the floor function decreases its argument by less than $1$, while the right inequality holds since $v^2/(2v-1) \leq v$. This completes the verification of Clause \ref{lem.depbound.c2}.\\

\underline{Verification of Clause \ref{lem.depbound.c5}}:\\

Let $a,b \in [2,\infty)_\zee$ be multiplicatively independent, let $\kappa \coloneqq \kappa_{a,b}$ be as provided by Lemma \ref{olem.bw} and define $\theta \coloneqq \log(a)/\log(b)$, so that the conclusion (\ref{olem.bw.conc1}) of Lemma \ref{olem.bw} verifies the hypothesis of Lemma \ref{lem.rotfill} for our current $\theta$ and $\kappa$.\\

We assign the constant $\mathtt{c}_{a,b}$ demanded by Clause \ref{lem.depbound.c5} as $\mathtt{c}_{a,b} \coloneqq 1/(8\kappa)$, and we choose the threshold $\mathtt{s}_{a,b} \in \enn$, depending only on $a$ and $b$, sufficiently large that the following two conditions hold for all $s \in [\mathtt{s}_{a,b},\infty)_\zee$.
\begin{align} \log(s) \geq 8\kappa && \sqrt{s} \geq \frac{2\log(a)}{\log(2)}\bigl(4\lceil\log(b)\rceil\bigr)^{\kappa} \label{lem.depbound.p5.1} \end{align}
Now let $s \in [\mathtt{s}_{a,b},\infty)_\zee$ and define $\ell, N \in \enn$ as follows.
\begin{align} \label{lem.depbound.p5.2} \ell \coloneqq \left\lfloor \frac{\log_2(s)}{4\kappa} \right\rfloor && N \coloneqq 2^{2\ell+2}\lceil\log(b)\rceil \end{align}
We compute:
\begin{align} \ell & \geq \frac{\log(s)}{4\kappa} - 1 \label{lem.depbound.p5.3}
\\ & \geq \frac{\log(s)}{8\kappa} \label{lem.depbound.p5.4}
\\ & = \mathtt{c}_{a,b}\log(s) \label{lem.depbound.p5.4.1} \end{align}
This computation may be justified as follows.
\begin{itemize}
\item The inequality in (\ref{lem.depbound.p5.3}) follows from the assignment of $\ell$ at the left of (\ref{lem.depbound.p5.2}), since we have $\log(y) \leq \log_2(y)$ for all $y \in [1,\infty)_\arr$.
\item (\ref{lem.depbound.p5.4}) follows from (\ref{lem.depbound.p5.3}) by the first condition in (\ref{lem.depbound.p5.1}).
\item (\ref{lem.depbound.p5.4.1}) follows from (\ref{lem.depbound.p5.4}) directly from our assignment of $\mathtt{c}_{a,b}$.
\end{itemize}
Thus it suffices to show the hypothesis of Lemma \ref{lem.gapshad} holds for $\mathbb{S}_{a,b}$ with this value of $\ell$. We also compute:
\begin{align} N^{\kappa} & = \bigl(4\lceil\log(b)\rceil\bigr)^{\kappa}\, 2^{2\kappa\ell} \label{lem.depbound.p5.5}
\\ & \leq \bigl(4\lceil\log(b)\rceil\bigr)^{\kappa} \sqrt{s} \label{lem.depbound.p5.6}
\\ & \leq \frac{s\log(2)}{2\log(a)} \label{lem.depbound.p5.7} \end{align}
This computation may be justified as follows.
\begin{itemize}
\item The equality in (\ref{lem.depbound.p5.5}) follows directly from the assignment of $N$ at the right of (\ref{lem.depbound.p5.2}).
\item (\ref{lem.depbound.p5.6}) follows from (\ref{lem.depbound.p5.5}) since $2\kappa\ell \leq \log_2(s)/2$ by the assignment of $\ell$ at the left of (\ref{lem.depbound.p5.2}).
\item (\ref{lem.depbound.p5.7}) follows from (\ref{lem.depbound.p5.6}) by the condition at the right of (\ref{lem.depbound.p5.1}).
\end{itemize}
Now let $x \in [2^s, 2^{s+\ell}]_\arr$. Noting that $N = 2^{2\ell+2}\lceil\log(b)\rceil \geq 4$ since $\ell \geq 0$ and $\lceil\log(b)\rceil \geq 1$, we instantiate Lemma \ref{lem.rotfill} with the data $(\theta,\kappa,N)$ and the target $y \coloneqq \log(x)/\log(b) + 1/N$.\\

The conclusion of Lemma \ref{lem.rotfill} returns $m \in [0,2N^{\kappa}]_\zee$ and $p \in \zee$ satisfying $\lvert m\theta - p - y \rvert \leq 1/N$. We compute:
\begin{align} p & \leq m\theta - y + \frac{1}{N} \label{lem.depbound.p5.8}
\\ & \leq 2N^{\kappa}\theta - \frac{\log(x)}{\log(b)} \label{lem.depbound.p5.9}
\\ & \leq \frac{s\log(2)}{\log(b)} - \frac{s\log(2)}{\log(b)} = 0 \label{lem.depbound.p5.10} \end{align}
This computation may be justified as follows.
\begin{itemize}
\item The inequality in (\ref{lem.depbound.p5.8}) follows from the conclusion of Lemma \ref{lem.rotfill} as instantiated above.
\item (\ref{lem.depbound.p5.9}) follows from (\ref{lem.depbound.p5.8}) since we have by construction that $m \leq 2N^{\kappa}$ and $\theta > 0$ and:
\begin{align}\label{lem.depbound.p5.10.1} y - \frac{1}{N} = \frac{\log(x)}{\log(b)} \end{align}
\item (\ref{lem.depbound.p5.10}) follows from (\ref{lem.depbound.p5.9}) by the computation ending in (\ref{lem.depbound.p5.7}), where we note that $\theta = \log(a)/\log(b)$ and that $\log(x) \geq s\log(2)$ since $x \geq 2^s$.
\end{itemize}
Thus $k \coloneqq -p$ is an element of $\enn$. Since $\lvert m\theta - p - y \rvert \leq 1/N$ we have $m\theta - p - y + 1/N \in [0,2/N]_\arr$. Thus multiplying by $\log(b)$ and again using (\ref{lem.depbound.p5.10.1}) we obtain:
\begin{align} 0 \leq m\log(a) + k\log(b) - \log(x) \leq  \frac{2\log(b)}{N} \label{lem.depbound.p5.11} \end{align}
We define $w \coloneqq a^m b^k \in \mathbb{S}_{a,b}$ and compute:
\begin{align} w & \leq x\, e^{2^{-2\ell-1}} \label{lem.depbound.p5.12}
\\ & \leq x\bigl(1 + 2^{-2\ell}\bigr) \label{lem.depbound.p5.13}
\\ & \leq x + 2^{s-\ell} \label{lem.depbound.p5.14} \end{align}
This computation may be justified as follows.
\begin{itemize}
\item The inequality in (\ref{lem.depbound.p5.12}) follows from the right inequality in (\ref{lem.depbound.p5.11}) since our assignment of $N$ at the right of (\ref{lem.depbound.p5.2}) implies $2\log(b)/N \leq 2^{-2\ell-1}$. 
\item (\ref{lem.depbound.p5.13}) follows from (\ref{lem.depbound.p5.12}) since $e^{t} \leq 1+2t$ for all $t \in [0,1]_\arr$ by convexity of the exponential function.
\item (\ref{lem.depbound.p5.14}) follows from (\ref{lem.depbound.p5.13}) since we have assumed $x \leq 2^{s+\ell}$.
\end{itemize}
Since the left inequality in (\ref{lem.depbound.p5.11}) provides $w \geq x$, the closed interval of length $2^{s-\ell}$ with left endpoint $x$ contains the element $w$ of $\mathbb{S}_{a,b}$. Since $x \in [2^s, 2^{s+\ell}]_\arr$ was arbitrary, Lemma \ref{lem.gapshad} provides $\mathtt{Dep}_{\mathbb{S}_{a,b}}(s) \geq \ell$. Combining this with the computation ending in (\ref{lem.depbound.p5.4.1}) completes the verification of Clause \ref{lem.depbound.c5} and thereby the proof of Lemma \ref{lem.depbound}. \end{Proof}

\subsection{Proof of Lemma \ref{lem.primedep}}

\begin{Proof}[of Lemma \ref{lem.primedep}] Assign $n_0 \coloneqq \bigl\lceil \exp(\exp(33.3)) \bigr\rceil$. For every $n \in [n_0,\infty)_{\zee}$ Lemma \ref{olem.primegap} provides $p \in \mathbb{P}$ with $n^3 < p < (n+1)^3$, and since $p \in \zee$ this gives $\mathbb{P} \cap [n^3+1,(n+1)^3-1]_{\zee} \neq \emptyset$. Thus the data $(n_0,\mathbb{P})$ satisfies the hypotheses of Clause \ref{lem.depbound.c1} of Lemma \ref{lem.depbound}.\\

Since $n_0 \leq 2\exp(\exp(33.3))$ we have $3\log_2(n_0) \leq 3 + 3\exp(33.3)/\log(2)$, and combining $\exp(33.3) \leq 3 \times 10^{14}$ with $\log(2) \geq 0.69$ provides $3\log_2(n_0) \leq 1.31 \times 10^{15}$, so that $3\log_2(n_0) \vee 12 \leq 2^{51}$. Instantiating Clause \ref{lem.depbound.c1} of Lemma \ref{lem.depbound} with the data $(n_0,\mathbb{P})$ therefore yields (\ref{lem.primedep.1}) for all $s \in [2^{51},\infty)_{\zee}$, and thereby completes the proof of Lemma \ref{lem.primedep}. \end{Proof}

\section{Proof of Lemma \ref{lem.ratio}} \label{sec.pfapp.ratio}

An outline of the proof of Lemma \ref{lem.ratio} is as follows.
\begin{itemize}
\item Clauses \ref{lem.ratio.c2} and \ref{lem.ratio.c3} are one assertion read through Clause \ref{rem.shadep.c1} of Remark \ref{rem.shadep}, and the proof spends no further effort on them.
\item The passage from Clause \ref{lem.ratio.c1} to Clause \ref{lem.ratio.c3} is a change of scale. Sublacunarity bounds the forward gap of $\mathscr{M}$ above every large $x$ by $2^{-2\ell}x$, and on the window $[2^s,2^{s+\ell}]_\arr$ this is at most $2^{s-\ell}$. Lemma \ref{lem.gapshad} converts the gap bound into shadowing at every large scale.
\item The passage from Clause \ref{lem.ratio.c3} to Clause \ref{lem.ratio.c1} shadows, for each large element $m_j$, the second multiple of $2^s$ above $m_j$ at the scale assigned in (\ref{lem.ratio.p3}). The shadow lies strictly above $m_j$ by (\ref{lem.ratio.p6.3}), hence at or above $m_{j+1}$, and within $2^{4-\ell}m_j$ of $m_j$ by (\ref{lem.ratio.p7.5}). The depth $\ell$ is then sent to infinity.
\end{itemize}

The dependencies in the proof of Lemma \ref{lem.ratio} are as follows.
\begin{deps}
\dep{Definitions}{\depitem{Progression-shadowing sets}{def.progshad}\depsep\depitem{Progression-shadowing depth}{def.shadep}\depsep\depitem{Sublacunary sets}{def.sublac}}
\dep{Remarks}{\depitemc{Shadowing depth and progression-shadowing}{rem.shadep}{rem.shadep.c1}}
\dep{Lemmas}{\depitem{Gap criterion for progression-shadowing}{lem.gapshad}}
\end{deps}

We now move on to the formal proof.

\begin{Proof}[of Lemma \ref{lem.ratio}] Let $\mathscr{M}$ be as in the statement of Lemma \ref{lem.ratio}, and write $(m_j)_{j \in \enn}$ for the increasing enumeration of $\mathscr{M} \setminus \{0\}$, which is an infinite set since $\mathscr{M}$ is. Thus $(m_j)_{j \in \enn}$ is a strictly increasing sequence of positive integers, and in particular it is unbounded.\\

\underline{Equivalence of Clauses \ref{lem.ratio.c2} and \ref{lem.ratio.c3}}:\\

Since Definition \ref{def.shadep} exhibits $\mathtt{Dep}_{\mathscr{M}}(s)$ as an element of $\enn \cup \{\infty\}$ for every $s \in \enn$, Clause \ref{lem.ratio.c2} asserts precisely that for every $\ell \in \ennplus$ there exists $s_0 \in \enn$ such that $\mathtt{Dep}_{\mathscr{M}}(s) \geq \ell$ for all $s \in [s_0,\infty)_\zee$. For all $\ell \in \ennplus$ and all $s \in \enn$, Clause \ref{rem.shadep.c1} of Remark \ref{rem.shadep} asserts that $\mathtt{Dep}_{\mathscr{M}}(s) \geq \ell$ if and only if $\mathscr{M}$ is $(s,\ell)$ progression-shadowing. Substituting this equivalence into the previous sentence exhibits Clause \ref{lem.ratio.c2} as a restatement of Clause \ref{lem.ratio.c3}.\\

\underline{Clause \ref{lem.ratio.c1} implies Clause \ref{lem.ratio.c3}}:\\

Assume that $\mathscr{M}$ is sublacunary, and let $\ell \in \ennplus$. Applying the hypothesized limit in (\ref{def.sublac.1}) we find there exists $J \in \enn$ such that the following holds for all $j \in [J,\infty)_\zee$.
\begin{align} m_{j+1} \leq \bigl(1 + 2^{-2\ell}\bigr) m_j \label{lem.ratio.p1} \end{align}
We define $s_0 \coloneqq \lfloor \log_2(m_J) \rfloor + 1$, so that $2^{s_0} > m_J$. Now let $s \in [s_0,\infty)_\zee$ and let $x \in [2^s,2^{s+\ell}]_\arr$. \\

Since $x \geq 2^{s_0} > m_J$ and the sequence $(m_j)_{j \in \enn}$ is unbounded, the set $\mathscr{J} \coloneqq \{ k \in \enn : m_k < x \}$ is finite and contains $J$. We define $j \coloneqq \max(\mathscr{J})$. Then $j \in [J,\infty)_\zee$ since $J \in \mathscr{J}$, while $m_j < x$ since $j \in \mathscr{J}$ and $x \leq m_{j+1}$ since $j + 1 \notin \mathscr{J}$. We compute:
\begin{align} m_{j+1} & \leq \bigl(1 + 2^{-2\ell}\bigr) m_j \label{lem.ratio.p2.1}
\\ & \leq x + 2^{-2\ell} x \label{lem.ratio.p2.2}
\\ & \leq x + 2^{s-\ell} \label{lem.ratio.p2.3} \end{align}
This computation may be justified as follows.
\begin{itemize}
\item The inequality in (\ref{lem.ratio.p2.1}) follows from (\ref{lem.ratio.p1}), which applies since $j \geq J$.
\item (\ref{lem.ratio.p2.2}) follows from (\ref{lem.ratio.p2.1}) since $m_j < x$.
\item (\ref{lem.ratio.p2.3}) follows from (\ref{lem.ratio.p2.2}) since $x \leq 2^{s+\ell}$ provides $2^{-2\ell}x \leq 2^{s-\ell}$.
\end{itemize}
Since also $m_{j+1} \geq x$, the closed interval of length $2^{s-\ell}$ with left endpoint $x$ contains the element $m_{j+1}$ of $\mathscr{M}$. Since $x$ was arbitrary in $[2^s,2^{s+\ell}]_\arr$, Lemma \ref{lem.gapshad} provides that $\mathscr{M}$ is $(s,\ell)$ progression-shadowing. Since $s$ was arbitrary in $[s_0,\infty)_\zee$ and $\ell$ was arbitrary in $\ennplus$, this verifies Clause \ref{lem.ratio.c3}.\\

\underline{Clause \ref{lem.ratio.c3} implies Clause \ref{lem.ratio.c1}}:\\

Assume Clause \ref{lem.ratio.c3}, and let $\ell \in [2,\infty)_\zee$. Clause \ref{lem.ratio.c3} returns $s_0 \in \enn$ such that $\mathscr{M}$ is $(s,\ell)$ progression-shadowing for all $s \in [s_0,\infty)_\zee$. Since the sequence $(m_j)_{j \in \enn}$ is unbounded there exists $J \in \enn$ such that $m_J \geq 2^{s_0+\ell-2}$. Let $j \in [J,\infty)_\zee$, and define $s \in \zee$ and $k \in \zee$ as follows.
\begin{align} s \coloneqq \lfloor \log_2(m_j) \rfloor - \ell + 2 && k \coloneqq \lfloor m_j 2^{-s} \rfloor + 2 \label{lem.ratio.p3} \end{align}
Since $m_j \geq m_J \geq 2^{s_0+\ell-2}$ we have $\lfloor \log_2(m_j) \rfloor \geq s_0 + \ell - 2$, whence $s \geq s_0$ and in particular $s \in \enn$. Since $2^{\lfloor \log_2(m_j) \rfloor} \leq m_j < 2^{\lfloor \log_2(m_j) \rfloor + 1}$, the assignment of $s$ provides:
\begin{align} 2^{s+\ell-2} \leq m_j < 2^{s+\ell-1} \label{lem.ratio.p4} \end{align}
Dividing (\ref{lem.ratio.p4}) by $2^s$ we find $m_j 2^{-s} \in [2^{\ell-2},2^{\ell-1})_\arr$, so that the assignment of $k$ provides $k \in [2^{\ell-2}+2,2^{\ell-1}+1]_\zee$, and in particular $k \in [1,2^\ell]_\zee$ since $\ell \geq 1$. Since $\lfloor y \rfloor + 2 - y \in (1,2]_\arr$ for all $y \in \arr$, the assignment of $k$ also provides:
\begin{align} m_j + 2^s < k2^s \leq m_j + 2^{s+1} \label{lem.ratio.p5} \end{align}
Since $s \geq s_0$ the set $\mathscr{M}$ is $(s,\ell)$ progression-shadowing, and instantiating Definition \ref{def.progshad} returns $q \in \mathscr{M}$ satisfying the following.
\begin{align} \lvert q - k2^s \rvert \leq 2^{s-\ell} \label{lem.ratio.p5.1} \end{align}
We compute:
\begin{align} q & \geq k2^s - 2^{s-\ell} \label{lem.ratio.p6.1}
\\ & > m_j + 2^s - 2^{s-\ell} \label{lem.ratio.p6.2}
\\ & \geq m_j \label{lem.ratio.p6.3} \end{align}
This computation may be justified as follows.
\begin{itemize}
\item The inequality in (\ref{lem.ratio.p6.1}) follows from (\ref{lem.ratio.p5.1}).
\item (\ref{lem.ratio.p6.2}) follows from (\ref{lem.ratio.p6.1}) by the left inequality in (\ref{lem.ratio.p5}).
\item (\ref{lem.ratio.p6.3}) follows from (\ref{lem.ratio.p6.2}) since $\ell \geq 0$.
\end{itemize}
Thus $q$ is an element of $\mathscr{M}$ with $q > m_j \geq 1$, so that $q \in \mathscr{M} \setminus \{0\}$, and since $m_{j+1}$ is the least element of $\mathscr{M} \setminus \{0\}$ exceeding $m_j$ we have $m_{j+1} \leq q$. We compute again:
\begin{align} m_{j+1} & \leq q \label{lem.ratio.p7.1}
\\ & \leq k2^s + 2^{s-\ell} \label{lem.ratio.p7.2}
\\ & \leq m_j + 2^{s+1} + 2^{s-\ell} \label{lem.ratio.p7.3}
\\ & \leq m_j + 2^{s+2} \label{lem.ratio.p7.4}
\\ & \leq \bigl(1 + 2^{4-\ell}\bigr) m_j \label{lem.ratio.p7.5} \end{align}
This computation may be justified as follows.
\begin{itemize}
\item The inequality in (\ref{lem.ratio.p7.1}) follows from the computation ending in (\ref{lem.ratio.p6.3}) together with the minimality of $m_{j+1}$ recorded immediately thereafter.
\item (\ref{lem.ratio.p7.2}) follows from (\ref{lem.ratio.p7.1}) by (\ref{lem.ratio.p5.1}).
\item (\ref{lem.ratio.p7.3}) follows from (\ref{lem.ratio.p7.2}) by the right inequality in (\ref{lem.ratio.p5}).
\item (\ref{lem.ratio.p7.4}) follows from (\ref{lem.ratio.p7.3}) since $2^{s-\ell} \leq 2^{s+1}$.
\item (\ref{lem.ratio.p7.5}) follows from (\ref{lem.ratio.p7.4}) since the left inequality in (\ref{lem.ratio.p4}) provides $2^{s+2} = 2^{4-\ell} \cdot 2^{s+\ell-2} \leq 2^{4-\ell} m_j$.
\end{itemize}
Since the sequence $(m_j)_{j \in \enn}$ is increasing, dividing the computation ending in (\ref{lem.ratio.p7.5}) by $m_j$ provides:
\begin{align*} 1 \leq \frac{m_{j+1}}{m_j} \leq 1 + 2^{4-\ell} \end{align*}
Since $j$ was arbitrary in $[J,\infty)_\zee$ and $\ell$ was arbitrary in $[2,\infty)_\zee$, we have verified that for every $\ell \in [2,\infty)_\zee$ there exists $J \in \enn$ such that $1 \leq m_{j+1}/m_j \leq 1 + 2^{4-\ell}$ for all $j \in [J,\infty)_\zee$. Letting $\ell \to \infty$ verifies (\ref{def.sublac.1}), and so $\mathscr{M}$ is sublacunary. This completes the proof of Lemma \ref{lem.ratio}. \end{Proof}

\section{Proof of Theorem \ref{thm.blprime}} \label{sec.pfapp.special.ramp}

An outline of the proof of Theorem \ref{thm.blprime} is as follows.
\begin{itemize}
\item We make an initial common reduction, noting that any multiplier set whose depth clears $\ell$ at every scale beyond a threshold $t$ feeds Theorem \ref{thm.exist} at every resolution $N \geq 2^{9(t+\ell)/\epsilon}$ at which $\mu$ carries entropy $\sigma$. The bandprofile conclusion (\ref{thm.exist.conc1}) returned there converts into a statement about the twisted measure $m_\circledcirc \mu$ by endomorphy of the exponential infimum. 
\item Each clause is then a threshold check against the matching depth bound. The value of $t$ named in Clause \ref{thm.blprime.c1} is engineered so that $s/5 - 3 \geq \ell$ for all $s \geq t$, the value in Clause \ref{thm.blprime.c2} so that $s/(2v-1) - v - 1$ clears the same bar, and the value in Clause \ref{thm.blprime.c3} so that $\mathtt{c}_{a,b}\log(s)$ does.
\item Here is where the asymmetry arises in the size of the three thresholds: linear in $\ell$ for the primes and the $v$-th powers, but exponential in $\ell/\mathtt{c}_{a,b}$ for the multiplicative semigroups.
\end{itemize}

The dependencies in the proof of Theorem \ref{thm.blprime} beyond those in the statement are as follows.
\begin{deps}
\dep{Remarks}{\depitem{Endomorphy of exponential infimum}{rem.einfendo}}
\dep{Lemmas}{\depitem{Bounds on progression-shadowing depth for canonical sets}{lem.depbound} \depsep \depitem{Progression-shadowing depth of the primes}{lem.primedep}}
\dep{Theorems}{\depitem{Existential endomorphic rigidity via entropy}{thm.exist}}
\end{deps}

We now move on to the formal proof.

\begin{Proof}[of Theorem \ref{thm.blprime}] Let $\epsilon, r, \sigma$ and $\mu$ be as in the statement of Theorem \ref{thm.blprime} and recall from that statement the notation $\ell$ assigned in (\ref{thm.blprime.h1}), which lies in $\ennplus$ and satisfies the inequality demanded by Clause \ref{thm.exist.c1}\ref{thm.exist.c1.2} of Theorem \ref{thm.exist}.\\

We first record a reduction common to the three clauses. Suppose $\mathscr{M}$ is a subset of $\enn$ and $t \in \enn$ satisfies $\mathtt{Dep}_{\mathscr{M}}(s) \geq \ell$ for all $s \in [t,\infty)_\zee$, and let $N \in [2^{9(t+\ell)/\epsilon},\infty)_\zee$ satisfy $\achn(\mu,N) \geq \sigma$. Then the data $(\epsilon,r,\sigma,\ell,\mathscr{M},t,N,\mu)$ satisfies the hypotheses gathered in Clause \ref{thm.exist.c1} of Theorem \ref{thm.exist}, and instantiating that theorem we are returned $m \in \mathscr{M} \cap [1,2N]_\zee$ for which the conclusion (\ref{thm.exist.conc1}) holds. \\

Since $m$ lies in $\ennplus$, Remark \ref{rem.einfendo} provides $\einf(\mu,\{m\}[0,r]_\zee) = \einf(m_\circledcirc\mu,[0,r]_\zee)$, and combining this with (\ref{thm.exist.conc1}) we obtain:
\begin{align*} \einf\bigl(m_\circledcirc\mu,[0,r]_\zee\bigr) \geq \sigma - \epsilon \end{align*}
Thus it suffices to verify, for each clause of the current theorem, that the multiplier set and the threshold $t$ specified by that clause satisfy the depth hypothesis of the above reduction.\\

\vspace{5 pt}

\underline{Verification of Clause \ref{thm.blprime.c1}}:\\

Let $s \in [t,\infty)_\zee$, where $t \coloneqq 2^{51} \vee 5(\ell+3)$ as in Clause \ref{thm.blprime.c1}. Since $s \geq 2^{51}$, Lemma \ref{lem.primedep} provides $\mathtt{Dep}_{\mathbb{P}}(s) \geq s/5 - 3$, while $s \geq 5(\ell+3)$ provides $s/5 - 3 \geq \ell$.\\

\underline{Verification of Clause \ref{thm.blprime.c2}}:\\

Let $v \in \ennplus$ and let $s \in [t,\infty)_\zee$, where $t \coloneqq (2v-1)(\ell+v+1)$ as in Clause \ref{thm.blprime.c2}. Clause \ref{lem.depbound.c2} of Lemma \ref{lem.depbound} provides $\mathtt{Dep}_{\mathbb{W}_v}(s) \geq s/(2v-1) - v - 1$, while $s \geq (2v-1)(\ell+v+1)$ provides $s/(2v-1) - v - 1 \geq \ell$.\\

\underline{Verification of Clause \ref{thm.blprime.c3}}:\\

Let $a,b \in [2,\infty)_\zee$ be multiplicatively independent, let $\mathtt{c}_{a,b}$ and $\mathtt{s}_{a,b}$ be as provided by Clause \ref{lem.depbound.c5} of Lemma \ref{lem.depbound}, and let $s \in [t,\infty)_\zee$, where $t \coloneqq \mathtt{s}_{a,b} \vee \lceil e^{\ell/\mathtt{c}_{a,b}} \rceil$ as in Clause \ref{thm.blprime.c3}, which lies in $\enn$ by the ceiling. Since $s \geq \mathtt{s}_{a,b}$, Clause \ref{lem.depbound.c5} of Lemma \ref{lem.depbound} provides $\mathtt{Dep}_{\mathbb{S}_{a,b}}(s) \geq \mathtt{c}_{a,b}\log(s)$, while $s \geq e^{\ell/\mathtt{c}_{a,b}}$ provides $\mathtt{c}_{a,b}\log(s) \geq \ell$. This completes the verification of the three clauses and thereby the proof of Theorem \ref{thm.blprime}. \end{Proof}

\section{Proof of Theorem \ref{thm.entdim}} \label{sec.pfapp.special.entdim}

\subsection{Proof of Theorem \ref{thm.entdim}: imports and remarks} \label{sec.pfapp.special.einfleb.import}

\subsubsection{Proof of Theorem \ref{thm.entdim}: imports}

We now record two basic results that will be needed in the proof of Theorem \ref{thm.entdim}. The first appears, for example, as Theorems 5.18 and 7.17 of \cite{MR1681462}.

\begin{oLemma}[Weak-star compactness of the probability measures] \label{olem.compact}
The set $\prob(\rzee)$ is a compact subset of $C(\rzee)^\dagger$ in the weak-star topology. \end{oLemma}

The next result appears, for example, in Ch.~I Theorem 3.1 of \cite{MR2039503}.

\begin{oLemma}[Uniform approximation by trigonometric polynomials] \label{olem.weier}
Let $f \in C(\rzee)$ and let $\epsilon \in (0,\infty)_{\arr}$. Then there exist $N \in \enn$ and $\alpha \in \cee^{[0,2N]_{\zee}}$ such that the following holds.
\begin{align} \sup_{x \in \rzee} \lvertt f(x) - \sum_{j=0}^{2N} \alpha(j)\, \mathfrak{e}^{j-N}(x) \rvertt \leq \epsilon \label{olem.weier.conc1} \end{align} \end{oLemma}

\subsubsection{Proof of Theorem \ref{thm.entdim}: remarks}

We now collect certain straightforward observations that will be needed in the proof of Theorem \ref{thm.entdim}.

\begin{Remark}[Lebesgue norm of a trigonometric polynomial] \label{rem.plancherel}
Let $\mathscr{F}$ be a nonempty finite subset of $\zee$, let $\alpha \in \cee^{\mathscr{F}}$ and define $\psi \coloneqq \sum_{j \in \mathscr{F}} \alpha(j) \mathfrak{e}^j$. We obtain from the elementary Plancherel identity that $\lambda(\lvert \psi \rvert^2) = \sum_{j \in \mathscr{F}} \lvert \alpha(j) \rvert^2$. \end{Remark}

\begin{Remark}[Submeasures and nonnegative continuous functions] \label{rem.submeaspos}
Let $\mu,\nu \in \operatorname{FMea}(\rzee)$. Combining \ref{conv.submeasure} with the identification of $\operatorname{FMea}(\rzee)$ inside $C(\rzee)^\dagger$ recorded in \ref{conv.dual}, we find that $\mu \preceq \nu$ if and only if $\mu(f) \leq \nu(f)$ for all nonnegative $f \in C(\rzee)$. \end{Remark}

\begin{Remark}[Closedness of bounds on the exponential infimum] \label{rem.einfclosed}
Let $\mathscr{F}$ be a nonempty finite subset of $\zee$ and let $c \in \arr$. Definition \ref{def.foustat} exhibits the set of $\zeta \in \prob(\rzee)$ satisfying $\einf(\zeta,\mathscr{F}) \geq c$ as the intersection, over all $\alpha \in \cee^{\mathscr{F}}$ with $\sum_{j \in \mathscr{F}} \lvert \alpha(j) \rvert^2 = 1$, of the sets of $\zeta \in \prob(\rzee)$ satisfying the inequality below.
\begin{align} \zeta\left( \lvertt \sum_{j \in \mathscr{F}} \alpha(j) \, \mathfrak{e}^j \rvertt^2 \right) \geq c \label{rem.einfclosed.1} \end{align}
Since each set appearing in this intersection is closed by \ref{conv.weakstar}, we find that the set of $\zeta \in \prob(\rzee)$ satisfying $\einf(\zeta,\mathscr{F}) \geq c$ is closed. \end{Remark}

\subsection{Proof of Theorem \ref{thm.entdim}: preliminary lemma on Lebesgue minorization} \label{sec.pfapp.special.entdim.einfleb}

Lemma \ref{lem.einfleb} below is a routine assertion that connects uniform Fourier analysis with measure theory. It serves as our bridge between the finitary apparatus and the classical conclusion. It is consumed exactly once, in the proof of Theorem \ref{thm.entdim}, where it converts the exponential infima of a weak-star limit into the minorization that theorem asserts.

\begin{Lemma}[Lebesgue minorization from exponential infima] \label{lem.einfleb}
Let $\nu \in \prob(\rzee)$ and let $c \in [0,1]_{\arr}$ be such that $\einf(\nu,n) \geq c$ for all $n \in \enn$. Then $\nu \succeq c \lambda$. \end{Lemma}

\begin{Proof}[of Lemma \ref{lem.einfleb}] Let $\nu$ and $c$ be as in the statement of Lemma \ref{lem.einfleb}. By Remark \ref{rem.submeaspos} it is enough to show that $c\lambda(f) \leq \nu(f)$ for every nonnegative $f \in C(\rzee)$. \\

We first claim that the following holds for every trigonometric polynomial $\psi$.
\begin{align} c\,\lambda\bigl( \lvert \psi \rvert^2 \bigr) \leq \nu\bigl( \lvert \psi \rvert^2 \bigr) \label{lem.einfleb.p1.1} \end{align}
Both sides of (\ref{lem.einfleb.p1.1}) are homogeneous of degree two in $\psi$, so we may assume that $\lambda(\lvert \psi \rvert^2) = 1$, the case $\psi = 0$ being trivial. Write $\psi = \sum_{j=0}^{2N} \alpha(j) \mathfrak{e}^{j-N}$ with $N \in \enn$ and $\alpha \in \cee^{[0,2N]_{\zee}}$, so that Remark \ref{rem.plancherel} gives $\sum_{j=0}^{2N} \lvert \alpha(j) \rvert^2 = 1$. Since $\mathfrak{e}^N$ is unimodular we have $\lvert \psi \rvert = \lvert \mathfrak{e}^N \psi \rvert$, and $\mathfrak{e}^N \psi$ is the linear combination of the $\mathfrak{e}^j$ with $j \in [0,2N]_{\zee}$ determined by $\alpha$.\\

We compute:
\begin{align} \nu\bigl( \lvert \psi \rvert^2 \bigr) & \geq \einf\bigl(\nu , [0,2N]_{\zee}\bigr) \label{lem.einfleb.p1.2}
\\ & = \einf(\nu , 2N) \label{lem.einfleb.p1.3}
\\ & \geq c \label{lem.einfleb.p1.4} \end{align}
This computation may be justified as follows.
\begin{itemize}
\item The inequality in (\ref{lem.einfleb.p1.2}) follows by instantiating Definition \ref{def.foustat} with the bandprofile $[0,2N]_{\zee}$ and the coefficient vector $\alpha$.
\item (\ref{lem.einfleb.p1.3}) follows from (\ref{lem.einfleb.p1.2}) by \ref{conv.einfint}.
\item (\ref{lem.einfleb.p1.4}) follows from (\ref{lem.einfleb.p1.3}) since $\einf(\nu,m) \geq c$ for all $m \in \enn$ by hypothesis.
\end{itemize}
Since $\lambda(\lvert \psi \rvert^2) = 1$, the computation ending in (\ref{lem.einfleb.p1.4}) verifies (\ref{lem.einfleb.p1.1}).\\

For an arbitrary nonnegative $f \in C(\rzee)$, Lemma \ref{olem.weier} guarantees that we may uniformly approximate the continuous function $\sqrt{f}$ by trigonometric polynomials $\psi$. We thus obtain a uniform approximation to $f$ by $\lvert \psi \rvert^2$. Combining this with the bound (\ref{lem.einfleb.p1.1}), which holds for every trigonometric polynomial, we find $c\lambda(f) \leq \nu(f)$ as required. This completes the proof of Lemma \ref{lem.einfleb}.\end{Proof}

\subsection{Proof of Theorem \ref{thm.entdim}: main} \label{sec.pfapp.special.entdim.main}

An outline of the proof of Theorem \ref{thm.entdim} is as follows.
\begin{itemize}
\item For each pair $(r,\epsilon)$ separately, the constraint set $\mathscr{C}(r,\epsilon)$ of (\ref{thm.entdim.p2.1}) is shown to be nonempty. In the nondegenerate case when $\sigma \coloneqq D - \epsilon >0$, Lemma \ref{lem.ratio} converts sublacunarity into a threshold $t$ beyond which the depth clears the $\ell$ of (\ref{thm.entdim.p0}).
\item Theorem \ref{thm.exist} then returns a multiplier, and the twisted measure $m_\circledcirc \mu$ witnesses the nonemptiness by way of the computation ending in (\ref{thm.entdim.p1.3}), where the two losses $\epsilon$ combine as $\sigma - \epsilon = D - 2\epsilon$.
\item Compactness assembles the separate witnesses into one. The constraint sets are weak-star closed and shrink as $r$ grows and $\epsilon$ falls, so the diagonal sequence $\mathscr{C}_k \coloneqq \mathscr{C}(k,1/k)$ is a nested family of nonempty closed subsets of the compact space $\prob(\rzee)$, and its intersection holds a single measure $\nu$ spread at every resolution simultaneously. This is the article's one deliberately soft step.
\item Membership of $\nu$ in every $\mathscr{C}_k$ forces $\einf(\nu,n) \geq D$ for all $n \in \enn$, and Lemma \ref{lem.einfleb} converts this family of exponential infima into the minorization $\nu \succeq D\lambda$ demanded by the conclusion.
\end{itemize}

The dependencies in the proof of Theorem \ref{thm.entdim} beyond those in the statement are as follows.
\begin{deps}
\dep{Conventions}{\depitem{$\einf(\zeta,n)$}{conv.einfint}}
\dep{Definitions}{\depitem{Sublacunary sets}{def.sublac}\depsep\depitem{Exponential infimum}{def.foustat}}
\dep{Remarks}\depitem{Normalization of resolution entropy}{rem.achnunit}\depcont{\depitem{Submeasures and endomorphs}{rem.pushsub}}
\depcont{\depitem{Endomorphy of exponential infimum}{rem.einfendo}} \depcont{\depitem{Closedness of bounds on the exponential infimum}{rem.einfclosed}}
\depcont{\depitem{Monotonicity of exponential infimum}{rem.statmono}}
\dep{Lemmas}{\depitem{Sublacunary sets and progression-shadowing depth}{lem.ratio}\depsep\depitem{Lebesgue minorization from exponential infima}{lem.einfleb}}
\dep{Theorems}{\depitem{Existential endomorphic rigidity via entropy}{thm.exist}}
\dep{Imports}{\depitem{Weak-star compactness of the probability measures}{olem.compact}}
\end{deps}

We now move on to the formal proof.

\begin{Proof}[of Theorem \ref{thm.entdim}] Let $\mathscr{M}$ and $\mu$ be as in the statement of Theorem \ref{thm.entdim} and define $D \in [0,1]_\arr$ by $D \coloneqq \ach_\ast^\uparrow(\mu)$, where membership in the unit interval follows from Remark \ref{rem.achnunit}. In keeping with the statement of the current theorem, all topological notions applied to sets of measures should be understood with respect to the weak-star topology.\\

Let $r \in \enn$ and let $\epsilon \in (0,1]_\arr$, and define a subset $\mathscr{C}(r,\epsilon)$ of $\prob(\rzee)$ as follows.
\begin{align} \mathscr{C}(r,\epsilon) \coloneqq \Bigl\{ \varsigma \in \overline{\Sigma(\mu,\mathscr{M})} \, : \, \einf(\varsigma,r) \geq D - 2\epsilon \Bigr\} \label{thm.entdim.p2.1} \end{align}
Since $\mathscr{M}$ is sublacunary, it follows immediately from Definition \ref{def.sublac} that $\mathscr{M}$ is infinite. Thus $\mathscr{M} \cap \ennplus$ is nonempty and hence so is $\overline{\Sigma(\mu,\mathscr{M})}$.\\

We now dispose of the degenerate case where $D \leq 2\epsilon$. Definition \ref{def.foustat} provides $\einf(\varsigma,r) \geq 0 \geq D - 2\epsilon$ for every $\varsigma \in \prob(\rzee)$, so that in this case $\mathscr{C}(r,\epsilon) = \overline{\Sigma(\mu,\mathscr{M})}$ and is therefore nonempty.\\

We may therefore assume $D > 2\epsilon$ and assign $\sigma \coloneqq D - \epsilon$, which lies in $(0,1]_\arr$ since $\epsilon < D \leq 1$. We further define $\ell \in \ennplus$ as follows.
\begin{align} \ell \coloneqq \left\lceil \frac{128(1-\sigma+\epsilon)}{\epsilon}\log_2\left(\frac{8(r+2)}{\epsilon}\right) \right\rceil \label{thm.entdim.p0} \end{align}
Since $\mathscr{M}$ is sublacunary, Clause \ref{lem.ratio.c2} of Lemma \ref{lem.ratio} provides $\lim_{s \to \infty} \mathtt{Dep}_{\mathscr{M}}(s) = \infty$, and hence there exists $t \in \enn$ such that $\mathtt{Dep}_{\mathscr{M}}(s) \geq \ell$ for all $s \in [t,\infty)_\zee$. Moreover $D$ is the limit superior of $\achn(\mu,N)$ as $N$ tends to infinity, so that the set of $N \in \ennplus$ satisfying $\achn(\mu,N) \geq \sigma$ is infinite, and we fix such an $N$ lying in $[2^{9(t+\ell)/\epsilon},\infty)_\zee$.\\

Thus the data $(\epsilon,r,\sigma,\ell,\mathscr{M},t,N,\mu)$ satisfies the hypotheses gathered in Clause \ref{thm.exist.c1} of Theorem \ref{thm.exist}. Instantiating that theorem we are returned $m \in \mathscr{M} \cap [1,2N]_{\zee}$ for which the conclusion (\ref{thm.exist.conc1}) holds. Since $m$ lies in $\ennplus$, Remark \ref{rem.pushsub} provides that $m_\circledcirc \mu$ is an element of $\prob(\rzee)$ and hence of $\Sigma(\mu,\mathscr{M})$.\\

We compute:
\begin{align} \einf\bigl( m_\circledcirc \mu , r \bigr) & = \einf\bigl( \mu , \{m\}[0,r]_{\zee} \bigr) \label{thm.entdim.p1.1}
\\ & \geq \sigma - \epsilon \label{thm.entdim.p1.2}
\\ & = D - 2\epsilon \label{thm.entdim.p1.3} \end{align}
This computation may be justified as follows.
\begin{itemize}
\item The equality in (\ref{thm.entdim.p1.1}) follows by combining \ref{conv.einfint} with Remark \ref{rem.einfendo}, where we note that $m$ is an element of $\ennplus$.
\item (\ref{thm.entdim.p1.2}) follows from (\ref{thm.entdim.p1.1}) by the conclusion (\ref{thm.exist.conc1}) of Theorem \ref{thm.exist}.
\item (\ref{thm.entdim.p1.3}) follows from (\ref{thm.entdim.p1.2}) by the assignment of $\sigma$.
\end{itemize}
The computation ending in (\ref{thm.entdim.p1.3}) shows that $m_\circledcirc \mu \in \mathscr{C}(r,\epsilon)$, so that $\mathscr{C}(r,\epsilon)$ is nonempty in the nondegenerate case as well. \\

We now allow $r \in \enn$ and $\epsilon \in (0,1]_\arr$ to vary. Remark \ref{rem.einfclosed} instantiated with the bandprofile $[0,r]_{\zee}$ provides that the set of $\varsigma \in \prob(\rzee)$ satisfying $\einf(\varsigma,r) \geq D-2\epsilon$ is closed. As $\overline{\Sigma(\mu,\mathscr{M})}$ is also closed, so is $\mathscr{C}(r,\epsilon)$. 
 Moreover, Remark \ref{rem.statmono} guarantees that the following holds for all $r_1,r_2 \in \enn$ and $\epsilon_1,\epsilon_2 \in (0,1]_\arr$ such that $r_1 \leq r_2$ and $\epsilon_1 \geq \epsilon_2$.
\begin{align*} \mathscr{C}(r_1,\epsilon_1) \supseteq  \mathscr{C}(r_2,\epsilon_2)\end{align*}
For $k \in \ennplus$ we write $\mathscr{C}_k \coloneqq \mathscr{C}(k,1/k)$, so that the previous display provides that $\mathscr{C}_{k+1}$ is a subset of $\mathscr{C}_k$. In summary, we find that $(\mathscr{C}_k)_{k \in \ennplus}$ is a downward-nested chain of nonempty closed subsets of $\overline{\Sigma(\mu,\mathscr{M})}$. Lemma \ref{olem.compact} provides that $\prob(\rzee)$ is compact, and so its closed subset $\overline{\Sigma(\mu,\mathscr{M})}$ is compact as well. Consequently:
\begin{align*} \emptyset \neq \bigcap_{k \in \ennplus} \mathscr{C}_k \subseteq \overline{\Sigma(\mu,\mathscr{M})} \end{align*}
We fix an element $\nu$ of the above intersection, and assign this as the parameter demanded by the current conclusion. Thus it remains to verify that $\nu \succeq D \lambda$. To that end, let $n \in \enn$ and let $k \in \ennplus$ satisfy $k \geq n$. We have:
\begin{align*} \einf(\nu,n) \geq \einf(\nu,k) \geq D - \frac{2}{k} \end{align*}
Here, the left inequality follows from Remark \ref{rem.statmono}, while the right inequality follows from the membership of $\nu$ in $\mathscr{C}_k$. Letting $k \to \infty$ over $[n,\infty)_\zee$ in the above display shows that $\einf(\nu,n) \geq D$, where $n \in \enn$ was arbitrary.\\

Instantiating Lemma \ref{lem.einfleb} with $\nu$ and $D$ now returns $\nu \succeq D \lambda$, and thereby completes the proof of Theorem \ref{thm.entdim}. \end{Proof}

\section{Proof of Theorem \ref{thm.rjcanon}}

\subsection{Proof of Theorem \ref{thm.rjcanon}: imports} \label{sec.pfapp.rjcanon.imports}

We now record two classical results which will be needed in the proof of Theorem \ref{thm.rjcanon}. The first is Furstenberg's lemma on nonlacunary semigroups, Lemma IV.1 of \cite{MR0213508}. Furstenberg defines a multiplicative semigroup of positive integers to be nonlacunary if it is not contained in the set of powers of a single integer, which is equivalent to the semigroup containing two multiplicatively independent elements, and we state the lemma in the latter form.

\begin{oLemma}[Furstenberg's nonlacunary semigroup lemma] \label{olem.furst}
Let $\mathscr{S}$ be a subset of $\ennplus$ which is closed under multiplication and contains two multiplicatively independent elements, and write $(m_j)_{j \in \enn}$ for the increasing enumeration of $\mathscr{S}$. Then $\lim_{j \to \infty} m_{j+1}/m_j = 1$. \end{oLemma}

The second is the asymptotic formula of Hardy and Ramanujan for the partition function $\mathtt{p}$ of Definition \ref{def.canonsets}, which appears in \cite{MR2280879}.

\begin{oLemma}[Hardy--Ramanujan asymptotic] \label{olem.hr}
We have:
\begin{align*} \lim_{n \to \infty} 4\sqrt{3} \, n \, \mathtt{p}(n) \exp\left(-\pi\sqrt{\frac{2n}{3}}\right) = 1 \end{align*} \end{oLemma}

\subsection{Proof of Theorem \ref{thm.rjcanon}: preliminary lemmas} \label{sec.pfapp.rjcanon.lemmas}

\subsubsection{Sublacunarity of the integer parts of a sequence}

Lemma \ref{lem.floorseq} is consumed twice in the proof of Theorem \ref{thm.rjcanon}, for the logarithmic powers and for the partition numbers, where it converts the growth of a defining sequence into sublacunarity of the set of its integer parts. The hypothesis on consecutive ratios passes to the integer parts at the cost of the additive error of the floor function, and this error is absorbed because the sequence is unbounded.\\

\begin{Lemma}[Sublacunarity of the integer parts of a sequence] \label{lem.floorseq}
Let $g: \ennplus \to [1,\infty)_\arr$ be nondecreasing and satisfy $\lim_{n \to \infty} g(n) = \infty$ and $\lim_{n \to \infty} g(n+1)/g(n) = 1$. Then the set $\{ \lfloor g(n) \rfloor : n \in \ennplus \}$ is sublacunary. \end{Lemma}

\begin{Proof}[of Lemma \ref{lem.floorseq}] Let $g$ be as in the statement of Lemma \ref{lem.floorseq}, and define $\mathscr{M} \coloneqq \{ \lfloor g(n) \rfloor : n \in \ennplus \}$. Since $g$ takes values in $[1,\infty)_\arr$ we have $\mathscr{M} \subseteq \ennplus$, so that $\mathscr{M} \setminus \{0\} = \mathscr{M}$, and since $g$ is unbounded so is $\mathscr{M}$. In particular $\mathscr{M}$ is infinite, and we write $(m_j)_{j \in \enn}$ for its increasing enumeration.\\

Let $\epsilon \in (0,1]_\arr$. Applying the hypothesized limit of the ratios $g(n+1)/g(n)$, we find there exists $N \in \ennplus$ such that the following holds for all $n \in [N,\infty)_\zee$.
\begin{align} g(n+1) \leq (1+\epsilon) g(n) \label{lem.floorseq.p1} \end{align}
Since $(m_j)_{j \in \enn}$ is unbounded, there exists $J \in \enn$ such that $m_J > g(N)$ and $m_J \geq 1/\epsilon$. Let $j \in [J,\infty)_\zee$. Since $m_j \in \mathscr{M}$, the set $\mathscr{N} \coloneqq \{ n \in \ennplus : \lfloor g(n) \rfloor = m_j \}$ is nonempty, and since $g$ is unbounded we have $\lfloor g(n) \rfloor \geq m_j + 1$ for all sufficiently large $n$, so that $\mathscr{N}$ is finite. We define $u \coloneqq \max(\mathscr{N})$, so that $\lfloor g(u) \rfloor = m_j$. We compute:
\begin{align} g(u) & \geq m_j \label{lem.floorseq.p2.2}
\\ & \geq m_J \label{lem.floorseq.p2.3}
\\ & > g(N) \label{lem.floorseq.p2.4} \end{align}
This computation may be justified as follows.
\begin{itemize}
\item The inequality in (\ref{lem.floorseq.p2.2}) follows from the assignment of $u$.
\item (\ref{lem.floorseq.p2.3}) follows from (\ref{lem.floorseq.p2.2}) since $j \geq J$ and the sequence $(m_j)_{j \in \enn}$ is increasing.
\item (\ref{lem.floorseq.p2.4}) follows from (\ref{lem.floorseq.p2.3}) by the choice of $J$.
\end{itemize}
Since $g$ is nondecreasing, the computation ending in (\ref{lem.floorseq.p2.4}) provides $u > N$, and in particular $u \in [N,\infty)_\zee$. Since $g$ is nondecreasing and $\lfloor g(u) \rfloor = m_j$ by the assignment of $u$, we have $\lfloor g(u+1) \rfloor \geq m_j$, while $u + 1 \notin \mathscr{N}$ by the assignment of $u$ provides $\lfloor g(u+1) \rfloor \neq m_j$. Thus:
\begin{align} \lfloor g(u+1) \rfloor \geq m_j + 1 \label{lem.floorseq.p2.5} \end{align}
Since $\lfloor g(u+1) \rfloor$ is an element of $\mathscr{M}$, and since $m_{j+1}$ is the least element of $\mathscr{M}$ exceeding $m_j$, the inequality (\ref{lem.floorseq.p2.5}) provides:
\begin{align} m_{j+1} \leq \lfloor g(u+1) \rfloor \label{lem.floorseq.p2.6} \end{align}
We compute again:
\begin{align} m_{j+1} & \leq \lfloor g(u+1) \rfloor \label{lem.floorseq.p3.1}
\\ & \leq (1+\epsilon) g(u) \label{lem.floorseq.p3.3}
\\ & < (1+\epsilon)(m_j + 1) \label{lem.floorseq.p3.4}
\\ & \leq (1+\epsilon)^2 m_j \label{lem.floorseq.p3.5}
\\ & \leq (1 + 3\epsilon) m_j \label{lem.floorseq.p3.6} \end{align}
This computation may be justified as follows.
\begin{itemize}
\item The inequality in (\ref{lem.floorseq.p3.1}) follows from (\ref{lem.floorseq.p2.6}).
\item (\ref{lem.floorseq.p3.3}) follows from (\ref{lem.floorseq.p3.1}) by instantiating (\ref{lem.floorseq.p1}) with $n \coloneqq u$, which is legitimate since $u \in [N,\infty)_\zee$.
\item (\ref{lem.floorseq.p3.4}) follows from (\ref{lem.floorseq.p3.3}) since $g(u) < \lfloor g(u) \rfloor + 1 = m_j + 1$.
\item (\ref{lem.floorseq.p3.5}) follows from (\ref{lem.floorseq.p3.4}) since $m_j \geq m_J \geq 1/\epsilon$ provides $m_j + 1 \leq (1+\epsilon) m_j$.
\item (\ref{lem.floorseq.p3.6}) follows from (\ref{lem.floorseq.p3.5}) since $\epsilon \leq 1$ provides $(1+\epsilon)^2 \leq 1 + 3\epsilon$.
\end{itemize}
Since the sequence $(m_j)_{j \in \enn}$ is increasing, dividing the computation ending in (\ref{lem.floorseq.p3.6}) by $m_j$ provides $1 \leq m_{j+1}/m_j \leq 1 + 3\epsilon$. Since $j$ was arbitrary in $[J,\infty)_\zee$ and $\epsilon$ was arbitrary in $(0,1]_\arr$, we have verified that for every $\epsilon \in (0,1]_\arr$ there exists $J \in \enn$ such that $1 \leq m_{j+1}/m_j \leq 1 + 3\epsilon$ for all $j \in [J,\infty)_\zee$. This is (\ref{def.sublac.1}), and so $\mathscr{M}$ is sublacunary. This completes the proof of Lemma \ref{lem.floorseq}. \end{Proof}

\subsection{Proof of Theorem \ref{thm.rjcanon}: main} \label{sec.pfapp.rjcanon.main}

An outline of the proof of Theorem \ref{thm.rjcanon} is as follows.
\begin{itemize}
\item Theorem \ref{thm.entdim} carries all the analysis. Each clause reduces to the assertion that its multiplier set is sublacunary, and this reduction is recorded once at the opening of the proof. Nothing below consumes the depth bounds of Lemma \ref{lem.depbound}, whose rates are spent in Theorem \ref{thm.blprime}.
\item Clause \ref{thm.rjcanon.c1} is verified directly on consecutive elements $m_j < m_{j+1}$. The cube interval above $m_j$ supplies an element within a factor tending to one, at (\ref{thm.rjcanon.p1.6}), and minimality of $m_{j+1}$ closes the gap. Clause \ref{thm.rjcanon.c4} is Furstenberg's lemma on nonlacunary semigroups, Lemma \ref{olem.furst}, applied to $\mathbb{S}_{a,b}$.
\item Clauses \ref{thm.rjcanon.c2} and \ref{thm.rjcanon.c3} are instances of Lemma \ref{lem.floorseq}. For the logarithmic powers the ratio hypothesis is the computation ending in (\ref{thm.rjcanon.p2.4}). For the partition numbers both the ratio hypothesis and unboundedness are read off the Hardy--Ramanujan asymptotic, at (\ref{thm.rjcanon.p3.3}) and (\ref{thm.rjcanon.p3.6}).
\end{itemize}

The dependencies in the proof of Theorem \ref{thm.rjcanon} are as follows.
\begin{deps}
\dep{Definitions}{\depitem{Multiplicative independence}{def.multindep}\depsep\depitem{Canonical multiplier sets}{def.canonsets}\depsep\depitem{Sublacunary sets}{def.sublac}}
\dep{Lemmas}{\depitem{Sublacunarity of the integer parts of a sequence}{lem.floorseq}}
\dep{Theorems}{\depitem{Limiting endomorphic rigidity via upper entropy dimension}{thm.entdim}}
\dep{Imports}{\depitem{Furstenberg's nonlacunary semigroup lemma}{olem.furst}\depsep\depitem{Hardy--Ramanujan asymptotic}{olem.hr}}
\end{deps}

We now move on to the formal proof.

\begin{Proof}[of Theorem \ref{thm.rjcanon}] Let $\mu$ be as in the statement of Theorem \ref{thm.rjcanon}. We first record a reduction common to the four clauses. Suppose $\mathscr{M}$ is a sublacunary subset of $\enn$. Instantiating Theorem \ref{thm.entdim} with $\mathscr{M}$ and $\mu$ returns $\nu \in \overline{\Sigma(\mu,\mathscr{M})}$ satisfying $\nu \succeq \ach_\ast^\uparrow(\mu)\lambda$. Thus it suffices to verify, for each clause of the current theorem, that the multiplier set specified by that clause is sublacunary.\\

\vspace{5 pt}

\underline{Verification of Clause \ref{thm.rjcanon.c1}}:\\

Let $n_0$ and $\mathscr{M}$ be as in Clause \ref{thm.rjcanon.c1}. The intervals $[n^3+1,(n+1)^3-1]_\zee$ for $n \in [n_0,\infty)_\zee$ are pairwise disjoint and each meets $\mathscr{M}$, so that $\mathscr{M}$ is infinite, and we write $(m_j)_{j \in \enn}$ for the increasing enumeration of $\mathscr{M} \setminus \{0\}$.\\

Now, let $\epsilon \in (0,1]_\arr$. Since $(m_j)_{j \in \enn}$ is unbounded, there exists $J \in \enn$ such that the following holds.
\begin{align*} m_J \geq n_0^3 \vee 8 \vee \left(\frac{14}{\epsilon} \right)^3 \end{align*} 
Let $j \in [J,\infty)_\zee$ and define $n \coloneqq \lfloor m_j^{1/3} \rfloor$, so that $n^3 \leq m_j < (n+1)^3$. Since $m_j \geq n_0^3$ we have $m_j^{1/3} \geq n_0$, whence $n \geq n_0$ since $n_0 \in \zee$. Since $n + 1 \in [n_0,\infty)_\zee$, the hypothesis of Clause \ref{thm.rjcanon.c1} returns $q \in \enn$ such that:
\begin{align} q \in \mathscr{M} \cap [(n+1)^3+1,(n+2)^3-1]_\zee \label{thm.rjcanon.p1.0} \end{align}
By (\ref{thm.rjcanon.p1.0}) we have $q \geq (n+1)^3 + 1 > m_j \geq 1$, so that $q$ is an element of $\mathscr{M} \setminus \{0\}$ exceeding $m_j$, and since $m_{j+1}$ is the least such element we obtain:
\begin{align} m_{j+1} \leq q \label{thm.rjcanon.p1.1} \end{align}
We compute:
\begin{align} q & \leq (n+2)^3 \label{thm.rjcanon.p1.2}
\\ & \leq \bigl( m_j^{1/3} + 2 \bigr)^3 \label{thm.rjcanon.p1.3}
\\ & = m_j \bigl( 1 + 2 m_j^{-1/3} \bigr)^3 \label{thm.rjcanon.p1.4}
\\ & \leq m_j \bigl( 1 + 14 m_j^{-1/3} \bigr) \label{thm.rjcanon.p1.5}
\\ & \leq (1 + \epsilon) m_j \label{thm.rjcanon.p1.6} \end{align}
This computation may be justified as follows.
\begin{itemize}
\item The inequality in (\ref{thm.rjcanon.p1.2}) follows from (\ref{thm.rjcanon.p1.0}), which provides $q \leq (n+2)^3 - 1$.
\item (\ref{thm.rjcanon.p1.3}) follows from (\ref{thm.rjcanon.p1.2}) since the assignment of $n$ provides $n \leq m_j^{1/3}$.
\item (\ref{thm.rjcanon.p1.4}) follows from (\ref{thm.rjcanon.p1.3}) by factoring out $m_j^{1/3}$ and cubing.
\item (\ref{thm.rjcanon.p1.5}) follows from (\ref{thm.rjcanon.p1.4}) by instantiating the elementary inequality $(1+y)^3 = 1 + 3y + 3y^2 + y^3 \leq 1 + 7y$, valid for all $y \in [0,1]_\arr$, with $y \coloneqq 2m_j^{-1/3}$, which lies in $[0,1]_\arr$ since $m_j \geq 8$.
\item (\ref{thm.rjcanon.p1.6}) follows from (\ref{thm.rjcanon.p1.5}) since $m_j \geq (14/\epsilon)^3$ provides $14 m_j^{-1/3} \leq \epsilon$.
\end{itemize}
Combining (\ref{thm.rjcanon.p1.1}) with the computation ending in (\ref{thm.rjcanon.p1.6}) we obtain:
\begin{align*} m_{j+1} \leq (1 + \epsilon) m_j \end{align*}
Since the sequence $(m_j)_{j \in \enn}$ is increasing, dividing the above display by $m_j$ provides $1 \leq m_{j+1}/m_j \leq 1 + \epsilon$. Since $j$ was arbitrary in $[J,\infty)_\zee$ and $\epsilon$ was arbitrary in $(0,1]_\arr$, this verifies (\ref{def.sublac.1}), and so $\mathscr{M}$ is sublacunary.\\
\vspace{5 pt}

\underline{Verification of Clause \ref{thm.rjcanon.c2}}:\\

We define $g: \ennplus \to \arr$ as follows.
\begin{align} g(n) \coloneqq n^{\log(n)} = \exp\bigl( \log(n)^2 \bigr) \label{thm.rjcanon.p2.0} \end{align}
Then (\ref{def.canonsets.c3.1}) provides $\mathbb{L} = \{ \lfloor g(n) \rfloor : n \in \ennplus \}$. Since $\log(n)^2$ is nonnegative and nondecreasing in $n \in \ennplus$ and tends to infinity with $n$, the assignment of $g$ in (\ref{thm.rjcanon.p2.0}) provides that $g$ takes values in $[1,\infty)_\arr$, is nondecreasing, and satisfies $\lim_{n \to \infty} g(n) = \infty$.\\

Now, let $n \in \ennplus$. We compute:
\begin{align} \log\left( \frac{g(n+1)}{g(n)} \right) & = \log(n+1)^2 - \log(n)^2 \label{thm.rjcanon.p2.1}
\\ & = \bigl( \log(n+1) - \log(n) \bigr) \bigl( \log(n+1) + \log(n) \bigr) \label{thm.rjcanon.p2.2}
\\ & \leq \frac{2\log(n+1)}{n} \label{thm.rjcanon.p2.3}
\\ & \leq \frac{2\sqrt{2}}{\sqrt{n}} \label{thm.rjcanon.p2.4} \end{align}
This computation may be justified as follows.
\begin{itemize}
\item The equality in (\ref{thm.rjcanon.p2.1}) follows from the assignment of $g$ in (\ref{thm.rjcanon.p2.0}).
\item (\ref{thm.rjcanon.p2.2}) follows from (\ref{thm.rjcanon.p2.1}) by factoring the difference of squares.
\item (\ref{thm.rjcanon.p2.3}) follows from (\ref{thm.rjcanon.p2.2}) since $\log(n+1) - \log(n) = \log(1 + 1/n) \leq 1/n$ and $\log(n) \leq \log(n+1)$.
\item (\ref{thm.rjcanon.p2.4}) follows from (\ref{thm.rjcanon.p2.3}) by instantiating the elementary inequality $\log(x) \leq \sqrt{x}$, valid for all $x \in [1,\infty)_\arr$, with $x \coloneqq n+1$, together with $\sqrt{n+1} \leq \sqrt{2n}$.
\end{itemize}
Since $g$ is nondecreasing, the left side of (\ref{thm.rjcanon.p2.1}) is nonnegative, and the right side of (\ref{thm.rjcanon.p2.4}) tends to $0$ as $n$ tends to infinity. Since $n$ was arbitrary in $\ennplus$, the computation ending in (\ref{thm.rjcanon.p2.4}) provides:
\begin{align} \lim_{n \to \infty} \frac{g(n+1)}{g(n)} = 1 \label{thm.rjcanon.p2.5} \end{align}
Thus $g$ satisfies the hypotheses of Lemma \ref{lem.floorseq}, and instantiating that lemma with $g$ provides that $\mathbb{L}$ is sublacunary.\\

\vspace{5 pt}

\underline{Verification of Clause \ref{thm.rjcanon.c3}}:\\

For $n \in \ennplus$ we define $A(n) \in (0,\infty)_\arr$ as follows.
\begin{align} A(n) \coloneqq 4\sqrt{3} \, n \, \mathtt{p}(n) \exp\left( -\pi\sqrt{\frac{2n}{3}} \right) \label{thm.rjcanon.p3.0} \end{align}
Lemma \ref{olem.hr} asserts the following.
\begin{align} \lim_{n \to \infty} A(n) = 1 \label{thm.rjcanon.p3.1} \end{align}
We define $g: \ennplus \to \arr$ by $g(n) \coloneqq \mathtt{p}(n)$. Since $\mathtt{p}$ is nondecreasing with values in $\ennplus$ and $\mathtt{p}(0) = 1 = \mathtt{p}(1)$, we have that $g$ is nondecreasing with values in $[1,\infty)_\arr$.\\

Now, let $n \in \ennplus$. We compute:
\begin{align} g(n) & = \frac{A(n)}{4\sqrt{3} \, n} \exp\left( \pi\sqrt{\frac{2n}{3}} \right) \label{thm.rjcanon.p3.2}
\\ & \geq \frac{A(n) \sqrt{n}}{4\sqrt{3}} \label{thm.rjcanon.p3.3} \end{align}
This computation may be justified as follows.
\begin{itemize}
\item The equality in (\ref{thm.rjcanon.p3.2}) follows from the assignment of $g$ and the assignment of $A$ in (\ref{thm.rjcanon.p3.0}).
\item (\ref{thm.rjcanon.p3.3}) follows from (\ref{thm.rjcanon.p3.2}) by instantiating the elementary inequality $\exp(x) \geq x^3/6$, valid for all $x \in [0,\infty)_\arr$, with $x \coloneqq \pi\sqrt{2n/3}$.
\end{itemize}
Since (\ref{thm.rjcanon.p3.1}) provides that $A(n)$ tends to $1$, and since $n$ was arbitrary in $\ennplus$, the computation ending in (\ref{thm.rjcanon.p3.3}) provides $\lim_{n \to \infty} g(n) = \infty$.\\

Now, let $n \in \ennplus$ again. We compute:
\begin{align} \frac{g(n+1)}{g(n)} & = \frac{A(n+1)}{A(n)} \cdot \frac{n}{n+1} \cdot \exp\left( \pi\sqrt{\frac{2}{3}} \bigl( \sqrt{n+1} - \sqrt{n} \bigr) \right) \label{thm.rjcanon.p3.4}
\\ & = \frac{A(n+1)}{A(n)} \cdot \frac{n}{n+1} \cdot \exp\left( \frac{\pi\sqrt{2/3}}{\sqrt{n+1} + \sqrt{n}} \right) \label{thm.rjcanon.p3.5} \end{align}
This computation may be justified as follows.
\begin{itemize}
\item The equality in (\ref{thm.rjcanon.p3.4}) follows from (\ref{thm.rjcanon.p3.2}), instantiated at $n$ and at $n+1$.
\item (\ref{thm.rjcanon.p3.5}) follows from (\ref{thm.rjcanon.p3.4}) since $\bigl(\sqrt{n+1} - \sqrt{n}\bigr)\bigl(\sqrt{n+1} + \sqrt{n}\bigr) = 1$.
\end{itemize}
Each of the three factors in (\ref{thm.rjcanon.p3.5}) tends to $1$ as $n$ tends to infinity: the first by (\ref{thm.rjcanon.p3.1}), the second since $n/(n+1)$ tends to $1$, and the third since $1/(\sqrt{n+1} + \sqrt{n})$ tends to $0$. Since $n$ was arbitrary in $\ennplus$, the computation ending in (\ref{thm.rjcanon.p3.5}) provides:
\begin{align} \lim_{n \to \infty} \frac{g(n+1)}{g(n)} = 1 \label{thm.rjcanon.p3.6} \end{align}
Thus $g$ satisfies the hypotheses of Lemma \ref{lem.floorseq}, and instantiating that lemma with $g$ provides that $\mathtt{p}(\enn)$ is sublacunary.\\

\vspace{5 pt}

\underline{Verification of Clause \ref{thm.rjcanon.c4}}:\\

Let $a,b \in [2,\infty)_\zee$ be multiplicatively independent. By (\ref{def.canonsets.c5.1}) the set $\mathbb{S}_{a,b}$ is a subset of $\ennplus$ which is closed under multiplication, contains $a^c$ for every $c \in \enn$ and hence is infinite, and contains the multiplicatively independent elements $a$ and $b$. Thus $\mathbb{S}_{a,b}$ satisfies the hypotheses of Lemma \ref{olem.furst}, and $\mathbb{S}_{a,b} \setminus \{0\} = \mathbb{S}_{a,b}$.\\

Writing $(m_j)_{j \in \enn}$ for the increasing enumeration of $\mathbb{S}_{a,b}$, Lemma \ref{olem.furst} provides $\lim_{j \to \infty} m_{j+1}/m_j = 1$, which is (\ref{def.sublac.1}), and so $\mathbb{S}_{a,b}$ is sublacunary. This completes the verification of the four clauses and thereby the proof of Theorem \ref{thm.rjcanon}.\end{Proof}

\section{Proof of Theorem \ref{thm.ad.cantor}} \label{sec.pfapp.cantor}

\subsection{Proof of Theorem \ref{thm.ad.cantor}: preliminaries}

We first make precise the measure named before the statement of Theorem \ref{thm.ad.cantor}, and record two classical facts about it.

\begin{Definition}[Bias-$p$ ternary Cantor measure] \label{def.cantor} Let $p \in (0,1/2]_\arr$. We define the \emph{bias-$p$ ternary Cantor measure} to be the self-similar measure of Theorem 4.4 in \cite{MR0625600} for the two similitudes of the line with ratio $1/3$ fixing $0$ and $1$ respectively, taken with weights $1-p$ and $p$ respectively, and we write $\nu_p$ for its image in $\prob(\rzee)$. \end{Definition}

Thus $\nu_p$ is a $\times 3$-invariant measure carried by the middle-thirds Cantor set, and it is the distribution of a point of that set whose ternary digits are chosen independently, each equal to $2$ with probability $p$ and to $0$ otherwise, as in Theorem 4.4(4) of the same reference.

\begin{Remark}[Fourier coefficients of the Cantor measure] \label{rem.cantorfou} For $p \in (0,1/2]_\arr$ we have:
\begin{align} \label{rem.cantorfou.1} \lvert \hat{\nu}_p(1) \rvert^2 =  \prod_{j=1}^{\infty} \left(1 - 4p(1-p)\sin^2\left(\frac{2\pi}{3^j}\right)\right) \end{align}
This is the classical product formula for the Fourier transform of an infinite convolution, which is due to Strichartz \cite{MR1078738} and for coin-tossing measures of the present kind appears as (2.1) and (2.2) of \cite{MR4042627}.\end{Remark}

\begin{Remark}[Resolution entropy of the Cantor measure] \label{rem.cantorent} Let $p \in (0,1/2]_\arr$ and $k \in \ennplus$. Then:
\begin{align} \label{rem.cantorent.1} \ach(\nu_p,3^k) = -k\bigl(p\log(p) + (1-p)\log(1-p)\bigr) \end{align}
This is the classical computation of the entropy of a Bernoulli shift, transferred to resolution cells as in Theorem 2.6.6 in \cite{MR2239987}. \end{Remark}

\subsection{Proof of Theorem \ref{thm.ad.cantor}: main}

The dependencies in the proof of Theorem \ref{thm.ad.cantor} beyond those in the statement are as follows.
\begin{deps}
\dep{Definitions}{\depitem{Bias-$p$ ternary Cantor measure}{def.cantor}}
\dep{Remarks}{\depitem{Fourier coefficients of the Cantor measure}{rem.cantorfou} \depsep \depitem{Resolution entropy of the Cantor measure}{rem.cantorent}}
\end{deps}

We now move on to the formal proof.

\begin{Proof}[of Theorem \ref{thm.ad.cantor}] Let $p$ be as in the statement of Theorem \ref{thm.ad.cantor}. We compute:
\begin{align} \lvert \hat{\nu}_p(1) \rvert^2 & \geq 1 - 4p(1-p)\sum_{j=1}^{\infty} \sin^2\left(\frac{2\pi}{3^j}\right) \label{thm.ad.cantor.1}
\\ & \geq 1 - 4p(1-p)\left(\frac{3}{4} + \sum_{j=2}^{\infty} \frac{4\pi^2}{9^j}\right) \label{thm.ad.cantor.2}
\\ & \geq 1 - 6p \label{thm.ad.cantor.3} \end{align}
This computation may be justified as follows.
\begin{itemize}
\item (\ref{thm.ad.cantor.1}) follows from (\ref{rem.cantorfou.1}) in  Remark \ref{rem.cantorfou} by applying the Weierstrass product inequality (see \cite{MR0274686}).
\item (\ref{thm.ad.cantor.2}) follows from (\ref{thm.ad.cantor.1}) since $\sin^2(2\pi/3) = 3/4$ and $\lvert \sin(y) \rvert \leq \lvert y \rvert$ for all $y \in \arr$.
\item (\ref{thm.ad.cantor.3}) follows from (\ref{thm.ad.cantor.2}) since we have:
\begin{align*} \frac{3}{4}+\sum_{j=2}^\infty \frac{4\pi^2}{9^j} = \frac{3}{4}+\frac{\pi^2}{18} \leq \frac{3}{2} \end{align*}
\end{itemize}
We obtain:
\begin{align} \label{thm.ad.cantor.x1} \inf_{j \in \enn} \lvert \hat{\nu}_p(3^j) \rvert = \lvert \hat{\nu}_p(1) \rvert \geq 1-6p  \end{align}
Here, the equality follows from the $\times 3$-invariance of $\nu_p$, while the right inequality follows by assuming without loss of generality that $1-6p \geq 0$ and using the elementary inequality $\sqrt{1-y} \geq 1-y$. On the other hand we compute:
\begin{align} \frac{8}{\log(1/p)}\cdot \ach_\ast^\uparrow(\nu_p) & \geq  \frac{8}{\log(1/p)} \cdot \limsup_{k \to \infty} \frac{\ach(\nu_p,3^k)}{k\log(3)} \label{thm.ad.cantor.5.0}
\\ & = -\frac{8}{\log(1/p)} \cdot \frac{p\log(p) + (1-p)\log(1-p)}{\log(3)} \label{thm.ad.cantor.5}
\\ & \geq \frac{8p}{\log(3)} \label{thm.ad.cantor.6}
\\ & \geq 6p \label{thm.ad.cantor.7} \end{align}
This computation may be justified as follows.
\begin{itemize}
\item The inequality in (\ref{thm.ad.cantor.5.0}) follows directly from the definition of $\ach_\ast^\uparrow$ in Definition \ref{def.entropy}.
\item (\ref{thm.ad.cantor.5}) follows from (\ref{rem.cantorent.1}) in Remark \ref{rem.cantorent}.
\item (\ref{thm.ad.cantor.6}) follows from (\ref{thm.ad.cantor.5}) since $-(1-p)\log(1-p) \geq 0$.
\item (\ref{thm.ad.cantor.7}) follows from (\ref{thm.ad.cantor.6}) since $\log(3) \leq 1.1$.
\end{itemize}
Combining (\ref{thm.ad.cantor.x1}) with (\ref{thm.ad.cantor.7}) completes the proof of Theorem \ref{thm.ad.cantor}. \end{Proof}

\appendix

\chapter{Dependency flowchart} \label{app.flowchart}

Figures \ref{fig.flowchart.a} and \ref{fig.flowchart.b} display every proved result of the article, every imported result, and the five theorems of the Invitation, arranged in bands by chapter.
\begin{itemize} \item An arrow from one box to another asserts that the first is consumed in the proof of the second, as recorded in the dependency box preceding that proof, or in the proof itself for the short proofs of Chapter \ref{cha.pfapp} which carry no box. A dotted arrow marks a corollary, or a theorem of the Invitation, derived without a separate proof.
\item The stacked small boxes beneath or beside a result list the imports and the results from the other figure consumed in its proof, each with a page reference, so that no arrow needs to cross from one figure to the other. A clause number after a result number indicates that only that clause is consumed.
\item Results whose proofs consume only definitions, conventions and remarks stand without inputs; the import at the right of the second band of Figure \ref{fig.flowchart.a} feeds a remark only, and the import at the left of its first band is consumed only in Chapter \ref{cha.pfapp}.
\item Box colors follow the rule colors of the environments: red for theorems, blue for lemmas, grey for corollaries, orange for imports; dashed boxes are the theorems of the Invitation restated from results of the body.\end{itemize}

\newcommand{\flowpg}[1]{{\color{downgrey}\fontsize{6}{7}\selectfont p.~\pageref{#1}}}
\tikzset{
  flowchart/.style={x=1in, y=1in, every node/.style={font=\sffamily\fontsize{7}{8.4}\selectfont}},
  res/.style={draw=#1, fill=#1!7, rectangle, minimum width=1.40in, minimum height=0.40in, text width=1.28in, align=center, inner sep=2pt, line width=0.5pt,
     append after command={\pgfextra{\draw[#1, line width=2.2pt] (\tikzlastnode.north west) -- (\tikzlastnode.south west);}}},
  imp/.style={res=gloryorange},
  inv/.style={res=gloryred, dashed, fill=white},
  stk/.style={draw=#1, fill=white, rectangle, minimum width=1.40in, minimum height=0.22in, text width=1.30in, align=center, inner sep=1pt, line width=0.4pt, font=\sffamily\fontsize{6}{7}\selectfont},
  stki/.style={stk=gloryorange, fill=gloryorange!7},
  stkg/.style={stk=#1, dashed},
  dep/.style={-stealth, line width=0.55pt, draw=downgrey!85!black},
  elabel/.style={font=\sffamily\fontsize{6}{7}\selectfont, fill=white, inner sep=1pt},
  band/.style={draw=topgrey, fill=topgrey!22, rounded corners=3pt, line width=0.4pt},
  bandlabel/.style={anchor=north west, font=\sffamily\bfseries\fontsize{7.5}{9}\selectfont, text=downgrey!75!black, inner sep=2pt, fill=none}}

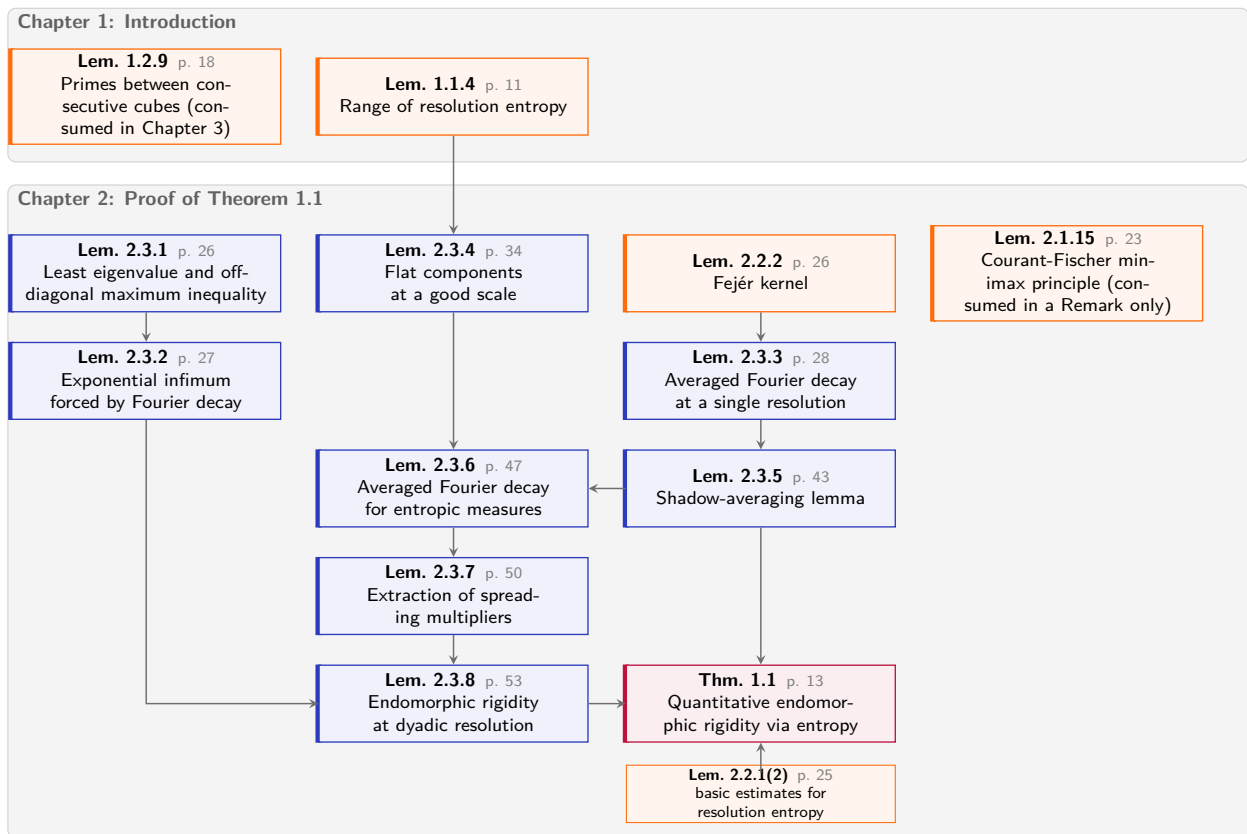
\begin{figure}[htbp]
\captionsetup{font=sf}
\centering
\begin{tikzpicture}[flowchart]
\draw[band] (-0.02,0.00) rectangle (6.44,-0.80);
\node[bandlabel] at (0.00,0.00) {Chapter~\ref{cha.intro}: Introduction};
\draw[band] (-0.02,-0.92) rectangle (6.44,-4.32);
\node[bandlabel] at (0.00,-0.92) {Chapter~\ref{cha.pfmain}: Proof of Theorem~\ref{thm.main}};
\node[imp] (nolem_primegap) at (0.70,-0.46) {\textbf{Lem.~\ref{olem.primegap}}\enspace\flowpg{olem.primegap}\\Primes between consecutive cubes (consumed in Chapter~\ref{cha.pfapp})};
\node[imp] (nolem_entrange) at (2.30,-0.46) {\textbf{Lem.~\ref{olem.entrange}}\enspace\flowpg{olem.entrange}\\Range of resolution entropy};
\node[res=gloryblue] (nlem_lsv) at (0.70,-1.38) {\textbf{Lem.~\ref{lem.lsv}}\enspace\flowpg{lem.lsv}\\Least eigenvalue and off-diagonal maximum inequality};
\node[res=gloryblue] (nlem_flatscale) at (2.30,-1.38) {\textbf{Lem.~\ref{lem.flatscale}}\enspace\flowpg{lem.flatscale}\\Flat components at a good scale};
\node[imp] (nolem_fejer) at (3.90,-1.38) {\textbf{Lem.~\ref{olem.fejer}}\enspace\flowpg{olem.fejer}\\Fej\'er kernel};
\node[imp] (nolem_courfish) at (5.50,-1.38) {\textbf{Lem.~\ref{olem.courfish}}\enspace\flowpg{olem.courfish}\\Courant-Fischer minimax principle (consumed in a Remark only)};
\node[res=gloryblue] (nlem_unifeinf) at (0.70,-1.94) {\textbf{Lem.~\ref{lem.unifeinf}}\enspace\flowpg{lem.unifeinf}\\Exponential infimum forced by Fourier decay};
\node[res=gloryblue] (nlem_fouravg) at (3.90,-1.94) {\textbf{Lem.~\ref{lem.fouravg}}\enspace\flowpg{lem.fouravg}\\Averaged Fourier decay at a single resolution};
\node[res=gloryblue] (nlem_decompave) at (2.30,-2.50) {\textbf{Lem.~\ref{lem.decompave}}\enspace\flowpg{lem.decompave}\\Averaged Fourier decay for entropic measures};
\node[res=gloryblue] (nlem_shadave) at (3.90,-2.50) {\textbf{Lem.~\ref{lem.shadave}}\enspace\flowpg{lem.shadave}\\Shadow-averaging lemma};
\node[res=gloryblue] (nlem_extract) at (2.30,-3.06) {\textbf{Lem.~\ref{lem.extract}}\enspace\flowpg{lem.extract}\\Extraction of spreading multipliers};
\node[res=gloryblue] (nlem_equildy) at (2.30,-3.62) {\textbf{Lem.~\ref{lem.equildy}}\enspace\flowpg{lem.equildy}\\Endomorphic rigidity at dyadic resolution};
\node[res=gloryred] (nthm_main) at (3.90,-3.62) {\textbf{Thm.~\ref{thm.main}}\enspace\flowpg{thm.main}\\Quantitative endomorphic rigidity via entropy};
\node[stki] at (3.90,-4.09) {\textbf{Lem.~\ref{olem.entropy}\ref{olem.entropy.c5}}\enspace\flowpg{olem.entropy}\\basic estimates for resolution entropy};
\draw[dep] (3.90,-3.98) -- (nthm_main.south);
\draw[dep] (nolem_entrange.south) -- (nlem_flatscale.north);
\draw[dep] (nolem_fejer.south) -- (nlem_fouravg.north);
\draw[dep] (nlem_lsv.south) -- (nlem_unifeinf.north);
\draw[dep] (nlem_fouravg.south) -- (nlem_shadave.north);
\draw[dep] (nlem_flatscale.south) -- (nlem_decompave.north);
\draw[dep] (nlem_shadave.west) -- (nlem_decompave.east);
\draw[dep] (nlem_decompave.south) -- (nlem_extract.north);
\draw[dep] (nlem_extract.south) -- (nlem_equildy.north);
\draw[dep] (nlem_unifeinf.south) -- (0.70,-3.62) -- (nlem_equildy.west);
\draw[dep] (nlem_equildy.east) -- (nthm_main.west);
\draw[dep] (nlem_shadave.south) -- (nthm_main.north);
\end{tikzpicture}
\caption{Dependency flowchart of the results of Chapters~\ref{cha.intro} and \ref{cha.pfmain}.}
\label{fig.flowchart.a}
\end{figure}

\begin{figure}[p]
\captionsetup{font=sf}
\centering
\begin{tikzpicture}[flowchart]
\draw[band] (-0.02,0.00) rectangle (6.44,-4.50);
\node[bandlabel] at (0.00,0.00) {Chapter~\ref{cha.pfapp}: Proofs of subsequent theorems};
\draw[band] (-0.02,-4.62) rectangle (6.44,-5.42);
\node[bandlabel] at (0.00,-4.62) {Section~\ref{sec.ad}: Invitation to article};
\node[imp] (nolem_dirichlet) at (0.70,-0.46) {\textbf{Lem.~\ref{olem.dirichlet}}\enspace\flowpg{olem.dirichlet}\\Dirichlet approximation};
\node[res=gloryblue] (nlem_gapshad) at (2.30,-0.46) {\textbf{Lem.~\ref{lem.gapshad}}\enspace\flowpg{lem.gapshad}\\Gap criterion for progression-shadowing};
\node[res=gloryblue] (nlem_ratio) at (3.90,-0.46) {\textbf{Lem.~\ref{lem.ratio}}\enspace\flowpg{lem.ratio}\\Sublacunary sets and progression-shadowing depth};
\node[res=gloryred] (nthm_exist) at (5.50,-0.46) {\textbf{Thm.~\ref{thm.exist}}\enspace\flowpg{thm.exist}\\Existential endomorphic rigidity via entropy};
\node[stkg=gloryred] at (5.50,-0.93) {\textbf{Thm.~\ref{thm.main}}\enspace\flowpg{thm.main}\\quantitative endomorphic rigidity};
\node[res=gloryblue] (nlem_rotfill) at (0.70,-1.02) {\textbf{Lem.~\ref{lem.rotfill}}\enspace\flowpg{lem.rotfill}\\Orbit filling from an irrationality measure};
\node[res=gloryblue] (nlem_depbound) at (2.30,-1.02) {\textbf{Lem.~\ref{lem.depbound}}\enspace\flowpg{lem.depbound}\\Bounds on progression-shadowing depth for canonical sets};
\node[imp] (nolem_bw) at (0.70,-1.58) {\textbf{Lem.~\ref{olem.bw}}\enspace\flowpg{olem.bw}\\Effective irrationality measure for ratios of logarithms};
\node[res=gloryblue] (nlem_primedep) at (2.30,-1.58) {\textbf{Lem.~\ref{lem.primedep}}\enspace\flowpg{lem.primedep}\\Progression-shadowing depth of the primes};
\node[stki] at (2.30,-2.05) {\textbf{Lem.~\ref{olem.primegap}}\enspace\flowpg{olem.primegap}\\primes between consecutive cubes};
\node[res=gloryred] (nthm_blprime) at (3.90,-1.58) {\textbf{Thm.~\ref{thm.blprime}}\enspace\flowpg{thm.blprime}\\Rigidity after a canonical endomorph};
\node[imp] (nolem_compact) at (5.50,-1.58) {\textbf{Lem.~\ref{olem.compact}}\enspace\flowpg{olem.compact}\\Weak-star compactness of the probability measures};
\node[imp] (nolem_weier) at (0.70,-2.50) {\textbf{Lem.~\ref{olem.weier}}\enspace\flowpg{olem.weier}\\Uniform approximation by trigonometric polynomials};
\node[res=gloryblue] (nlem_einfleb) at (2.30,-2.50) {\textbf{Lem.~\ref{lem.einfleb}}\enspace\flowpg{lem.einfleb}\\Lebesgue minorization from exponential infima};
\node[res=gloryred] (nthm_entdim) at (5.50,-2.50) {\textbf{Thm.~\ref{thm.entdim}}\enspace\flowpg{thm.entdim}\\Limiting endomorphic rigidity via upper entropy dimension};
\node[res=gloryblue] (nlem_floorseq) at (0.70,-3.06) {\textbf{Lem.~\ref{lem.floorseq}}\enspace\flowpg{lem.floorseq}\\Sublacunarity of the integer parts of a sequence};
\node[stki] at (3.90,-3.06) {\textbf{Lem.~\ref{olem.furst}}\enspace\flowpg{olem.furst}\\Furstenberg's nonlacunary semigroup lemma};
\node[stki] at (3.90,-3.30) {\textbf{Lem.~\ref{olem.hr}}\enspace\flowpg{olem.hr}\\Hardy--Ramanujan asymptotic};
\node[res=gloryred] (nthm_rjcanon) at (5.50,-3.06) {\textbf{Thm.~\ref{thm.rjcanon}}\enspace\flowpg{thm.rjcanon}\\Rudolph-Johnson for canonical multiplier sets};
\node[res=gloryred] (nthm_ad_cantor) at (0.70,-3.80) {\textbf{Thm.~\ref{thm.ad.cantor}}\enspace\flowpg{thm.ad.cantor}\\Rigidity failure for a lacunary semigroup (stated in the Invitation)};
\node[res=midgrey] (ncor_rjprime) at (2.30,-3.80) {\textbf{Cor.~\ref{cor.rjprime}}\enspace\flowpg{cor.rjprime}\\Rudolph-Johnson for cube-succeeding primes};
\node[stki] at (2.30,-4.27) {\textbf{Lem.~\ref{olem.primegap}}\enspace\flowpg{olem.primegap}\\primes between consecutive cubes};
\node[res=midgrey] (ncor_rjclass) at (3.90,-3.80) {\textbf{Cor.~\ref{cor.rjclass}}\enspace\flowpg{cor.rjclass}\\Classical Rudolph-Johnson};
\node[inv] (nthm_ad_sublac) at (0.70,-5.08) {\textbf{Thm.~\ref{thm.ad.sublac}}\enspace\flowpg{thm.ad.sublac}\\Invitation: sequential form of Theorem~\ref{thm.entdim}};
\node[inv] (nthm_equisparse) at (2.30,-5.08) {\textbf{Thm.~\ref{thm.equisparse}}\enspace\flowpg{thm.equisparse}\\Invitation: clauses of Corollary~\ref{cor.rjprime} and Theorem~\ref{thm.rjcanon}};
\node[inv] (nthm_ad_primeff) at (3.90,-5.08) {\textbf{Thm.~\ref{thm.ad.primeff}}\enspace\flowpg{thm.ad.primeff}\\Invitation: clauses of Theorem~\ref{thm.blprime}};
\node[res=gloryred] (nthm_lyons) at (5.50,-5.08) {\textbf{Thm.~\ref{thm.lyons}}\enspace\flowpg{thm.lyons}\\Lyons' conjecture at full entropy dimension};
\draw[dep] (5.50,-0.82) -- (nthm_exist.south);
\draw[dep] (2.30,-1.94) -- (nlem_primedep.south);
\draw[dep] (2.30,-4.16) -- (ncor_rjprime.south);
\draw[dep] (4.60,-3.06) -- (nthm_rjcanon.west);
\draw[dep] (nolem_dirichlet.south) -- (nlem_rotfill.north);
\draw[dep] (nlem_gapshad.south) -- (nlem_depbound.north);
\draw[dep] (nlem_gapshad.east) -- (nlem_ratio.west);
\draw[dep] (nlem_rotfill.east) -- (nlem_depbound.west);
\draw[dep] ([xshift=0.45in]nolem_bw.north) -- (1.15,-1.30) -- (1.85,-1.30) -- ([xshift=-0.45in]nlem_depbound.south);
\draw[dep] (nlem_depbound.south) -- (nlem_primedep.north);
\draw[dep] (nlem_depbound.east) -- (3.90,-1.02) -- (nthm_blprime.north);
\draw[dep] (nlem_primedep.east) -- (nthm_blprime.west);
\draw[dep] ([yshift=-0.10in]nthm_exist.west) -- (4.70,-0.56) -- (4.70,-1.58) -- (nthm_blprime.east);
\draw[dep] ([yshift=-0.10in]nthm_exist.west) -- (4.70,-0.56) -- (4.70,-2.40) -- ([yshift=0.10in]nthm_entdim.west);
\draw[dep] (nlem_ratio.north) -- (3.90,-0.13) -- (6.30,-0.13) -- (6.30,-2.50) -- (nthm_entdim.east);
\draw[dep] (nolem_compact.south) -- (nthm_entdim.north);
\draw[dep] (nolem_weier.east) -- (nlem_einfleb.west);
\draw[dep] ([yshift=-0.10in]nlem_einfleb.east) -- ([yshift=-0.10in]nthm_entdim.west);
\draw[dep] (nthm_entdim.south) -- (nthm_rjcanon.north);
\draw[dep] (nlem_floorseq.north) -- (0.70,-2.78) -- (5.05,-2.78) -- ([xshift=-0.45in]nthm_rjcanon.north);
\draw[dep, dotted] (nthm_rjcanon.south) -- (5.50,-3.50) -- (2.30,-3.50) -- (ncor_rjprime.north);
\draw[dep, dotted] (nthm_rjcanon.south) -- (5.50,-3.50) -- (3.90,-3.50) -- (ncor_rjclass.north);
\draw[dep, dotted] ([yshift=-0.10in]nthm_entdim.east) -- (6.30,-2.60) -- (6.30,-5.08) -- (nthm_lyons.east);
\end{tikzpicture}
\caption{Dependency flowchart of the results of Chapter~\ref{cha.pfapp} and the Invitation.}
\label{fig.flowchart.b}
\end{figure}
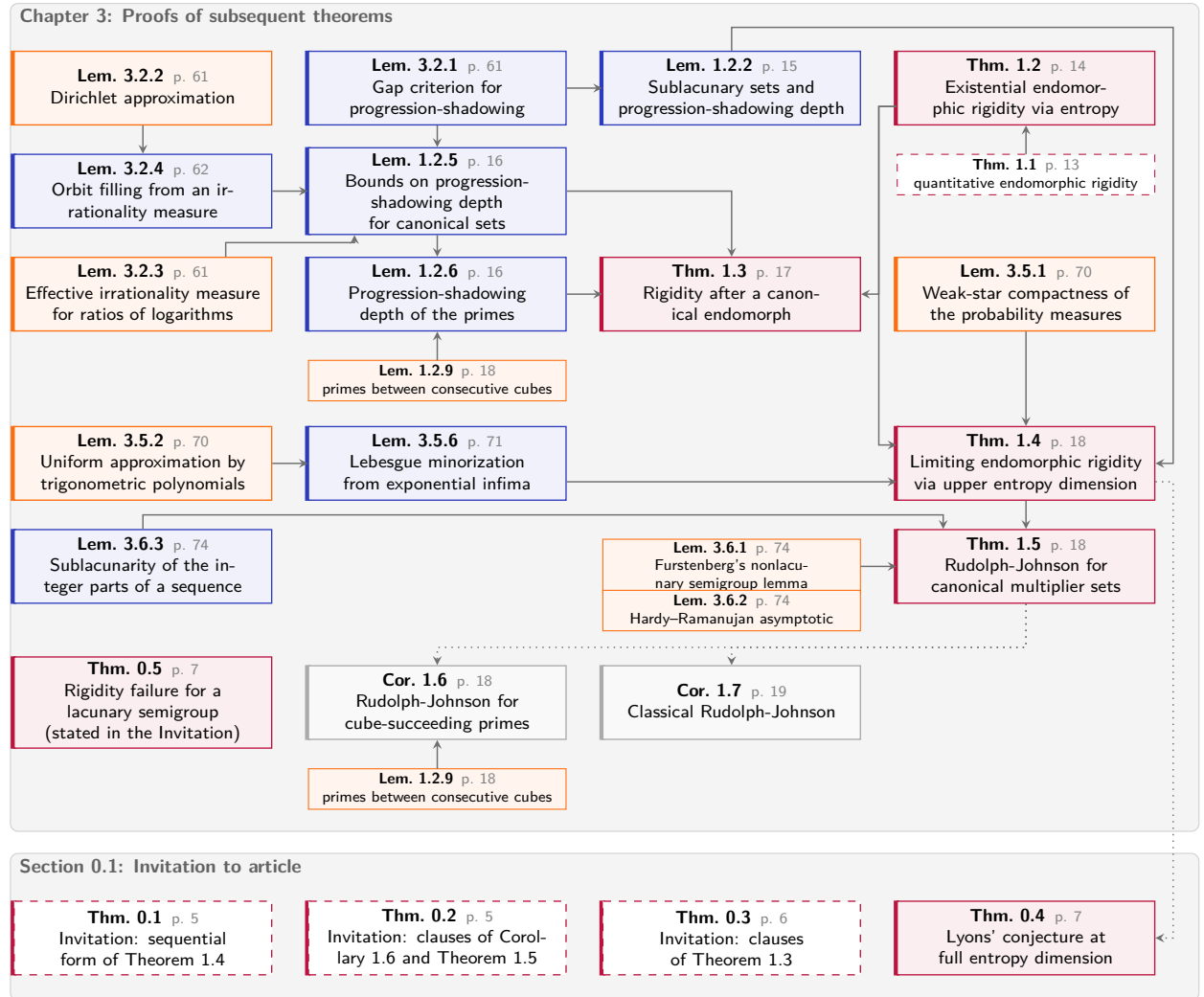

\textbf{Peter Burton:} \\
Department of Mathematics and Statistics\\
University of Wyoming, Laramie WY \\
\url{peterburton1728@gmail.com}\\
\\
\textbf{Kate Juschenko:} \\
Department of Mathematics\\
University of Texas at Austin, Austin TX \\
\url{k.juschenko@austin.utexas.edu}


\begin{thebibliography}{99}

\bibitem{MR4082895} C.~Badea and S.~Grivaux, Kazhdan constants, continuous probability measures with large Fourier coefficients and rigidity sequences, \emph{Comment. Math. Helv.} \textbf{95} (2020), 99--127.

\bibitem{MR4822284} C.~Badea and S.~Grivaux, Around Furstenberg's times $p$, times $q$ conjecture: times $p$-invariant measures with some large Fourier coefficients, \emph{Discrete Anal.} (2024), Paper No.~10, 31 pp.

\bibitem{MR1234835} A.~Baker and G.~W\"ustholz, Logarithmic forms and group varieties, \emph{J. Reine Angew. Math.} \textbf{442} (1993), 19--62.

\bibitem{MR858400} A.~Balog and A.~Perelli, Exponential sums over primes in short intervals, \emph{Acta Math. Hungar.} \textbf{48} (1986), 223--228.

\bibitem{BdJLR14} V.~Bergelson, A.~del Junco, M.~Lema\'nczyk and J.~Rosenblatt, Rigidity and non-recurrence along sequences, \emph{Ergodic Theory Dynam. Systems} \textbf{34} (2014), 1464--1502.

\bibitem{Bosh94} M.~D. Boshernitzan, Density modulo $1$ of dilations of sublacunary sequences, \emph{Adv. Math.} \textbf{108} (1994), 104--117.

\bibitem{MR2018901} N.~Bourbaki, \emph{Integration I. Chapters 1--6}, Elements of Mathematics, Springer-Verlag, Berlin, 2004.

\bibitem{MR2563089} J.~Bourgain, E.~Lindenstrauss, P.~Michel and A.~Venkatesh, Some effective results for $\times a \times b$, \emph{Ergodic Theory Dynam. Systems} \textbf{29} (2009), 1705--1722.

\bibitem{EntRates} P.~Burton and K.~Juschenko, Entropic rates from eigenfunction structure for $\times 2 \times 3$, preprint (2026).

\bibitem{MR2239987} T.~M. Cover and J.~A. Thomas, \emph{Elements of Information Theory}, second ed., Wiley-Interscience, Hoboken, NJ, 2006.

\bibitem{MR3584567} A.~W. Dudek, An explicit result for primes between cubes, \emph{Funct. Approx. Comment. Math.} \textbf{55} (2016), 177--197.

\bibitem{MR2247967} M.~Einsiedler, A.~Katok and E.~Lindenstrauss, Invariant measures and the set of exceptions to Littlewood's conjecture, \emph{Ann. of Math. (2)} \textbf{164} (2006), 513--560.

\bibitem{MR2648695} M.~Einsiedler and E.~Lindenstrauss, Diagonal actions on locally homogeneous spaces, in \emph{Homogeneous Flows, Moduli Spaces and Arithmetic}, Clay Math. Proc.~\textbf{10}, Amer. Math. Soc., Providence, RI, 2010, pp.~155--241.

\bibitem{MR2723325} M.~Einsiedler and T.~Ward, \emph{Ergodic Theory with a view towards Number Theory}, Graduate Texts in Mathematics, vol.~259, Springer-Verlag London, Ltd., London, 2011.

\bibitem{MR3410920} T.~Eisner, B.~Farkas, M.~Haase and R.~Nagel, \emph{Operator Theoretic Aspects of Ergodic Theory}, Graduate Texts in Mathematics, vol.~272, Springer, Cham, 2015.

\bibitem{EG11} T.~Eisner and S.~Grivaux, Hilbertian Jamison sequences and rigid dynamical systems, \emph{J. Funct. Anal.} \textbf{261} (2011), 2013--2052.

\bibitem{MR4746037} A.~Fan, H.~Queff\'elec and M.~Queff\'elec, The Furstenberg set and its random version, \emph{Enseign. Math. (2)} \textbf{70} (2024), 61--120.

\bibitem{FMK25} A.~Ferr\'e Moragues and A.~Koutsogiannis, Furstenberg systems of certain sequences of superpolynomial growth, arXiv:2510.11957 (2025).

\bibitem{MR1681462} G.~B. Folland, \emph{Real Analysis: Modern Techniques and Their Applications}, second ed., John Wiley \& Sons, Inc., New York, 1999.

\bibitem{MR2585398} N.~Frantzikinakis, Equidistribution of sparse sequences on nilmanifolds, \emph{J. Anal. Math.} \textbf{109} (2009), 353--395.

\bibitem{MR3613710} N.~Frantzikinakis, Some open problems on multiple ergodic averages, \emph{Bull. Hellenic Math. Soc.} \textbf{60} (2016), 41--90.

\bibitem{MR0213508} H.~Furstenberg, Disjointness in ergodic theory, minimal sets, and a problem in Diophantine approximation, \emph{Math. Systems Theory} \textbf{1} (1967), 1--49.

\bibitem{MR4042627} X.~Gao, J.~Ma, K.~Song and Y.~Zhang, On the Fourier transform of coin-tossing type measures, \emph{J. Math. Anal. Appl.} \textbf{484} (2020), 123706.

\bibitem{MR4176691} X.~Gao and S.~Wen, Note on Fourier--Stieltjes coefficients of coin-tossing measures, \emph{Bull. Aust. Math. Soc.} \textbf{102} (2020), 479--489.

\bibitem{MR2280879} G.~H. Hardy and S.~Ramanujan, Asymptotic formulae in combinatory analysis, \emph{Proc. London Math. Soc. (2)} \textbf{17} (1918), 75--115.

\bibitem{MR2445243} G.~H. Hardy and E.~M. Wright, \emph{An Introduction to the Theory of Numbers}, sixth ed., Oxford University Press, Oxford, 2008.

\bibitem{MR2978290} R.~A. Horn and C.~R. Johnson, \emph{Matrix Analysis}, second ed., Cambridge University Press, Cambridge, 2013.

\bibitem{MR0625600} J.~E.~Hutchinson, Fractals and self-similarity, \emph{Indiana Univ. Math. J.} \textbf{30} (1981), 713--747.

\bibitem{MR1194793} A.~S.~A. Johnson, Measures on the circle invariant under multiplication by a nonlacunary subsemigroup of the integers, \emph{Israel J. Math.} \textbf{77} (1992), 211--240.

\bibitem{MR1351328} A.~Johnson and D.~J. Rudolph, Convergence under $\times_q$ of $\times_p$ invariant measures on the circle, \emph{Adv. Math.} \textbf{115} (1995), 117--140.

\bibitem{MR3642325} O.~Kallenberg, \emph{Random Measures, Theory and Applications}, Probability Theory and Stochastic Modelling, vol.~77, Springer, Cham, 2017.

\bibitem{MR330068} A.~A. Karatsuba, Estimates for trigonometric sums by Vinogradov's method, and some applications, \emph{Trudy Mat. Inst. Steklov.} \textbf{112} (1971), 241--255; English transl.\ \emph{Proc. Steklov Inst. Math.} \textbf{112} (1973), 251--265.

\bibitem{MR2039503} Y.~Katznelson, \emph{An Introduction to Harmonic Analysis}, third ed., Cambridge University Press, Cambridge, 2004.

\bibitem{MR419394} L.~Kuipers and H.~Niederreiter, \emph{Uniform Distribution of Sequences}, Pure and Applied Mathematics, Wiley-Interscience, New York--London--Sydney, 1974.

\bibitem{MR2827904} E.~Lindenstrauss, Equidistribution in homogeneous spaces and number theory, in \emph{Proceedings of the International Congress of Mathematicians, Volume~I}, Hindustan Book Agency, New Delhi, 2010, pp.~531--557.

\bibitem{MR941238} R.~Lyons, On measures simultaneously $2$- and $3$-invariant, \emph{Israel J. Math.} \textbf{61} (1988), 219--224.

\bibitem{MR0274686} D.~S.~Mitrinović, \emph{Analytic inequalities}, Die Grundlehren der mathematischen Wissenschaften \textbf{165}, Springer-Verlag, New York-Berlin, 1970.

\bibitem{MR0107575} A.~R\'enyi, On the dimension and entropy of probability distributions, \emph{Acta Math. Acad. Sci. Hungar.} \textbf{10} (1959), 193--215.

\bibitem{MR0132570} A.~R\'enyi, On measures of entropy and information, in \emph{Proc. 4th Berkeley Sympos. Math. Statist. and Prob., Vol.~I}, Univ. California Press, Berkeley--Los Angeles, Calif., 1961, pp.~547--561.

\bibitem{JFM40.0388.03} F.~Riesz, Sur les op\'erations fonctionnelles lin\'eaires, \emph{C. R. Acad. Sci. Paris} \textbf{149} (1909), 974--977.

\bibitem{MR1157815} W.~Rudin, \emph{Functional Analysis}, second ed., International Series in Pure and Applied Mathematics, McGraw-Hill, New York, 1991.

\bibitem{MR1062766} D.~J. Rudolph, $\times 2$ and $\times 3$ invariant measures and entropy, \emph{Ergodic Theory Dynam. Systems} \textbf{10} (1990), 395--406.

\bibitem{MR1078738} R.~S.~Strichartz, Self-similar measures and their Fourier transforms. I, \emph{Indiana Univ. Math. J.} \textbf{39} (1990), 797--817.

\bibitem{MR2760403} T.~Tao, \emph{An Epsilon of Room, I: Real Analysis}, Graduate Studies in Mathematics, vol.~117, American Mathematical Society, Providence, RI, 2010.

\bibitem{MR4741329} S.~Usuki, $\times a$ and $\times b$ empirical measures, the irregular set and entropy, \emph{Ergodic Theory Dynam. Systems} \textbf{44} (2024), 1673--1692.

\bibitem{MR2358379} A.~Venkatesh, The work of Einsiedler, Katok and Lindenstrauss on the Littlewood conjecture, \emph{Bull. Amer. Math. Soc. (N.S.)} \textbf{45} (2008), 117--134.

\bibitem{MR62183} I.~M. Vinogradov, \emph{The Method of Trigonometrical Sums in the Theory of Numbers}, translated and revised by K.~F. Roth and A.~Davenport, Interscience, London--New York, 1954.

\bibitem{MR0648108} P.~Walters, \emph{An Introduction to Ergodic Theory}, Graduate Texts in Mathematics, vol.~79, Springer-Verlag, New York--Berlin, 1982.

\bibitem{MR4720212} A.~Weingartner, Uniform distribution of $\alpha n$ modulo one for a family of integer sequences, \emph{Unif. Distrib. Theory} \textbf{18} (2023), no.~2, 19--30.

\bibitem{MR1511862} H.~Weyl, \"Uber die Gleichverteilung von Zahlen mod.~Eins, \emph{Math. Ann.} \textbf{77} (1916), 313--352.

\end{thebibliography}
\end{document}